\documentclass[article, a4paper, 11pt,reqno]{amsart}
\usepackage{tikz}
\usetikzlibrary{calc,fit, shapes.misc, intersections}
\usepackage{subcaption}
\usepackage{adjustbox} 

\usepackage[top = 1in, bottom = 1in, left=1.25in, right=1.25in]{geometry}

\usepackage[utf8]{inputenc}

\usepackage{amsthm,amsmath,amssymb}
\usepackage{amsmath}
\usepackage{amsfonts}
\usepackage{amssymb}
\usepackage[font=small]{caption}
\usepackage{ifdraft}
\usepackage{graphicx}
\usepackage[shortlabels]{enumitem} 
\usepackage{dsfont}
\usepackage{thmtools}
\usepackage{thm-restate}

\usepackage{verbatim}
\usepackage{centernot}
\setlist[enumerate]{label={\upshape{(\roman*)}}}
\usepackage{esvect}
\usepackage{relsize}
\usepackage{url}

\usepackage[draft=false,breaklinks,colorlinks]{hyperref} 
\hypersetup{colorlinks, linkcolor={red!50!black}, citecolor={green!50!black}, urlcolor={blue!50!black}}
\hypersetup{hypertexnames=false}

\usepackage[
backend=biber,
natbib=true,
style=alphabetic,
isbn=false,
maxbibnames=9,
url=false,
doi=false
]{biblatex}
\usepackage{csquotes}
\usepackage{xcolor}

\definecolor{ForestGreen}{RGB}{34,139,34}

\renewbibmacro{in:}{}
\DeclareFieldFormat{pages}{#1}

\DeclareFieldFormat[article,misc]{volume}{\mkbibbold{#1}\nopunct}
\DeclareFieldFormat[article,misc]{number}{(#1)}

 \DeclareFieldFormat{title}{\mkbibquote{#1}}
 \DeclareFieldFormat[misc]{date}{Preprint (#1)}
 
\newbibmacro{string+doi}[1]{%
	\iffieldundef{doi}{\iffieldundef{url}{#1}{\href{\thefield{url}}{#1}}}{\href{http://dx.doi.org/\thefield{doi}}{#1}}}	
\DeclareFieldFormat*{title}{\usebibmacro{string+doi}{#1}}

\ifdraft{
	\usepackage[colorinlistoftodos]{todonotes}
	\usepackage{letltxmacro}
	\LetLtxMacro{\oldtodo}{\todo}
	\usepackage{xspace}
	\renewcommand{\todo}[2][]{\oldtodo[#1]{#2}\xspace}%
	\usepackage{xcolor}
	\usepackage[mark,dirty={Modified!}]{gitinfo2}
	
}{
	\usepackage[disable]{todonotes}
}

\usepackage[nameinlink, capitalise]{cleveref}

\declaretheorem[style=definition, numberwithin=section]{definition}
\declaretheorem[sibling=definition]{theorem, lemma, corollary, proposition, fact, claim}
\declaretheorem[style=definition, sibling=definition]{remark, example, conjecture, problem, question, observation, assumption}

\newcommand{\averagedegree}{d}
\newcommand{\eps}{\varepsilon}
\newcommand{\dred}{\mathfrak{d}}
\newcommand{\AB}{L}
\newcommand{\forbidden}{Z^{\mathrm{forb}}}
\newcommand{\oriented}[1]{\vv{#1}}
\newcommand{\fmat}[1]{#1}
\DeclareMathOperator{\dist}{dist}
\DeclareMathOperator{\probability}{Pr}
\DeclareMathOperator{\expectation}{\mathbf{E}}
\DeclareMathOperator{\im}{im}

\newcounter{propcounter}

\allowdisplaybreaks

\newenvironment{proofclaim}[1][Proof of the claim]{\begin{proof}[#1]\renewcommand*{\qedsymbol}{\(\square\)}}{\end{proof}}
\renewcommand*{\qedsymbol}{\(\blacksquare\)}

\crefname{theorem}{theorem}{theorems}
\crefname{corollary}{corollary}{corollaries}
\crefname{example}{example}{examples}
\crefname{lemma}{lemma}{lemmas}
\crefname{proposition}{proposition}{propositions}
\crefname{definition}{definition}{definitions}
\crefname{observation}{observation}{observations}
\crefname{claim}{Claim}{Claims}
\crefname{conjecture}{conjecture}{conjectures}

\crefformat{section}{\S#2#1#3}
\crefformat{subsection}{\S#2#1#3}
\crefformat{subsubsection}{\S#2#1#3}
\crefrangeformat{section}{\S\S#3#1#4 to~#5#2#6}
\crefmultiformat{section}{\S\S#2#1#3}{ and~#2#1#3}{, #2#1#3}{ and~#2#1#3}

\title{
	The asymptotic version of the Erd\H os-S\'os conjecture and beyond}
\author[A.~Davoodi]{Akbar Davoodi}
\author[D.~Piguet]{Diana Piguet}
\address{The Czech Academy of Sciences, Institute of Computer Science, Pod Vod\'{a}renskou v\v{e}\v{z}\'{\i} 2, Prague,182 07, Czech Republic}
\email{davoodi@cs.cas.cz, piguet@cs.cas.cz}
\author[H.~Řada]{Hanka Řada}
\address{Czech Technical University in Prague, Faculty of Information Technology, Department of Applied Mathematics, Thákurova 9, 160 00 Prague 6, Czech Republic}
\email{hanka.rada@fit.cvut.cz}
\author[N.~Sanhueza-Matamala]{Nicolás Sanhueza-Matamala}
\address{Departamento de Ingeniería Matemática and CI$^2$MA, Universidad de Concepción, Chile}
\email{nsanhuezam@udec.cl}

\begin{document}

\begin{abstract}
Klimo\v{s}ov\'a, Piguet, and Rozho\v{n} conjectured that any graph with minimum degree $k/2$ and sufficiently many vertices of degree $k$ should contain all trees with $k$ edges.
We prove an asymptotic version of this conjecture for dense host graphs.
We obtain interesting corollaries: the first is an asymptotic version of the Erd\H{o}s--S\'os conjecture for dense host graphs, which works without any bounded-degree restriction on the guest trees.
Secondly, by leveraging recent results by Pokrovsky, we can translate our results to sparse host graphs in the case of bounded-degree guest trees.
\end{abstract}

\maketitle

\tableofcontents

\section{Introduction}

One of the most classical questions in graph theory is to determine the number of edges in a host graph $G$ that forces the existence of a copy of another guest graph~$H$.
Turán's theorem~\cite{Turan1941} gives a complete answer whenever $H$ is a clique,
and the Erd\H{o}s--Stone--Simonovits theorem \cite{ErdosStone1946, ErdosSimonovits1966} gives a satisfactory answer (up to lower order terms) for every graph with chromatic number at least $3$.
However, much less is known for bipartite $H$, and even the particular case of trees remains widely open.
A seminal conjecture by Erd\H{o}s and S\'{o}s~\cite{Erdos1964} says that graphs with average degree larger than $k-1$ should contain all $k$-edge trees.

\begin{conjecture}[Erd\H os--S\'os conjecture]\label{conj:E-S}
	Every graph $G$ with average degree $\averagedegree(G) > k-1$ contains every tree with $k$ edges.
\end{conjecture}

The conjecture has been verified for specific families of trees
(e.g. paths~\cite{ErdosGallai1959},
trees with diameter at most four~\cite{McLennan2005}, 
spiders~\cite{FanHongLiu2018}, etc.
See also~\cite{GKL2016,DGMT2018} for further variations in hypergraphs, and the survey by Stein~\cite{Stein2020} for more results).
Besomi, Stein and Pavez-Signé~\cite{BPS2021} verified the Erd\H{o}s--So\'s conjecture for bounded-degree trees in dense host graphs, which then was extended to sparse host graphs by Pokrovskiy~\cite{Pokrovskiy2024b}.
A solution of the Erd\H os--S\'os conjecture for all large enough trees was announced in the early 1990s by Ajtai, Komlós, Simonovits and Szemerédi,
though it has not yet been made available
(a sketch of the proof can be found in~\cite{AKSS2015}).

For some classes of trees (such as stars, paths, trees with diameter at most three) even a weaker condition than that of \Cref{conj:E-S} on the host graph suffices: specifically, it is enough to assume that $\Delta(G) \geq k$ and $\delta(G) \geq k/2$ (this can be seen, e.g. by following the proofs by Erd\H{o}s and Gallai~\cite{ErdosGallai1959}).
This is indeed weaker, since if a graph satisfies $\averagedegree(G) > k-1$ then a well-known argument implies that $G$ contains a subgraph $G'$ with $\Delta(G') \geq k$ and $\delta(G') \geq k/2$.

However, this new condition on the host graph is not enough to ensure the containment of all trees.
It fails for trees of diameter four, as shown by examples of Havet, Reed, Stein and Wood~\cite[\S 1]{HRSW2020}.
In those examples, the host graphs consist either of two cliques or a complete bipartite graph together with an extra universal vertex, which turns out to be the unique vertex of degree at least $k$.
In view of this situation, it is natural to think that perhaps the following is true: if we have a substantial number of vertices of degree at least $k$ in $G$ (instead of just one), and we also have the condition $\delta(G) \geq k/2$, then we should find all $k$-edge trees as subgraphs of~$G$.
A conjecture along these lines was proposed by Klimo\v sov\'a, Piguet, and Rozho\v n {\cite[Conjecture 1.4]{Rozhon2019}}.

\begin{conjecture}[Klimo\v sov\'a, Piguet, Rozho\v n]\label{conj:KPR}
	Every $n$-vertex graph~$G$ with $\delta(G) \geq k/2$ and at least $n/(2 \sqrt{k})$ vertices of degree at least~$k$ contains all $k$-edge trees.
\end{conjecture}

Our main result is an approximate version of \Cref{conj:KPR} for dense graphs, i.e., with quadratically many edges. 

\begin{theorem}[Main result]\label{thm:approDense-MinMax}
	For any $\eta, q>0$, there exists an $n_0\in \mathbb N$ such that for every $n\ge n_0$ and all $k\ge qn$, any $n$-vertex graph $G$ with minimal degree $\delta (G)\ge(1+\eta)k/2$ and with at least $\eta n$ vertices of degree at least $(1+\eta)k$ contains all $k$-edge trees.
\end{theorem}

\Cref{thm:approDense-MinMax} implies an approximate dense version of the Erd\H os-S\'os conjecture.

\begin{corollary}\label{cor:dense-approx E-S}
	For any $\eta, q > 0$ there exists an $n_0 \in \mathbb N$ such that for every $n\ge n_0$ and all $k\ge qn$ any $n$-vertex graph with average degree more than $(1+\eta)k$ contains any tree on at most $k$ edges as its subgraph. 
\end{corollary}

\Cref{cor:dense-approx E-S} strengthens similar results by Rozho\v n~\cite{Rozhon2019} and Besomi, Pavez-Sign\'e and Stein~\cite{BPS2019,BPS2021}, which had the extra requirement that the trees $T$ to be found as subgraphs satisfy $\Delta(T) = o(k)$.
It also gives a proof independent of the one proposed by Ajtai, Koml\'os, Simonovits, and Szemer\'edi~\cite{AKSS2015}
 in the case the host graph is dense, under the very mild strengthening that its average degree is required to be slightly larger than $k$ (see \cref{ssec:comparision} for a comparison between their approach and ours).
 Recently, Reed and Stein~\cite{BruceTalk} announced a resolution of the Erd\H{o}s--Sós conjecture for the case of `dense trees', i.e. that one can replace $(1 + \eta)k$ with $k-1$ in \Cref{cor:dense-approx E-S}.

Both of our results mentioned above work only in the setting of dense host graphs and linear-sized trees.
However, a recent key structural result by Pokrovskiy~\cite{Pokrovskiy2024a} proves that the dense case in fact encapsulates most of the difficulty of the general tree-embedding problem; his result reduces problems about embedding bounded-degree trees in host graphs (non-necessarily dense) to the case of dense host graphs.
In our situation, these tools allow us to translate our \Cref{thm:approDense-MinMax} to the sparse setting, in the case of bounded-degree trees.

\begin{corollary}\label{cor:approDense-MinMax-sparse}
For $\Delta\in \mathbb N$, and $\eta>0$, let $k_0\in \mathbb N$ be sufficiently large. 
Then for any $k\ge k_0$ and any graph $G$ with minimal degree $\delta(G)\ge (1+\eta)k/2$ and at least $\eta |V(G)|$ vertices of degree at least $(1+\eta) k$ contains as its subgraph  any tree $T$ on at most $k$ edges  
with bounded maximal degree $\Delta(T)\le \Delta$.
\end{corollary}

There is a chance that our methods can be used to attack other tree-embedding conjectures which combine minimum and maximum degree conditions; see \Cref{ssec:stateoftheart} for more discussion.

As a consequence of \Cref{cor:dense-approx E-S}, we can quickly obtain bounds for the multicolour Ramsey numbers of trees.
Given graphs $T_1, \dotsc, T_r$, the \emph{$r$-colour Ramsey number of $T_1, \dotsc, T_r$}, denoted by $R_r(T_1, \dotsc, T_r)$, is the least $N$ so that every complete graph on $N$ vertices which is edge-coloured with $\{1, \dotsc, r\}$, contains a monochromatic copy of $T_i$ in the $i$th colour, for some $i$.
We write $R_r(T)$ if $T = T_1 = \dotsb = T_r$.
Erd\H{o}s and Graham~\cite{Erdos-1981} conjectured that for every $r \geq 2$ and every $n$-vertex tree $T$, the bound $R_r(T) \leq rn + O(1)$ should hold.
This would follow from the validity of \Cref{conj:E-S}.
The case $r = 2$ was proven by Zhao~\cite{Zhao-2011} for all large $n$.
The following result gives an asymptotically tight bound for $R_r(T_1, \dotsc, T_r)$ for arbitrary trees, generalising a result of Piguet and Stein~\cite{PiguetStein} (for two colours) and Klimo\v sov\'a, Piguet, and Rozho\v n~\cite{KlimosovaPiguetRozhon2020} (who assumed further `skew' properties of the trees).
Also if $T = T_1 = \dotsb = T_r$ note that this gives an approximate version of the conjecture of Erd\H{o}s and Graham.

\begin{corollary} \label{corollary:ramsey}
	For $r \geq 2$ and $\eps > 0$, there exists $n_0$ such that for any trees $T_1, \dotsc, T_r$ with $\sum_{i=1}^r |V(T_i)| \geq n_0$, we have $R_r(T_1, \dotsc, T_r) \leq (1 + \eps) \sum_{i=1}^r |V(T_i)|$.
\end{corollary}

\subsection{Rough sketch of the proof of our main result}
The proof follows the well-known regularity method.
First, we prepare both the host graph and the tree for embedding (i.e., finding a copy of the tree within the graph).
We prepare the host graph by applying Szemerédi’s Regularity Lemma, which yields the so-called \emph{reduced graph}.
For the tree, we find a negligible set of cut vertices that partition the tree into smaller subtrees with specific properties.

Next, in the reduced graph, we identify a suitable structure to assist in embedding the tree.
This involves selecting two adjacent clusters, each with a sufficiently large degree relative to a structure that matches the shape of the tree.
One of these clusters will have a total (weighted) degree slightly greater than $k$ (after re-scaling by the cluster size), a condition guaranteed by the degree-inheritance principle of the cluster graph. The second cluster will be a carefully chosen neighbour.
These two adjacent clusters will host the cut vertices of the tree.
The tree-specific structure is designed to accommodate the small subtrees.
While many previous tree-embedding problems have used a simple matching for this structure, we have developed a significantly more sophisticated approach -- a \emph{skew-matching pair}. This structure builds on the idea of fractional matching but assigns different weights to each end-vertex. This refinement enables a much more precise replication of the global structure of the tree, offering a substantial improvement over previous methods and underscoring its key role in proving our result.

Finally, we embed the tree using the structure established in the cluster graph during the previous phase. The embedding process combines standard properties of regular pairs with probabilistic arguments.

\subsection{Organisation of the paper}
In \Cref{sec:notaion}, we set some very general notation.
In \Cref{section:corollaryproof}, we derive from our main result (\Cref{thm:approDense-MinMax}) a proof of the announced corollaries, \Cref{cor:dense-approx E-S} and \Cref{cor:approDense-MinMax-sparse}.
 In \Cref{sec:main-proofs}, we first give some rough overview of the proof, and we prove our main result (\Cref{thm:approDense-MinMax}) in \Cref{ssec:mainproof}, assuming the three main technical propositions: the Tree-Coating Lemma (\Cref{prop:coatoftree}), the Structural Proposition (\Cref{prop:weighted-structural}), and the Tree-Embedding Lemma (\Cref{lemma:treeembedding}).
 In \Cref{sec:coating} we prove the Tree-Coating Lemma.

 The proof of the Structural Proposition (\Cref{prop:weighted-structural}) spans the next four sections.
 In \Cref{sec:struct} we give some additional notation and definitions of the objects we will use in that proof, including the key `skew-matchings'.
 In \Cref{section:fractionalGE} we investigate the structure of skew-matchings in general graphs, based on known results for matchings in graphs.
 In \Cref{sec:Matching-lemmas} we prove various `building blocks' associated with skew-matchings that serve as lemmas in our proof of \Cref{prop:weighted-structural} (Structural Proposition); then we finally give the proof in~\Cref{sec:new-proof-structural}.
 Finally, in~\Cref{sec:embed} we prove \Cref{lemma:treeembedding} (Tree Embedding Lemma).
 We conclude with closing remarks in \Cref{sec:conclusion}.
 
 We remark that the bulk of the length of the paper is the proof of the two main technical lemmas: the proof of the Structural Proposition (\Cref{prop:weighted-structural}) (which is in \Cref{sec:struct}--\Cref{sec:new-proof-structural}) and the proof of the Tree-Embedding Lemma (\Cref{lemma:treeembedding}) (\Cref{sec:embed}).
 The latter can be read independently of the previous sections (using the definition of skew-matchings pair as a black-box).
\addtocontents{toc}{\protect\setcounter{tocdepth}{1}}
\section{Notation}\label{sec:notaion}

\subsection*{Basic notation}

We sometimes use the hierarchy symbol $\ll$, where $a\ll b$ informally means \emph{``$a$ is much smaller than $b$''},
and formally translates to ``there is a monotone increasing function $f:(0,1)\rightarrow (0,1)$ such that for any $a, b$ satisfying $a\le f(b)$, the following holds''. 

\subsection*{Graphs}
Given a graph $G$ and a subset $S \subseteq V(G)$, we write $G-S$ to refer to the graph obtained from $G$ after removing $S$ and any edges incident to $S$.

For two disjoint subsets $X$ and $ Y $ of $ V(G) $, we denote by $ d(X,Y) $ the \emph{bipartite density} of the pair $ (X,Y) $,  defined by \[ d(X,Y) :=\frac{|E(X,Y)|}{|X| |Y|}, \] in which $|E(X,Y)|$ denotes the number of edges between $ X $ and $ Y $. 
For a graph $G$, we denote by $d(G)$ the \emph{average degree} of $G$, i.e. $d(G):= 2|E(G)|/|V(G)|$. 
For a vertex $v\in V(G)$, let $N_G(v)$ denote the set of neighbours of $v$ in $G$. We will omit $G$ from the notation if the graph is clear from context. A \emph{vertex-cover} of a graph $G$ (or just \emph{cover} for short) is a set $C\subseteq V(G)$ such that for every edge $\fmat{xy}\in E(G)$ we have $\{x,y\}\cap C\neq \emptyset$.

A graph $H$ \emph{embeds} in a graph $G$ if there is a subgraph $H'\subseteq G$ such that $H'$ is isomorphic to $H$. 

\subsection*{Digraphs}
In \emph{digraphs} every edge is oriented, meaning that it consists of an ordered pair of vertices, which we shall denote by a pair with an arrow on top (e.g. $\oriented{uv}$ to mean the ordered pair $(u,v)$) to differentiate well with the non-oriented pairs.
In the digraphs we will use during our proofs we admit cycles of length 2 (where the pairs of edges $\oriented{uv}$ and $\oriented{vu}$ are both present), but no digraph will have parallel edges in the same direction, and we also forbid loops.
If every time the pair $\oriented{uv}$ is present we also have $\oriented{vu}$, we say that the digraph is \emph{symmetric}. 
Given a vertex $u \in V(G)$, the sets $N^+_G(u) = \{ v \in V(G) : \oriented{uv} \in E(G) \}$ and $N^-_G(u) = \{ v \in V(G) : \oriented{vu} \in E(G) \}$ are the \emph{outneighbours} and \emph{inneighbours} of $u$, respectively.
Again, we will omit $G$ from the notation if the digraph is clear from context.

Given a graph $G$, its \emph{associated digraph $G^\leftrightarrow$} is the digraph with the same vertex set as $G$ and  both  $\oriented{uv}$ and $\oriented{vu}$ are present for each undirected $\fmat{uv} \in E(G)$.
Observe that this digraph is symmetric.
Moreover, we use $V(G)$ instead of $V(G^{\leftrightarrow})$ and $N_G(u)$ instead of $N_G^{+}(u)$ or $N_G^{-}(u)$ when we deal with associated digraphs. (It is because in this case, we have $V(G) = V(G^{\leftrightarrow})$ and $N_G(u) = N_G^{+}(u) = N_G^{-}(u)$.)

\subsection*{Weighted graphs and digraphs}
A \emph{weighted graph} is a pair $(G,w)$ where $G$ is a graph and $w: E(G) \rightarrow \mathbb{R}^+$ is a function.
We define the \emph{weighted degree} of a vertex~$v$ as $\deg_w(v):= \sum_{u\in N(v)}w(\fmat{vu})$.
We can implicitly assume that $w(vu) = 0$ for any $uv \notin E(G)$, so it makes sense to evaluate $w$ on non-edges of $G$. 
Given this assumption, we also have $\deg_w(v)= \sum_{u\in V(G)}w(\fmat{vu})$.
Given a set $S \subseteq V(G)$, we also define $\deg_w(v, S):= \sum_{u\in S} w(\fmat{vu})$.
Also, the \emph{minimum} and \emph{maximum weighted degree} are defined as
$\delta_w(G):=\min\{\deg_w(v): v\in V(G)\}$,
and $\Delta_w(G):= \max\{\deg_w(v): v\in V(G)\}$, respectively. Similarly, for $v\in V(G)$ we define $N_w(v):= \{u\in V(G):\:  w(\fmat{vu})>0\}$.

Similarly, we can define a weighted digraph putting a weight on each directed edge.
We will mostly consider (weighted) \emph{digraphs} $G^{\leftrightarrow}$ associated to weighted graphs $(G,w)$ as well with their weights inherited from $G$, i.e. $w(\oriented{uv})=w(\oriented{vu})=w(uv)$.
However, in general the weights in digraphs do not need to be symmetric, and we may have symmetric digraph with non-symmetric weights\footnote{For example, in associated weighted digraphs $G^\leftrightarrow$ after altering its inherited weights; we may lose symmetry when using the concept of truncated weighted digraphs (\Cref{def:truncated-weighted-graph}).}. Similarly as above, we can define degrees in weighted directed graphs by $\deg_w(u):=\sum_{u\in V(G)}w(\oriented{vu})$,  $\deg_w(u, S):= \sum_{u\in S} w(\oriented{vu})$, and $N_w(G):= \{u\in V(G):\:  w(\oriented{vu})>0\}$.

\addtocontents{toc}{\protect\setcounter{tocdepth}{2}}
\section{Proof of the corollaries} \label{section:corollaryproof}

\subsection{Proof of \Cref{cor:dense-approx E-S}}

We  prove that our main result (\Cref{thm:approDense-MinMax}) implies the asymptotic version of the Erd\H os-S\'os Conjecture in the setting of dense graphs
(\Cref{cor:dense-approx E-S}).

\begin{proof}[Proof of~\Cref{cor:dense-approx E-S}]
	Let $k = rn$ with $r \geq q > 0$,
	and let $G$ be a graph on $n$ vertices with average degree at least $(1 + \eta)k$.
	It is well-known~\cite[Proposition 1.2.2]{Diestel} that $G$ contains an induced subgraph $H$ such that $\delta(H) \ge  \averagedegree(H)/2 \geq \averagedegree(G)/2 \geq (1 + \eta)k/2$.
	Let $m$ be the number of vertices of $H$. We clearly have $(1 + \eta)k/2 \leq \delta(H) < m \leq n$.
	
	For a given $\lambda > 0$, let $X_\lambda$ be the set of vertices of $H$ whose degree in $H$ is at least $(1 + \lambda)k$.
	Then we have
	\begin{align*}
		(1 + \eta) k m & \leq m \averagedegree(H) = \sum_{v \in V(H)} \deg_H(v) \leq |X_\lambda| m + (m - |X_\lambda|) (1 + \lambda) k,
	\end{align*}
	which, by rearranging, gives
	\begin{align*}
		|X_\lambda| & \geq \frac{(\eta - \lambda)km}{m - (1 + \lambda)k} \geq (\eta - \lambda) k.
	\end{align*}
	From now on, fix the choice $\lambda := \eta k /(m+k)$.
	This choice satisfies $\eta \geq \lambda \geq \eta r/(1+r)$,
	and from the previous calculations we deduce that $H$ satisfies $\delta(H) > (1 + \lambda)k/2$,
	and has at least $(\eta - \lambda)k = \lambda m$ vertices of degree at least $(1 + \lambda)k$.
	Thus, the statement follows by applying \Cref{thm:approDense-MinMax} to $H$,
	with $\lambda$ playing the role of $\eta$.
\end{proof}

\subsection{Proof of \Cref{cor:approDense-MinMax-sparse}}\label{sec:corollary-sparse}

Now we prove \Cref{cor:approDense-MinMax-sparse}, which is a sparse version of \Cref{thm:approDense-MinMax} for bounded-degree trees.
We shall use the following very recent `hyperstability' result of Pokrovskiy~\cite[Theorem 5.2]{Pokrovskiy2024a}.

\begin{theorem}\label{thm:Pokrovskiy}
For any $\Delta \geq 1$ and $\varepsilon > 0$, there exists $d > 0$ such that the following holds.
Let~$T$ be a tree with $d$ edges and $\Delta(T)\le \Delta$. For any $n$-vertex graph $G$ having no copies of $T$, it is possible to delete $\varepsilon dn$ edges to get a graph $H$ each of whose components has a vertex-cover of order at most $(2+\varepsilon)d$. 
\end{theorem}

We shall also use the following result of Piguet and Stein~\cite[Theorem 2]{PiguetStein}, which is a dense approximate version of the Loebl--Komlós--Sós conjecture.
\begin{theorem}\label{thm:LKS}
Let $\eta,q \gg 1/n$.
Let $G$ be an $n$-vertex graph $G$ and $k\geq qn$.
If at least half of the vertices in $G$ have degree at least $(1 + \eta)k$, then $G$ contains every tree with at most $k$ edges.
\end{theorem}

\begin{proof}[Proof of \Cref{cor:approDense-MinMax-sparse}]
The proof consists of several steps.
First, we apply \Cref{thm:Pokrovskiy} to obtain a subgraph with nearly the same number of edges, but composed of a union of components, each with a small vertex-cover.
In the second step, we select a suitable component to work with and consider two cases: either the component is sufficiently small, or it is still quite large.
In the first case, we refine the component to obtain a subgraph $H_1$ which retains similar properties to the original graph. Exploiting the small size of the component, we reduce the problem to its dense instance and apply \Cref{thm:approDense-MinMax} (our main result) to $H_1$ to find a copy of the tree.
In the second case, we carefully select a small subgraph $H_2$ of the component that includes the entire vertex-cover. This reduction transforms the problem into the dense instance of the Loebl-Koml\'os-S\'os case, so we can apply \Cref{thm:LKS} to $H_2$ to obtain a copy of the tree.
\medskip

\noindent \emph{Step 1: Decomposition via hyperstability.}
Let $\Delta\in \mathbb N$ and $\eta>0$ be given.
Set \[\varepsilon:= \left(\frac{\eta}{8}\right)^4 , \qquad \eta' := \frac{\eta}{6^2},\] and let $k_0$ be large enough such that it satisfies the requirement on $d$ in \Cref{thm:Pokrovskiy} with input $\Delta$ and $\varepsilon$, 
and ensures $k_0 \geq qn$ for any $n \ge n_0$, where $q = \max\left\{ \varepsilon^{-1}, \Delta \right\}$ and $n_0$ is the largest threshold required by Theorems~\ref{thm:LKS} and \ref{thm:approDense-MinMax}, so that the conditions of both theorems are satisfied.

Let $T$ be a tree on at most $k$ edges with $\Delta(T) \leq \Delta$.
Assume that $G$ does not contain any copy of $T$, as otherwise we are done.
By \Cref{thm:Pokrovskiy}, we can erase a set of edges $E'\subseteq E(G)$ with $|E'|\le \varepsilon kn$ so that the resulting graph $H=(V(G), E(G)\setminus E')$ consists of components $C_1, \ldots, C_t$, each of which has a cover of order at most $3k$.  
\medskip

\noindent \emph{Step 2: Component selection.}
Let $L:=\{v\in V(G)\::\: \deg_G(v)\ge (1+\eta)k\}$ and $L_i=L\cap C_i$. Let $I^-:=\{i\in [t]\::\: |L_i|<\frac{\eta}{2}|C_i|\}$ and $I^+:=[t]\setminus I^-$. Let $L^+:=\bigcup_{i\in I^+}L_i$. We have 
\[|L^+|=|L|-|L^-|\ge \eta n-\sum_{i\in I^-}|L_i|>\eta n-\frac{\eta}{2}\sum_{i\in I^-}|C_i|\geq \frac{\eta}{2}n\Big(2-\dfrac{\sum_{i\in I^-}|C_i|}{n}\Big)\geq \frac{\eta}{2}n.\]
We shall discard $\bigcup_{i\in I^-}C_i$ and consider only $H^+:=\bigcup_{i\in I^+}C_i$, where each component contains a substantial portion of its vertices of large degree (in $G$).

Let $B_i:= \{v\in C_i\::\: \deg_{H}(v)<\deg_G(v)-\frac{\eta}{8}k\}$.
In words, $B_i$ corresponds to those vertices in $C_i$ whose degree dropped substantially when erasing the edges $E'$.
We have 
\[ \sum_{i\in I^+}|B_i| = \left|\bigcup_{i\in I^+}B_i\right|\le \left|\bigcup_{i\in [t]}B_i\right|\leq \frac{2 |E'|}{\eta k/8}\leq \frac{16\varepsilon kn}{\eta k}=\frac{16\varepsilon}{\eta}n,\]
which implies that
$$\sum_{i\in I^+}|C_i|\geq \sum_{i\in I^+}|L_i|=|L^+|\geq \frac{\eta n}{2}\geq \frac{\eta}{2}\cdot \frac{\eta}{16\varepsilon}\sum_{i\in I^+}|B_i|.$$
Hence, there is an index $i_0\in I^+$ such that
\begin{equation}
	|C_{i_0}|\geq \frac{\eta^2}{32\varepsilon}|B_{i_0}|.
	\label{equation:pokrovsky-ci0}
\end{equation}
We select this component and show there is a copy of $T$ in $C_{i_0}$. 
\medskip

\noindent \emph{Step 3: Small component case - reduction to the dense case.}
First, assume that $|C_{i_0}|<(\eta')^{-1} k$. Then observe that, by \eqref{equation:pokrovsky-ci0},
\[|B_{i_0}|\leq \frac{32\varepsilon}{\eta^2}|C_{i_0}|< \frac{32 \varepsilon}{\eta'\eta^2}k<\frac{\eta}{4}k.\]
Set $H_1:= H[C_{i_0}\setminus B_{i_0}]$. Trivially $|V(H_1)| \leq |C_{i_0}|$ holds. 
Then, as $i_0\in I^+$, we have  $|L_{i_0}\cap V(H_1)|\ge |L_{i_0}|-|B_{i_0}|\geq  \frac{\eta}{2}|C_{i_0}|-\frac{32\varepsilon }{\eta^2}|C_{i_0}|\geq\frac{\eta }{4}|C_{i_0}|\geq \eta'|V(H_1)|$. 
Moreover, for every $v\in V(H_1)$, we have 
\[\deg_{H_1}(v)\geq \deg_{H}(v)-|B_{i_0}|\geq \deg_G(v)-\frac{\eta}{8}k-\frac{\eta}{4}k.\]
Therefore, $\delta(H_1)\geq (1+\eta)k/2-\frac{3\eta}{4}k/2\geq (1+\eta')k/2$ and for any $v\in L_{i_0}\cap V(H_1)$, we obtain $\deg_{H_1}(v)\ge (1+\eta)k-\frac{3\eta}{8}k\geq (1+\eta')k$. 
We thus can apply \Cref{thm:approDense-MinMax} to $H_1$ (in which we have at least $k$ vertices), with $\eta'$ playing the role of $\eta$, $k/|V(H_1)|$ implicitly satisfying the role of $q$, and $k$ corresponding to itself. This allows us to obtain a copy of $T$ in $H_1 \subseteq G$.
\medskip

\noindent \emph{Step 4: Large component case - reduction to the Loebl-Koml\'os-S\'os dense case.}
We are left to consider the case when $|C_{i_0}|\geq (\eta')^{-1}k$.
Let $D_{i_0}$ be the vertex-cover of~$C_{i_0}$ of size at most $3k$.
This means that every edge in $C_{i_0}$ is incident to at least one vertex in $D_{i_0}$. 
Now, we note that
\[|L_{i_0}|\geq \frac{\eta}{2}|C_{i_0}|>\frac{32\varepsilon}{\eta^2}|C_{i_0}|+\frac{\eta}{4}|C_{i_0}|\geq |B_{i_0}|+\frac{\eta}{4\eta'}k\geq|B_{i_0}|+3|D_{i_0}|,\] where the first inequality holds since $i_0 \in I^+$, and the remaining inequalities follow from \eqref{equation:pokrovsky-ci0}, the choices of parameters $\varepsilon$ and $\eta'$, and $|D_{i_0}| \leq 3k$. 
Hence, we can choose $ U_{i_0} $ as a subset of $ L_{i_0} \setminus (D_{i_0} \cup B_{i_0}) $ with size $ 2|D_{i_0}| $.

Now, we set $H_2:= H[U_{i_0}\cup D_{i_0}]$. Observe that as $D_{i_0}$ is a vertex-cover and disjoint from $U_{i_0}$, so every vertex in $U_{i_0}$ has all of its neighbors in $H_2$ in $D_{i_0}$.
For any vertex $v \in U_{i_0} \cap V(H_2)$, we have
\[\deg_{H_2}(v)= \deg_{H}(v)\ge \deg_G(v)-\frac{\eta}{8}k\ge (1+\eta')k.\]
This is because $\deg_G(v) \geq (1 + \eta)k$ for $v \in U_{i_0} \subseteq L_{i_0}$,
and since $U_{i_0} \cap B_{i_0} = \emptyset$, the vertex $v$ lost at most $\frac{\eta}{8}k$ of its degree in passing from $G$ to $H$.

As ensuring a sufficiently high minimum degree in $H_2$ is challenging, our strategy is to apply \Cref{thm:LKS} instead of \Cref{thm:approDense-MinMax} to $H_2$.
To apply \Cref{thm:LKS}, we note that
$|V(H_2)| = |U_{i_0}| + |D_{i_0}| = 3|D_{i_0}| \leq 9k$.
Since $|U_{i_0}| = 2|D_{i_0}| > \frac{|V(H_2)|}{2}$, the conditions of \Cref{thm:LKS} are satisfied.
Thus, we apply Theorem~\ref{thm:LKS} to $H_2$ (in which we have at least $k$ vertices), with $k/|V(H_2)|$ playing the role of $q$, $\eta'$ playing the role of $\eta$, and $|V(H_2)|$ playing the role of $n$. These choices ensure that the conditions of the theorem are satisfied, allowing us to find a copy of $T$ in $H_2 \subseteq G$.
This finishes this case, and since there are no more cases, this finishes the proof.
\end{proof}

We remark that, using the stability result of Pokrovsky (\Cref{thm:Pokrovskiy}) in conjunction with our dense approximate version of the Erd\H{o}s-S\'os conjecture (\Cref{cor:dense-approx E-S}), it is possible to deduce an approximate version of the Erd\H{o}s-S\'os conjecture for bounded-degree trees in sparse graphs.
However, this was already deduced by Pokrosvky \cite[Theorem 1.6]{Pokrovskiy2024a} from the results of Rozho\v n~\cite{Rozhon2019} and Besomi, Pavez-Sign\'e and Stein~\cite{BPS2019,BPS2021} mentioned in the introduction.
In fact, this strategy was used to obtain an \emph{exact} version of the Erd\H{o}s-S\'os conjecture for large, bounded-degree trees~\cite[Theorem 1.16]{Pokrovskiy2024a}.
This strengthening additionally requires further stability analysis based on ideas of Besomi, Pavez-Sign\'e and Stein~\cite{BPS2021}, see~\cite{Pokrovskiy2024b} for full details.

\subsection{Proof of \Cref{corollary:ramsey}}

\begin{proof}[Proof of \Cref{corollary:ramsey}]
	We prove it with $2 r \eps$ in place of $\eps$.
	Let $n := \sum_{i=1}^r |V(T_i)|$, which we assume to be sufficiently large with respect to $r$ and $\eps$; and let $N := (1 + 2r\eps) n$.
	Let $K_N$ be $r$-edge-coloured.
	By averaging, there must exist $i \in \{1, \dotsc, r\}$ such that the subgraph $G_i$ consisting of the $i$th coloured edges has average degree at least $(1 + \eps)\max\{ |V(T_i)|, \eps n\}$.
	By \Cref{cor:dense-approx E-S}, $G_i$ must contain each tree with at most $\max\{ |V(T_i)|, \eps n\} \geq |E(T_i)|$ edges, so in particular there is an $i$-coloured copy of $T_i$ in $K_N$, as desired.
\end{proof}

\section{Proof of the main result}\label{sec:main-proofs}

This section focuses on proving \Cref{thm:approDense-MinMax}. The proof relies on the regularity method, which consists of three standard steps: pre-processing the host graph $G$ using Szemer\'edi's Regularity Lemma (\Cref{thm:SzemLemma}) which yields a so-called \emph{reduced graph} that captures the large-scale structure of $G$.
Then, we find a suitable structure (\Cref{prop:weighted-structural}) in the reduced graph.
Finally, we use the structure in the reduced graph to embed the guest tree $T$ in $G$ (\Cref{lemma:treeembedding}).

Before delving into the proof, we provide some insights and outline three key statements necessary for this proof.
(Readers unfamiliar with the regularity method and its terminology can check \Cref{ssec:regularity} for the main concepts and definitions.)
We will identify a specific structure in the reduced graph of the host graph, which we call a \emph{skew-matching pair} (see \Cref{def:weighted-anchored-pair}).
This structure helps in embedding the tree.
Similar to how a `connected matching' in the reduced graph serves as a guide to embed a long path or a cycle in a host graph, a skew-matching pair in the reduced graph will correspond to a specific part of the host graph where there is enough space to embed the tree.
However, the shape of the skew-matching pair is heavily dependent on the structure of the tree we want to embed.
To understand what kind of shape we need, we will observe that any tree can be broken down into a negligible-sized set of cut vertices and two forests $\mathcal{F}_A, \mathcal{F}_B$ that consist of tiny trees.
What is most relevant here is the colour classes (i.e. the natural bipartition) of the forests $\mathcal{F}_A$ and $\mathcal{F}_B$.
Trees with forests with colour classes of the same sizes are grouped into the same tree-class~$\mathcal{T}$.
Then the skew-matching pair encodes and represents some tree-class $\mathcal{T}$.

With this structure in mind, we can expand our initial discussion of the three main necessary lemmas for our approach.
We have a host graph $G$ and a guest tree $T$ which we need to embed into $G$.
The first element of the proof is the Tree-Coating Lemma (\Cref{prop:coatoftree}), which guarantees that $T$ belongs to a tree-class $\mathcal{T}$, according to the partition into cut vertices and forests which we mentioned before.
This tree-class~$\mathcal{T}$ defines the skew-matching pair we will look for in the reduced graph of $G$.
The second element we need is the Structural Proposition (\Cref{prop:weighted-structural}), that ensures that the required skew-matching pair exists in the reduced graph.
Finally, we need the Tree-Embedding Lemma (\Cref{lemma:treeembedding}). This lemma states that if the reduced graph contains a specific skew-matching pair, we can embed in the host graph $G$ any tree belonging to the tree-class this pairs represents.
This gives an embedding of $T$ into $G$, as desired.

To keep this section from becoming too complex due to heavy notation and definitions, we are intentionally leaving out the precise definitions of the skew-matching pair and the corresponding tree-classes for now, treating these concepts as black-box structures. We will provide exact definitions when needed, in \Cref{sec:struct} and~\Cref{sec:coating}, respectively.

\subsection{Tree-classes and the Tree-Coating Lemma}\label{ssec:tree-clases}
We begin by stating our Tree-Coating Lemma.
This lemma ensures that any tree $T$ belongs to a tree-class, depending on the outcome of some process which decomposes the tree into small parts.
This decomposition depends on a parameter $\rho>0$ which we are free to choose.
As explained above, this gives a set of cut-vertices of negligible size and two forests $\mathcal F_A$ and $\mathcal F_B$ of tiny trees.
Let $a_1, a_2$ be the size of the colour classes of $\mathcal F_A$, and $b_1,b_2$ the colour classes of $\mathcal F_B$.
We then include $T$ in the tree-class $\mathcal T^\rho_{a_1,a_2,b_1,b_2}$ (see \Cref{def:tree-class}).
For technical reasons, it is not actually $T$ which is going to be decomposed but instead some slightly larger supertree $T' \supset T$, which is given by the following lemma.

\begin{restatable}[Tree-Coating Lemma]{lemma}{restatetreecoating} \label{prop:coatoftree}
	For any $1/16 \geq \eta \geq \rho>0$ and any tree $T$ of size $|V(T)|\ge 1000/\rho$, there are natural numbers $a_i, b_i\ge \eta |V(T)|$ for $i\in [2]$, and there is a tree $T'\supset T$ belonging to $\mathcal T^\rho_{a_1,a_2,b_1,b_2}$ with $|V(T')|\le (1+4\eta)|V(T)|$. 
\end{restatable}

The proof of \Cref{prop:coatoftree}, together with the definition of $\mathcal T^\rho_{a_1,a_2,b_1,b_2}$, are deferred to \Cref{sec:coating}. 

\subsection{Skew-matching pairs and the Structural Proposition}

A skew-matching~$\sigma$ is a substructure of a weighted graph, inspired by the concept of fractional matching.
We have a fixed `skew' value $\gamma > 0$ and for every edge $uv$ we attribute some non-negative weight which is distributed unequally between $u$ and $v$, in a ratio of $\gamma$.
The weight $W(\sigma)$ of a skew-matching $\sigma$ represents its size, i.e. the sum of all weights assigned to all edges.
A $(\gamma_A, \gamma_B)$-skew-matching pair $(\sigma_A, \sigma_B)$ consists of two skew-matchings $\sigma_A$ and $\sigma_B$, each with respective skew parameters $\gamma_A$ and $\gamma_B$.
Additionally, these skew-matchings will be well-positioned relative to each other.

\begin{restatable}[Structural Proposition]{proposition}{restatestructural} \label{prop:weighted-structural}
	Let $a_1,a_2,b_1,b_2\in \mathbb N$ be such that $a_1+a_2+b_1+b_2=k$ and let $(H,w)$ be a weighted graph with $w:E(H)\rightarrow (0,1]$ such that $\delta_w(H)\ge \frac{k}{2}$ and $\Delta_w(H)\ge k$.
	Let $\gamma_A:= \frac{a_2}{a_1}$ and $\gamma_B:=\frac{b_2}{b_1}$. 	
	Then~$H$ admits a $(\gamma_A,\gamma_B)$-skew-matching pair $(\sigma_A,\sigma_B)$
	with weights $W(\sigma_A)=a_1+a_2$ and $W(\sigma_B)=b_1+b_2$.   
\end{restatable}

In \Cref{sec:struct}, one can find the definition of a skew-matching pair, as well as all the notions needed to prove \Cref{prop:weighted-structural}.
Then the proof of \Cref{prop:weighted-structural} is given in \Cref{sec:new-proof-structural}.

\subsection{Reduced graph and the Embedding Lemma} \label{ssec:regularity}
As already mentioned, we use the Szemerédi Regularity Lemma~\cite{Szemeredi1978} on the host graph to obtain a regular partition.
The necessary definitions to work with graph regularity are as follows.
Given a graph $G$ and $\varepsilon > 0$, a pair $(X,Y)$ with $X, Y \subseteq V(G)$ disjoint, is said to be \emph{$\varepsilon$-regular}, if for any sets $X' \subseteq X$ and $Y' \subseteq Y$ with $|X'| \geq \varepsilon |X|$ and $|Y'| \geq \varepsilon |Y|$ it holds that $|d(G[X', Y']) - d(G[X, Y])| < \varepsilon$.
We say that a partition $\{V_0, \dotsc, V_t\}$ of $V(G)$ is an \emph{$\varepsilon$-regular partition} if $|V_0| \leq \varepsilon |V(G)|$, and for every $1 \leq i \leq t$, all but at most $\varepsilon t$ values of $1 \leq j \leq t$ are such that the pair $(V_i, V_j)$ is not $\varepsilon$-regular.
We call an $\varepsilon$-regular partition \emph{equitable} if $|V_i| = |V_j|$ for every $1 \leq i < j \leq t$.
The following version of the Regularity Lemma follows by standard arguments from its `degree form'~\cite[Theorem 1.10]{KomlosSimonovits1996}.

\begin{theorem}[Szemer\'edi's Regularity Lemma]\label{thm:SzemLemma}
	For every $\varepsilon>0$ there is $n_0$ and $M_0$ such that every graph of size at least $n_0$ admits an $\varepsilon$-regular equitable partition $\{V_0, V_1, \dotsc, V_t\}$ with $1/\varepsilon \leq t \leq M_0$.
\end{theorem}

This regular partition is then captured in a \emph{weighted reduced graph}, as defined below. 

\begin{definition}[Weighted $\dred$-reduced graph] \label{definition:reducedgraph}
	Given a graph $G$, an $\varepsilon$-regular equitable partition $\mathcal{P} = \{V_0, \dotsc, V_t\}$ of $V(G)$ and $\dred > 0$, we define the \emph{$\dred$-reduced graph} $\Gamma_{\dred, \varepsilon}$ as follows.
	The vertex set of $\Gamma_{\dred,\varepsilon}$ is $\{1, \dotsc, t\}$, and there is an edge $ij \in E(\Gamma_{\dred, \eps})$ if and only if the pair $(V_i, V_j)$ is $\varepsilon$-regular and $d(V_i, V_j) \geq \dred$.
	The \emph{weighted $\dred$-reduced graph} consists of $\Gamma_{\dred, \varepsilon}$ endowed with a natural weight function $w: E(\Gamma_{\dred, \eps}) \rightarrow [0,1]$ defined by $w(\fmat{ij}) = d(V_i, V_j)$.
	Note that $\Gamma_{\dred, \varepsilon}$ depends on the choice of the partition $\mathcal{P}$, but since it will always be clear from context which partition is being used, we omit it from the notation.	
\end{definition}

We shall apply the Structural Proposition (\Cref{prop:weighted-structural}) in the weighted $\dred$-reduced graph $\Gamma_{\dred, \varepsilon}$ to obtain a skew-matching pair there.
That is precisely the input which the next lemma requires, and it provides an embedding of a tree in the original host graph.

\begin{restatable}[Tree Embedding Lemma]{lemma}{restatetreembedding} \label{lemma:treeembedding}
	For any $\eta,  \dred, q>0$, and $t\in \mathbb N$, there are $ \varepsilon = \varepsilon(\eta,  \dred, q), \rho = \rho(\eta,  \dred, q, t)>0$ and $n_0 = n_0(\eta,  \dred, q, t) \in\mathbb N$  such that for $n\ge n_0$ the following holds.
	Suppose $G$ is an $n$-vertex graph, and $\mathcal{P}=\{V_0, V_1, \ldots, V_N\}$ is an $\varepsilon$-regular equitable partition of $G$ with $N\le t$. 
	Suppose $k\ge qn$ and that we have natural numbers $a_1, b_1, a_2, b_2 \geq \eta k$ 
	 such that
	 \[k= a_1 + a_2 + b_1 + b_2\,.\] 
	Suppose the weighted $\dred$-reduced graph $\Gamma_{\dred,\varepsilon}$ corresponding to $G$ admits a $(a_2/a_1, b_2/b_1)$-skew-matching pair $(\sigma_A, \sigma_B)$,
		with weights
	 $W(\sigma_A) \ge (1+\eta)(a_1 + a_2) N/n$ and $W(\sigma_B)\ge (1+\eta)(b_1 + b_2) N/n$.
	Then $G$ contains any tree $T\in \mathcal{T}_{a_1, a_2 ,b_1, b_2}^{\rho}$. 
\end{restatable}

\Cref{lemma:treeembedding} is proven in \Cref{sec:embed}.

\subsection{The proof of \Cref{thm:approDense-MinMax}}\label{ssec:mainproof}

Now, we give the proof of \Cref{thm:approDense-MinMax}, assuming the validity of the main statements introduced before (\Cref{prop:coatoftree}, \Cref{prop:weighted-structural}, \Cref{lemma:treeembedding}).

\begin{proof}[Proof of~\Cref{thm:approDense-MinMax}]
The proof splits naturally into four steps.

\medskip  \noindent \emph{Step 1: Setting up the parameters.}
Suppose we are given input parameters $\eta>0$ for the approximation factor, and $q>0$ for the ratio of the size of the trees with respect to the size of the host graph. We may assume that $q,\eta\ll 1$, or we just replace them with smaller values.  We set the following parameters to satisfy
\begin{equation}\label{eq:parameters}
0< 1/n_0 \ll \rho \ll 1/M_0 \ll \varepsilon\ll \dred \ll \eta, q \ll 1,
\end{equation}
as follows. We set $\dred:= \frac{\eta q }{100 }$, $\varepsilon:= \min\{\frac{\dred}{2}, \varepsilon'\}$, where $\varepsilon'$ is the output of the Tree Embedding Lemma (\Cref{lemma:treeembedding}) given the input $\eta/40$, $\dred$, and $q$ playing the roles of $\eta$, $\dred$, and $q$ respectively.
Let $M_0, n'_0$ be the outputs of the Szemer\'edi Regularity Lemma (\Cref{thm:SzemLemma}) with input $\varepsilon$. 
Let $\rho, n''_0$ be the output of the Tree Embedding Lemma (\Cref{lemma:treeembedding}) with input $\eta/40, \dred, q, M_0$ in place of $\eta, \dred, q, t$.
Finally, let $n_0 = \max\{n'_0, n''_0, 1000/\rho\}$.

\medskip \noindent \emph{Step 2: Processing the tree.}
We apply the Tree-Coating Lemma (\Cref{prop:coatoftree}) with input $\eta/12$, $\rho$, and~$T$, playing the roles of $\eta, \rho$, and $T$. We obtain $T'\in \mathcal T^\rho_{a_1,a_2, b_1, b_2}$ with $a_i, b_i\ge \frac{\eta k}{12}\ge \frac{\eta}{20}|V(T')|$ and $|V(T')|\le (1+\frac{\eta}{3})k$.
Recall that, by assumption, the host graph $G$ satisfies $\delta(G)\ge (1+\eta) \frac{k}{2}\ge (1+\frac{\eta}{20})\frac{|V(T')|}{2}$ and that at least $\frac{\eta n}{20}$ vertices of $G$ have degree at least $(1+\eta)k\ge (1+\frac{\eta}{20})|V(T')|$.

\medskip \noindent \emph{Step 3: Preparing the host graph.}
We apply the Szemer\'edi Regularity Lemma (\Cref{thm:SzemLemma}) on $G$ with parameter $\varepsilon$.
This application yields an $\varepsilon$-regular equitable partition $\{V_0, V_1, \ldots, V_t\}$, with $1/\eps \leq t\le M_0$.
Given this partition, we define the weighted ${\dred}$-reduced graph $\Gamma_{\dred, \varepsilon}$ endowed with the weight function $w:E(\Gamma_{\dred, \varepsilon})\rightarrow [0,1]$, defined by $w(\fmat{ij}):=d(V_i,V_j)$.
Recall that by the definition of $\Gamma_{\dred, \varepsilon}$,
if $(V_i, V_j)$ is not an $\varepsilon$-regular pair,
or if $d(V_i,V_j) <\dred$, then $ij \not\in E( \Gamma_{\dred, \varepsilon} )$.

From now on, let $r = k/n$.
Standard calculations show that $\Gamma_{\dred, \varepsilon}$ inherits minimum and maximum degree conditions from $G$.
Concretely, we claim that $(\Gamma_{\dred, \varepsilon}, w)$ satisfies $\delta_w(\Gamma_{\dred, \eps})\ge (1+\frac{\eta}{40})\frac{rt}{2}$ and $\Delta_w(\Gamma_{\dred, \varepsilon})\ge (1+\frac{\eta}{40})rt$.
We give those arguments for completeness.
Indeed, for  $i\in V(\Gamma_{\dred, \varepsilon})$ we have 
\begin{align*}
\deg_w(i)
	&=\sum_{j\in N_{\Gamma_{\dred, \varepsilon}}(i)}w(\fmat{ij})
	= \sum_{j\in N_{\Gamma_{\dred, \varepsilon}}(i)} d(V_i, V_j) \\
	& \ge \left( \sum_{j\in [t] \setminus \{i\}} d(V_i,V_j) \right) -(\varepsilon+\dred)t
	 = \left( \sum_{j\in [t] \setminus \{i\}} \sum_{v \in V_i} \frac{\deg_G(v, V_j)}{|V_i|^2} \right) -(\varepsilon+\dred)t\\
	& \geq \sum_{v \in V_i} \frac{\deg_G(v) - |V_i|}{|V_i|^2} -(\varepsilon+\dred)t
	\geq \frac{\delta(G)}{|V_i|} - 1 -(\varepsilon+\dred)t\\
	& \geq \frac{(1 + \eta)k/2}{|V_i|} -(2\varepsilon+\dred)t
	= \frac{(1 + \eta)rn/2}{|V_i|} - (2\varepsilon+\dred)t\\
	& \geq (1 + \eta)\frac{rt}{2} - (2\varepsilon+\dred)t
	\geq \left(1 + \frac{\eta}{40} \right)\frac{rt}{2},
\end{align*}
where we used $|V_i| \leq n/t$ in the second to last inequality,
and $r \geq q$ and the choice of $\dred, \varepsilon$ in the last inequality.

To calculate $\Delta_w(\Gamma_{\dred, \varepsilon})$, denote by $L$ the set of  $v\in V(G)$ such that $\deg_w(v)\ge (1+\eta)k = (1 + \eta) rn$. By the pigeonhole principle, there is an index $i\in[t]$ with $|V_i\cap L|\ge \frac{\eta n-\varepsilon n}{t} \geq \eta(1 - \varepsilon)n/t \geq \eta(1 - \varepsilon) |V_i| > \varepsilon |V_i|$.
Then 
\begin{align*}
\deg_w (i)
&=\sum_{j\in N_{\Gamma_{\dred, \varepsilon}}(i)}w(\fmat{ij})
= \sum_{j\in N_{\Gamma_{\dred, \varepsilon}}(i)} d(V_i, V_j)
\ge \left( \sum_{j\in [t] \setminus \{i\}} d(V_i,V_j) \right) -(\varepsilon+\dred)t \\
& \ge \left( \sum_{j\in [t] \setminus \{i\}} d(V_i \cap L,V_j) - \varepsilon \right) -(\varepsilon+\dred)t \\
& \ge \left( \sum_{v \in V_i \cap L} \frac{\deg_G(v) - |V_i|}{|V_i \cap L||V_i|} \right) -(2\varepsilon+\dred)t \\
& \ge \left( \sum_{v \in V_i \cap L} \frac{(1 + \eta) k}{|V_i \cap L||V_i|} \right) -1 -(2\varepsilon+\dred)t 
 \ge (1 + \eta)\frac{k}{|V_i|} -(3\varepsilon+\dred)t \\
&  = (1 + \eta)\frac{rn}{|V_i|} -(3\varepsilon+\dred)t
\geq (1 + \eta)rt -(3\varepsilon+\dred)t
\geq \left(1 + \frac{\eta}{40}\right)rt,
\end{align*} 
where we used $|V_i| \leq n/t$ in the second to last inequality,
and $r \geq q$ and the choice of $\dred, \varepsilon$ in the last inequality.

These estimates on $\delta_w(\Gamma_{\dred, \eps})$ and $\Delta_w(\Gamma_{\dred, \eps})$ allow us to apply the Structural Proposition (\Cref{prop:weighted-structural}) with $(\Gamma_{\dred, \varepsilon}, w)$, $(1+\frac{\eta}{40})rt$, $(1+\frac{\eta}{40})a_i t / n$ and $(1+\frac{\eta}{40})b_i t / n$ playing the roles of $(H,w)$, $k$ and $a_i, b_i$, for $i\in [2]$, respectively.
We obtain an $(\frac{a_2}{a_1},\frac{b_2}{b_1})$-skew-matching $(\sigma_A,\sigma_B)$ with $W(\sigma_A)=(1+\frac{\eta}{40})(a_1+a_2)t/n$ and $W(\sigma_B)=(1+\frac{\eta}{40})(b_1+b_2)t/n$.

\medskip \noindent \emph{Step 4: Embedding the tree.}
Having found the skew-matching pair in the reduced graph, we finalise by applying the Tree Embedding Lemma (\Cref{lemma:treeembedding}).
First, we note that $\min\{a_1,a_2\}\ge \frac{\eta}{20}|V(T')|\ge \frac{\eta}{40}\left((1+\frac{\eta}{40})|V(T')|\right)$ and similarly for $b_1,b_2$, we have that $\min\{b_1,b_2\}\ge \frac{\eta}{40}((1+\frac{\eta}{40})|V(T')|)$. 
Thus, we can apply the Tree Embedding Lemma (\Cref{lemma:treeembedding}) with
$t$ and  $\eta/40$,
playing the roles of
$N$ and $\eta$, respectively.
We obtain that $G$ contains any tree in $\mathcal{T}^\rho_{a_1, b_1, a_2, b_2}$,
so in particular $T' \subseteq G$.
Since $T \subseteq T'$, this proves \Cref{thm:approDense-MinMax}.
\end{proof}

\section{Coating the Tree}\label{sec:coating}

To prepare the embedding, we will first partition the trees into suitably small parts.
We will use the following handy concept from Hladk\'y, Komlós, Piguet, Simonovits, Stein, and Szemerédi~\cite[Definition 3.3]{HKPSSS2017d}.
It gives a partition of a tree into smaller trees that also satisfy several additional useful properties.

If $T$ is a tree rooted at $r$, and $\widetilde{T} \subseteq T$ is a subtree with $r \notin V(\widetilde{T})$, the \emph{seed of $\widetilde{T}$} is the unique vertex $x \in V(T) \setminus V(\widetilde{T})$ which is farthest from $r$ and also belongs to every $(r, v)$-path in $T$, for every $v \in V(\widetilde{T})$.
We emphasize that the seed of $\widetilde{T}$ is \emph{not} contained in $V(\widetilde{T})$, and that $\widetilde{T}$ does not necessarily contain all descendants of its seed.

\begin{definition}[$\ell$-fine partition]\label{def:ell-fine}
	
	Let $T$ be a tree on $k$ vertices rooted at a vertex $r$.
	An \emph{$\ell$-fine partition of $T$} is a quadruple $(W_A, W_B, \mathcal{F}_A, \mathcal{F}_B)$, where $W_A, W_B \subseteq V(T)$ and $\mathcal{F}_A, \mathcal{F}_B$ are families of subtrees of $T$ such that
	
	\stepcounter{propcounter}
	\begin{enumerate}[(\Alph{propcounter}\arabic*),topsep=0.7em, itemsep=0.5em]
		\item \label{item:vertex partition} the three sets $W_A$, $W_B$, and $\{ V(T^\ast) \}_{T^\ast \in \mathcal{F}_A \cup \mathcal{F}_B}$ partition $V(T)$ (in particular, the trees in $\mathcal{F}_A \cup \mathcal{F}_B$ are pairwise vertex-disjoint),
		\item $r \in W_A \cup W_B$,
		\item \label{item:fp-sizeWAWB} $\max \{ |W_A|, |W_B| \} \leq 336 k/\ell$,
		\item for $w_1, w_2 \in W_A \cup W_B$, the distance between $w_1$ and $w_2$ in $T$ is odd if and only if one of them lies in $W_A$ and the other one in $W_B$,
		\item \label{item:fp-smallshrubs} $|V(T^\ast)| \leq \ell$ for every $T^\ast \in \mathcal{F}_A \cup \mathcal{F}_B$,
		\item \label{item:fp-locationroots} $V(T^\ast) \cap N(W_B) = \emptyset$ for every $T^\ast \in \mathcal{F}_A$, and 
		$V(T^\ast) \cap N(W_A) = \emptyset$ for every $T^\ast \in \mathcal{F}_B$;
		\item \label{item:fp-seedsinWAWB} each tree of $\mathcal{F}_A \cup \mathcal{F}_B$ has its seeds in $W_A \cup W_B$,
		\item \label{item:fp-twoseeds} $|N(V(T^\ast)) \cap (W_A \cup W_B)| \leq 2$ for each $T^\ast \in \mathcal{F}_A \cup \mathcal{F}_B$,
		\item \label{item:fp-distance} if $ N(V(T^\ast)) \cap (W_A \cup W_B)$ contains two distinct vertices $w_1, w_2$ for some $T^\ast \in \mathcal{F}_A \cup \mathcal{F}_B$, then $\dist_T(w_1, w_2) \geq 6$.
	\end{enumerate}
\end{definition}

\begin{figure}[h]
	\centering
	\includegraphics[width=0.5\linewidth]{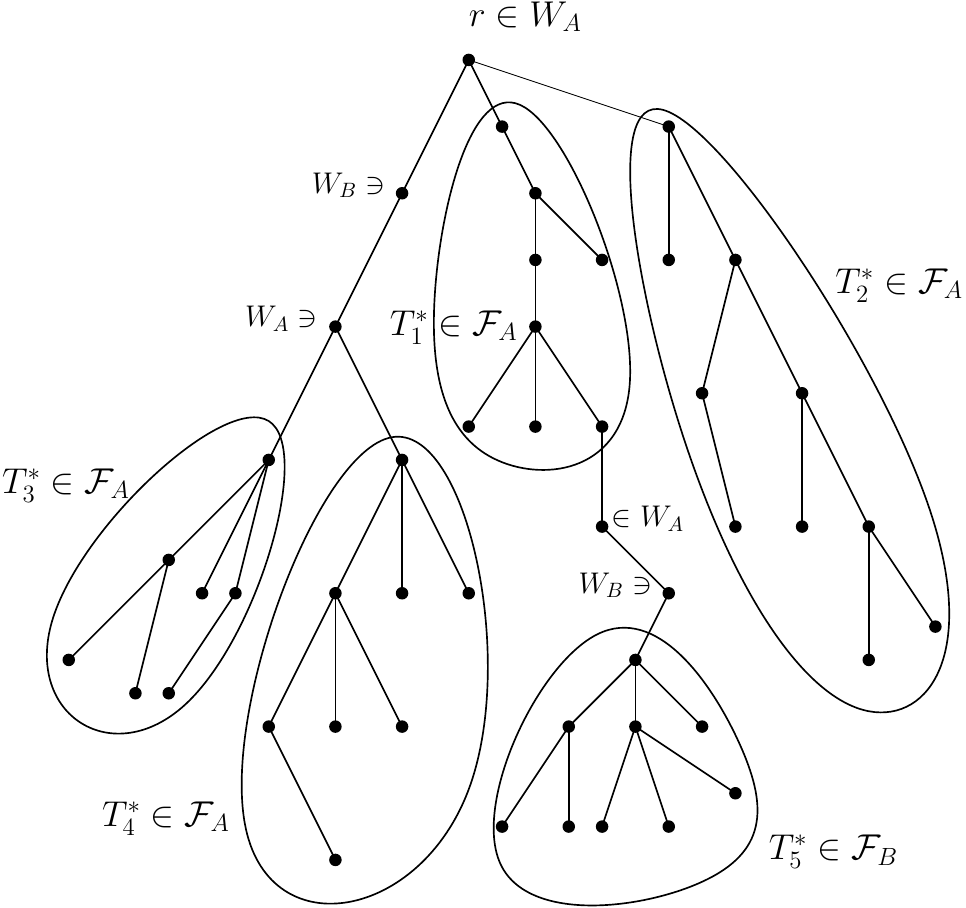}
	\caption{A schematic view of an $ \ell$-fine partition of a tree with $ 45 $ vertices rooted at $ r $ into $W_A, W_B, \mathcal{F}_A, \mathcal{F}_B$, satisfying \ref{item:vertex partition}--\ref{item:fp-distance}.}
	\label{lfine}
\end{figure}

The trees $T^\ast \in \mathcal{F}_A \cup \mathcal{F}_B$ will be called \emph{shrubs}.
We remark that $\ell$-fine partitions can be defined differently so they satisfy even more properties, but we are only citing those we need.
See \Cref{lfine} for a visual representation of a tree together with an $\ell$-fine partition.

The next lemma~\cite[Lemma 3.5]{HKPSSS2017d} says that any tree has an $\ell$-fine partition.

\begin{lemma}\label{lem:ell-fine}
	Let $T$ be a tree on $k$ vertices rooted at $r$, and let $1 \leq \ell \leq k$.
	Then $T$ has an $\ell$-fine partition.
\end{lemma}

\begin{definition}[Tree-class $\mathcal T^\rho_{a_1,a_2,b_1,b_2}$]\label{def:tree-class}
We denote by $\mathcal T^{\rho}_{a_1,a_2,b_1,b_2}$ the set of trees $T$, such that there is a $(\rho|V(T)|)$-fine partition $(W_A, W_B, \mathcal{F}_A,  \mathcal{F}_B)$ of $T$ such that $|V_i(\mathcal{F}_A)|=a_i $, $|V_i(\mathcal{F}_B)|=b_i$, for $i \in \{ 1,2 \}$,
 where $V_1(\mathcal {F}_A)$  (resp. $V_2(\mathcal F_A)$) is the set of vertices of $\mathcal F_A$ that are at odd (resp. even) distance from $W_A$,  and $V_i(\mathcal F_B)$ are defined analogously with respect to $W_B$.
 \end{definition}

As mentioned earlier, for technical reasons we shall work with a tree $T'\supset T$ which is slightly larger than $T$, but that belongs to the class $\mathcal T^{\rho}_{a_1, a_2, b_1, b_2}$ for some suitable values of $a_1, a_2, b_1, b_2$.
In our case `suitable' will mean that none of these values are very small; this fact will be useful during the embedding phase later.
The Tree-Coating Lemma (\Cref{prop:coatoftree}) describes precisely the process of finding this supertree $T'$.
We now have the tools and vocabulary to prove it.
We restate the statement for the convenience of the reader.

\restatetreecoating*

\begin{proof}[Proof of \Cref{prop:coatoftree}]
	Pick any vertex $r\in V(T)$ as the root of~$T$ and apply \Cref{lem:ell-fine} with input~$T$,~$|V(T)|$,~$r$,  and $\rho |V(T)|$, playing the roles of~$T$,~$k$,~$r$ and $\ell$ to obtain a 
	$\rho|V(T)|$-fine partition of~$T$. First assume that both $W_A$ and $W_B$ are non-empty. Pick a vertex  $x_A\in W_A$ and a vertex $x_B\in W_B$ and attach ${\eta}|V(T)|$ paths of length~$2$ to each of the two vertices~$x_A$ and~$x_B$. This adds ${\eta}|V(T)|$ vertices to $\mathcal F_\AB\cap V_i(T)$
	for $\AB\in \{A,B\}$ and $i\in[2]$, where $V_i(T), i\in [2]$ are the two colour classes of $T$,  and increases the size of the tree by $4\eta |V(T)|$. As $2\le \rho |V(T)|$, all the newly formed shrubs in $\mathcal F_A\cup \mathcal F_B$ have size at most $ \rho |V(T')|$ and thus, the newly obtained tree $T'$ belongs to the class $\mathcal{T}^\rho_{a_1, a_2,b_1,b_2}$, for suitable $a_i$ and $b_i$, for $ i\in [2]$. 
	
	If (w.l.o.g.) $W_A=\emptyset$, by \Cref{def:ell-fine}\ref{item:fp-locationroots}, we have that $\mathcal F_A=\emptyset$.
	By \ref{item:fp-sizeWAWB}, we have $|W_B| \leq 336/\rho \leq |V(T)|/2$, and therefore
	 $|\mathcal F_B\cap V_i(T)|\ge |V(T)|/4 \ge 4\eta |V(T)|$, for some $i\in[2]$.
	 Add a vertex $x_A$ to $W_A$, connecting it to any vertex in $W_B$, and as previously attach ${\eta}|V(T)|$ paths of length~$2$ to $x_A$. For the other side, if the vertices in $\mathcal F_B$ at odd distance to $W_B$ are less than ${\eta}{|V(T)|}$, then 
	place a star with $2{\eta}|V(T)|-1$ leaves centred on an arbitrary vertex $x_B\in W_B$. Otherwise attach  ${\eta}|V(T)|-2$ paths of length~$2$ to some vertex $x_B\in W_B$ and a star with $2$ leaves and connect the centre to $x_B$ as well. 
	In this way, we ensure (as $3 \leq \rho |V(T)|$) that the newly created tree belongs to $\mathcal T^\rho_{a_1,a_2,b_1,b_2}$ for suitable $a_i, b_i\ge {\eta |V(T)|}$, $i\in[2]$, and that we have increased the tree by exactly $4\eta |V(T)|$ vertices.
\end{proof}

\section{Skew-matchings and other matching structures}\label{sec:struct}

In the next four sections, we prove the Structural Proposition (\Cref{prop:weighted-structural}).
In this section, we will introduce the basic definitions that will allow us to describe the structure we will seek in the reduced graph.
The basic building block of our argument is a ``skew'' oriented fractional matching, defined in Section \ref{ssec:skew}.
Its definition is inspired by the standard fractional matching, which (we recall) is a weight function over the edges of a graph such that the sum of weights over the edges incident to any given vertex is at most one (\Cref{def:fractional_matching}).
As we shall see, a skew-matching satisfies a similar property, but instead is a weight function defined over \emph{oriented} edges of the graph, which allows us to distribute the weight in every edge in an unbalanced way.
We remark that essentially equivalent concepts were used before (e.g.~\cite[\S 2.4]{Balogh-McDowell-Molla-Mycroft-2018}) in the context of graph tilings.

After introducing the main definitions in \Cref{ssection:fractionalmatching} and \Cref{ssec:skew}, we introduce many extra auxiliary concepts.
In \Cref{subsection:comparingmatchings} we introduce notation to compare different fractional matchings and skew-matchings; and in \Cref{ssection:truncated} we introduce the auxiliary notion of the \emph{truncated weighted graph}.

\subsection{Fractional matchings} \label{ssection:fractionalmatching}
The following standard definition generalises the notion of matchings to a fractional setting.
\begin{definition}[Fractional matching]\label{def:fractional_matching}
	A \emph{fractional matching} is a weight function $\mu:E(G)\rightarrow [0,1]$ such that for any vertex $v\in V(G)$, we have $\sum_{u\in N(v)}\mu(\fmat{uv})\le 1$.
\end{definition}

We shall abuse notation here and use the sign $\mu$ on vertices as well and we call \[\mu(v):=\sum_{u\in N(v)}\mu(\fmat{uv})\] the \emph{weight of  $\mu$ on $v$}. 
The \emph{weight} of a fractional matching $\mu$ is 
\[W(\mu):= \sum_{e\in E(G)}\mu(e).\]

We say that two fractional matchings $\mu'$ and $\mu''$ are \emph{disjoint} if $\mu'(x)+\mu''(x)\le 1$ for all $x\in V(G)$. We denote by $V(\mu)$ the set of vertices $v\in V(G)$ such that $\mu(v)>0$.  	

\subsection{Skew-matchings}\label{ssec:skew}

Now we present the definition of skew-matchings.
Note that, according to our convention on $G^{\leftrightarrow}$, both edges $\oriented{uv}$ and $\oriented{vu}$ are present in the following definition.

\begin{definition}[$\gamma$-skew oriented fractional matching]\label{def:sigma}
	Let $G$ be a  graph, $G^{\leftrightarrow}$ its associated digraph, and $\gamma \geq 0$.
	A \emph{$\gamma$-skew oriented fractional matching} (or just \emph{$\gamma$-skew-matching}) is a function $\sigma: E(G^{\leftrightarrow}) \rightarrow [0,1]$ such that for any vertex $u \in V(G)$,
	\[\sum_{v \in N(u)} \left( \frac{1}{1+\gamma}\sigma(\oriented{uv}) + \frac{\gamma}{1+\gamma} \sigma(\oriented{vu}) \right)\leq 1.  \]
\end{definition}
As explained before, $\gamma$-skew-matchings can be understood as fractional matchings (\Cref{def:fractional_matching}) in graphs where the weight of the edge is distributed in an unbalanced way, meaning that one end of the edge gets $\gamma$ times the weight of the other end, and the direction of this imbalance is given by the direction of the edge in the digraph $G^\leftrightarrow$.

With a slight abuse of notation, given a $\gamma$-skew-matching $\sigma$ we shall use the symbol $\sigma$ on vertices as well.
We define $$\sigma(u):= \sigma^1(u)+\sigma^2(u),$$
where 
\[\sigma^1(u):=\frac{1}{1+\gamma}\sum_{v \in N(u)} \sigma(\oriented{uv}) \qquad \text{and} \qquad \sigma^2(u):=\frac{\gamma}{1+\gamma}\sum_{v \in N(u)} \sigma(\oriented{vu}). \]

Now we define two crucial concepts associated with a $\gamma$-skew-matching $\sigma$.
The \emph{anchor} of $\sigma$, denoted by $\mathcal{A}(\sigma)$, is the set \[\mathcal A(\sigma):=\{u\in V(G)\::\:\sigma^1(u)>0\}.\]
The \emph{weight} of $\sigma$ is \[W(\sigma) := \sum_{e \in E(G^{\leftrightarrow})} \sigma(e).\]

For our proof, we will usually need to work with several skew-matchings, possibly with different values of $\gamma$.
If $\sigma, \sigma'$ are $\gamma$-skew and $\gamma'$-skew-matchings respectively,
we say $\sigma, \sigma'$ are \emph{disjoint} if, for every $u \in V(G)$,
\[ \sigma(u) + \sigma'(u) \leq 1. \]

In what follows, we often work with a weighted graph $(G,w)$ and its associated digraph $G^{\leftrightarrow}$  and we do not always explicitly distinguish between $G$ and $G^{\leftrightarrow}$, as they are in one-to-one correspondence.
The following definition is key: it essentially says that a skew-matching $\sigma$ does not place too much weight in the neighbourhood of a vertex $u$, meaning that it respects the weight dictated by some edge weighting $w$.
This will allow us to embed trees using the space, allowed by $\sigma$, that is also joined correctly with $u$.

\begin{definition}
	[Fitting in the $w$-neighbourhood]
	\label{def:weighted-anchored}
	Let $(G,w)$ be a weighted digraph, 
	$\gamma \geq 0$, and $u\in  V(G)$.
	Given a $\gamma$-skew oriented fractional matching $\sigma$, we say that its anchor \emph{fits in the $w$-neighbourhood of $u$}
	if $\sigma^1(v)\le w(\oriented{uv})$, for every $v\in V(G)$.
\end{definition}

Note that, in particular, if $\mathcal{A}(\sigma)$ fits in the $w$-neighbourhood of $u$, then $\mathcal{A}(\sigma) \subseteq N_w(u)$.

The following definition is also key.
It corresponds to two skew-matchings, with different skews, which are disjoint from each other and they fit in the neighborhoods of two vertices forming an edge $cd$.
This will be the structure we desire in the reduced graph.
Referencing the $\ell$-fine partitions of a tree $T$ discussed before: the (few) vertices of the sets $W_A, W_B$ will be mapped to the clusters corresponding to $cd$, while the shrubs in $\mathcal{F}_A, \mathcal{F}_B$ will be mapped in the space ensured by $\sigma_A, \sigma_B$, respectively.

\begin{definition}[Edge-anchored skew-matching pair]\label{def:weighted-anchored-pair}
	Let $(G,w)$ be a weighted digraph, 
	and $\gamma_A, \gamma_B\ge 0$.
	A \emph{$(\gamma_A, \gamma_B)$-skew-matching pair anchored in $\oriented{cd}\in E(G)$} is a pair $(\sigma_A,\sigma_B)$ such that
	\stepcounter{propcounter}
	\begin{enumerate}[(\Alph{propcounter}\arabic*),topsep=0.7em, itemsep=0.5em]
		\item \label{itdef:disjoint}$\sigma_A$ and $\sigma_B$ are disjoint,
		\item \label{itdef:anchorc}$\sigma_A$ is a $\gamma_A$-skew-matching, whose  anchor fits  in the $w$-neighbourhood of $c$, 
		\item \label{itdef:anchord} $\sigma_B$ is a $\gamma_B$-skew-matching, whose  anchor fits in the $w$-neighbourhood of  $d$, and
		\item \label{itdef:anchorpartition} for every $x \in N(c) \cap N(d)$, we have
		\[ \max\{ w(\oriented{cx}), w(\oriented{dx}) \} \geq \sigma_A^1(x)+\sigma_B^1(x). \]
	\end{enumerate}
	Observe that there is no weight requirement on $\oriented {cd}$.
	
	We say that a weighted (unoriented) graph $(G,w)$ \emph{admits} a $(\gamma_A,\gamma_B)$-skew-matching pair $(\sigma_A, \sigma_B)$, if $(\sigma_A, \sigma_B)$ is defined on its associated weighted digraph  $G^\leftrightarrow$ and is anchored in some edge $\oriented{cd}\in E(G^\leftrightarrow)$, and $w(\fmat{cd})>0$.  
\end{definition}

\begin{remark}
	We can express \ref{itdef:anchorpartition} equivalently by saying that there is a partition $\{X_c, X_d\}$ of $N(c)\cap N(d)$ such that for every $x\in X_c$, we have $w(\oriented{cx})\ge \sigma_A^1(x)+\sigma_B^1(x)$, and for every $x\in X_d$, we have $w(\oriented{dx})\ge \sigma_A^1(x)+\sigma_B^1(x)$.
\end{remark}

We say that a vertex $v\in V(G)$ is \emph{covered} by a fractional matching $\sigma$ (or a skew oriented fractional matching) if $\sigma(v)=1$. 
Analogously, we say that a set $U$ is \emph{covered}, if each vertex $v\in U$ is covered.
Moreover, we define $V(\sigma):=\{v\in V(G) \ : \ \sigma(v)>0\}$.

The value 
 \[\deg_w(u, \sigma)=\sum_{v\in V(G)}\min\{\sigma(v), w(\oriented{uv})\}\]
is called the \emph{saturation} of $N(u)$ by $\sigma$. We say that $\sigma$ \emph{saturates} $N(u)\cap U$ for some $U\subseteq V(G)$, if for every vertex $v\in U$, we have that $\sigma(v)\ge w(\oriented{uv})$. If $U=V(G)$, we simply say that $\sigma$ \emph{saturates} the neighbourhood of $u$ 		(or \emph{saturates $N(u)$} for short).
If this happens, then  $\deg_w(u)=\deg_w(u, \sigma)$.

\subsection{Comparing matchings} \label{subsection:comparingmatchings}
Now we define a few relations between fractional matchings and skew-matchings, which will allow us to compare them to each other.
Notice that both edges $\oriented{xy}$ and $\oriented{yx}$ are often used.

We begin by comparing fractional matchings between them.
We write $\mu\le \mu'$ for two fractional matchings $\mu$ and $\mu'$,  whenever $\mu(\fmat{xy})\le \mu'(\fmat{xy})$ for every $xy\in E(G)$.
Moreover, $\mu = \mu'$ if and only if $\mu(\fmat{xy}) = \mu'(\fmat{xy})$ for every $xy \in E(G)$.

Similarly, we can compare skew-matchings between them, even if they differ in their skews.
A $\gamma'$-skew-matching $\sigma'$ is a  \emph{skew sub-matching} of a $\gamma$-skew  matching $\sigma$, if for every $\oriented{uv}\in E(G^{\leftrightarrow})$, we have $\frac{\sigma'(\oriented{uv})}{1+\gamma'}\le \frac{\sigma(\oriented{uv})}{1+\gamma}$ and $\frac{\gamma'}{1+\gamma'}\sigma'(\oriented{uv})\le \frac{\gamma}{1+\gamma}\sigma(\oriented{uv})$. This is denoted by $\sigma'\le \sigma$. 

In the next definition we introduce the language and notation to compare fractional matchings with skew-matchings, and vice versa.

\begin{definition}[$\preceq, \trianglelefteq$]\label{def:<=}
	Let $G$ be a graph with associated digraph $G^{\leftrightarrow}$, let $\sigma$ be a $\gamma$-skew oriented fractional matching in $G^{\leftrightarrow}$, and let $\mu$ be a fractional matching in $G$. We write $\mu\preceq \sigma$ if for every $\oriented{xy} \in E(G^{\leftrightarrow})$ with $\mu(\fmat{xy})>0$, we have 
\[
\mu(\fmat{xy})\le  \frac{1}{1+\gamma}\sigma(\oriented{xy})+\frac{\gamma}{1+\gamma}\sigma(\oriented{yx}).
\] 

Similarly, we write $\sigma\trianglelefteq \mu$ 
if for every $\oriented{xy} \in E(G^{\leftrightarrow})$ with $\sigma(\oriented{xy})+\sigma(\oriented{yx})>0$, we have 
\[
\frac{\sigma(\oriented{xy})+\gamma\sigma(\oriented{yx})}{1+\gamma}\le \mu(\fmat{xy}).
\]
\end{definition}

\begin{remark}
	We emphasize that in both definitions, that of $\mu\preceq \sigma$ and $\sigma\trianglelefteq \mu$, we require the inequality to hold for each oriented edge in $G^{\leftrightarrow}$, which means that for each (unoriented) edge $xy \in E(G)$ we need to verify the inequality both for $\oriented{xy}$ and  $\oriented{yx}$.
	
	For instance, in the definition of $\mu \preceq \sigma$ we need \[\mu(\fmat{xy})\le  \frac{1}{1+\gamma}\sigma(\oriented{yx})+\frac{\gamma}{1+\gamma}\sigma(\oriented{xy})\] to hold as well; and in the definition of $\sigma \trianglelefteq \mu$ we also need that \[\frac{\sigma(\oriented{yx})+\gamma\sigma(\oriented{xy})}{1+\gamma}\le \mu(\fmat{xy}).\]
\end{remark}

We also introduce notation to compare multiple skew-matchings with a fractional matching.
Let $\sigma_i$ be a $\gamma_i$-skew-matching for every $i\in [k]$. 
We write $\sum_{i=1}^{k}\sigma_i\trianglelefteq \mu$ if for every $\oriented{xy}\in E(G^{\leftrightarrow})$ with $\sum_{i=1}^{k}\big(\sigma_i(\oriented{xy})+\sigma_i(\oriented{yx})\big)>0$, we have
\[\sum_{i=1}^{k}\frac{\sigma_i(\oriented{xy})+\gamma_i\sigma_i(\oriented{yx})}{1+\gamma_i}\le \mu(\fmat{xy}).\]

\begin{remark}\label{rem:adding-skew<=mu}
In the last definition, we compute a sum of skew-matchings with different skews. This sum does not define a new skew-matching. We only compare the weight of these skew-matchings on every oriented edge with the weight of the fractional matching $\mu$ on the corresponding unoriented edge (and we do it for both possible orientations). In this way we investigate if all these skew-matchings together ``fit'' in a fractional matching~$\mu$. In other words, if $\sum_{i=1}^{k}\sigma_i\trianglelefteq \mu$, we have $\sum_{i=1}^k\sigma_i(u)\le \mu(u)\le 1$ for all $u\in V(G)$. 
\end{remark}

\begin{observation}
	Let $G$ be a graph, $G^{\leftrightarrow}$ be its associated digraph, $\mu,\widehat{\mu}$ be fractional matchings in $G$ and $\sigma$, $\widehat{\sigma}$ be skew fractional matchings in $G^{\leftrightarrow}$.

All the values $\mu(\fmat{uv}), \sigma(\oriented{uv}), \mu(u),\widehat{\mu}(u)$ are non-negative real numbers.
Therefore, expressions such as $\mu(\fmat{uv}) \leq \sigma(\oriented{uv})$, $\mu(u)\leq \widehat{\mu}(u)$, $\sigma(\oriented{uv}) + \mu(\fmat{uv})$ are well-defined. 
	
We can also gather some straightforward observations and consequences of our previous definitions, that will allow us to work more mechanically with different objects.	
\begin{enumerate}
\item Having $\mu(u) \leq \widehat{\mu}(u)$ for every $u\in V(G)$ does not imply that $\mu\leq \widehat{\mu}$.
\item Similarly, having $\sigma  \trianglelefteq \mu$ does not imply that $\sigma(\oriented{uv})  \leq \mu(\fmat{uv})$ for every $uv$.
\item On the other hand, having $\sigma \trianglelefteq \mu$ implies that $\sigma(u) \leq \mu(u)$ for every $u\in V(G)$.
\item If $\mu \preceq \sigma$ and $\sigma \trianglelefteq \widehat{\mu}$, then $\mu \leq \widehat{\mu}$.
\item Having $\sigma \trianglelefteq \mu$ and $\mu \preceq \widehat{\sigma}$ does not imply that $\sigma\leq \widehat{\sigma}$.
\item If $
	\sigma \trianglelefteq \mu$ and $W(\sigma) = 2 W(\mu)$, then $\sigma(u) = \mu(u)$ for every $u \in V(G)$.
\end{enumerate}	 
Note that, in general, we have $\sigma(u) \neq \sum_{v \in N(u)}\sigma(\oriented{uv}) + \sigma(\oriented{vu})$.
\end{observation}

While fractional matchings are defined on  a non-oriented  graph ~$G$ and skew-matchings are defined on $G^\leftrightarrow$, we can easily obtain a $1$-skew-matching $\sigma$ from a fractional matching~$\mu$.
This can be done, for instance, by choosing an orientation for each edge $xy \in E(G)$ and setting $\sigma(\oriented{xy})=2\mu(\fmat{xy})$ and $\sigma(\oriented{yx})=0$.
We observe, however, that the definition of ``weight'' in both cases differs, because then \[W(\sigma)=\sum_{\oriented{xy}\in E(G^\leftrightarrow)}\sigma(\oriented{xy})=\sum_{\fmat{xy}\in E(G)}\left(2\mu(\fmat{xy})+0\right)=2W(\mu).\] 

There are many different $1$-skew-matchings $\sigma$ that we can define using a fractional matching $\mu$. For example  $\sigma(\oriented{xy})=\sigma(\oriented{yx})=\mu(\fmat{xy})$ is also a possibility. In this case, we again obtain that $W(\sigma) = 2W(\mu)$.
Conversely, we can obtain a fractional matching $\mu$ from a $1$-skew-matching $\sigma$ just by ``forgetting'' the orientation of the edges.
Then, we get
\[W(\mu)=\sum_{\fmat{xy}\in E(G)}\mu(\fmat{xy})=\sum_{\fmat{xy}\in E(G)}\frac{\sigma(\oriented{xy})+\sigma(\oriented{yx})}{2}=\frac{W(\sigma)}{2}.\]
The next lemma encapsulates the above discussion for future reference.

\begin{lemma} \label{lemma:fractionalfrom1skew}
	Let $G$ be a graph and $G^{\leftrightarrow}$ its associated digraph.
	Let $\sigma$ be a $1$-skew-matching in $G^{\leftrightarrow}$.
	Then $G$ has a fractional matching $\mu$ such that
	\begin{enumerate}
		\item $W(\mu) = \frac{1}{2} W(\sigma)$,
		\item for all $x \in V(G)$, $\sigma(x) = \mu(x)$, and
		
		\item if $\sigma'$ is a $\gamma$-skew-matching in $G^{\leftrightarrow}$ such that $\sigma \leq \sigma'$, then $\mu \preceq \sigma'$.
	\end{enumerate} 
\end{lemma}

In the following definition, we extend the concept of \emph{covering} and \emph{saturation} from fractional matchings or skew-matchings, to sum of skew-matchings, or sum of a fractional matching and skew-matchings.
This is done in a natural way: the notions of covering and saturation refer to the weights that the given objects induce on the vertices and do not depend on the type of object.

\begin{definition}[Saturation, Covering]\label{def:ext-saturation}
	Let $(G,w)$ be a weighted graph, $(G^\leftrightarrow, w)$ its associated weighted digraph,  $\nu_1, \nu_2$ represent fractional matchings and/or skew-matchings\footnote{The values of the skew are not relevant in this definition.} 
	(including the option of everywhere $0$-valued functions), and $u\in V(G)$. Then the value 
	\[ \deg_w(u,\nu_1+\nu_2):=\sum_{v\in N_w(u)}\min\{\nu_1(v)+\nu_2(v),w(\oriented{uv})\}\]
	is  called the \emph{saturation} of $N_w(u)$ by $\nu_1+\nu_2$.  		
		 We say that $\nu_1+\nu_2$ \emph{saturates the neighbourhood of $u$} (or \emph{saturates} $N_w(u)$ for short) if $\deg_w(u)$ equals $\deg_w(u, \nu_1+\nu_2)$. 
		 
	Analogously, $\nu_1+\nu_2$ \emph{saturates} $N_w(u)\cap U$, for some $U\subseteq V(G)$, if for every vertex $v\in N_w(u)\cap U$ we have $\nu_1(v)+\nu_2(v)\ge w(\oriented{uv})$. 
	We say that $\nu_1+\nu_2$ \emph{covers} a vertex $v\in V(G)$ if $\nu_1(v)+\nu_2(v)=1$. We say that it covers a set $S \subseteq V(G)$, if it covers every vertex $v\in S$. 
\end{definition}
 
%
%
%
%

\subsection{Truncated weighted graphs} \label{ssection:truncated}
Now we shall define another important concept used in our proofs, that of a \emph{truncated weighted graph}.
To motivate this definition, consider the following scenario.
Suppose we need to find a fractional matching in a weighted graph $(G,w)$ of sufficiently large weight, but so far we have only found a fractional matching~$\mu$ of insufficient weight.
In our settings, we always need to make sure our fractional objects (matchings or skew-matchings) `fit' in the $w$-neighbourhood of a given vertex, say $c$.
Thus, if we want to build another matching $\bar \mu$, disjoint from $\mu$, and we want to consider the sum $\mu + \bar \mu$, we need to keep in mind that we can only allocate weight in a vertex as determined by the weight function $w$.
More precisely, for every vertex $x \in N(c)$ we need to ensure $\bar \mu(x) + \mu(x) \leq w(\oriented{cx})$.
To work with these kinds of restrictions, we will define a new weight function $\bar w$, which is essentially $w(\oriented{cx}) - \mu(x)$ for every $x \in V(G)$, so the restriction we just discussed is expressed more succinctly as $\bar \mu(x) \leq \bar w(\oriented{cx})$.
\begin{definition}[Truncated weighted digraph]\label{def:truncated-weighted-graph}
Let $(G,w)$ be a weighted symmetric digraph and $\mu$ be a fractional matching in $G$ (forgetting the orientation of the edges).
For every directed $ux \in E(G)$, let $\bar{w}(\oriented{ux}):= \max\{0, w(\oriented{ux})-\mu(x)\}$.
We call $(G, \bar w)$  the \emph{$\mu$-truncated weighted digraph obtained from $(G, w)$}.
\end{definition}

\begin{remark}
	In the previous definition we need $G$ to be a symmetric digraph, to be able to define a fractional matching in a meaningful way; in fact we will use it essentially only for associated digraphs obtained from undirected graphs.
	But we do not require the weights themselves in $(G,w)$ to be symmetric, and generally the truncated weighted digraphs $(G, \bar w)$ will also not have symmetric weights.
\end{remark}

The next proposition ensures that we can indeed correctly combine skew-matchings as long as they respect the weights from a truncated weighted digraph; so we achieved what we set out to do with the definition.

\begin{proposition}\label{prop:adding-skew-matchings}
Let $(G,w)$ be a weighted symmetric digraph, and $\fmat{uv}\in E(G)$. 
Let $\mu$ and $\bar\mu$ be disjoint fractional matchings in (unoriented) $G$ and let $(G, \bar w)$ be the $\mu$-truncated weighted digraph obtained from $(G,w)$. 
Suppose that
\stepcounter{propcounter}
\begin{enumerate}[\upshape{(\Alph{propcounter}\arabic*)},topsep=0.7em, itemsep=0.5em]
	\item \label{item:addingskew-in-1} $(\sigma_A, \sigma_B)$ is a $(\gamma_A, \gamma_B)$-skew-matching pair in $(G, w)$ anchored in $\oriented{uv}$ with $\sigma_A + \sigma_B \trianglelefteq \mu$; and
	\item \label{item:addingskew-in-2} $(\bar \sigma_A, \bar \sigma_B)$ is a $(\gamma_A, \gamma_B)$-skew-matching pair in $(G, \bar w)$ anchored in $\oriented{uv}$ with $\bar \sigma_A + \bar \sigma_B \trianglelefteq \bar \mu$.
\end{enumerate}
Then $(\sigma_A + \bar \sigma_A, \sigma_B + \bar \sigma_B)$ is a $(\gamma_A, \gamma_B)$-skew-matching pair in $(G, w)$, anchored in $\oriented{uv}$, with $\sigma_A + \bar \sigma_A + \sigma_B + \bar \sigma_B \trianglelefteq \mu + \bar \mu$.
\end{proposition}

\begin{proof}
	We need to verify the properties \ref{itdef:disjoint}--\ref{itdef:anchorpartition} for $(\sigma_A + \bar \sigma_A, \sigma_B + \bar \sigma_B)$.
	The disjointness of $\mu$ and $\bar\mu$, along with $\sigma_A+\sigma_B\trianglelefteq \mu$ and $\bar \sigma_A+\bar \sigma_B\trianglelefteq \bar \mu$, ensure that $\sigma_A+\bar\sigma_A+\sigma_B+\bar\sigma_B\trianglelefteq \mu+\bar\mu$.
	Thus, by \Cref{rem:adding-skew<=mu}, we have that $\sigma_A+\bar\sigma_A$ is disjoint from $\sigma_B+\bar\sigma_B$, which gives  \ref{itdef:disjoint}.
	This also gives that $\sigma_A+\bar\sigma_A$ is a $\gamma_A$-skew-matching and $\sigma_B+\bar\sigma_B$ is a $\gamma_B$-skew-matching. 

	Now we verify \ref{itdef:anchorc}--\ref{itdef:anchord}.
	For all $z\in N_{\bar w}(u)$, we have $\bar w(\oriented{uz}) > 0$, and therefore $\bar w (\oriented{uz}) = w(\oriented{uz}) - \mu(z)$ by definition.
	Then $\sigma_A^1(z)+\bar\sigma_A^1(z)\le \mu(z)+\bar w(\oriented{uz}) = w(\oriented{uz})$, where in the first inequality we used that $\sigma_A \trianglelefteq \mu$ and that $\bar \sigma_A$ fits in the $\bar w$-neighbourhood of $u$. 
	On the other hand, if $\bar w(\oriented{uz})=0$, then $\bar\sigma_A^1(z)=0$.
	Thus, $\sigma_A^1(z)+\bar\sigma_A^1(z)=\sigma_A^1(z)\le w(\oriented{uz})$.
	Hence, we have that $\mathcal A(\sigma_A+\bar\sigma_A)$ fits in the $w$-neighbourhood of $u$, which gives  \ref{itdef:anchorc}.
	A symmetric argument gives that $\mathcal{A}(\sigma_B + \bar \sigma_B)$ fits in the $w$-neighbourhood of $v$, thus giving \ref{itdef:anchord}.

	Finally, consider any $z\in N_w(u)\cap N_w(v)$.
	Using that $\mathcal{A}(\bar \sigma_A), \mathcal{A}(\bar \sigma_B)$ fit in the $\bar w$-neighbourhood of $u$ and $v$ respectively, we get 
	\begin{align}
	\bar\sigma_A^1(z)+\sigma_A^1(z)+\bar\sigma_B^1(z)+\sigma_B^1(z) \le \max\{\bar w(\oriented{uz}), \bar w(\oriented{vz})\}+\sigma^1_A(z)+\sigma^1_B(z).
	\label{equation:addingskewspartitioncheck}
	\end{align} 
	Now we consider three cases.
	If $z\in N_{\bar w}(u)\cap N_{\bar w}(v)$, then $\bar{w}(\oriented{uz}) = w(\oriented{uz}) - \mu(z)$  and $\bar{w}(\oriented{vz}) = w(\oriented{vz}) - \mu(z)$.
	Using this in \eqref{equation:addingskewspartitioncheck} gives
	\[ \bar\sigma_A^1(z)+\sigma_A^1(z)+\bar\sigma_B^1(z)+\sigma_B^1(z) \le \max\{w(\oriented{uz}), w(\oriented{vz})\}+\sigma^1_A(z)+\sigma^1_B(z) - \mu(z), \]
	and the last term is at most $\max\{w(\oriented{uz}), w(\oriented{vz})\}$ since $\sigma_A + \sigma_B \trianglelefteq \mu$, thus giving \ref{itdef:anchorpartition} in this case.
	Now assume that $z\in N_{\bar w}(u) \setminus N_{\bar w}(v)$.
	Then $0 < \bar{w}(\oriented{uz}) = w(\oriented{uz}) - \mu(z)$, but $\bar{w}(\oriented{vz}) = 0$, so $\max\{\bar w(\oriented{uz}), \bar w(\oriented{vz})\} = w(\oriented{uz}) - \mu(z)$.
	Using this in \eqref{equation:addingskewspartitioncheck} gives
	\begin{align*}
		\bar\sigma_A^1(z)+\sigma_A^1(z)+\bar\sigma_B^1(z)+\sigma_B^1(z)  \le w(\oriented{uz}) - \mu(z) +\sigma^1_A(z)+\sigma^1_B(z) \leq w(\oriented{uz}),
	\end{align*}
	where again we used $\sigma_A + \sigma_B \trianglelefteq \mu$ in the last inequality.
	The case $z\in N_{\bar w}(v) \setminus N_{\bar w}(u)$ follows by a symmetric argument, so it only remains to check the case where $z \notin  N_{\bar w}(u) \cup N_{\bar w}(v)$.
	In this case, $\bar{w}(\oriented{uz}) = \bar{w}(\oriented{vz}) = 0$, so $\bar\sigma_A^1(z) + \bar\sigma_B^1(z) = 0$.
	Hence, left-hand side of \eqref{equation:addingskewspartitioncheck} becomes
	\[ \sigma_A^1(z)+\sigma_B^1(z) \le \max\{w(\oriented{uz}), w(\oriented{vz})\}, \]
	where we used that \ref{itdef:anchorpartition} holds for $(\sigma_A, \sigma_B)$.
	This gives \ref{itdef:anchorpartition}, and we are done.
\end{proof}

%

\begin{remark} \label{crem:adding-skew-matchings}
	\Cref{prop:adding-skew-matchings} is meaningful in the case where we take some of the skew-matchings to be identically equal to zero. We denote such a zero-valued skew-matching by~$\sigma_\emptyset$.
	Indeed, any $\gamma$-skew-matching $\sigma$ with its anchor 
	fitting in the $w$-neighbourhood of a vertex $v$ can be paired-up with an empty skew-matching to form a $(\gamma, \gamma')$-skew pair $(\sigma, \sigma_\emptyset)$ in $(G,w)$ for any $\gamma'\ge 0$.
	This will allow us to apply \Cref{prop:adding-skew-matchings} replacing the input of a skew pair with just a single skew-matching.
\end{remark}

The next proposition will be useful to estimate degrees of a vertex in different combinations of matchings and weights obtained after truncations.

\begin{proposition}\label{prop:adding-degree-truncated}
Let $(G,w)$ be a weighted symmetric digraph.
Let $\mu'\le \mu$ be fractional matchings in (unoriented) $G$ and let $(G, w')$ be the $\mu'$-truncated weighted digraph obtained from $(G,w)$.
Then 
\begin{equation*}
\deg_w(v, \mu')+\deg_{w'}(v, \mu-\mu')=\deg_w(v, \mu).
\end{equation*}
\end{proposition}

\begin{proof}
We have 
\[
	\deg_{w}(v, \mu')+\deg_{w'}(v, \mu-\mu')
	=\sum_{x\in V(G)}\left(\min\{w(\oriented{vx}),\mu'(x)\}+\min\{w'(\oriented{vx}), \mu(x)-\mu'(x)\}\right),\]
so it suffices to show that, for every $x \in V(G)$, we have
\[ \min\{w(\oriented{vx}),\mu'(x)\}+\min\{w'(\oriented{vx}), \mu(x)-\mu'(x)\} = \min\{ w(\oriented{vx}), \mu(x) \}. \]

Suppose first that $\mu'(x) < w(\oriented{vx})$ holds.
Then we have $w'(\oriented{vx}) = w(\oriented{vx}) - \mu'(x)$, and the left-hand side of the desired equality becomes $\mu'(x) + \min\{w(\oriented{vx}) - \mu'(x), \mu(x) - \mu'(x)\} = \min\{w(\oriented{vx}), \mu(x)\}$, as we wanted to show.
Hence, we can suppose that $w(\oriented{vx}) \leq \mu'(x)$.
In this case, we have $w'(\oriented{vx}) = 0$.
Then the left-hand side of the desired equality becomes $w(\oriented{vx})$, which is also equal to $\min\{w(\oriented{vx}), \mu(x)\}$ since $w(\oriented{vx}) \leq \mu'(x) \leq \mu(x)$.
\end{proof}

\section{Fractional structure of Gallai--Edmonds decompositions} \label{section:fractionalGE}

In this section, we investigate the structure of skew-matchings in general graphs, relying on known structural results about matchings in graphs.
Our starting point is the classical Gallai--Edmonds decomposition theorem~\cite{Gallai1964, Edmonds1965}, which provides a vertex-partition which exhibits crucial information about maximal matchings in a graph.
Recall that a graph $G$ is said to be \emph{factor-critical}
if for any vertex $v\in V(G)$ there is a perfect matching $M\subseteq E(G)$ covering $G-\{v\}$.
We say that a component of a graph is \emph{factor-critical} if the graph induced by this component is factor-critical.
The following statement appears in \cite[Theorem 2.2.3]{Diestel}.

\begin{theorem}[Gallai--Edmonds theorem]\label{thm:Gallai-Edmond}
	For any graph $G$,
	there is a set $S\subseteq V(G)$,
	called a \emph{separator},
	such that each component of $G-S$ is factor-critical
	and there is a matching $M_S\subseteq E(G)$ between $S$ and $V(G)\setminus S$,
	that covers $S$ and matches the elements of $S$ into different components of $G-S$;
	i.e., for any  component $K$ of $G-S$, we have $|V(K)\cap V(M_S)|\le 1$. 
\end{theorem}

\begin{figure}

\begin{center}
    \begin{tikzpicture}[scale=1, every node/.style={draw, circle, fill=white, inner sep=2pt}]
        \draw[rounded corners] (-4,-2) rectangle (6,-3);
        
        \node[draw=none] at (6.5,-2.5) {\(S\)}; 
        \node[draw=none] at (5.5,-1.55) {$M_S$};  
        
        \node (S1) at (-3,-2.5) {};
        \node (S2) at (-2,-2.5) {};
        \node (S3) at (-1,-2.5) {}; 
        \node (S4) at (0,-2.5) {};
        \node (S5) at (1,-2.5) {};
        \node (S6) at (2,-2.5) {};
        \node (S7) at (3,-2.5) {};
        \node (S8) at (4,-2.5) {};
        \node (S9) at (5,-2.5) {};
        
        \draw[thick] (-4,1) ellipse (0.4 and 1.5);
        \node (F1) at (-4,0) {};           
        \node (F2) at (-3.85,0.3) {};
        \node (F3) at (-4.15,0.5) {};
        \node (F4) at (-3.85,0.8) {};
        \node (F5) at (-4.15,1.0) {};
        \node (F6) at (-3.85,1.3) {};
        \node (F7) at (-4.15,1.5) {};
        \node (F8) at (-3.85,1.8) {};
        \node (F9) at (-4.15,2.0) {};
        
        \draw[thick] (-3,1) ellipse (0.4 and 2.5);
        \node (A1) at (-3,-1) {};
        \node (A2) at (-3.15,-0.5) {};
        \node (A3) at (-2.85,-0.3) {};
        \node (A4) at (-3.15,0) {};
        \node (A5) at (-2.85,0.2) {};
        \node (A6) at (-3.15,0.5) {};
        \node (A7) at (-2.85,0.7) {};
        \node (A8) at (-3.15,1) {};
        \node (A9) at (-2.85,1.2) {};
        \node (A10) at (-3.15,1.5) {}; 
        \node (A11) at (-2.85,1.7) {}; 
        \node (A12) at (-3.15,2) {}; 
        \node (A13) at (-2.85,2.2) {}; 
        \node (A14) at (-3.15,2.5) {}; 
        \node (A15) at (-2.85,2.7) {}; 
        
        \draw[thick] (-2,1) ellipse (0.4 and 2);
        \node (B1) at (-2,-0.7) {};
        \node (B2) at (-2.15,-0.4) {};
        \node (B3) at (-1.85,-0.2) {};
        \node (B4) at (-2.15,0.1) {};
        \node (B5) at (-1.85,0.3) {};
        \node (B6) at (-2.15,0.6) {};
        \node (B7) at (-1.85,0.8) {};
        \node (B8) at (-2.15,1.1) {};
        \node (B9) at (-1.85,1.3) {};
        \node (B10) at (-2.15,1.6) {};
        \node (B11) at (-1.85,1.8) {};
        \node (B12) at (-2.15,2.1) {};
        \node (B13) at (-1.85,2.3) {};
        
        \draw[thick] (-1,0) ellipse (0.4 and 1.5);
        \node (C1) at (-1,-1) {};
        \node (C2) at (-1.15,-0.7) {};
        \node (C3) at (-0.85,-0.5) {};
        \node (C4) at (-1.15,-0.2) {};
        \node (C5) at (-0.85,0) {};
        \node (C6) at (-1.15,0.3) {};
        \node (C7) at (-0.85,0.5) {};
        \node (C8) at (-1.15,0.8) {};
        \node (C9) at (-0.85,1.0) {};
        
        \draw[thick] (0,1) ellipse (0.4 and 1.5);
        \node (E1) at (0,0) {}; 
        \node (E2) at (0.15,0.3) {};
        \node (E3) at (-0.15,0.5) {};
        \node (E4) at (0.15,0.8) {};
        \node (E5) at (-0.15,1.0) {};
        \node (E6) at (0.15,1.3) {};
        \node (E7) at (-0.15,1.5) {};
        \node (E8) at (0.15,1.8) {};
        \node (E9) at (-0.15,2.0) {};
        
        \node[draw=none] at (1,1) {$\mathcal{K}^*_S$};
        
        \node (D1) at (1,-1) {};
        \node (D2) at (2,-1) {};
        \node (D3) at (3,-1) {};
        \node (D4) at (4,-1) {};
        \node (D5) at (5,-1) {};
        \node (D6) at (6,-1) {};
        \node (D7) at (7,-1) {};
        \node (D8) at (8,-1) {};
        
        
        \draw (F2) -- (F3);
        \draw (F4) -- (F5);
        \draw (F6) -- (F7);
        \draw (F8) -- (F9);
        
        \draw (A2) -- (A3);
        \draw (A4) -- (A5);
        \draw (A6) -- (A7);
        \draw (A8) -- (A9);
        \draw (A10) -- (A11);
        \draw (A12) -- (A13);
        \draw (A14) -- (A15);
        
        \draw (B2) -- (B3);
        \draw (B4) -- (B5);
        \draw (B6) -- (B7);
        \draw (B8) -- (B9);
        \draw (B10) -- (B11);
        \draw (B12) -- (B13);
        
        \draw (C2) -- (C3);
        \draw (C4) -- (C5);
        \draw (C6) -- (C7);
        \draw (C8) -- (C9);
        
        \draw (E2) -- (E3);
        \draw (E4) -- (E5);
        \draw (E6) -- (E7);
        \draw (E8) -- (E9);
        
        \draw[dashed] (S1) -- (A1);
        \draw[dashed] (S2) -- (B1);
        \draw[dashed] (S3) -- (C1);
        \draw[dashed] (S4) -- (E1);
        \draw[dashed] (S5) -- (D1);
        \draw[dashed] (S6) -- (D2);
        \draw[dashed] (S7) -- (D3);
        \draw[dashed] (S8) -- (D4);
        \draw[dashed] (S9) -- (D5);
    \end{tikzpicture}
\end{center}
\caption{A view of a $(G, S, M)$ Gallai--Edmonds triple. Note that however there might be more edges, here we just draw the matching $ M_S $ between $S$ and components of $G-S$ by dashed edges. We also illustrate the notion of a factor-critical component by completing $M_S$ with a matching (full edges) in each non-singleton component of $\mathcal K^*-V(M_S)$.}
	\label{fig:edmond-gallai}
\end{figure}
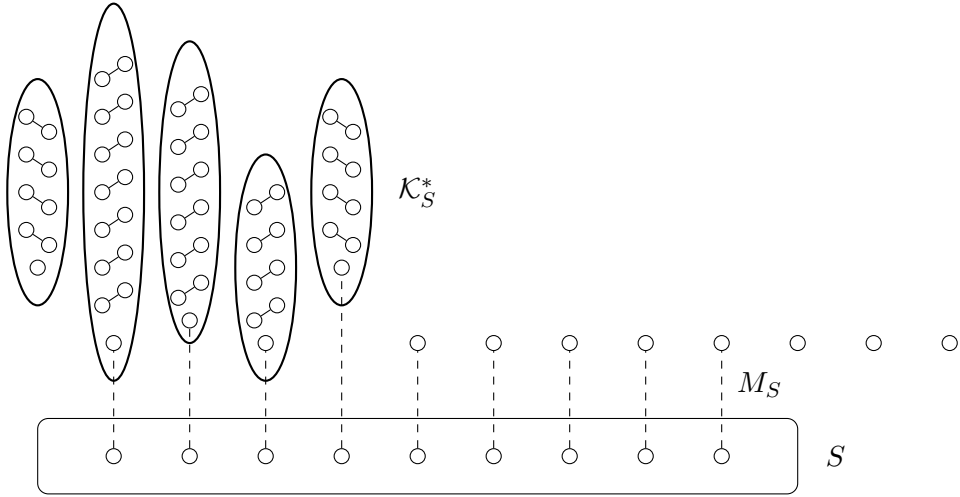

Given a graph $G$ and $S \subseteq V(G)$, $M_S \subseteq E(G)$ as in \Cref{thm:Gallai-Edmond}, we say $(G, S, M_S)$ is a \emph{Gallai--Edmonds triple}.
It will be important to us to consider the components in $G-S$ and to distinguish whether they correspond to single vertices or not.
We let $\mathcal{K}_S$ be the set of components of $G-S$, and we let $\mathcal{K}^\ast_S \subseteq \mathcal{K}_S$ correspond to the non-singleton components of $G-S$.
We also let $U_{S} = \{ v \in V(G) : \{v\} \in \mathcal{K}_S \setminus \mathcal{K}^\ast_S \}$ be the set of vertices corresponding to singleton components in $\mathcal{K}_S$.
See \Cref{fig:edmond-gallai} for an illustration of all of these objects.

In our proofs we will work with fractional and skew-fractional matchings which are supported on edges beyond those of $M_S$.
The following definition captures the set of edges that we will use.

\begin{definition}[Gallai--Edmonds support]\label{def:GE-support}
	Let $(G, S, M)$ be a Gallai--Edmonds triple.
	We say that $E_{S,M} \subseteq E(G)$ is the \emph{Gallai--Edmonds support} of $(G,S,M)$ if 
	\[ E_{S,M} := M \cup E(G[S, U_S]) \cup \bigcup_{K \in \mathcal{K}^\ast_S} E(G[K]). \]
\end{definition}

\begin{remark} \label{remark:GE-support}
	Observe that $E_{S,M}$ does not include any edge with two endpoints in $S$. 
\end{remark}


\subsection{Fractional matchings and $c$-optimal fractional matchings}
We now consider fractional matchings associated with a Gallai--Edmonds triple $(G,S,M)$.
The fractional matchings we consider will be supported by the Gallai--Edmonds support, and cover both the separator $S$ and the non-singleton components of $G - S$.

\begin{definition}[Fractional Gallai--Edmonds triple] \label{definition:gallaiedmonds-fractional}
Let $\mu$ be a fractional matching in $G$. We say that $(G, S, \mu)$ is a \emph{fractional Gallai--Edmonds triple} if there is a matching $M\subseteq E(G)$, such that $(G,S,M)$ is a Gallai--Edmonds triple, and
\stepcounter{propcounter}
\begin{enumerate}[\upshape{(\Alph{propcounter}\arabic*)},topsep=0.7em, itemsep=0.5em]
		\item \label{item:fractionalgallaiedmonds-support} $\mu$ is supported in the Gallai--Edmonds support $E_{S,M}$,
		\item \label{item:fractionalgallaiedmonds-cover} $S \cup \bigcup_{K \in \mathcal{K}^\ast_S} V(K)$ is covered by $\mu$. 
	\end{enumerate} 
\end{definition}

Among all fractional Gallai--Edmond triples $(G, S, \mu)$, we will work with those that are in some sense ``optimal'' with respect to a vertex $c$: they put as much weight as possible in the neighbourhood of $c$; with respect to a weight function $w$ on the edges of $G$.
To formalise this definition we will need a fractional notion of an \emph{alternating path}.

\begin{definition}[Alternating path]\label{def:alternatingpath}
	For a fractional matching $\mu$ in a graph $G$ and a vertex $u\in V(G)$, a \emph{$\mu$-alternating path starting at $u$}  is a path $P=(u,v_1,v_2,\ldots, v_\ell)$ in $G$ starting at $u$, such that $\mu(\fmat{v_{2i-1}v_{2i}})>0$ for all $i \in \{1, \ldots, \lfloor\ell/2\rfloor\}$.
	The \emph{thickness} of $P$ is  $\theta(P):=\min\{\mu(\fmat{v_{2i-1}v_{2i}})\::\: i \in \{1, \ldots, \lfloor\ell/2\rfloor \} \}$. 
\end{definition}

Given this, we can precisely define what are the fractional matchings we want to use in our proofs.

\begin{definition}[$c$-optimal fractional matchings]\label{def:coptimal}
	Let $(G,w)$ be a weighted graph, $G^\leftrightarrow$ be its associated digraph (with inherited weights),  $(G, S, M)$ be a Gallai--Edmonds triple, and $c \in S$.
	We say that a fractional matching $\mu$ is a \emph{$c$-optimal fractional matching} if $(G, S, \mu)$ is a fractional Gallai--Edmonds triple, and moreover the following holds.
	For each singleton component $\{u\}$ of $G - S$ with $\mu(u)<w(\oriented{cu})$,
	for each $\mu$-alternating path $P_u$ starting at $u$, and
	for each $v \in V(P_u) \cap  S$, 
	let
	\[N_\mu(v) = \{ x\in N(v): \mu(\fmat{vx})>0 \}.\] Then we have
	\stepcounter{propcounter}
	\begin{enumerate}[\upshape{(\Alph{propcounter}\arabic*)},topsep=0.7em, itemsep=0.5em]
		\item  \label{it:weighted-R}$N_\mu(v)$ consists of singleton components of $G -S$,
		\item \label{it:weighted-Scovered}$\mu(v)=1$,
		\item \label{it:mu-coveredbyN(c)} for each $x\in N_\mu(v)$ we have 
		$\mu(x)\le  w(\oriented{cx})$
	\end{enumerate}
\end{definition}

\begin{remark} \label{remark:coptimal-neighbourhoodofR}
	If $\{u\}$ is a singleton component of $G-S$ with $\mu(u) < w(\oriented{cu})$, then we have $N_w(u) \subseteq S$.
	Also, for each $v \in N_w(u)$ the edge $uv$ forms a $\mu$-alternating path of length 1.
	Then in particular properties \ref{it:weighted-R}--\ref{it:mu-coveredbyN(c)} apply to all $v \in N_w(u)$.
\end{remark}

\begin{remark} \label{remark:coptimal-neighbourhoodc}
	If $v$ is as above, and $x \in N_{\mu}(v)$, property \ref{it:mu-coveredbyN(c)} implies that $w(\oriented{cx}) \geq \mu(x) > 0$, so in particular $x \in N_w(c)$.
	In short, $N_\mu(v) \subseteq N_w(c)$.
\end{remark}

The following property is key, and shows that indeed $c$-optimal fractional matchings exist given a Gallai--Edmonds triple $(G, S, M)$, as long as we pick $c$ in the separator $S$.

\begin{proposition} \label{proposition:weighted-GEmatching-c}
	Let $(G,w)$ be a weighted graph, $(G, S, M)$ be a Gallai--Edmonds triple, and $c \in S$.
	Then there exists a fractional Gallai--Edmonds triple $(G,S,\mu)$,
	such that $\mu$ is a $c$-optimal fractional matching.
\end{proposition}

\begin{figure}
\begin{center}
\begin{tikzpicture}[scale=1, every node/.style={draw, circle, fill=white, inner sep=2pt}]
    \draw[rounded corners] (0.8,-2) -- (5.5,-2) -- (5.5,-3) -- (0.8,-3);
    \node[draw=none] at (0.5,-2.5) {\(\ldots\)};
     \node[draw=none] at (3.3,-0.5) {\(x\)};
      \node[draw=none, text=blue] at (4.5,0) {\(N_\mu(v)\)};
        \node[draw=none, text=red] at (1.6,-1.4) {\(P_u\)};

    \node (D1) at (1,-0.5) {};
    \node (D2) at (2,-0.5) {};
    \node (D3) at (3,-0.5) {};
    \node (D4) at (4,-0.5) {};
    \node (D5) at (5,-0.5) {};
    \node (D6) at (6,-0.5) {};
    \node (D7) at (7,-0.5) {};
    \node (D8) at (8,-0.5) {};

    \draw[thick, blue] 
  ($ (D3)!0.5!(D4) $)        
  ellipse [x radius=0.8, y radius=0.4];
    \node[draw=none, anchor=west] at (D8.east) {\(u\)};
    
    \node (S5) at (1,-2.5) {};
    \node (S6) at (2,-2.5) {};
    \node (S7) at (3,-2.5) {};
    \node (S8) at (4,-2.5) {};
    \node (S9) at (5,-2.5) {};
    
    \node[draw=none, anchor=north] at ($(S7.south)+(0.22,0.1)$) {\(v\)};

    \node[draw=none] at (6.5,-2.5) {\(S\)};
    \node[draw=none] at (5.9,-1.4) {\(\mu\)};

    \draw[dashed] (S5) -- (D1);
    \draw[dashed] (S6) -- (D2);
    \draw[dashed] (S7) -- (D3);
    \draw[dashed] (S8) -- (D4);
    \draw[dashed] (S9) -- (D5);
    
    \draw[line width=2pt, red]         (D8) -- (S9);
    \draw[line width=2pt, red, dashed] (S9) -- (D5);
    \draw[line width=2pt, dashed] (S9) -- (D6);
    \draw[line width=2pt, red]         (D5) -- (S7);
    \draw[line width=2pt, red, dashed] (S7) -- (D3);
    \draw[line width=2pt, dashed] (S7) -- (D4);
    \draw[line width=2pt, red]         (D3) -- (S6);
    \draw[line width=2pt, red,dashed] (S6) -- (D2);

\end{tikzpicture}
\end{center}
	\caption{An alternating path $ P_u $ in bold red edges starting at $ u $ in a $(G, S, M)$ Gallai--Edmonds triple structure. The support of the fractional matching $\mu$ is illustrated in dashed edges. For any vertex $ v $ in $ V(P_u)\cap S $, the set
$N_\mu(v) $ (shown in blue) consists of the neighbours of $v $ connected by support edges of $\mu$.}
	\label{fig:alternating}
\end{figure}
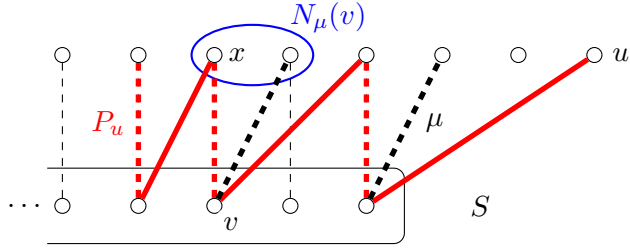

To prove \Cref{proposition:weighted-GEmatching-c} we will need the following simple result.

\begin{proposition}\label{prop:factor-critical->fractional matching}
	Every factor-critical graph on more than one vertex 
	has a perfect fractional matching.
\end{proposition}

\begin{proof}
	Since $G$ is factor-critical, for each vertex $v \in V(G)$ there exists a matching~$M_v$ which covers $G - \{v\}$.
	Let $\mu_v:E(G)\rightarrow \{0,1\}$ be such that $\mu_v(M_v)\equiv 1$ and $\mu_v(E(G)\setminus M_v)\equiv 0$.
	Then $\mu = \sum_{v \in V(G)} \mu_v / (|V(G)|-1)$ is a perfect fractional matching of~$G$.
\end{proof}

\begin{proof}[Proof of \Cref{proposition:weighted-GEmatching-c}]
	Among all possible Gallai--Edmonds triples $(G,S,M)$,
	choose one with matching $M_S$ that minimizes the number of vertices in $N(c)$ which lie in singleton components of $G-S$ and are not contained in $M_S$.
	Let $E_{S,M_S}$ be the corresponding Gallai--Edmonds support of $(G, S, M_S)$.
	The matching $M_S$ can be understood as a fractional matching $\mu_S$ with weights $1$ on $M_S$ and $0$ on other edges.
	This fractional matching can be extended to cover all non-trivial components~$K$ of $G-S$; either by using \Cref{prop:factor-critical->fractional matching} if $V(K)\cap M_S=\emptyset$, or using the factor-criticality of~$K$ to find a perfect matching in $K-V(M_S)$. 
	By doing this, we have obtained a fractional matching $\mu'$ such that
	\begin{enumerate}
		\item $\mu'$ covers~$S$,
		\item the only uncovered vertices by $\mu'$ form isolated components in $G-S$,
		\item $\mu'$ is supported in $E_{S, M_S}$, and
		\item \label{item:weigthtedGE-item3} for each isolated component $\{u\}$ of $G - S$, we have $\mu'(u) \in \{0,1\}$, and $\mu'(u) = 1$ if and only if $u$ was covered by $M_S$.
	\end{enumerate}  
	
	\begin{claim}
		$\mu'$ satisfies \ref{it:weighted-Scovered} and \ref{it:weighted-R}.
	\end{claim}

	\begin{proofclaim}
	Let $\{u\}$ be any singleton component of $G-S$ with $\mu'(u) < w(\oriented{cu})$,
	let $P_u$ be any $\mu'$-alternating path $P_u$ starting at $u$,
	and let $v \in V(P_u) \cap S$ be arbitrary.
	Since $\mu'$ covers $S$, indeed~\ref{it:weighted-Scovered} holds.
	
	It only remains to see that \ref{it:weighted-R} holds.
	Note that $N_{\mu'}(v)\cap S = \emptyset$ because $\mu'$ is supported in $E_{S,M_S}$.
	Moreover, it follows from $0\leq \mu'(u)< w(\oriented{cu})\leq 1$ that $u \in N(c)$.
	By \ref{item:weigthtedGE-item3}, we have that $\mu'(u) = 0$, and that $u$ was not covered by $M_S$.
	For a contradiction, suppose that $x \in N_{\mu'}(v)$ is not a singleton component of $G - S$.
	By definition of $N_{\mu'}(v)$, we have that $\mu'(\fmat{vx}) > 0$.
	We also know that $x \notin S$.
	Since $\mu'$ is supported in $E_{S, M_S}$ in fact we deduce that $vx \in M_S$.
	Since $v \in V(P_u) \cap S$, we deduce that there exists $i \geq 1$ and a $\mu'$-alternating path $u v_1 v_2 \dotsb v_{2i-1} v_{2i}$ with $\mu'(\fmat{v_{2j-1} v_{2j}}) = 1$ for all $1 \leq j \leq i$, and $v_{2i-1} = v$ and $v_{2i} = x$.
	This implies that $v_1 v_2, \dotsc, v_{2i-1} v_{2i} \subseteq M_S$.
	Define $M'_S := (M_S \setminus \{v_1 v_2, \dotsc, v_{2i - 1} v_{2i} \}) \cup \{ u v_1, v_2 v_3, \dotsc, v_{2i-2} v_{2i-1}\}$.
	In summary, by passing from $M_S$ to $M'_S$ we have uncovered $x$ but covered $u$.

	But $(G, S, M'_S)$ is a Gallai--Edmonds triple,
	and since $u \in N(c)$ and $x$ does not form a singleton component in $G - S$, by passing from $M_S$ to $M'_S$ we have decreased the number of vertices in $N(c)$ which lie in singleton components and are not covered.
	This contradicts the minimality of $M_S$.
	Thus, \ref{it:weighted-R} indeed holds.
	\end{proofclaim}

	Let $\mu$ be the fractional matching that maximizes $\deg_w(c, \mu)$ among all fractional matchings $\mu'$ whose support is a subset of $E_{S,M_S}$ and which satisfy \ref{it:weighted-R}, and~\ref{it:weighted-Scovered}.
	We claim that $\mu$ satisfies~\ref{it:mu-coveredbyN(c)}.
	Aiming at a contradiction, suppose this is not the case.
	Thus, there exist a singleton component $\{u\}$ of $G-S$ with $\mu(u) < w(\oriented{cu})$,
	a $\mu$-alternating path $P_{uv}$ starting at $u$ and finishing at $v\in S$,
	and $x\in N_{\mu}(v)$ such that $\mu(x)>w(\oriented{cx})$.
	We can assume that among all the possible such paths we choose one of shortest length.
	Using~\ref{it:weighted-R} we obtain that $x\in G-S$.
	Note that $x \notin V(P_{uv})$, as otherwise we can pass to a shorter path.
	Then the path $P$ consisting of $P_{uv}$, extended with $vx$ at the end is a $\mu$-alternating path starting from~$u$.
	Indeed, $x \in N_\mu(v)$ and therefore $\mu(\fmat{vx})>0$, with $x \notin S$. 

	Suppose $P= v_0 v_1 \dotsb v_{2\ell}$ with $v_0=u$ and $v_{2\ell}=x$.
	By definition of $P$ and $N_\mu(v)$, we have that the thickness $\theta(P)$ of $P$ (see \Cref{def:alternatingpath}) satisfies $\theta(P)>0$.
	By \ref{it:weighted-R}, we have that  $E(P)\subseteq E_{S,M_S}$.

	Let $\delta:= \min\{\mu(x)-w(\oriented{cx}), w(\oriented{cu})-\mu(u), \theta(P)\}>0$.
	We define a new fractional matching $\mu_{P}$ such that
	 $\mu_P(e)= \mu(e)$ if $e\not\in E(P)$; and $
	\mu_P(\fmat{v_{2i+1}v_{2i+2}}):= \mu(\fmat{v_{2i+1}v_{2i+2}}) -\delta$ and $\mu_P(\fmat{v_{2i}v_{2i+1}}):= \mu(\fmat{v_{2i}v_{2i+1}})+\delta$ for every $i \in \{0, \ldots, \ell-1\}$. 
	Then $\mu_P(z)=\mu(z)$ for every $z\in V(G)\setminus \{u,x\}$, and $\mu_P(x)\ge w(\oriented{cx})$.
	We also have $w(\oriented{cu})\ge \mu_P(u)>\mu(u)$. Since $u\in N(c)$, this implies $\deg_w(c,\mu_P)>\deg_w(c,\mu)$, contradicting the maximality of $\deg_w(c,\mu)$ (or equivalently, our assumption that $\mu$ fails to satisfy~\ref{it:mu-coveredbyN(c)}).
\end{proof}

\subsection{Reachable vertices}
Now we take a closer look at the structure of $c$-optimal fractional matchings.
Of particular interest are the vertices outside the separator $S$ that belong to alternating paths, which we call \emph{reachable}. 

\begin{definition}[Reachable vertices]\label{def:weighted-R}
	Let $(G, S, M)$ be a Gallai--Edmonds triple, $c\in S$ and $w:E(G)\rightarrow [0,1]$ be a weight function. 
	Let $(G, S,  \mu)$ be a fractional Gallai--Edmonds triple such that $\mu$ is $c$-optimal.
	We define the set $\mathcal{R}$ of \emph{reachable vertices with respect to $c$} as the set of vertices in $V(G)\setminus S$ contained in some $\mu$-alternating path $P_u$ starting at a singleton component $u$ of $G- S$ with $\mu(u)<w(\oriented{cu})$.
	Also, let $S_{\mathcal{R}}=\bigcup_{x\in {\mathcal{R}}}N(x)$.
\end{definition}

If $x \in {\mathcal{R}}$, then either $x$ is a singleton component in $G-S$ for which $w(\oriented{cx}) > 0$, or there exists a path $P_u$ as in the definition such that $x \in V(P_u)$.
By taking the neighbour of $x$ in $P_u$ that is closer to $u$, we obtain some $v \in S$ such that $x \in N_{\mu}(v)$, so by \Cref{remark:coptimal-neighbourhoodc} we have that $w(\oriented{cx}) > 0$ as well.
Hence, in any case, we have verified the following three consequences.
\begin{observation}\label{ob:R=singletons}
	We have ${\mathcal{R}} \subseteq U_S$, i.e. $\mathcal{R}$ consists only of singleton components of $G-S$.
	In particular, $S_{\mathcal{R}}\subseteq S$. 
\end{observation}

\begin{observation}\label{ob:R=neighbours}
	If $x \in \mathcal{R}$, then $w(\oriented{cx}) > 0$. In particular, $\mathcal{R} \subseteq N_w(c)$.
\end{observation}

\begin{observation} \label{observation:reachable-wversusmu}
	For every $y \in \mathcal{R}$, we have that $w(\oriented{cy}) \geq \mu(y)$.
\end{observation}

The following observation allows us to compare the sizes of $\mathcal{R}$ and $S_\mathcal{R}$.

\begin{observation} \label{observation:reachable-weightsR-SR}
	For every $y \in S_\mathcal{R}$ and $z \in V(H)$ such that $\mu(yz) > 0$, we have $z \in \mathcal{R}$.
	In particular, we have
	\[ \sum_{x \in \mathcal{R}} \mu(x) = \sum_{y \in S_{\mathcal{R}}} \mu(y), \]
	and $|S_\mathcal{R}| \leq |\mathcal{R}|$.
\end{observation}

\begin{proof}
	Let $y \in S_\mathcal{R}$ and $z \in V(H)$ such that $\mu(yz) > 0$.
	We have $yz$ is in the Gallai--Edmonds support, hence $z \notin S$.
	Since $S_\mathcal{R}$ is the set of neighbourhoods of vertices in $\mathcal{R}$, there exists $r \in \mathcal{R}$ such that $ry \in E(H)$.
	By the definition of $\mathcal{R}$, $y$ is the endpoint of a suitable $\mu$-alternating path which we can extend by adding $z$; this implies that $z \in \mathcal{R}$, as desired.
	
	This implies that for every edge such that $\mu(xy) > 0$, $x \in \mathcal{R}$ if and only if $y \in S_\mathcal{R}$.
	This yields the claimed equality.
	Since $\mu$ covers $S$ (by property \ref{item:fractionalgallaiedmonds-cover}), the right-hand side equals $|S_\mathcal{R}|$.
	On the other hand, the left-hand side is at most $|\mathcal{R}|$ because $\mu$ is a fractional matching. This proves the claimed inequality.
\end{proof}

\subsection{GE pairs}
In the proof of \Cref{prop:weighted-structural} we shall work with the concept of \emph{GE pair}.
This will be a pair $(\tilde \sigma, \tilde \mu)$ consisting of a $\gamma$-skew-matching $\tilde \sigma$ and a fractional matching $\tilde \mu$, which are disjoint.
Essentially, we want to generalise the concept of a fractional Gallai--Edmonds triple by allowing some of the weight to be provided by a skew-matching, not only by a fractional matching.

\begin{definition}\label{def:GE-weighted-pair}
	Let $(G,w)$ be a weighted digraph, $(G, S, M)$ be a Gallai--Edmonds triple, $c \in S$, and $\gamma>1$.
	Let $\mu$ be a $c$-optimal fractional matching and let $\mathcal R$ be the set of reachable vertices with respect to $c$ and $\mu$. 
	We say that a pair $(\tilde \sigma, \tilde \mu)$ is a \emph{$\gamma$-GE pair for $c$ with respect to $w$ and $\mu$}, if 
	\stepcounter{propcounter}
	\begin{enumerate}[(\Alph{propcounter}\arabic*),topsep=0.7em, itemsep=0.5em]
		\item \label{item:gepair-1} $\tilde \sigma$ is a $\gamma$-skew-matching,
		\item \label{item:gepair-2} $\tilde\mu$ is a fractional matching disjoint from $\tilde \sigma$,
		\item \label{item:gepair-3} $\mathcal A(\tilde \sigma)$ is contained in $S_{\mathcal R}$,
		\item \label{item:gepair-4} $\mathcal A(\tilde \sigma)$ fits in the $w$-neighbourhood of $c$,
		\item \label{item:gepair-5} 
		$V(\tilde \sigma)\setminus\mathcal A(\tilde \sigma)\subseteq \mathcal R$. 
		\item \label{item:gepair-6} for all $y\in \mathcal R$ we have $ \tilde \mu(y)+\tilde \sigma(y)\le w(\oriented{cy})$,
		\item \label{item:gepair-7} for all $y\in V(G)\setminus \mathcal R$ we have $\tilde \mu(y)+\tilde \sigma(y)\ge w(\oriented{cy})$,
		\item \label{item:gepair-8} $\tilde\mu+\tilde \sigma$ covers $S$, and
		\item \label{item:gepair-9} $\tilde \mu$ equals to $ \mu$ when restricted to the graph $G-(\mathcal R\cup S_{\mathcal R})$ and any supporting edge of $\tilde \mu$ intersecting the set  $\mathcal R\cup S_{\mathcal R}$ lies in the bipartite graph $G[\mathcal R,S_{\mathcal R}]$.
	\end{enumerate} 
%
	
For brevity, we will say that $(\tilde \sigma, \tilde \mu)$ is a \emph{$(G, w, S, M, c,\mu, \gamma)$-GE pair}.
\end{definition}

\begin{remark}
	Note that in our definition we only consider the case $\gamma > 1$.
	While it is in principle possible to define $\gamma$-GE pairs with $\gamma \leq 1$, this will never be used in our proofs.
	The main reason are properties \ref{item:gepair-3} and \ref{item:gepair-5}. For any the directed edges $\oriented{uv}$ with $\tilde \sigma(\oriented{uv}) > 0$ must satisfy $u \in S_{\mathcal R}$ and $v \in \mathcal R$.
	Hence, if $\gamma \leq 1$ then in such a situation the contribution of $\oriented{uv}$ to the weight of $u$ is at least the contribution of that edge in $v$.
	This means that we could replace the weight in $\tilde \sigma$ with some weight in $\mu$ (which covers $u$ and $v$ equally) without decreasing the total weight.
	Thus the definition only represents a real gain over fractional Gallai--Edmonds triples in the case $\gamma > 1$.
\end{remark}

\begin{remark}\label{lemma:gepairs-exist}
	A $\gamma$-GE pair always exists as, for any $c$-optimal fractional matching $\mu$, the pair $(\mu,\emptyset)$ satisfies \Cref{def:GE-weighted-pair}.
	Formally, this follows from the fact that $(G, S, \mu)$ is a fractional Gallai--Edmonds triple, together with the definition of $\mathcal{R}$ and \Cref{observation:reachable-wversusmu,observation:reachable-weightsR-SR}.
\end{remark}

\subsection{Separation with GE pairs}

Let $(\tilde \sigma, \tilde \mu)$ be a $(G, w, S, M, c, \mu, \gamma)$-GE pair.
Among all such pairs, we will work with pairs which saturate the neighbourhood of~$c$ as much as possible, that is, we want that
$\deg_w(c, \tilde \mu+\tilde \sigma)$
is maximum over all  choices of $(\tilde \sigma, \tilde \mu)$. 
In such a case, we say that $(\tilde \sigma, \tilde \mu)$ is \emph{optimal}.

We adapt a generalization of an alternating path argument to obtain some structural information about optimal GE pairs; that is the content of the next two lemmas.

\begin{lemma}[Separating Lemma I] \label{lemma:separatinglemma-1}
	Let $(\tilde \sigma, \tilde \mu)$ be an optimal $(G, w, S, M, c, \mu, \gamma)$-GE pair.
	Suppose there exists $d \in \mathcal{R}$ such that $\tilde \sigma(d) + \tilde \mu(d) < w(\oriented{cd})$.
	Then we have $\tilde \sigma(x)= w(\oriented{cx})$  for all $x\in N_w(c)\cap N_w(d)$. 
\end{lemma}

\begin{figure}[ht]
\centering
\begin{subfigure}{0.45\linewidth}
\centering
\begin{tikzpicture}[
  every node/.style={draw, circle, fill=white, minimum size=1.4cm, inner sep=2pt}
]
  \pgfmathsetmacro{\r}{0.7}

  \node[draw=none] at (-1.4,0) {$\ldots$};
  \node[draw=none] at (1.9,0) {$S$};
  \node[draw=none] at (0.8,-0.5) {$x$};
  \node[draw=none] at (2.7,3.6) {$d$};
  \node[draw=none] at (0.8,3.6) {$y$};
    \node[draw=none] at (-1,1.5) {$\textcolor{red}{\tilde\mu}$};

  \node (S1) [draw=red, thick, fill=red!30] at (0,0) {};
  \node (D1) [draw=red, thick, fill=red!30] at (0,3) {};
  \node (D2)  at (2,3) {};

  \draw[rounded corners]
    (-1.4,0.9) -- (1.4,0.9) -- (1.4,-0.9) -- (-1.4,-0.9);

  \draw[red, thick] ($ (S1)+(\r,0) $) -- ($ (D1)+(\r,0) $);
  \draw[red, thick] ($ (S1)+(-\r,0) $) -- ($ (D1)+(-\r,0) $);

\path[fill=black!22, draw=none]
  (D2) ++(-90:\r) arc (-90:90:\r) -- ($(D2)+(0,-\r)$) -- cycle;  
  
  \draw[black] ($(D2)+(0,-\r)$) -- ($(D2)+(0,\r)$);

  \coordinate (T1a) at (tangent cs:node=S1, point={(D2)}, solution=1);
  \coordinate (T2a) at (tangent cs:node=D2, point={(S1)}, solution=2);
  \draw (T1a) -- (T2a);

  \coordinate (T1b) at (tangent cs:node=S1, point={(D2)}, solution=2);
  \coordinate (T2b) at (tangent cs:node=D2, point={(S1)}, solution=1);
  \draw (T1b) -- (T2b);
\end{tikzpicture}
\caption{For simplicity, we illustrate here the easy case, when 
${w(\oriented{cd})}=1/2$, ${w(\oriented{cx})}=w(\oriented{cy})=1$,  $\textcolor{red}{\tilde\mu(\fmat{xy})}=1$ thus $ \textcolor{blue}{\tilde \sigma(x)}=0$, and $\textcolor{blue}{\tilde\sigma(d)}+\textcolor{red}{\tilde\mu(d)}=0$. This situation contradicts the conclusion of \Cref{lemma:separatinglemma-1}.
}
\end{subfigure}
\hfill
\begin{subfigure}{0.45\linewidth}
\centering
\begin{tikzpicture}[
  every node/.style={draw, circle, fill=white, minimum size=1.4cm, inner sep=2pt}
]
  \pgfmathsetmacro{\r}{0.7}
  \pgfmathsetmacro{\eps}{0.06} 

  \node[draw=none] at (-1.4,0) {$\ldots$};
  \node[draw=none] at (1.9,0) {$S$};
  \node[draw=none] at (0.8,-0.5) {$x$};
  \node[draw=none] at (2.7,3.6) {$d$};
  \node[draw=none] at (0.8,3.6) {$y$};
    \node[draw=none] at (-1,1.5) {$\textcolor{red}{\tilde\mu'}$};
      \node[draw=none] at (2,1.5) {$\textcolor{blue}{\tilde\sigma'}$};

  \node (S1)  at (0,0) {};
  \node (D1)  at (0,3) {};
  \node (D2)  at (2,3) {};
  
\path[fill=black!22, draw=none]
  (D2) ++(-90:\r) arc (-90:90:\r) -- ($(D2)+(0,-\r)$) -- cycle;  
  
\path[fill=blue!22, draw=none]
  (D1) ++(-90:\r) arc (-90:90:\r) -- ($(D1)+(0,-\r)$) -- cycle;
 \path[fill=blue!22, draw=none]
  (D2) ++(90:\r) arc[start angle=90, end angle=270, radius=\r]
  -- ($(D2)+(0,-\r)$) -- cycle;
  \path[fill=blue!18, draw=none]
  (S1) ++(-90:\r) arc (-90:90:\r) -- ($(S1)+(0,-\r)$) -- cycle;

\path[fill=red!22, draw=none]
  (S1) ++(90:\r) arc[start angle=90, end angle=270, radius=\r]
  -- ($(S1)+(0,-\r)$) -- cycle;

\path[fill=red!22, draw=none]
  (D1) ++(90:\r) arc[start angle=90, end angle=270, radius=\r]
  -- ($(D1)+(0,-\r)$) -- cycle;

  \draw[rounded corners]
    (-1.4,0.9) -- (1.4,0.9) -- (1.4,-0.9) -- (-1.4,-0.9);

 \draw[red, thick] (S1) ++(90:\r) arc[start angle=90, end angle=270, radius=\r];

  \draw[red, thick] (D1) ++(90:\r) arc[start angle=90, end angle=270, radius=\r];
  \draw[red, thick] ($ (D1)+(0,\r) $) -- ($ (S1)+(0,-\r) $);
  \draw[blue, thick] ($(D1)+(\eps,\r)$) -- ($(S1)+(\eps,-\r)$);

  \draw[blue, thick] (D1) ++(-90:\r) arc[start angle=-90, end angle=90, radius=\r];
  
  \draw[blue, thick] (S1) ++(-90:\r) arc[start angle=-90, end angle=90, radius=\r];
  
\pgfmathsetmacro{\cx}{\r*cos(42)}
\pgfmathsetmacro{\cy}{\r*sin(42)}

\coordinate (X45) at ($(S1)+(\cx,\cy)$);      
\coordinate (Xepsbot) at ($(S1)+(\eps,-\r)$); 
  
  \pgfmathsetmacro{\cdx}{\r*cos(153)} 
\pgfmathsetmacro{\cdy}{\r*sin(153)} 
\coordinate (D2left) at ($(D2)+(\cdx,\cdy)$);
  
    \pgfmathsetmacro{\cxx}{\r*cos(325)} 
\pgfmathsetmacro{\cxy}{\r*sin(325)} 
\coordinate (S1Bot) at ($(S1)+(\cxx,\cxy)$);

\draw[blue, thick] (Xepsbot) -- (X45);

\draw[blue, thick] (X45) -- ($ (D1)+(\r, 0) $);
\draw[blue, thick] (X45) --  (D2left);
\draw[blue, thick] (S1Bot) --  ($ (D2)+(0, -\r) $);

  \draw[blue, thick] (D2) ++(90:\r) arc[start angle=90, end angle=270, radius=\r];
  \draw[blue, thick] ($(D2)+(0,-\r)$) -- ($(D2)+(0,\r)$);
  
  \draw[red, thick] ($ (S1)+(-\r,0) $) -- ($ (D1)+(-\r,0) $);

  \draw ($ (S1)+(\r,0) $) -- ($ (D1)+(\r,0) $);
  \coordinate (T1a) at (tangent cs:node=S1, point={(D2)}, solution=1);
  \coordinate (T2a) at (tangent cs:node=D2, point={(S1)}, solution=2);
  \draw (T1a) -- (T2a);
  
   \coordinate (T1b) at (tangent cs:node=S1, point={(D2)}, solution=2);
  \coordinate (T2b) at (tangent cs:node=D2, point={(S1)}, solution=1);
  \draw (T1b) -- (T2b);
  
\end{tikzpicture}
\caption{Assuming $\gamma=2$, we obtain $\delta = 1/4$, and after  modification, we have a new pair $(\textcolor{red}{\tilde \mu'}, \textcolor{blue}{\tilde \sigma'})$ with  $\textcolor{red}{\tilde\mu'(\fmat{xy})}=1/2$,  $\textcolor{blue}{\tilde\sigma'(\oriented{xy})}=3/4$, $\textcolor{blue}{\tilde\sigma'(\oriented{xd})}=3/4$. Observe that we have $\textcolor{blue}{\tilde \sigma'(d)}+\textcolor{red}{\tilde \mu'(d)}=\textcolor{blue}{\tilde \sigma'(d)}=1/2=w(\oriented{cd})
$.}
\end{subfigure}

\caption{Situations in the proof of \Cref{lemma:separatinglemma-1}. In the second figure $\textcolor{red}{\tilde \mu'}+\textcolor{blue}{\tilde \sigma'}$ covers $x$ and $y$ as in the first figure, but additionally it covers also half of $d$. Hence $\deg_w(c, \textcolor{red}{\tilde \mu'}+\textcolor{blue}{\tilde \sigma'})>\deg_w(c, \textcolor{red}{\tilde \mu}+\textcolor{blue}{\tilde \sigma})$. Therefore, $(\textcolor{red}{\tilde \mu}, \textcolor{blue}{\tilde \sigma})$ was not an optimal GE pair, as assumed.
}
\label{fig:sep1}
\end{figure}
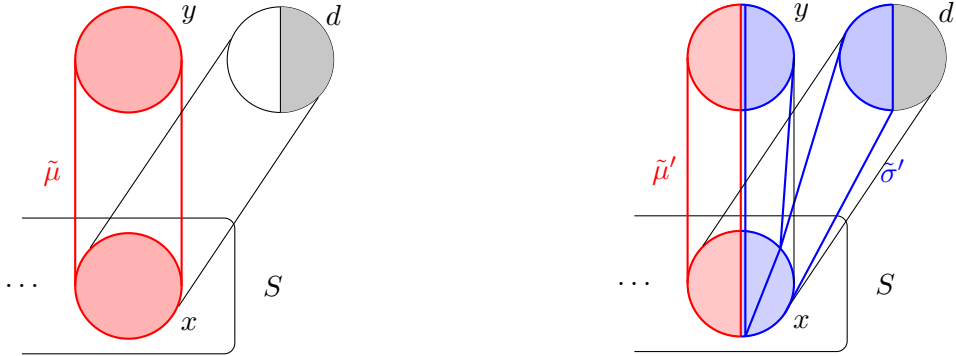

\begin{proof}
	Let $d \in \mathcal{R}$ as in the statement, and let $x \in N_w(c) \cap N_w(d)$ be arbitrary.
	Since $d \in \mathcal{R}$, we have $x \in S_\mathcal{R} $ by \Cref{ob:R=singletons}.
	 We also have $d \in N_w(c)$ by \Cref{ob:R=neighbours}.
	Now, we leverage the fact that if there is 
	any possibility to increase the total saturation (w.r.t. $c$); then clearly the GE pair $(\tilde \sigma, \tilde \mu)$ is not maximizing the saturation of $N_w(c)$, which will be a contradiction.
	
	By \ref{item:gepair-4}, we have $\tilde \sigma(x)\le w(\oriented{cx})$.
	Aiming at a contradiction, suppose that $\tilde \sigma(x)<w(\oriented{cx})$.
	By \ref{item:gepair-8}, we have that $\tilde \sigma(x) + \tilde \mu(x) = 1$, so we must have $\tilde \mu(x) > 0$.
	Hence, there exists $y \in N_{\tilde \mu}(x)$ such that $\tilde\mu(\fmat{xy}) > 0$.  See Picture~(A) of \Cref{fig:sep1}  for an example.
	Note that it could happen that $y = d$.
	Define
	\begin{equation}\label{def:sepI-delta} \delta :=  \min\left\lbrace \frac{\gamma-1}{\gamma}(w(\oriented{cx}) - \tilde \sigma(x)),
	\frac{1}{\gamma}(w(\oriented{cd}) - \tilde \sigma(d) - \tilde \mu(d)),\frac{\gamma-1}{\gamma}
	\tilde \mu(\fmat{xy}) \right\rbrace. 
	\end{equation}
	The first and second terms of this minimum are strictly positive by assumption, and the last term is positive by the choice of $y$.
	Hence, $\delta > 0$.
	
	We define $\tilde \mu'$ by 
	\begin{align*}
		\tilde \mu'(\fmat{xy}) &= \tilde \mu(\fmat{xy}) - \frac{\gamma}{\gamma-1}\delta\overset{\eqref{def:sepI-delta}}{\ge} 0, \mbox{ and}\\
		\tilde \mu'(\fmat{uv}) &= \tilde \mu(\fmat{uv})\mbox{, in any other case.}
	\end{align*}
	Now, we suppose first that $y \neq d$.
	We define $\tilde \sigma'$ by 
	\begin{align*}\tilde \sigma'(\oriented{xd}) &= \tilde \sigma(\oriented{xd}) + (1+\gamma) \delta,\\
	\tilde \sigma'(\oriented{xy}) &= \tilde \sigma(\oriented{xy}) +\frac{1+\gamma}{\gamma-1}\delta, \mbox{ and}\\
	\tilde \sigma'(\oriented{uv}) &= \tilde \sigma(\oriented{uv}) \mbox{, in any other case.}
	\end{align*}

	See Picture~(B) of \Cref{fig:sep1} for an illustration. Now, we have 
	\begin{align}\label{eq:sepI-1}
		\tilde\mu'(d)+\tilde \sigma'(d)&=\tilde \mu (d)+\tilde \sigma(d)+\gamma\delta\overset{\eqref{def:sepI-delta}}{\le} w(\oriented{cd})\le 1,\\
		\label{eq:sepI-2}\tilde \mu'(y)+\tilde \sigma'(y)&=\tilde \mu(y)-\frac{\gamma}{\gamma-1} \delta+\tilde \sigma(y)+\frac{\gamma}{\gamma-1} \delta=\tilde \mu(y)+\tilde\sigma(y)\overset{\text{\ref{item:gepair-6}}}{\le} w(\oriented{cy})\le 1,\\
		\label{eq:sepI-3}\tilde\sigma'(x)+\tilde \mu'(x)&=\tilde \sigma(x)+\delta+\frac{\delta}{\gamma-1}+\tilde \mu(x)-\frac{\gamma}{\gamma-1}\delta=\tilde \sigma(x)+\tilde \mu(x)\overset{\text{\ref{item:gepair-2}}}{=} 1,\\
		\label{eq:sepI-4}\tilde \sigma'(x)&=\tilde \sigma(x)+\delta+ \frac{\delta}{\gamma-1}\overset{\eqref{def:sepI-delta}}{\le} w(\oriented{cx}).
	\end{align}
	
	In the case when $y=d$, we define $\sigma'$ as before for every edge except $\oriented{xd}$, where we set
	\[ \tilde \sigma'(\oriented{xd}) = \tilde \sigma (\oriented{xd}) + (1 + \gamma)\delta + \frac{1 + \gamma}{\gamma - 1} \delta. \]
	Observe that \eqref{eq:sepI-3}--\eqref{eq:sepI-4} still hold with this choice; and instead of \eqref{eq:sepI-1} and \eqref{eq:sepI-2} we have
	\begin{equation}
		\tilde\mu'(d)+\tilde \sigma'(d)=\tilde \mu (d) - \frac{\gamma}{\gamma-1} \delta + \tilde \sigma(d) + \gamma \delta+\frac{\gamma}{\gamma-1} \delta\overset{\eqref{def:sepI-delta}}{\le} w(\oriented{cd}).
		\label{eq:sepI-5}
	\end{equation}
	
	We claim that, in any case, $(\tilde \sigma', \tilde \mu')$ is a $(G, w, S, M, c, \mu, \gamma)$-GE pair.
	First, observe that by the assumption that $d\in \mathcal R$, we get that $x\in S_{\mathcal R}$.
	Together with~\ref{item:gepair-9} on the pair $(\tilde \sigma, \tilde \mu)$, we deduce that $y\in \mathcal R$.
	Property~\ref{item:gepair-1} follows directly from the definition of $\tilde\sigma'$ and from \eqref{eq:sepI-1}-\eqref{eq:sepI-2} and \eqref{eq:sepI-4}-\eqref{eq:sepI-5}; and \ref{item:gepair-2} follows from~\eqref{eq:sepI-1}-\eqref{eq:sepI-3} and \eqref{eq:sepI-5}.
	Properties~\ref{item:gepair-3} and \ref{item:gepair-5} are trivial, since $x \in S_\mathcal R$ and $d,y \in \mathcal{R}$.
	Property~\ref{item:gepair-4} follows from~\eqref{eq:sepI-4}; and~\ref{item:gepair-6} follows from \eqref{eq:sepI-1}-\eqref{eq:sepI-2} and~\eqref{eq:sepI-5}.
	We also have that $\tilde \sigma'(z) + \tilde \mu'(z) = \tilde \sigma(z) + \tilde \mu(z)$ 	for any $z \notin \{d, x, y\}$,	so to get \ref{item:gepair-7} and \ref{item:gepair-8}, it suffices to check it for~$x$, and this follows from~\eqref{eq:sepI-3}. 
	The last property~\ref{item:gepair-9} follows from the fact that the pair $(\tilde \sigma', \tilde \mu')$ satisfies~\ref{item:gepair-9}, together with the fact that $x \in S_\mathcal{R}$ and $d, y \in \mathcal{R}$.
	
	Finally, we note that the change in $\deg_w(c, \tilde \sigma'+\tilde \mu')$ with respect to $\deg_w(c, \tilde \sigma+\tilde \mu)$ depends only (possibly) on the vertices $d$, $y$ and $x$, and we have already checked that the change in $x$ and in $y$ is zero.
	From \eqref{eq:sepI-1} or \eqref{eq:sepI-5}, we get
	\[ \deg_w(c, \tilde \sigma'+\tilde \mu') - \deg_w(c,\tilde \sigma+\tilde \mu) = \gamma\delta>0\]
	where we used $\delta > 0$ in the last step.
	This contradicts the optimality of $(\tilde \sigma, \tilde \mu)$.
\end{proof}

\begin{lemma}[Separating Lemma II] \label{lemma:separatinglemma-2}
	Let $(\tilde \sigma, \tilde \mu)$ be an optimal $(G, w, S, M, c,\mu, \gamma)$-GE pair.
	Suppose there exists $d \in \mathcal{R}$ such that $\tilde \sigma(d) + \tilde \mu(d) < w(\oriented{cd})$.
	
	Let $\sigma_d\le \tilde \sigma$ and $\mu_d\le \tilde \mu$ be the $\gamma$-skew-matching and fractional matching, respectively, obtained from $\tilde \sigma$ and $\tilde \mu$, respectively, by considering only the edges of their support that intersect $N_w(d)$.
	Then there is no edge between the sets \[\{x\in N_w(c)\cap S\::\: \tilde \sigma(x)<w(\oriented{cx})\}\]
	and 
	\[\{y\in N_w(c)\setminus S\::\: \sigma_d(y)+\mu_d(y)>0\}.\]
\end{lemma}

\begin{figure}[ht]
\centering
\begin{subfigure}{0.48\linewidth}
\centering
\begin{adjustbox}{max width=\linewidth}
\begin{tikzpicture}[
  every node/.style={draw, circle, fill=white, minimum size=1.4cm, inner sep=2pt}
]
  \pgfmathsetmacro{\r}{0.7}
  \pgfmathsetmacro{\eps}{0.06} 
  
\pgfmathsetmacro{\Lx}{-2.0}  
\pgfmathsetmacro{\Ly}{0.0}   

\pgfmathsetmacro{\rectL}{\Lx - \r - 0.3}  
\pgfmathsetmacro{\dotsx}{\rectL - 0.6}    
\pgfmathsetmacro{\cx}{-\r*cos(42)}
\pgfmathsetmacro{\cy}{\r*sin(42)}

\coordinate (X45) at ($(S1)+(\cx,\cy)$);      

  \pgfmathsetmacro{\cdx}{\r*cos(153)} 
\pgfmathsetmacro{\cdy}{\r*sin(153)} 
\coordinate (D2left) at ($(D2)+(\cdx,\cdy)$);
  
    \pgfmathsetmacro{\cxx}{\r*cos(325)} 
\pgfmathsetmacro{\cxy}{\r*sin(325)} 
\coordinate (S1Bot) at ($(S1)+(\cxx,\cxy)$);

  \node[draw=none] at (1.9,0) {$S$};
   \node[draw=none] at (-1.2,-0.5) {$x$};
     \node[draw=none] at (0.8,-0.5) {$z$};
  \node[draw=none] at (2.7,3.6) {$d$};
  \node[draw=none] at (0.8,3.6) {$y$};
    \node[draw=none] at (1,2) {$\textcolor{red}{\mu_d}$};
       \node[draw=none] at  (-0.4,2) {$\textcolor{blue}{\sigma_d}$};
        \node[draw=none] at (-2,2) {$\textcolor{red}{\tilde\mu}$};

 \node (D3)[draw=red, thick,  fill=red!22]   at (\Lx,3) {};   
  
  \node (L) [draw=red, thick,  fill=red!22]  at (-2,0) {};   

  \draw [red, thick]($ (L)+(\r,0) $) -- ($ (D3)+(\r,0) $);
    \draw [red, thick]($ (L)+(-\r,0) $) -- ($ (D3)+(-\r,0) $);

  \node  (S1)  at (0,0) {};
  \node (D1)  at (0,3) {};
    
  \draw[red, thick] ($ (S1)+(\r,0) $) -- ($ (D1)+(\r,0) $);
 
\path[fill=blue!22, draw=none]
  (D1) ++(90:\r) arc[start angle=90, end angle=270, radius=\r] -- ($(D1)+(0,-\r)$) -- cycle;
 \path[fill=blue!22, draw=none]
  (S1) ++(90:\r) arc[start angle=90, end angle=270, radius=\r] -- ($(S1)+(0,-\r)$) -- cycle;
  
\path[fill=red!22, draw=none]
  (S1) ++(-90:\r) arc[start angle=-90, end angle=90, radius=\r]
  -- ($(S1)+(0,-\r)$) -- cycle;
\path[fill=red!22, draw=none]
  (D1) ++(-90:\r) arc[start angle=-90, end angle=90, radius=\r]
  -- ($(D1)+(0,-\r)$) -- cycle;

 \draw[blue, thick] (S1) ++(90:\r) arc[start angle=90, end angle=270, radius=\r];
  \draw[blue, thick] (D1) ++(90:\r) arc[start angle=90, end angle=270, radius=\r];
  \draw[red, thick] (D1) ++(-90:\r) arc[start angle=-90, end angle=90, radius=\r];
  \draw[red, thick] (S1) ++(-90:\r) arc[start angle=-90, end angle=90, radius=\r];

  \draw[blue, thick] ($ (D1)+(0,\r) $) -- ($ (D1)+(0,-\r) $);
   \draw[blue, thick] ($ (S1)+(0,\r) $) -- ($ (S1)+(0,-\r) $); 
  \draw[red, thick] ($(D1)+(\eps,\r)$) -- ($(S1)+(\eps,-\r)$);
  \draw [blue, thick]($ (S1)+(-\r,0) $) -- ($ (D1)+(-\r,0) $);

 \draw[blue, thick] ($ (S1)+(0,-\r) $) -- (X45);

\coordinate (LYa) at (tangent cs:node=L,  point={(D1)}, solution=1);
\coordinate (YLa) at (tangent cs:node=D1, point={(L)},  solution=2);
\draw (LYa) -- (YLa);

\coordinate (LYb) at (tangent cs:node=L,  point={(D1)}, solution=2);
\coordinate (YLb) at (tangent cs:node=D1, point={(L)},  solution=1);
\draw (LYb) -- (YLb);

  \node (D2)    at (2,3) {};
 \path[fill=blue!22, draw=none]
  (D2) ++(90:\r) arc[start angle=90, end angle=270, radius=\r]
  -- ($(D2)+(0,-\r)$) -- cycle;
  \draw[blue, thick] (D2) ++(90:\r) arc[start angle=90, end angle=270, radius=\r];
  \draw[blue, thick] ($(D2)+(0,-\r)$) -- ($(D2)+(0,\r)$);
\draw[blue, thick] (X45) --  (D2left);

  \draw[blue, thick] (X45) -- ($ (D1)+(0, -\r) $);
  \draw[blue, thick] ($ (S1)+(0, \r) $) --  ($ (D2)+(0, -\r) $);

  \coordinate (T1b)  at (tangent cs:node=S1,  point={(D2)}, solution=2);
\coordinate (T2b) at (tangent cs:node=D2, point={(S1)},  solution=1);
\draw (T1b)  -- (T2b) ;

\node[draw=none] at (\dotsx,0) {$\ldots$};

\draw[rounded corners]
  (\rectL,0.9) -- (1.4,0.9)
  -- (1.4,-0.9)
  -- (\rectL,-0.9);

\end{tikzpicture}
\end{adjustbox}
\caption{For simplicity, we illustrate here the trivial case, when 
${w(\oriented{cd})}={w(\oriented{cx})}=w(\oriented{cy})=1$,  $\textcolor{red}{\tilde\mu(x)}=1$ thus $ \textcolor{blue}{\tilde \sigma(x)}=0$, and $\textcolor{blue}{\tilde\sigma(d)}+\textcolor{red}{\tilde\mu(d)}=1/2<w(\oriented{cd})$. The vertices $x$ and $y$ are connected by an edge, contradicting the conclusion of \Cref{lemma:separatinglemma-2}.}
\end{subfigure}
\hfill
\begin{subfigure}{0.48\linewidth}
\centering
\begin{adjustbox}{max width=\linewidth}
\begin{tikzpicture}[
  every node/.style={draw, circle, fill=white, minimum size=1.4cm, inner sep=2pt}
]
  \pgfmathsetmacro{\r}{0.7}
  \pgfmathsetmacro{\eps}{0.06} 
  
\pgfmathsetmacro{\Lx}{-2.0}  
\pgfmathsetmacro{\Ly}{0.0}   

\pgfmathsetmacro{\rectL}{\Lx - \r - 0.3}  
\pgfmathsetmacro{\dotsx}{\rectL - 0.6}    
\pgfmathsetmacro{\cx}{-\r*cos(42)}
\pgfmathsetmacro{\cy}{\r*sin(42)}

\coordinate (X45) at ($(S1)+(\cx,\cy)$);      

  \pgfmathsetmacro{\cdx}{\r*cos(153)} 
\pgfmathsetmacro{\cdy}{\r*sin(153)} 
\coordinate (D2left) at ($(D2)+(\cdx,\cdy)$);
  
    \pgfmathsetmacro{\cxx}{\r*cos(325)} 
\pgfmathsetmacro{\cxy}{\r*sin(325)} 
\coordinate (S1Bot) at ($(S1)+(\cxx,\cxy)$);

  \node[draw=none] at (1.9,0) {$S$};
   \node[draw=none] at (-1.2,-0.5) {$x$};
     \node[draw=none] at (0.8,-0.5) {$z$};
  \node[draw=none] at (2.7,3.6) {$d$};
  \node[draw=none] at (0.8,3.6) {$y$};
    \node[draw=none] at (1,1.1) {$\textcolor{red}{\mu_\delta}$};
       \node[draw=none] at  (-0.4,2) {$\textcolor{blue}{\sigma_d}$};
        \node[draw=none] at (-2,2) {$\textcolor{red}{\tilde\mu}$};

 \node (D3)[draw=red, thick,  fill=red!22]   at (\Lx,3) {};   
  
  \node (L) [draw=red, thick,  fill=red!22]  at (-2,0) {};   

  \draw [red, thick]($ (L)+(\r,0) $) -- ($ (D3)+(\r,0) $);
    \draw [red, thick]($ (L)+(-\r,0) $) -- ($ (D3)+(-\r,0) $);

  \node  (S1)  at (0,0) {};
  \node (D1)  at (0,3) {};
    
  \draw ($ (S1)+(\r,0) $) -- ($ (D1)+(\r,0) $);
 
\path[fill=blue!22, draw=none]
  (D1) ++(90:\r) arc[start angle=90, end angle=270, radius=\r] -- ($(D1)+(0,-\r)$) -- cycle;
 \path[fill=blue!22, draw=none]
  (S1) ++(90:\r) arc[start angle=90, end angle=270, radius=\r] -- ($(S1)+(0,-\r)$) -- cycle;
  
\path[fill=red!22, draw=none]
  (S1) ++(-90:\r) arc[start angle=-90, end angle=90, radius=\r]
  -- ($(S1)+(0,-\r)$) -- cycle;

 \draw[blue, thick] (S1) ++(90:\r) arc[start angle=90, end angle=270, radius=\r];
  \draw[blue, thick] (D1) ++(90:\r) arc[start angle=90, end angle=270, radius=\r];
  \draw[red, thick] (S1) ++(-90:\r) arc[start angle=-90, end angle=90, radius=\r];

  \draw[blue, thick] ($ (D1)+(0,\r) $) -- ($ (D1)+(0,-\r) $);
   \draw[blue, thick] ($ (S1)+(0,\r) $) -- ($ (S1)+(0,-\r) $); 
  \draw[red, thick] ($(S1)+(\eps,\r)$) -- ($(S1)+(\eps,-\r)$);

 \draw[blue, thick] ($ (S1)+(0,-\r) $) -- (X45);


\coordinate (LYa) at (tangent cs:node=L,  point={(D1)}, solution=1);
\coordinate (YLa) at (tangent cs:node=D1, point={(L)},  solution=2);
\draw (LYa) -- (YLa);

\coordinate (LYb) at (tangent cs:node=L,  point={(D1)}, solution=2);
\coordinate (YLb) at (tangent cs:node=D1, point={(L)},  solution=1);
\draw (LYb) -- (YLb);
  
 \draw[blue, thick] (X45) -- ($ (D1)+(0, -\r) $);
  \draw [blue, thick]($ (S1)+(-\r,0) $) -- ($ (D1)+(-\r,0) $);
  \node (D2)    at (2,3) {};

 \path[fill=blue!22, draw=none]
  (D2) ++(90:\r) arc[start angle=90, end angle=270, radius=\r]
  -- ($(D2)+(0,-\r)$) -- cycle;
\path[fill=red!22, draw=none]
  (D2) ++(-90:\r) arc[start angle=-90, end angle=90, radius=\r]
  -- ($(D2)+(0,-\r)$) -- cycle;

  \draw[blue, thick] (D2) ++(90:\r) arc[start angle=90, end angle=270, radius=\r];

  \draw[red, thick] (D2) ++(-90:\r) arc[start angle=-90, end angle=90, radius=\r];
 
  \draw[blue, thick] ($(D2)+(0,-\r)$) -- ($(D2)+(0,\r)$);
    \draw[red, thick] ($ (D2)+(\eps,\r) $) -- ($ (D2)+(\eps,-\r) $);

 \draw[blue, thick] (X45) --  (D2left);
  \draw[blue, thick] ($ (S1)+(0, \r) $) --  ($ (D2)+(0, -\r) $);
  \coordinate (T1b)  at (tangent cs:node=S1,  point={(D2)}, solution=2);
\coordinate (T2b) at (tangent cs:node=D2, point={(S1)},  solution=1);
\draw [red, thick](T1b)  -- (T2b) ;
  \draw[red, thick] ($ (S1)+(\eps, \r) $) --  ($ (D2)+(\eps0, -\r) $);

\node[draw=none] at (\dotsx,0) {$\ldots$};

\draw[rounded corners]
  (\rectL,0.9) -- (1.4,0.9)
  -- (1.4,-0.9)
  -- (\rectL,-0.9);

\end{tikzpicture}
\end{adjustbox}
\caption{The parameter $\delta=1/2$ in this situation. This amount of weight in $\mu_d$ is transferred from $\fmat{zy}$ to $\fmat{zd}$. The skew matching $\tilde \sigma$ doesn't change. Now $\textcolor{blue}{\tilde\sigma(y)}+\textcolor{red}{\mu_\delta(y)}=1/2<w(\oriented{cy})$.}
\end{subfigure}

\caption{The new GE pair created in the right picture saturates the neighbourhood of $c$ by the same amount as the one on the left. As $\fmat{xy}$ is an edge the GE pair $(\textcolor{blue}{\tilde \sigma}, \textcolor{red}{\mu_\delta})$ cannot be optimal, and thus neither is $(\textcolor{blue}{\tilde \sigma}, \textcolor{red}{\tilde \mu})$, contradicting the assumption of \Cref{lemma:separatinglemma-2}.}
\label{fig:sep2}
\end{figure}

\begin{proof}
	Arguing by contradiction, we suppose that there is a vertex $x\in N_w(c)\cap S$ with $w(\oriented{cx})>\tilde\sigma(x)$ that is adjacent to a vertex $y\in N_w(c)\setminus S$  with $\sigma_d(y)+\mu_d(y)>0$. 	
	By Separating Lemma I (\Cref{lemma:separatinglemma-1}), we know that $x\not\in N_w(d)$.
	This implies, in particular, that $y \neq d$.
	It also implies that $\mu_d(x)=0$. 
	Also observe that, because $\sigma_d(y)+\mu_d(y)>0$, and the way $\sigma_d$ and $\mu_d$ is defined, we have that $y\in \mathcal R$ by~\ref{item:gepair-5} and~\ref{item:gepair-9}.
	This straightforwardly leads to $x\in S_{\mathcal R}$.
	
	We now split the proof into two cases. \medskip
	
	\noindent \emph{Case~1: $y \in V(\mu_d)$}.
	First consider the case when $y \in V(\mu_d)$, i.e. that $\mu_d(y) > 0$.
	Then there must exist some $z \in N_w(d) \cap S_{\mathcal R}$ with $0 < \mu_d(zy)$, and thus $\tilde \mu(zy) > 0$.
	Let
	\begin{equation} \label{def:sepII-delta-1}
	\delta =  \min \left\lbrace \tilde \mu(zy), w(\oriented{cd}) - \tilde \sigma(d) - \tilde \mu(d) \right\rbrace. 
	\end{equation}
	The first term is non-zero by the choice of $z$; and the second term is non-zero by our assumption on $d$.
	Hence, $\delta > 0$.
	
	We define a skew matching $\sigma_\delta$ and a fractional matching $\mu_\delta$ by setting $\sigma_\delta \equiv \tilde \sigma$ on all edges, and 
	\begin{align*}
	\mu_\delta (zd)&= \tilde \mu(zd)+\delta,\\
		\mu_\delta(zy) & := \tilde \mu(zy) -  \delta,
	\end{align*}
	and $\mu_\delta = \tilde \mu$ on any other edge. See Figure \ref{fig:sep2} for illustration.
	
Observe that the construction also implies that we have not modified the weights of vertices outside of $\{z, d, y\}$. Also observe that
\begin{align}
\label{eq:sepII-1}
\sigma_\delta(z)+\mu_\delta(z)&=\tilde\sigma(z)+\tilde \mu(z)+\delta-\delta=\tilde \sigma(z)+\tilde \mu(z)\overset{\text{\ref{item:gepair-8}}}{=}1,\\
\label{eq:sepII-2}\sigma_\delta(y)+\mu_\delta(y)&=\tilde\sigma(y)+\tilde \mu(y)-\delta\overset{\text{\ref{item:gepair-6}}}{<}w(\oriented{cy})\le 1,\\
\label{eq:sepII-3}\sigma_\delta(d)+\mu_\delta(d)&=\tilde\sigma(d)+\tilde \mu(d)+\delta\overset{\eqref{def:sepII-delta-1}}{\le}w(\oriented{cd})\le 1.
\end{align}	
	
	We will show that $(\sigma_\delta, \mu_\delta)$ is a $(G, w, S, M, c,\mu, \gamma)$-GE pair.
	Properties \ref{item:gepair-1} and \ref{item:gepair-3}--\ref{item:gepair-5} are automatic as $\sigma_\delta\equiv \tilde \sigma$. Property \ref{item:gepair-2} follows from \eqref{eq:sepII-1}--\eqref{eq:sepII-3}. Property \ref{item:gepair-6} follows from \eqref{eq:sepII-2}--\eqref{eq:sepII-3}. As $d,y\in \mathcal R$, we have that \ref{item:gepair-7}--\ref{item:gepair-8} follow from \eqref{eq:sepII-1}. Finally, property \ref{item:gepair-9} follows from the fact that $\fmat{zd}\in G[S_{\mathcal R}, \mathcal R]$.
	
	Now, the fact that there is an edge between $x$ and $y$ and from \eqref{eq:sepII-2}, we deduce from \Cref{lemma:separatinglemma-1} that $(\sigma_\delta, \mu_\delta)$ is not an optimal  $(G, w, S, M, c,\mu, \gamma)$-GE pair. As 
	\begin{equation*}
	\deg_w(c, \sigma_\delta+\mu_\delta)=\deg_w(c, \tilde \sigma+\tilde \mu)+\delta-\delta=\deg_w(c, \tilde \sigma+\tilde \mu),
	\end{equation*}
	we deduce that $(\tilde \sigma, \tilde \mu)$ is not an optimal  $(G, w, S, M, c,\mu, \gamma)$-GE pair, either.
	This is a contradiction.

	\medskip
	
	\noindent \emph{Case~2: $y \notin V(\mu_d)$}.
	Since $\sigma_d(y) + \mu_d(y) > 0$, this implies that $y\in V(\sigma_d)\setminus V(\mu_d)$, i.e., there is a $z\in N_w(d)$ such that $\tilde \sigma(\oriented{zy})>0$ and $\tilde \mu(\fmat{zy})=0$.
	By \ref{item:gepair-3}--\ref{item:gepair-4}, we have that in fact $z \in N_w(c) \cap N_w(d)$.

\begin{figure}[ht]
\centering
\begin{subfigure}{0.48\linewidth}
\centering
\begin{adjustbox}{max width=\linewidth}
\begin{tikzpicture}[
  every node/.style={draw, circle, fill=white, minimum size=1.4cm, inner sep=2pt}
]
  \pgfmathsetmacro{\r}{0.7}
  \pgfmathsetmacro{\eps}{0.06}

  \pgfmathsetmacro{\Lx}{-2.0}
  \pgfmathsetmacro{\Ly}{0.0}
  \pgfmathsetmacro{\rectL}{\Lx - \r - 0.3}
  \pgfmathsetmacro{\dotsx}{\rectL - 0.6}

  \node[draw=none] at (1.9,0) {$S$};
  \node[draw=none] at (-1.2,-0.5) {$x$};
  \node[draw=none] at (0.8,-0.5) {$z$};
  \node[draw=none] at (2.7,3.6) {$d$};
  \node[draw=none] at (0.8,3.6) {$y$};
  \node[draw=none] at (-2.3,2) {$\textcolor{red}{\tilde\mu}$};
  \node[draw=none] at (-0.15,2.1){$\textcolor{blue}{ \sigma_d(\oriented{zy})}$};
  \node[draw=none] at (1.25,2)   {$\textcolor{blue}{ \sigma_d(\oriented{zd})}$};

  \node (L)  [draw=red, thick, fill=red!22] at (\Lx,\Ly) {};
  \node (D3) [draw=red, thick, fill=red!22] at (\Lx,3)   {};

  \node (S1) [draw=blue, thick, fill=blue!18] at (0,0) {};
  \node (D1) [draw=blue, thick, fill=blue!18] at (0,3) {};

  \node (D2) at (2,3) {};

  \pgfmathsetmacro{\cx}{\r*cos(42)}
  \pgfmathsetmacro{\cy}{\r*sin(42)}
  \coordinate (X45)    at ($(S1)+(\cx,\cy)$);
  \coordinate (Xepsbot) at ($(S1)+(\eps,-\r)$);

  \coordinate (L45)     at ($(L)+(\cx,\cy)$);
  \coordinate (Lepsbot) at ($(L)+(\eps,-\r)$);

  \pgfmathsetmacro{\cdx}{\r*cos(153)}
  \pgfmathsetmacro{\cdy}{\r*sin(153)}
  \coordinate (D1left) at ($(D1)+(\cdx,\cdy)$);
  \coordinate (D2left) at ($(D2)+(\cdx,\cdy)$);

  \pgfmathsetmacro{\sxSW}{\r*cos(335)}
  \pgfmathsetmacro{\sySW}{\r*sin(335)}
  \coordinate (D1SW) at ($(D1)+(\sxSW,\sySW)$);

  \draw[red, thick] ($(L)+(\r,0)$) -- ($(D3)+(\r,0)$);
  \draw[red, thick] ($(L)+(-\r,0)$) -- ($(D3)+(-\r,0)$);

  \draw[blue, thick] ($(S1)+(0,\r)$) -- ($(S1)+(0,-\eps/2)$);
  \draw[blue, thick] ($(S1)+(\eps,\r)$) -- ($(S1)+(\eps,0)$);


  \draw[blue, thick] (D1SW) -- ($(S1)+(0,\r)$);
  \draw[blue, thick] ($(S1)+(-\r,0)$) -- ($(D1)+(-\r,0)$);

  \coordinate (T1a) at (tangent cs:node=S1, point={(D2)}, solution=1);
  \coordinate (T2a) at (tangent cs:node=D2, point={(S1)}, solution=2);
  \draw (T1a) -- (T2a);
  
    \coordinate (T1b)  at (tangent cs:node=S1,  point={(D2)}, solution=2);
\coordinate (T2b) at (tangent cs:node=D2, point={(S1)},  solution=1);
\draw (T1b)  -- (T2b) ;

  \draw ($ (S1)+(\r,0) $) -- ($ (D1)+(\r,0) $);

  \node[draw=none] at (\dotsx,0) {$\ldots$};

  \draw[rounded corners]
    (\rectL,0.9) -- (1.4,0.9)
    -- (1.4,-0.9) -- (\rectL,-0.9);

 \path[fill=blue!22, draw=none]
  (D2) ++(90:\r) arc[start angle=90, end angle=270, radius=\r]
  -- ($(D2)+(0,-\r)$) -- cycle;
    
  \draw[blue, thick] (D2) ++(90:\r) arc[start angle=90, end angle=270, radius=\r];
  \draw[blue, thick] ($(D2)+(0,-\r)$) -- ($(D2)+(0,\r)$);
\draw[blue, thick] ($(S1)+(\eps,\r)$) --  (D2left); 
\draw[blue, thick] (S1Bot) --  ($ (D2)+(0, -\r) $);
  
\coordinate (LYa) at (tangent cs:node=L,  point={(D1)}, solution=1);
\coordinate (YLa) at (tangent cs:node=D1, point={(L)},  solution=2);
\draw (LYa) -- (YLa);

\coordinate (LYb) at (tangent cs:node=L,  point={(D1)}, solution=2);
\coordinate (YLb) at (tangent cs:node=D1, point={(L)},  solution=1);
\draw (LYb) -- (YLb);


 \pgfmathsetmacro{\angB}{-60}
\coordinate (S1B) at ($(S1)+(\angB:\r)$);

\pgfmathsetmacro{\angBshift}{-65}
\coordinate (S1Bshift) at ($(S1)+(\angBshift:\r)$);

\draw[blue, thick] ($(S1)+(\eps,0)$) -- (S1B);
\draw[blue, thick] ($(S1)+(0,-\eps/2)$) -- (S1Bshift); 

\end{tikzpicture}
\end{adjustbox}
\caption{The vertices $x$ and $y$ are connected, contradicting the conclusion of \Cref{lemma:separatinglemma-2}. In this picture $\textcolor{blue}{\sigma_d(y)}>0$. However, $y\not\in V(\mu_d)$ (i.e., we are in Case~2)}
\end{subfigure}
\hfill
\begin{subfigure}{0.48\linewidth}
\centering
\begin{adjustbox}{max width=\linewidth}
\begin{tikzpicture}[
  every node/.style={draw, circle, fill=white, minimum size=1.4cm, inner sep=2pt}
]
  \pgfmathsetmacro{\r}{0.7}
  \pgfmathsetmacro{\eps}{0.06}

  \pgfmathsetmacro{\Lx}{-2.0}
  \pgfmathsetmacro{\Ly}{0.0}
  \pgfmathsetmacro{\rectL}{\Lx - \r - 0.3}
  \pgfmathsetmacro{\dotsx}{\rectL - 0.6}

  \node[draw=none] at (1.9,0) {$S$};
  \node[draw=none] at (-1.2,-0.5) {$x$};
  \node[draw=none] at (0.8,-0.5) {$z$};
  \node[draw=none] at (2.7,3.6) {$d$};
  \node[draw=none] at (0.8,3.6) {$y$};
  \node[draw=none] at (-2.3,2) {$\textcolor{red}{\tilde\mu}$};
  \node[draw=none] at (-0.55,2.1){$\textcolor{blue}{ \sigma_\delta(\oriented{zy})}$};
  \node[draw=none] at (1.4,2)   {$\textcolor{blue}{ \sigma_\delta(\oriented{zd})}$};

  \node (L)  [draw=red, thick, fill=red!22] at (\Lx,\Ly) {};
  \node (D3) [draw=red, thick, fill=red!22] at (\Lx,3)   {};

  \node (S1) [draw=blue, thick, fill=blue!18] at (0,0) {};
  \node (D1) at (0,3) {};
 \path[fill=blue!22, draw=none]
  (D1) ++(90:\r) arc[start angle=90, end angle=270, radius=\r]
  -- ($(D1)+(0,-\r)$) -- cycle;
    \draw[blue, thick] ($(D1)+(0,-\r)$) -- ($(D1)+(0,\r)$);
    \draw[blue, thick] (D1) ++(90:\r) arc[start angle=90, end angle=270, radius=\r];

  \node (D2)  [draw=blue, thick, fill=blue!18] at (2,3) {};

  \pgfmathsetmacro{\cx}{\r*cos(42)}
  \pgfmathsetmacro{\cy}{\r*sin(42)}
  \coordinate (X45)    at ($(S1)+(\cx,\cy)$);
  \coordinate (Xepsbot) at ($(S1)+(\eps,-\r)$);

  \coordinate (L45)     at ($(L)+(\cx,\cy)$);
  \coordinate (Lepsbot) at ($(L)+(\eps,-\r)$);

  \pgfmathsetmacro{\cdx}{\r*cos(153)}
  \pgfmathsetmacro{\cdy}{\r*sin(153)}
  \coordinate (D1left) at ($(D1)+(\cdx,\cdy)$);
  \coordinate (D2left) at ($(D2)+(\cdx,\cdy)$);

  \pgfmathsetmacro{\sxSW}{\r*cos(335)}
  \pgfmathsetmacro{\sySW}{\r*sin(335)}
  \coordinate (D1SW) at ($(D1)+(\sxSW,\sySW)$);

  \draw[red, thick] ($(L)+(\r,0)$) -- ($(D3)+(\r,0)$);
  \draw[red, thick] ($(L)+(-\r,0)$) -- ($(D3)+(-\r,0)$);

  \draw[blue, thick] ($(S1)+(0,\r)$) -- ($(S1)+(0,0)$);
  \draw[blue, thick] ($(S1)+(\eps,\r)$) -- ($(S1)+(\eps,-\eps/2)$);


  \draw[blue, thick] ($(D1)+(0,\eps)$) -- ($(S1)+(0,\r)$);
  \draw[blue, thick] ($(S1)+(-\r,0)$) -- ($(D1)+(-\r,0)$);

  \coordinate (T1a) at (tangent cs:node=S1, point={(D2)}, solution=1);
  \coordinate (T2a) at (tangent cs:node=D2, point={(S1)}, solution=2);
  \draw (T1a) -- (T2a);
  
    \coordinate (T1b)  at (tangent cs:node=S1,  point={(D2)}, solution=2);
\coordinate (T2b) at (tangent cs:node=D2, point={(S1)},  solution=1);
\draw [blue, thick] (T1b)  -- (T2b) ;

  \draw ($ (S1)+(\r,0) $) -- ($ (D1)+(\r,0) $);

  \node[draw=none] at (\dotsx,0) {$\ldots$};

  \draw[rounded corners]
    (\rectL,0.9) -- (1.4,0.9)
    -- (1.4,-0.9) -- (\rectL,-0.9);


\draw[blue, thick] ($(S1)+(\eps,\r)$) --  (D2left); 
  
\coordinate (LYa) at (tangent cs:node=L,  point={(D1)}, solution=1);
\coordinate (YLa) at (tangent cs:node=D1, point={(L)},  solution=2);
\draw (LYa) -- (YLa);

\coordinate (LYb) at (tangent cs:node=L,  point={(D1)}, solution=2);
\coordinate (YLb) at (tangent cs:node=D1, point={(L)},  solution=1);
\draw (LYb) -- (YLb);


 \pgfmathsetmacro{\angB}{-120}
\coordinate (S1B) at ($(S1)+(\angB:\r)$);

\pgfmathsetmacro{\angBshift}{-125}
\coordinate (S1Bshift) at ($(S1)+(\angBshift:\r)$);

\draw[blue, thick] ($(S1)+(\eps,-\eps/2)$) -- (S1B);
\draw[blue, thick] ($(S1)+(0,0)$) -- (S1Bshift); 

  \end{tikzpicture}
\end{adjustbox}
\caption{We obtain a new skew-matching $\textcolor{blue}{\sigma_\delta}$ by increasing the weight of $\textcolor{blue}{\sigma_d}$ on the edge $\fmat{zd}$ and decreasing it by the same amount on the edge $\fmat{zy}$. This ``transfer the empty space'' from $d$ to $y$, which is directly connected to $x$.}
\end{subfigure}
\caption{In both pictures above, we have a GE pair with the same saturation of $N_w(c)$. If the fist one is optimal, so is the second one. However, this is a contradiction with \Cref{lemma:separatinglemma-1}.}
\label{fig:sep2'}
\end{figure}
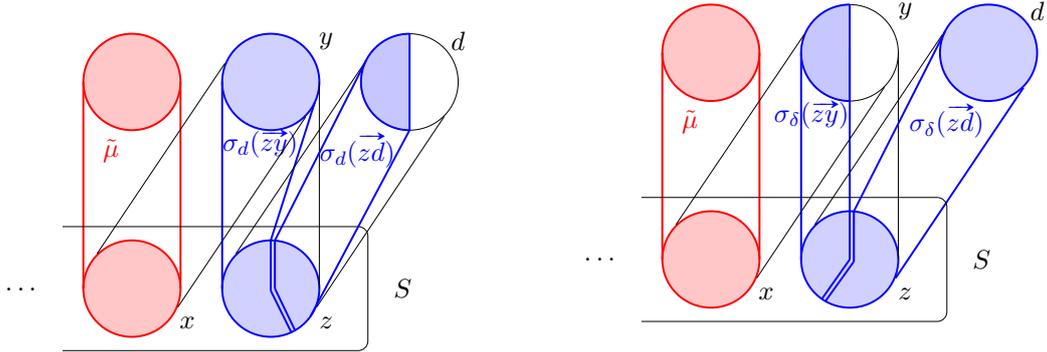
	
	Let
	\begin{equation}\label{def:sepII-delta-2}
	\delta:= \min\left\{\frac{ \sigma(\oriented{zy})}{\gamma}, \frac{w(\oriented{cd})-\tilde \sigma(d)-\tilde \mu(d)}{\gamma}\right\},
	\end{equation}
	From the choice of $z$ and our assumption, we have that $\delta>0$.
	
	Define a skew-matching $\sigma_\delta$ by setting
	\[\sigma_\delta(\oriented{zd}):=\tilde \sigma(\oriented{zd})+ (1+\gamma)\delta,\] 
	\[\sigma_\delta(\oriented{zy}):= \tilde \sigma(\oriented{zy})-(1+\gamma)\delta,\]
	and $\sigma_\delta\equiv\tilde \sigma$ on all other edges.
	Also let $\mu_\delta \equiv \tilde \mu$.
	See \Cref{fig:sep2'} for illustration.
	
	Observe that the construction also implies that we do not modify the weights of vertices outside $\{z, d, y\}$. Also note that we have
	\begin{align}
\label{eq:sepII-1'}	\sigma_\delta(d)+\mu_\delta(d) &= \tilde \sigma(d) +\tilde \mu(d)+ \gamma \delta\overset{\eqref{def:sepII-delta-2}}{\le} w(\oriented{cd})\le 1,\\
	\label{eq:sepII-2'}\sigma_\delta(y)+\mu_\delta(y) &= \tilde \sigma(y) +\tilde \mu(y) - \gamma \delta\overset{\text{\ref{item:gepair-6}}}{<} w(\oriented{cy})\le 1,	\\
	\label{eq:sepII-3'}\sigma_\delta(z)+\mu_\delta(z) &=\tilde \sigma(z) +\tilde \mu(z)+ \delta -  \delta = \tilde \sigma(z) +\tilde \mu(z) \overset{\text{\ref{item:gepair-8}}}{=}1,\\
	\label{eq:sepII-4'}	\sigma_\delta(z)&= \tilde \sigma(z)+ \delta -  \delta =  \tilde \sigma(z)\le w(\oriented{cz})\le 1.
	\end{align}
	
	We will show that $(\sigma_\delta, \mu_\delta)$ is a $(G, w, S, M, c,\mu, \gamma)$-GE pair.
	Properties \ref{item:gepair-1} and \ref{item:gepair-2} follow from \eqref{eq:sepII-1'}--\eqref{eq:sepII-3'}.
	Properties \ref{item:gepair-3} and \ref{item:gepair-5} follow from the definition of $\sigma_\delta$ and that $d\in \mathcal R$.
	Property \ref{item:gepair-4} follows from~\eqref{eq:sepII-4'}.
	Property \ref{item:gepair-6} follows from \eqref{eq:sepII-1'} and \eqref{eq:sepII-2'}.
	Properties \ref{item:gepair-7} and \ref{item:gepair-8} follow from~\eqref{eq:sepII-3'}.
	Finally, \ref{item:gepair-9} is trivial since $\mu_\delta \equiv \tilde \mu$.
		 
	Since $d, y \in N_w(c)$, we have that
	\begin{align}\label{eq:sepII-sameweight}
		\deg_w(c, \sigma_\delta+\tilde \mu)=
		 \deg_w(c, \tilde \sigma+\tilde \mu) +\gamma \delta -\gamma\delta = \deg_w(c, \tilde \sigma+\tilde \mu).
	\end{align}
	By \eqref{eq:sepII-2'}, using that $\fmat{xy}$ is an edge, Separating Lemma I (\Cref{lemma:separatinglemma-1}) implies that $(\sigma_\delta, \mu_\delta)$ is not an optimal$(G, w, S, M, c,\mu, \gamma)$-GE pair; but then \eqref{eq:sepII-sameweight} implies that neither is $(\tilde \sigma, \tilde \mu)$.
	This is contradiction, and it finishes the proof.
\end{proof}

\begin{remark}
After a careful moment of reflection, one can observe that the statements of Separating Lemma I (\Cref{lemma:separatinglemma-1}) and Separating Lemma II (\Cref{lemma:separatinglemma-2}) are about the same underlying property, just one step further in an ``$(\tilde \sigma+\tilde \mu)$-alternating path'' from $d$ in \Cref{lemma:separatinglemma-2}.
Similar properties hold for longer ``alternating paths''.
However, for the sake of simplification, here we choose to focus on just what we actually need in the proof.
\end{remark}

\section{The Matching Lemmas}\label{sec:Matching-lemmas}

To streamline the proof of \Cref{prop:weighted-structural}, we have extracted and formalized some recurring arguments into stand-alone lemmas.
These lemmas not only simplify the current proof but also offer versatile, black-box tools that can be applied in future work involving skew-matching pairs.

This section is organised as follows.
In \Cref{ssec:BML} we introduce basic lemmas that serve as foundational building blocks.
These lemmas are used directly in the proof of \Cref{prop:weighted-structural} or as components of the more intricate `advanced' matching lemmas presented in \Cref{ssec:AML}.
There are four advanced matching lemmas: the Improved Balancing Lemma, the Completion Lemma, the  Greedy Lemma, and the $(k, k/2)$ Lemma.
The proof of the basic lemmas are presented in \Cref{subsection:matching-basic-proofs}, and the next subsections present each a proof of an advanced matching lemma.

\subsection{Basic Matching Lemmas}\label{ssec:BML}
The first lemma determines how large a skew-matching can be fitted inside a given fractional matching.
The set $U$ represents the neighbourhood of a vertex where we want the anchor to fit in.

\begin{lemma}[Extending-out]\label{lem:new-extending-out}
Let $H$ be a graph and $U,V\subseteq V(H)$ be disjoint sets. Suppose there is a fractional matching $\mu$ between $U$ and $V$. Then there is a $\gamma$-skew-matching $\sigma$ in $H^\leftrightarrow$ such that
\stepcounter{propcounter}
\begin{enumerate}[\upshape{(\Alph{propcounter}\arabic*)},topsep=0.7em, itemsep=0.5em]
	\item \label{item:extendingout-1} $\sigma\trianglelefteq \mu$,
	\item \label{item:extendingout-2} $W(\sigma)=(1+\min\{\gamma, \gamma^{-1}\})W(\mu)$, and
	\item \label{item:extendingout-3} $\mathcal A(\sigma) \subseteq U$. 
\end{enumerate} 
\end{lemma}

The second lemma will be used in situations where we have ``access'' to the fractional matching from both sides (i.e., the fractional matching lies within the neighbourhood of a vertex).
In this situation, we can pack the skew-matching more efficiently.

\begin{lemma}[Balancing-out]\label{lem:new-balancing}
Let $H$ be a graph, and let $U\subseteq V(H)$, such that there is a fractional matching $\mu$ in $H[U]$.
Then there is a $\gamma$-skew-matching $\sigma \in H^\leftrightarrow$ such that
\stepcounter{propcounter}
\begin{enumerate}[\upshape{(\Alph{propcounter}\arabic*)},topsep=0.7em, itemsep=0.5em]
	\item \label{item:balancingout-out1} $\sigma\trianglelefteq \mu$,
	\item \label{item:balancingout-out2} $W(\sigma)=2W(\mu)$, and
	\item \label{item:balancingout-out3} $\mathcal A(\sigma)\subseteq U$. 
\end{enumerate} 
\end{lemma}

The next lemma is a more specific version of the previous two, stated directly for weighted graphs.
It will provide a lower bound on the weight of a skew-matching that can be built while respecting the saturation of the neighbourhood of a vertex by a fractional matching $\mu$.

\begin{lemma}[Combination]
	\label{lem:new-combination}
	Let $(H,w)$ be a weighted graph, $v\in V(H)$, $\mu$ a fractional matching in $H$ and $\gamma \geq 0$.
	Then there is a $\gamma$-skew-matching $\sigma$ in $H^\leftrightarrow$ such that
	\stepcounter{propcounter}
	\begin{enumerate}[\emph{(\Alph{propcounter}\arabic*)},topsep=0.7em, itemsep=0.5em]
		\item $\sigma\trianglelefteq\mu$,
		\item $W(\sigma)\ge \deg_w(v, \mu)$, and
		\item $\mathcal A(\sigma)$ fits in the $w$-neighbourhood of $v$. 
	\end{enumerate} 
\end{lemma} 

In contrast with the previous lemmas, the next lemma will be used to build a fractional matching from the existence of a skew-matching.
This lemma is presented in a more general form, where one $\gamma_A$-skew-matching is found inside another $\gamma_B$-skew-matching.
Applying this with $\gamma_A = 1$, this can be combined with \Cref{lemma:fractionalfrom1skew} to obtain fractional matchings.

\begin{lemma}[Extending-out skew-matching]\label{prop:weighted-extending-skew}
	Let $(H,w)$ be a weighted graph, $\gamma_A > 0$ and $\gamma_B\ge 1$.
	Let $u\in V(H)$, and $\sigma_B$ be a $\gamma_B$-skew-matching in $H^\leftrightarrow$ such that $\mathcal A(\sigma_B)$ fits in the $w$-neighbourhood of $u$.
	Then there is a $\gamma_A$-skew-matching in $H^\leftrightarrow$ such that
	\stepcounter{propcounter}
	\begin{enumerate}[\upshape{(\Alph{propcounter}\arabic*)},topsep=0.7em, itemsep=0.5em]
		\item \label{item:extending-out-skew-1} $\sigma_A\le \sigma_B$,
		\item \label{item:extending-out-skew-2} $W(\sigma_A)=\frac{1+\min\{\gamma_A, \gamma_A^{-1}\}}{1+\gamma_B}W(\sigma_B)$, and
		\item \label{item:extending-out-skew-3} $\mathcal A(\sigma_A)$ fits in the $w$-neighbourhood of $u$.
	\end{enumerate}
\end{lemma}

\subsection{Advanced Matching Lemmas}\label{ssec:AML}
The first of the more complex lemmas  gives a condition under which we can completely fill up a fractional matching with a skew-matching pair.
 
\begin{lemma}[Improved balancing]\label{lem:new-improved}
Let $H$ be a graph and $U,V\subseteq V(H)$ be disjoint sets. Suppose there is a fractional matching $\mu$ running between $U$ and $V$. Let $\alpha_1,  \beta_1>0$ and $\alpha_2, \beta_2\ge 0$ be such that 
\stepcounter{propcounter}
\begin{enumerate}[\upshape{(\Alph{propcounter}\arabic*)},topsep=0.7em, itemsep=0.5em]
\item \label{itnew:fills}$\alpha_1+\alpha_2+\beta_1+\beta_2=2W(\mu)$,
\item \label{itnew:fits}$\max\{\alpha_1,\alpha_2\}+\min\{\beta_1, \beta_2\}\le W(\mu)$.
\end{enumerate}
Set $\gamma_A:= \alpha_2 / \alpha_1$ and $\gamma_B:= \beta_2 / \beta_1$. 
Then there is a $\gamma_A$-skew-matching $\sigma_A$ and a $\gamma_B$-skew-matching $\sigma_B$ in $H^\leftrightarrow$, such that
\begin{enumerate}[\upshape{(\Alph{propcounter}\arabic*)},resume,topsep=0.7em, itemsep=0.5em]
	\item \label{item:improvedbalancing-out1} $\sigma_A+\sigma_B\trianglelefteq \mu$,
	\item \label{item:improvedbalancing-out2} $W(\sigma_A)=\alpha_1+\alpha_2$, $W(\sigma_B)=\beta_1+\beta_2$,
	\item \label{item:improvedbalancing-out3} $\mathcal A(\sigma_A)\subseteq U$, and $\mathcal A(\sigma_B)\subseteq U\cup V$. 
\end{enumerate} 
\end{lemma}

The following lemma is the most complex one. As in the previous lemma, we combine two different skew-matchings and place them disjointly within a fractional matching. However, in this case, we impose significantly weaker conditions for the graph to satisfy, leading to a much more challenging scenario.

\begin{lemma}[Completion]\label{lem:new-completion}
	Let $(H, w)$ be a weighted graph, $U, V\subseteq V(H)$ be disjoint sets, $u\in V(H)$,  and $\mu$ be a fractional matching running between $U$ and $V$.
	 Let $\alpha_1, \beta_1>0$ and $\alpha_2, \beta_2\ge 0$ be such that 
	\stepcounter{propcounter}
	\begin{enumerate}[\upshape{(\Alph{propcounter}\arabic*)},topsep=0.7em, itemsep=0.5em]
		\item \label{item:completion-1} for all $y\in V$, we have $\mu(y)\leq w(\oriented{uy})$,
		\item \label{eq:neq-fit}
		$\max\{\alpha_1, \alpha_2\}+\min\{\beta_1, \beta_2\}\le W(\mu)$, and
		\item \label{eq:VybalancovaniMozne}
		$\min\{\alpha_1,\alpha_2\}+\max\{\beta_1,\beta_2\}\ge W(\mu)$.
	\end{enumerate}
	Set $\gamma_A:= \alpha_2 / \alpha_1$ and $\gamma_B:= \beta_2 / \beta_1$.
	Then $H^\leftrightarrow$ admits a $\gamma_A$-skew-matching $\sigma_A$ and a $\gamma_B$-skew-matching $\sigma_B$ such that
	\begin{enumerate}[\upshape{(\Alph{propcounter}\arabic*)},resume,topsep=0.7em, itemsep=0.5em]
		\item \label{item:completion-out1} $\sigma_A+\sigma_B\trianglelefteq\mu$,
		\item \label{item:completion-out2} $W(\sigma_A)=\alpha_1+\alpha_2$,
		\item \label{item:completion-out3} $W(\sigma_B)\ge \max\{0, \deg_w(u, \mu)-W(\sigma_A)\}$,
		\item \label{item:completion-out4} $\mathcal A(\sigma_A)\subseteq U$,
		\item \label{item:completion-out5} $\mathcal A(\sigma_B)$ fits in the $w$-neighbourhood of $u$, and moreover,
		\item \label{item:completion-moreover} if $\gamma_B\le 1$, we can also ensure that $\sigma_A(y)+\sigma_B(y)=\mu(y)$ for all $y\in V$.
	\end{enumerate}
\end{lemma}

During the proof of \Cref{prop:weighted-structural}, we repeatedly encounter situations where we rely on a minimum degree condition to place our skew-matching.
We present three lemmas, which have similar, but slightly different, scenarios and outcomes.
We recall the notation $\deg_w(v, S):= \sum_{u\in S} w(\oriented{vu})$ for $v \in V(H)$ and $S \subseteq V(H)$.
%
%
%
\begin{lemma}[First Greedy Lemma] \label{prop:weighted-greedy-1}
		Let $(H,w)$ be a weighted graph, $u, v \in V(H)$, $\kappa \geq 0$, $(\sigma_A,\sigma_B)$ be a $(\gamma_A,\gamma_B)$-skew-matching pair in $H^\leftrightarrow$ anchored in $\oriented{uv}\in E(H^\leftrightarrow)$, and let $U, V\subseteq V(H)$ be disjoint sets.
		Suppose that
		\stepcounter{propcounter}
		\begin{enumerate}[\upshape{(\Alph{propcounter}\arabic*)},topsep=0.7em, itemsep=0.5em]
			\item \label{item:greedy-1-condition-1}	$\deg_w(u, V)\ge  \kappa + \sum_{x\in V}(\sigma_A (x)+\sigma_B(x))$, and
			\item \label{item:greedy-1-condition-2} $|N_w(x)\cap U|\ge \gamma_A \kappa+\sum_{y\in U}(\sigma_A (y)+\sigma_B(y))$, for every $x\in V$.
		\end{enumerate}
		Then there is a $\gamma_A$-skew oriented fractional matching $\sigma_A'$ in $H^\leftrightarrow$ with
		\begin{enumerate}[\upshape{(\Alph{propcounter}\arabic*)},resume,topsep=0.7em, itemsep=0.5em]
			\item \label{item:greedy-1-out-1} $\mathcal A(\sigma_A')\subseteq V$,
			\item \label{item:greedy-1-out-2} $W(\sigma_A')\ge (1+\gamma_A)\kappa$,
			\item \label{item:greedy-1-out-3} $\sigma'_A$ is supported in $H[V, U]$, and such that
			\item \label{item:greedy-1-out-4} $(\sigma_A + \sigma'_A, \sigma_B)$ is a $(\gamma_A, \gamma_B)$-skew-matching anchored in $\oriented{uv}$.
		\end{enumerate}
\end{lemma}

\begin{lemma}[Second Greedy Lemma] \label{prop:weighted-greedy-2}
	Let $(H,w)$ be a weighted graph, $u, v \in V(H)$, $\kappa \geq 0$, $(\sigma_A,\sigma_B)$ be a $(\gamma_A,\gamma_B)$-skew-matching pair in $H^\leftrightarrow$ anchored in $\oriented{uv}\in E(H^\leftrightarrow)$, and let $U, V\subseteq V(H)$ be disjoint sets.
	Suppose that
	\stepcounter{propcounter}
	\begin{enumerate}[\upshape{(\Alph{propcounter}\arabic*)},topsep=0.7em, itemsep=0.5em]
		\item \label{item:greedy-2-condition-1} $\deg_w(u, V)\ge  \kappa + \sum_{x\in V}(\sigma_A(x)+\sigma_B(x))$,
		\item \label{item:greedy-2-condition-2} 
		$|N_w(x)\cap (U\cup V)|\ge (1+\gamma_A)\kappa+\sum_{y\in U\cup V}(\sigma_A(y)+\sigma_B(y))$, for every $x\in V$, 
	\end{enumerate}
	Then there is a $\gamma_A$-skew oriented fractional matching $\sigma_A'$ in $H^\leftrightarrow$ with
	\begin{enumerate}[\upshape{(\Alph{propcounter}\arabic*)},resume,topsep=0.7em, itemsep=0.5em]
		\item \label{item:greedy-2-out-1} $\mathcal A(\sigma_A')\subseteq V$,
		\item \label{item:greedy-2-out-2} $W(\sigma_A')\ge (1+\gamma_A)\kappa$,
		\item \label{item:greedy-2-out-3} $\sigma'_A$ is supported in $\{ xy \in E(H) : x \in V, y \in V \cup U \}$, and such that
		\item \label{item:greedy-2-out-4} $(\sigma_A + \sigma'_A, \sigma_B)$ is a $(\gamma_A, \gamma_B)$-skew-matching anchored in $\oriented{uv}$.
	\end{enumerate}
\end{lemma}

\begin{lemma}[Third Greedy Lemma] \label{prop:weighted-greedy-3}
	Let $(H,w)$ be a weighted graph, $u, v \in V(H)$, $\kappa \geq 0$, $(\sigma_A,\sigma_B)$ be a $(\gamma_A,\gamma_B)$-skew-matching pair in $H^\leftrightarrow$ anchored in $\oriented{uv}\in E(H^\leftrightarrow)$, and let $U, V\subseteq V(H)$ be disjoint sets.
	Suppose that
	\stepcounter{propcounter}
	\begin{enumerate}[\upshape{(\Alph{propcounter}\arabic*)},topsep=0.7em, itemsep=0.5em]
		\item \label{item:greedy-3-condition-1} $|U|\ge \gamma_A \kappa+\sum_{y\in U}(\sigma_A (y)+\sigma_B(y))$, and
		\item \label{item:greedy-3-condition-2} $ \deg_w(u, N_w(y)\cap V)\ge \kappa + \sum_{x\in V}(\sigma_A (x)+\sigma_B(x))$, for every $y\in U$.
	\end{enumerate}
	Then there is a $\gamma_A$-skew oriented fractional matching $\sigma_A'$ in $H^\leftrightarrow$ with
	\begin{enumerate}[\upshape{(\Alph{propcounter}\arabic*)},resume,topsep=0.7em, itemsep=0.5em]
		\item \label{item:greedy-3-out-1} $\mathcal A(\sigma_A')\subseteq V$,
		\item \label{item:greedy-3-out-2} $W(\sigma_A')\ge (1+\gamma_A)\kappa$,
		\item \label{item:greedy-3-out-3} $\sigma'_A$ is supported in $H[V,U]$, and such that
		\item \label{item:greedy-3-out-4} $(\sigma_A + \sigma'_A, \sigma_B)$ is a $(\gamma_A, \gamma_B)$-skew-matching anchored in $\oriented{uv}$.
	\end{enumerate}
\end{lemma}

The following proposition will allow us to conclude if there is an edge $cd$ and a skew-matching which saturates sufficiently each of the neighbourhoods of $c$ and $d$.
This lemma (in a less general version, using only matchings) appeared in the work of Ajtai, Koml\'os, Simonovits, and Szemer\'edi~\cite{AKSS2015}; we adapt their strategy to our situation.

\begin{lemma}[The $(k, k/2)$-lemma]
	\label{lem:new-k-k/2}Let $k\ge 2$. 
	Let $(H,w)$ be a weighted graph, $\fmat{cd}\in E(H)$, and  $\mu$ a fractional matching in $H$, such that
	\stepcounter{propcounter}
	\begin{enumerate}[\upshape{(\Alph{propcounter}\arabic*)},topsep=0.7em, itemsep=0.5em]
		\item $\deg_w(c,\mu)\ge k$, and
		\item $\deg_w(d,\mu)\ge k/2$.
	\end{enumerate}
	Let $\alpha_1, \beta_1>0$ and $\alpha_2, \beta_2\ge 0$ be such that $\alpha_1+\alpha_2+\beta_1+\beta_2=k$.
	Set $\gamma_A:= \alpha_2 / \alpha_1$ and $\gamma_B:= \beta_2 / \beta_1$. 
	Then $H^\leftrightarrow$ admits a $(\gamma_A,\gamma_B)$-skew-matching pair $(\sigma_A,\sigma_B)$, anchored in $\oriented{cd}$ or in $\oriented{dc}$, such that
	\begin{enumerate}[\upshape{(\Alph{propcounter}\arabic*)}, resume,topsep=0.7em, itemsep=0.5em]
		\item $W(\sigma_A)=\alpha_1+\alpha_2$, 
		\item $W(\sigma_B)=\beta_1+\beta_2$, and
		\item $\sigma_A+\sigma_B\trianglelefteq \mu$.
	\end{enumerate}
\end{lemma}

\subsection{The proofs of the Basic Matching Lemmas} \label{subsection:matching-basic-proofs}

We give the proof of the Basic Matching Lemmas (Lemmas \ref{lem:new-extending-out}--\ref{prop:weighted-extending-skew}).

\begin{proof}[Proof of \Cref{lem:new-extending-out} (Extending-out)]
For every $xy$ with $x\in U$ and $y\in V$, define \[\sigma(\oriented{xy}):= (1+\gamma)\frac{\mu(\fmat{xy})}{\max\{1, \gamma\}},\] and set $\sigma$ to $0$ on all other edges.
It is straightforward to verify that $\sigma$ is a $\gamma$-skew fractional matching.
Indeed, for $x\in U$ we have $\sum_{y}\bigl(\sigma(\overrightarrow{xy})/(1+\gamma)\bigr)=\sum_y \mu(xy)/\max\{1,\gamma\}\le \mu(x)\le 1$,
and for $y\in V$ we have $\sum_{x}\bigl(\gamma\,\sigma(\overrightarrow{xy})/(1+\gamma)\bigr)=\gamma\sum_x \mu(xy)/\max\{1,\gamma\}\le \mu(y)\le 1$.

Now we verify the required properties.
Since the only edges $\oriented{xy}$ with non-zero weight have $x \in U$, this readily implies that $\mathcal A(\sigma)\subseteq U$, which gives \ref{item:extendingout-3}.
Also, we have 
\begin{align*}
W(\sigma)=& \sum_{x\in U,y\in V}\sigma(\oriented{xy})=\frac{1+\gamma}{\max\{1,\gamma \}}\sum_{x\in U, y\in V}\mu(\fmat{xy})= (1+\min\{\gamma, \gamma^{-1}\})W(\mu),\end{align*}
which gives \ref{item:extendingout-2}.
Finally, for $x\in U$ and $y\in V$, we have 
\begin{align*}
\frac{\sigma(\oriented{xy})+\gamma\sigma(\oriented{yx})}{1+\gamma}= \frac{\sigma(\oriented{xy})}{1+\gamma}=\frac{\mu(\fmat{xy})}{\max\{1, \gamma\}}\le \mu(\fmat{xy}),
\end{align*}and 
\begin{align*}
\frac{\sigma(\oriented{yx})+\gamma\sigma(\oriented{xy})}{1+\gamma}=
\frac{\gamma\sigma(\oriented{xy})}{1+\gamma}=\frac{\gamma\mu(\fmat{xy})}{\max\{1, \gamma\}}\le \mu(\fmat{xy}).
\end{align*}
This means that  $\sigma\trianglelefteq \mu$, so \ref{item:extendingout-1} holds.
\end{proof}

\begin{proof}[Proof of \Cref{lem:new-balancing} (Balancing-out)]
	For every $\fmat{xy}$ with $x,y\in U$, we define $\sigma(\oriented{xy})=\sigma(\oriented{yx}):= \mu(\fmat{xy})$; and we set $\sigma$ to $0$ on all other edges.
	We have 
	\begin{align*}
	W(\sigma)=\sum_{\{\fmat{xy} \in E(G) \::\: x,y\in U\}}\sigma(\oriented{xy})+\sigma(\oriented{yx})=\sum_{\{\fmat{xy} \in E(G) \::\: x,y\in U\}} 2\mu(\fmat{xy})=2W(\mu).
	\end{align*}
	Therefore, for $x,y\in U$  we obtain
	\begin{align*}
	\frac{\sigma(\oriented{xy})+\gamma\sigma(\oriented{yx})}{1+\gamma}= \frac{(1+\gamma)\mu(\fmat{xy})}{1+\gamma}=\mu(\fmat{xy}),
	\end{align*}and 
	\begin{align*}
	\frac{\sigma(\oriented{yx})+\gamma\sigma(\oriented{xy})}{1+\gamma}=\frac{(1+\gamma)\mu(\fmat{xy})}{1+\gamma}=\mu(\fmat{xy}).
	\end{align*}
	Hence, $\sigma\trianglelefteq \mu$.
	Finally, $\mathcal A(\sigma)\subseteq U$ holds by construction (because no directed edge with tail outside $U$ receives weight).
	We thus have \ref{item:balancingout-out1}--\ref{item:balancingout-out3}.
\end{proof}

%

\begin{proof}[Proof of \Cref{lem:new-combination} (Combination)]
	Let $U:=V(\mu)$, and let $\tilde \mu\le \mu$ be a maximal fractional matching so that $\tilde \mu(x)\le w(\oriented{vx})$ for all $x\in U$. 
	By \Cref{lem:new-balancing} there is a $\gamma$-skew-matching $\tilde \sigma\trianglelefteq \tilde \mu$ of weight $W(\tilde \sigma)=2W(\tilde \mu)$ with its anchor $\mathcal A(\tilde\sigma)$ contained in $U$.
	By the choice of $\tilde \mu$ and $\tilde \sigma\trianglelefteq \tilde \mu$, we have that the anchor $\mathcal A(\tilde \sigma)$ fits in the $w$-neighbourhood of $v$.

	Let $(H,w')$ be the $\tilde \mu$-truncated weighted graph obtained from $(H,w)$. 
	We define $U':= N_{w'}(v) = \{x\in V(H)\::\: w'(\oriented{vx})>0\}$ and $V':= U\setminus U'$.
	Let $\mu'\le  \mu-\tilde \mu$ be maximal so that $\mu'(x)\le w'(\oriented{vx})$ for all $x\in U'$ and has only support edges intersecting $U'$.
	
	Observe that we have $w'(\oriented{vx})=0$ or $w'(\oriented{vy})=0$ for every $xy$ with $\mu'(\fmat{xy})>0$.
	Indeed, if not, then we could increase $\tilde\mu(\fmat{xy})$ by a small $\varepsilon>0$, staying within the edge-wise bounds $\tilde\mu\le \mu$ and the vertex bounds $\tilde\mu(x)\le w(\oriented{vx})$ and $\tilde\mu(y)\leq w(\oriented{vy})$, contradicting the maximality of $\tilde\mu$.
	Hence, $\mu'$ runs between $U'$ and $V'$.

	Because of this, we have $W(\mu') = \sum_{x \in U'} \mu'(x)$, and since $\mu'(x) \leq w'(\oriented{vx})$ for all $x \in U'$ then we also have $W(\mu')=\deg_{w'}(v, \mu')$. 

	We apply \Cref{lem:new-extending-out} with $U', V', H, \mu'$ playing the roles of $U, V, H, \mu$, respectively.
	By doing so, we obtain a $\gamma$-skew-matching $\sigma' \in H^\leftrightarrow$ such that $\sigma' \trianglelefteq\mu'$, with weight $W(\sigma')=(1+\min\{\gamma, \gamma^{-1}\})W(\mu')\ge \deg_{w'}(v, \mu')$, and such that its anchor $\mathcal A(\sigma')$ is contained in $U'$.
	By the definition of $U'$ and $\mu'$, this implies that $\mathcal A(\sigma')$ fits in the $w'$-neighbourhood of $v$. 
	
	We apply \Cref{prop:adding-skew-matchings} (with $(H,w), \fmat{vx}, \tilde \mu, \mu', w',\tilde\sigma, \sigma_\emptyset, \sigma'$, and $\sigma_\emptyset$ , with an arbitrary $x$, playing the role of $(G,w),\fmat{uv}, \mu, \bar\mu, \sigma_A, \sigma_B, \bar\sigma_A$, and $\bar\sigma_B$ respectively) to consider the sum of $\tilde \sigma$ and $\sigma'$ 
	Let $\sigma := \tilde \sigma+\sigma'$.
	Thanks to the definition of~$\tilde \mu$, we have $\sigma \trianglelefteq\mu$.
	Moreover, $\sigma$ has weight 
	$W(\sigma)\ge \deg_w(v,\tilde\mu)+\deg_{w'}(v,\mu-\tilde\mu)=\deg_w(v,\mu)$, where we used $\deg_w(v,\tilde\mu)=2W(\tilde\mu)$ and \Cref{prop:adding-degree-truncated} with $H, w, \tilde \mu, \mu$, and $w'$ paying the role of $G, w, \mu', \mu$ and $w'$, respectively. 
	Also, its anchor $\mathcal A(\sigma)$ fits in the $w$-neighbourhood of $v$.
	This terminates the proof of \Cref{lem:new-combination}.
\end{proof}

\begin{proof}[Proof of \Cref{prop:weighted-extending-skew} (Extending-out Skew-matching)]
	Define a $\gamma_A$-skew-matching $\sigma_A$ so that for each 
	$\oriented{xy}\in E(H^\leftrightarrow)$ with $\sigma_B(\oriented{xy})>0$, 
	we have \[\sigma_A(\oriented{xy}):= \frac{1+\min\{\gamma_A,\gamma_A^{-1}\}}{1+\gamma_B}\sigma_B(\oriented{xy}),\] and $\sigma_A(\oriented{xy}):=0$ for all other $\oriented{xy} \in E(H^\leftrightarrow)$.
	Note that since $\gamma_B \geq 1$, we have $1 + \min\{\gamma_A, \gamma_{A}^{-1}\} \leq 2 \leq 1 + \gamma_B$, which ensures this indeed defines a $\gamma_A$-skew-matching.
	
	Then we have
	 \[W(\sigma_A)=\frac{1+\min\{\gamma_A,\gamma_A^{-1}\}}{1+\gamma_B}W(\sigma_B),\] which gives \ref{item:extending-out-skew-2}.
	 For every 
	$\oriented{xy}\in E(H^\leftrightarrow)$ we have 
	 \[\frac{\sigma_A(\oriented{xy})}{1+\gamma_A}= \frac{1+\min\{\gamma_A,\gamma_A^{-1}\}}{1+\gamma_A}\cdot  \frac{\sigma_B(\oriented{xy})}{1+\gamma_B}\le \frac{\sigma_B(\oriented{xy})}{1+\gamma_B}. \] 
	 Thus, $\sigma^1_A(x)\le \sigma_B^1(x)$ for all $x\in V(H)$. 
	 Further,  for all $\oriented{xy}\in E(H^\leftrightarrow)$, we have
	 \begin{align*}
		\frac{\gamma_A\sigma_A(\oriented{xy})}{1+\gamma_A}=& \frac{\gamma_A+\min\{(\gamma_A)^2, 1\}}{1+\gamma_A}\cdot \frac{\sigma_B(\oriented{xy})}{1+\gamma_B} \le \frac{\gamma_B\sigma_B(\oriented{xy})}{1+\gamma_B},
	\end{align*}
	since $\frac{\gamma_A+\min\{\gamma_A^2,1\}}{1+\gamma_A}\le 1$ (it equals $1$ if $\gamma_A\ge 1$ and equals $\gamma_A$ if $\gamma_A\le 1$), and $\gamma_B\ge 1$.
	This means that
	$\sigma_A\le \sigma_B$, so \ref{item:extending-out-skew-1} holds.
	Finally, from $\sigma^1_A(x)\le \sigma_B^1(x)$ we see that $\mathcal A(\sigma_A)$ fits in the $w$-neighbourhood of $u$, because $\mathcal A(\sigma_B)$ also fits in the $w$-neighbourhood of $u$, so we have \ref{item:extending-out-skew-3}.
	\end{proof}

\subsection{Proof of the Improved Balancing Lemma}

We shall need the following auxiliary `allocation' lemma.

\begin{lemma}\label{lem:weighted-filling-matching}
	Let $\alpha_1, \alpha_2, \beta_1, \beta_2, \gamma \geq 0$ such that
	\stepcounter{propcounter}
	\begin{enumerate}[\upshape{(\Alph{propcounter}\arabic*)},topsep=0.7em, itemsep=0.5em]
		\item \label{item:aloc-in1} $\alpha_1+\alpha_2+\beta_1+\beta_2\le 2\gamma$, and
		\item \label{item:aloc-in2} $\min\{\beta_1, \beta_2\}+\max\{\alpha_1,\alpha_2 \}\le \gamma$,
	\end{enumerate} then there exist $\hat \beta_1\le \beta_1$ and $\hat \beta_2\le \beta_2$ with $\hat \beta_1\cdot \beta_2=\hat \beta_2\cdot \beta_1$ such that
	\begin{enumerate}[\upshape{(\Alph{propcounter}\arabic*)}, resume, topsep=0.7em, itemsep=0.5em]
		\item \label{item:aloc-out2} $\hat\beta_1 +(\beta_2-\hat\beta_2)+\alpha_1\le \gamma$ and
		\item \label{item:aloc-out1}  $\hat\beta_2 +(\beta_1-\hat\beta_1)+\alpha_2\le \gamma$.
	\end{enumerate}
%
\end{lemma}

\begin{proof}
	Without loss of generality we can assume that $\beta_2\ge\beta_1$.
	We will also assume that $\alpha_2\ge \alpha_1$, and we explain how to remove this assumption later.
	Those two assumptions imply that $\beta_1+\alpha_2\le \beta_1+\alpha_1\le \gamma$.
	If $\beta_2+\alpha_2\le \gamma$, then we set $\hat\beta_1=\beta_1$ and $\hat\beta_2=\beta_2$, from which it is straightforward to verify \ref{item:aloc-out2}--\ref{item:aloc-out1}.
	Hence, we may assume that $\beta_2 + \alpha_2 >\gamma$.
	Together with $\gamma \geq \alpha_2 + \beta_1$ we get \[\beta_2>\gamma-\alpha_2 \geq \beta_1.\]
	This implies that there exists $\lambda\in [0,1]$ such that \[\lambda \beta_2+(1-\lambda)\beta_1=\gamma-\alpha_2.\]
	Set $\hat \beta_2:= \lambda \beta_2$ and $\hat \beta_1=\lambda \beta_1$.
	This choice gives $\hat \beta_2+(\beta_1-\hat\beta_1)+\alpha_2=\gamma$, so \ref{item:aloc-out1} holds.
	To see \ref{item:aloc-out2}, we note that
	\begin{align*}
		\hat \beta_1 + (\beta_2 - \hat \beta_2) + \alpha_1
		& = \alpha_1 - (\lambda \beta_2 + (1 - \lambda)\beta_1) + \beta_2 + \beta_1 \\
		& = \alpha_1 + \alpha_2 + \beta_2 + \beta_1 - \gamma \leq \gamma,
	\end{align*}
	where in the last inequality we used~\ref{item:aloc-in1}.
	Finally, note that by defining instead $\hat{\beta}_2 = (1 - \lambda) \beta_2$ and $\hat{\beta}_1 = (1 - \lambda) \beta_1$ we get \ref{item:aloc-out2}--\ref{item:aloc-out1} with the roles of $\alpha_1, \alpha_2$ swapped, which removes the assumption on $\alpha_1, \alpha_2$.
\end{proof}

\begin{proof}[Proof of \Cref{lem:new-improved} (Improved balancing)]
	Our assumptions \ref{itnew:fills} and \ref{itnew:fits} imply that we can apply \Cref{lem:weighted-filling-matching} with $\alpha_1, \alpha_2, \beta_1, \beta_2, W(\mu)$ playing the roles of $\alpha_1, \alpha_2, \beta_1, \beta_2, \gamma$.
	We obtain $\hat\beta_1\le \beta_1$ and $\hat \beta_2 \leq \beta_2$ such that $\hat\beta_2 / \hat\beta_1 = \beta_2 / \beta_1 = \gamma_B$ and such that
	\begin{equation}\label{eq:lem-improved-hatbeta_1+}
		\hat\beta_1+(\beta_2-\hat\beta_2)+\alpha_1 \leq W(\mu),
		\end{equation} and 
		\begin{equation}\label{eq:lem-improved-hatbeta_2+}
		\hat\beta_2+(\beta_1-\hat\beta_1)+\alpha_2 \leq W(\mu).
	\end{equation}
	
	For every $x\in U, y\in V$,  we set 
	\begin{align*}
			\sigma_A(\oriented{xy}):= \frac{\alpha_1+\alpha_2}{W(\mu)} \mu(\fmat{xy}),
			\qquad 
		\sigma_B(\oriented{xy}):= \frac{\hat\beta_1+\hat\beta_2}{W(\mu)} \mu(\fmat{xy}),
	\end{align*}and 
\begin{align*}
	\sigma_B(\oriented{yx}):= \frac{(\beta_1+\beta_2)-(\hat\beta_1+\hat\beta_2)}{W(\mu)}\mu(\fmat{xy}),
\end{align*} and we put $\sigma_A$ and $\sigma_B$ equal to $0$ on all other edges.
Note that $\mathcal A(\sigma_A)\subseteq U$ and $\mathcal A(\sigma_B)\subseteq U\cup V$, so \ref{item:improvedbalancing-out3} holds.
The weights are
	\begin{align*}
		W(\sigma_A)=&\sum_{x\in U, y\in V}\sigma_A(\oriented{xy})=\frac{\alpha_1+\alpha_2}{W(\mu)}\sum_{x\in U, y\in V}\mu(\fmat{xy})
		= \alpha_1+\alpha_2
	\end{align*}
	and
	\begin{align*}
		W(\sigma_B)& =\sum_{x\in U, y\in V}\Big(\sigma_B(\oriented{xy})+\sigma_B(\oriented{yx})\Big)=\frac{\beta_1+\beta_2}{W(\mu)}\sum_{x\in U, y\in V}\mu(\fmat{xy}) = \beta_1+\beta_2,
		\end{align*}
		so \ref{item:improvedbalancing-out2} holds.
		Finally, we observe that $\hat{\beta}_1 = t \beta_1$ and $\hat{\beta}_2 = t \beta_2$, for some $t \in [0,1]$.
		From this, it follows that for every $x\in U, y\in V$, we have 
		\begin{align*}
		\frac{\sigma_A(\oriented{xy})+\gamma_A\sigma_A(\oriented{yx})}{1+\gamma_A}+\frac{\sigma_B(\oriented{xy})+\gamma_B\sigma_B(\oriented{yx})}{1+\gamma_B}
		& =  \frac{\alpha_1+\hat\beta_1+(\beta_2-\hat\beta_2)}{W(\mu)} \mu(\fmat{xy}) \overset{\eqref{eq:lem-improved-hatbeta_1+}}{\leq} \mu(\fmat{xy}),
		\end{align*}and
			\begin{align*}
			\frac{\gamma_A\sigma_A(\oriented{xy})+\sigma_A(\oriented{yx})}{1+\gamma_A}+\frac{\gamma_B\sigma_B(\oriented{xy})+\sigma_B(\oriented{yx})}{1+\gamma_B}
			& = \frac{\alpha_2+\hat\beta_2+(\beta_1-\hat\beta_1)}{W(\mu)} \mu(\fmat{xy}) \overset{\eqref{eq:lem-improved-hatbeta_2+}}{\leq }\mu(\fmat{xy}),
		\end{align*}
		so $\sigma_A+\sigma_B\trianglelefteq \mu$ holds, giving \ref{item:improvedbalancing-out1}.
\end{proof}

\subsection{Proof of the Completion Lemma}
The proof is technical, so we divide it into six main steps. 

\begin{enumerate}
	\item \emph{Step 1: Defining $\sigma_A$.}
	The fractional matching $\mu$ is sufficiently large to accommodate the entire $\gamma_A$-skew-matching $\sigma_A$. Since the anchor must be contained in $U$, there is no flexibility in the choice of orientation. We define $\sigma_A$ by distributing its weight proportionally according to the weight of the fractional matching $\mu$.

\item \emph{Step 2: Partitioning $\mu$ and $\sigma_A$.} We partition $\mu$ into $\mu'$ and $\widehat{\mu}$ such that all of $\mu'$ can be ``reached from $u$ from both sides'' (i.e., both endpoints are in the $w$-neighborhood of $u$). This partition makes $\mu'$ easier to handle, as we can choose on which side to place the anchor. Similarly, we partition $\sigma_A$ into $\sigma_A'$ and $\widehat{\sigma_A}$ according to the proportions of $\mu'$ and $\widehat{\mu}$ on each edge.

\item \emph{Step 3: Perfectly filling $\mu'$.} We complete $\sigma_A'$ by adding $\sigma_B'$ to perfectly fill $\mu'$. This is achieved in two steps, formulated in \Cref{cl:new-balancing} and \Cref{cl:upTomu'}. In \Cref{cl:new-balancing}, we identify a minimum $\bar{\sigma}_B$ to compensate for the skew in $\sigma_A'$. This requires carefully selecting the orientation of $\bar{\sigma}_B$ on each support edge of $\mu'$. After compensating for the skew in $\sigma_A'$, we proceed to \Cref{cl:upTomu'}, where $\bar{\sigma}_B$ is complemented by a $\gamma_B$-skew-matching to obtain $\sigma_B'$, balancing its orientations within the remainder of $\mu'$.

\item \emph{Step 4: Fitting the remainder into $\widehat{\mu}$.} The fractional matching $\widehat{\mu}$ can only be accessed from $u$ on one side, meaning there is no choice regarding the placement of the anchor for a $\gamma_B$-skew-matching. After removing a fractional sub-matching $\tilde \mu$ with $\widehat \sigma_A\trianglelefteq \tilde \mu$ (i.e., where $\widehat \sigma_A$ ``lives''), we fit a $\gamma_B$-skew-matching into the remaining portion $\widehat\mu-\tilde \mu$.
After these four steps, we will have proven properties~\ref{item:completion-1}--\ref{item:completion-out5}.

\item \emph{Steps 5 and 6: The case $\gamma_B \leq 1$.} We are left to refine our approach to prove property~\ref{item:completion-moreover} under the additional assumption that $\gamma_B\le 1$. This is treated in two additional steps, depending on the case if $\gamma_A\ge 1$ and on the case when $\gamma_A<1$. When $\gamma_A\ge 1$, the skews $\gamma_A, \gamma_B$ work in our favour and the proof of property~\ref{item:completion-moreover} is straightforward using the objects we have defined so far.
However, if $\gamma_A<1$, we have to define things differently. Instead of placing a $\gamma_B$-skew-matching in the fractional matching $\widehat\mu-\tilde \mu$ as previously, we first need to complement $\widehat\sigma_A$ with a $\gamma_B$-skew-matching $\sigma$ that together with $\widehat\sigma_A$ uses the fractional matching $\widehat\mu$ equally on both of its sides. Only then we can complete it with a $\gamma_B$-skew-matching $\sigma'$ to ensure the set $V$ is fully covered. 
\end{enumerate}

\begin{proof}[Proof of \Cref{lem:new-completion} (Completion)]
We follow the six steps outlined above. \medskip

\noindent \emph{Step 1: Defining $\sigma_A$}.
For all $x\in U$ and $y\in V$, we set 
	\begin{align}\label{eq:new-completion-def-sigma_A}
		\sigma_A(\oriented{xy}):= \frac{\alpha_1+\alpha_2}{W(\mu)} \mu(\fmat{xy}),
	\end{align} and we put $\sigma_A$ equal to $0$ on all other edges. 
	Observe that \begin{align}\label{eq:lem-compl-WeightOfSigma_A}
		W(\sigma_A)=\sum_{x\in U, y\in V}\sigma_A(\oriented{xy})=\frac{\alpha_1+\alpha_2}{W(\mu)}\sum_{x\in U, y\in V}\mu(\fmat{xy})=\alpha_1+\alpha_2.
	\end{align}
	By construction, we have that
	\begin{equation}
		\label{equation:completion-anchorofsAisinU}
		\mathcal A(\sigma_A) \subseteq U,
	\end{equation}
	i.e. the anchor of $\sigma_A$ is contained in $U$.
	Using \ref{eq:neq-fit}, we also have
	\[
	\frac{\max\{1,\gamma_A\}}{1+\gamma_A}\sigma_A(\oriented{xy}) = \frac{\max\{\alpha_1,\alpha_2\}}{\alpha_1+\alpha_2}\sigma_A(\oriented{xy}) \leq \frac{W(\mu)}{\alpha_1+\alpha_2}\sigma_A(\oriented{xy}) = \mu(\fmat{xy}),
	\]
	and therefore
	$\sigma_A\trianglelefteq\mu$. \medskip
	
	\noindent \emph {Step 2: Partitioning $\mu$ and $\sigma_A$.}
	Let $\mu'\le \mu$ be a fractional matching such that for all $x\in U$ and $y\in V$ we have 
	\begin{align}\label{eq:def-mu'}
		\mu'(\fmat{xy}):= \frac{ \min\{w(\oriented{ux}), \mu(x)\}}{\mu(x)}\mu(\fmat{xy}).
	\end{align}
	By \ref{item:completion-1}, we have
	\begin{equation}\label{eq:mu'_y_leq_w_uy}
		\mu'(y)\leq\mu(y)\leq w(\oriented{uy}), \text{ for all } y \in V.
	\end{equation}
 	Directly from the definition, we also have
	\begin{align}\label{eq:for_mu'}
	\mu'(x)=\min\{w(\oriented{ux}), \mu(x)\}, \text{ for all } x\in U.
	\end{align}  and thus 
	\begin{align}
	\nonumber 	\deg_w(u, \mu')&=\sum_{x\in U}\min\{w(\oriented{ux}), \mu'(x)\}+\sum_{y\in V}\min\{w(\oriented{uy}), \mu'(y)\}\\
		&\overset{\eqref{eq:mu'_y_leq_w_uy},\eqref{eq:for_mu'}}{\le} \sum_{x\in U}\mu'(x)+\sum_{y\in V}\mu'(y)=2W(\mu')\label{eq:deg_u_mu}.
	\end{align}
	
	Now set
	\begin{equation} \label{eq:widehat_mu}
		\widehat \mu:= \mu-\mu'.
	\end{equation} 
	We partition $\sigma_A$ into the part that will fit in $\mu'$ and the part that will fit in $\widehat\mu$.
	For all $x \in U, y \in V$, set
	\begin{align}\label{eq:def-simga_A'}
		\sigma_A'(\oriented{xy}):= \frac{\mu'(\fmat{xy})}{\mu(\fmat{xy})}\sigma_A(\oriented{xy})=\frac{\alpha_1+\alpha_2}{W(\mu)}\cdot \mu'(\fmat{xy}),
	\end{align} 
	and let $\sigma_A'$ be equal to $0$ on all other edges.
	Let 
	\begin{equation} \label{eq:widehat_sigma_A}
	\widehat\sigma_A:= \sigma_A-\sigma_A'.
	\end{equation}
	Concretely, we have for $x \in U, y \in V$,
	\begin{equation}\label{eq:def-widehatsimga_A}
		\widehat\sigma_A(\oriented{xy})=\frac{\widehat\mu(\fmat{xy})}{\mu(\fmat{xy})}\sigma_A(\oriented{xy})=\frac{\alpha_1+\alpha_2}{W(\mu)}\cdot \widehat\mu(\fmat{xy}),
	\end{equation}
	and $\widehat\sigma_A$ equal to $0$ on all other edges.	
	Observe that, by our construction and the fact that $\sigma_A \trianglelefteq \mu$, we have that $\sigma_A'\trianglelefteq\mu'$ and
	\begin{equation}
		\widehat\sigma_A\trianglelefteq\widehat\mu.
		\label{equation:completion-hatsigmaAvshatmu}
	\end{equation} 
	
	\noindent \emph{Step 3: Filling $\mu'$ perfectly.}
	In the first claim, our goal is to counterbalance the skew of $\sigma_A'$ by carefully positioning a minimal $\gamma_B$-skew-matching, ensuring that each support edge of the fractional matching carries equal weight at both endpoints.
	\begin{claim}\label{cl:new-balancing}
		There are a fractional matching $\bar\mu\le \mu'$ and  a $\gamma_B$-skew-matching $\bar\sigma_B$ such that $\sigma_A'+\bar \sigma_B\trianglelefteq \bar \mu$, and $W(\sigma_A')+W(\bar \sigma_B)=2W( \bar\mu)$. 
	\end{claim}
	
	Before proving the claim, we gather two useful facts for the rest of the proof.
	First, we claim that
	\begin{equation}\label{eq:max-min}
		\mbox{if }\max\{1,\gamma_A\}>\min\{1, \gamma_A\}\mbox{, then }\max\{1,\gamma_B\}>\min\{1,\gamma_B\}.
	\end{equation}	
	Indeed, suppose otherwise.
	Then $\gamma_A \neq 1$ and $\gamma_B = 1$, so $\alpha_1 \neq \alpha_2$ and $\beta_1 = \beta_2$.
	We have 
	\begin{align*}
		\max\{\alpha_1,\alpha_2\}+\min\{\beta_1,\beta_2\}
		& = \max\{\alpha_1,\alpha_2\}+\max\{\beta_1,\beta_2\}\\
		& > \min\{\alpha_1,\alpha_2\}+\max\{\beta_1,\beta_2\}\\
		& \overset{\text{\ref{eq:VybalancovaniMozne}}}{\ge} W(\mu) \overset{\text{\ref{eq:neq-fit}}}{\ge}\max\{\alpha_1,\alpha_2\}+\min\{\beta_1,\beta_2\},
	\end{align*}a contradiction.
	This proves \eqref{eq:max-min}.
	
	Next, we gather a useful inequality.
	Observe that
	\begin{align*}
		\nonumber 
		\frac{\max\{1,\gamma_A\}}{1+\gamma_A}\sigma_A'(\oriented{xy}) & +\frac{\min\{\beta_1,\beta_2\}}{W(\mu)}\mu'(\fmat{xy})
		\\ \nonumber 
		& \overset{\eqref{eq:def-simga_A'}}{=}\frac{\max\{\alpha_1,\alpha_2\}+\min\{\beta_1,\beta_2\}}{W(\mu)}\mu'(\fmat{xy})\\ \nonumber 
		& \overset{\text{\ref{eq:VybalancovaniMozne}-\ref{eq:neq-fit}}}{\le } \frac{\min\{\alpha_1,\alpha_2\}+\max\{\beta_1,\beta_2\}}{W(\mu)}\mu'(\fmat{xy})\\
		& = \frac{\min\{1,\gamma_A\}}{1+\gamma_A}\sigma_A'(\oriented{xy})+\frac{\max\{\beta_1,\beta_2\}}{W(\mu)}\mu'(\fmat{xy}). 
	\end{align*}
	and therefore
	\begin{equation}
		\frac{\max\{1,\gamma_A\} - \min\{1, \gamma_A\}}{1+\gamma_A}\sigma_A'(\oriented{xy}) \leq \frac{\max\{\beta_1,\beta_2\} - \min\{\beta_1,\beta_2\}}{W(\mu)}\mu'(\fmat{xy}).
		\label{equation:completion-maxsprimavsmin}
	\end{equation}

	\begin{proofclaim}[Proof of \Cref{cl:new-balancing}]
		We will separate the proof into two cases, depending if $(1-\gamma_A)(1-\gamma_B)\le 0$, or $(1-\gamma_A)(1-\gamma_B)> 0$. \medskip
		
		\noindent \emph{Case 1: $(1-\gamma_A)(1-\gamma_B)\le 0$}.
		We define a $\gamma_B$-skew matching $\bar \sigma_B$ as follows.
		If $\gamma_A = 1$, then $\bar \sigma_B$ is identically zero.
		Otherwise, if $\gamma_A\neq 1$, for every $x\in U$ and every $y\in  V$ we choose $\bar \sigma_B(\oriented{xy})$ so that
		\begin{equation}\label{eq:completion-s3-claim-case1}
	\bar	\sigma_B(\oriented{xy}):=	\frac{1+\gamma_B}{1+\gamma_A}\cdot \frac{\max\{1,\gamma_A\}-\min\{1,\gamma_A\}}{\max\{1, \gamma_B\}-\min\{1, \gamma_B\}}\cdot \sigma_A'(\oriented{xy}),
		\end{equation}
		and $\bar \sigma_B(\oriented{xy}) = 0$ in every other edge.
		Observe that if $\gamma_A\neq 1$, by~\eqref{eq:max-min}, we also have $\gamma_B \neq 1$.
		Together with $\gamma_B \leq 1$, this yields $\gamma_B < 1$, and in particular, $\max\{1,\gamma_B\}>\min\{1,\gamma_B\}$ and thus $\bar \sigma_B(\oriented{xy})$ is correctly defined.
		
		We claim that
		\begin{equation} \label{eq:completion-s3-claim-case1-trianglemuprime}
			\sigma'_A + \bar \sigma_B \trianglelefteq \mu'
		\end{equation} holds (note this also verifies that $\bar \sigma_B$ is  indeed a $\gamma_B$-skew matching).
		
		If $\gamma_A = 1$ then \eqref{eq:completion-s3-claim-case1-trianglemuprime} is immediate, because $\sigma'_A \trianglelefteq \mu'$ is true.
		Hence, we assume $\gamma_A\neq 1$.
		We need to verify, for every $x \in U$, $y \in V$, that
		\begin{align}\label{eq:onTheU-side}
			\frac{\sigma_A'(\oriented{xy})}{1+\gamma_A}+\frac{\bar \sigma_B(\oriented{xy})}{1+\gamma_B}\le \mu'(\fmat{xy}), 
		\end{align}and 
		\begin{align}\label{eq:onTheV-side}
			\frac{\gamma_A\sigma_A'(\oriented{xy})}{1+\gamma_A}+\frac{\gamma_B\bar \sigma_B(\oriented{xy})}{1+\gamma_B}\le \mu'(\fmat{xy}),
		\end{align}
		hold, as required by the desired condition \eqref{eq:completion-s3-claim-case1-trianglemuprime}.
		
		Inequalities \eqref{equation:completion-maxsprimavsmin} and~\eqref{eq:completion-s3-claim-case1} together, imply, for each $x \in U$, $y \in V$, that
		\[ \bar \sigma_B(\oriented{xy})\le \frac{\beta_1+\beta_2}{W(\mu)}\mu'(\fmat{xy}),\] and thus 
		\begin{align}
			\nonumber \frac{\max\{1,\gamma_A\}}{1+\gamma_A}\sigma_A'(\oriented{xy})+\frac{\min\{1,\gamma_B\}}{1+\gamma_B}\bar \sigma_B(\oriented{xy}) & \overset{\eqref{eq:def-simga_A'}}{\leq} \frac{\max\{\alpha_1,\alpha_2\}+\min\{\beta_1,\beta_2\}}{W(\mu)}\mu'(\fmat{xy})\\
			& \overset{\text{\ref{eq:neq-fit}}}{\le} \mu'(\fmat{xy}).\label{eq:firstinequ}
		\end{align}
		Now, from~\eqref{eq:completion-s3-claim-case1}, we obtain
		\begin{align}\nonumber 
			\frac{\min\{1,\gamma_A\}}{1+\gamma_A}\sigma_A'(\oriented{xy})+\frac{\max\{1,\gamma_B\}}{1+\gamma_B}\bar \sigma_B(\oriented{xy})=&\frac{\max\{1,\gamma_A\}}{1+\gamma_A}\sigma_A'(\oriented{xy})+\frac{\min\{1,\gamma_B\}}{1+\gamma_B}\bar \sigma_B(\oriented{xy})\\
			\overset{\eqref{eq:firstinequ}}{\le}&\mu'(\fmat{xy}). \label{eq:secondineq}
		\end{align}
		Inequalities~\eqref{eq:firstinequ} and~\eqref{eq:secondineq} imply~\eqref{eq:onTheU-side} and~\eqref{eq:onTheV-side}, as the condition $(1-\gamma_A)(1-\gamma_B) \leq 0$ and~\eqref{eq:max-min} imply that $\gamma_A > 1$ if and only if $ \gamma_B < 1$.
		Thus indeed \eqref{eq:completion-s3-claim-case1-trianglemuprime} holds.
		
		Finally, for every $x\in U$ and every $y\in N_{\mu'}(x)$, set 
		\begin{align*}
			\bar\mu(\fmat{xy})& := \frac{\min\{1,\gamma_A\}}{1+\gamma_A}\sigma_A'(\oriented{xy})+\frac{\max\{1,\gamma_B\}}{1+\gamma_B}\bar \sigma_B(\oriented{xy})\\
			& = \frac{\max\{1,\gamma_A\}}{1+\gamma_A}\sigma_A'(\oriented{xy})+\frac{\min\{1,\gamma_B\}}{1+\gamma_B}\bar \sigma_B(\oriented{xy}),
		\end{align*}
		where the equality is due to~\eqref{eq:completion-s3-claim-case1}.
		From the definition we obtain $W(\sigma_A')+W(\bar \sigma_B)=2W(\bar\mu)$, and that $\sigma_A'+\bar \sigma_B\trianglelefteq\bar\mu$, as desired.\medskip
		
		\noindent \emph{Case 2: $(1-\gamma_A)(1-\gamma_B) > 0$}
		In this case, the proof goes similarly as above, but the support edges of $\bar \sigma_B$ will go in the opposite direction as the support edges of $\sigma_A'$. This means that for every $x\in U$ and every $y\in V$ we set
		\begin{equation}\label{eq:nas-stary-b')}
		\bar \sigma_B(\oriented{yx}):= 	\frac{1+\gamma_B}{1+\gamma_A}\cdot  \frac{\max\{1,\gamma_A\}-\min\{1,\gamma_A\}}{\max\{1, \gamma_B\}-\min\{1, \gamma_B\}}\cdot \sigma_A'(\oriented{xy}),
		\end{equation}and $\bar \sigma_B= 0$ on all other edges. 
		
		Analogously as above, we claim that~\eqref{eq:completion-s3-claim-case1-trianglemuprime} holds. We need to verify  for every $x\in U$ and every $y\in V$ that
		\begin{align}\label{eq:onTheU-side'}
			\frac{\sigma_A'(\oriented{xy})}{1+\gamma_A}+\frac{\gamma_B\bar \sigma_B(\oriented{yx})}{1+\gamma_B}\le \mu'(\fmat{xy}), 
		\end{align}and 
		\begin{align}\label{eq:onTheV-side'}
			\frac{\gamma_A\sigma_A'(\oriented{xy})}{1+\gamma_A}+\frac{\bar \sigma_B(\oriented{yx})}{1+\gamma_B}\le \mu'(\fmat{xy}).
		\end{align} 
		The calculations go verbatim,
		simply by switching the orientation of the support edges of $\bar\sigma_B$ and realising that $(1-\gamma_A)(1-\gamma_B)>0$ implies that $\gamma_A=\max\{1,\gamma_A\}$ if and only if $1=\min\{1,\gamma_B\}$. 
		
		Then, we set 
		\begin{align*}
			\bar\mu(\fmat{xy}):=& \frac{\min\{1,\gamma_A\}}{1+\gamma_A}\sigma_A'(\oriented{xy})+\frac{\max\{1,\gamma_B\}}{1+\gamma_B}\bar \sigma_B(\oriented{yx})\\
			=& \frac{\max\{1,\gamma_A\}}{1+\gamma_A}\sigma_A'(\oriented{xy})+\frac{\min\{1,\gamma_B\}}{1+\gamma_B}\bar \sigma_B(\oriented{yx})
		\end{align*}
		for every $x\in U$ and every $y\in N_\mu(x)$.
		We obtain $W(\sigma_A')+W(\bar \sigma_B)=2W(\bar\mu)$, and that $\sigma_A'+\bar \sigma_B\trianglelefteq\bar\mu$, as desired.
	\end{proofclaim}
	
	In this second claim, we complete the remaining unmatched portion of the fractional matching using a $\gamma_B$-skew-matching, alternating its orientation to ensure that the weight on each support edge is equally balanced between its endpoints.
	\begin{claim}\label{cl:upTomu'}
		There is a $\gamma_B$-skew-matching $ \sigma_B'$ with
		\begin{equation}
			\label{equation:completion:wsigmab'} W(\sigma_A')+W( \sigma_B')=2W(\mu')
		\end{equation}
		and
		\begin{equation}
			\label{equation:completion-sigmaABtrianglemu}
			\sigma_A'+\sigma_B'\trianglelefteq \mu'.
		\end{equation}
		Moreover, $\sigma_B' \ge \bar \sigma_B$,
		and the anchor of $\sigma'_B$ fits in the $w$-neighbourhood of $u$.
	\end{claim}	

	\begin{proofclaim} 
		Let $(H, \bar w)$ be the $\bar\mu$-truncated weighted graph obtained from $(H,w)$.
		Recall that $\mu' \leq \mu$ and that $V(\mu) \subseteq U \cup V$.
		We apply \Cref{lem:new-balancing} (Balancing-out) with $\mu' - \bar \mu$ in place of $\mu$.
		We obtain a $\gamma_B$-skew-matching $\tilde \sigma_B\trianglelefteq\mu'-\bar\mu$ in $ H^{\leftrightarrow}$
		of weight $W(\tilde \sigma_B)=2W(\mu'-\bar\mu)$ and with its anchor $\mathcal A(\tilde\sigma_B)$ contained in $U\cup V$.
		
		Set $\sigma_B':= \bar\sigma_B+\tilde \sigma_B$.
		As $\tilde \sigma_B\trianglelefteq \mu'-\bar\mu$ and $\sigma_A'+\bar \sigma_B\trianglelefteq\bar\mu$, we have $\sigma_A'+\sigma_B'\trianglelefteq\mu'$, as required, and this also shows that $\sigma'_B$ is well-defined.
		Observe that $\sigma_B' \ge \bar \sigma_B$ holds immediately.
		From \Cref{cl:new-balancing}, we have  $W(\sigma'_A)+W(\bar\sigma_B) = 2 W(\bar \mu)$. Thus we obtain $W(\sigma_A')+W(\sigma_B')=2W(\mu')$, as required. 
		
		Finally, using $\sigma_B'\trianglelefteq\mu'$ and both \eqref{eq:mu'_y_leq_w_uy} and \eqref{eq:for_mu'}, we obtain that the anchor  $\mathcal A(\sigma_B')$ fits in the $w$-neighbourhood of $u$. 
\end{proofclaim}

\noindent \emph{Step 4:  Fitting the remainder into $\widehat{\mu}$.} Recall that $\widehat\mu:= \mu-\mu'$ and $\widehat \sigma_A:= \sigma_A - \sigma_A'$.
By \eqref{equation:completion-hatsigmaAvshatmu}, we can let $\tilde \mu\le \widehat \mu$ be a minimal fractional matching such that $\widehat\sigma_A\trianglelefteq \tilde \mu$.
Recall that $\sigma_A$ (and therefore, $\widehat \sigma_A$) is supported only in edges $\oriented{xy}$ with $x \in U$ and $y \in V$.
This implies that, for every such $\oriented{xy}$, we have
\[ \tilde \mu(\fmat{xy}) = \frac{\max\{1, \gamma_A\}}{1 + \gamma_A} \sigma_A(\oriented{xy}), \]
and we have $\tilde \mu (xy) = 0$ for any other edge.
Hence, we have
	\begin{equation}\label{eq:widthtildemu}
		W(\tilde \mu)=\max\{1,\gamma_A\}\frac{W(\widehat\sigma_A)}{1+\gamma_A}.
	\end{equation}
	
	Now we obtain a $\gamma_B$-skew-matching $\widehat \sigma_B$ such that $\widehat \sigma_B\trianglelefteq\widehat{\mu}-\tilde{\mu}$.

	\begin{claim}\label{cl:I_need_a_reference}
	There is a $\gamma_B$-skew-matching $\widehat \sigma_B\trianglelefteq\widehat{\mu}-\tilde{\mu}$ in $H^{\leftrightarrow}$ of weight 
	\begin{equation}
		W(\widehat \sigma_B)=(1+\min\{\gamma_B, \gamma_B^{-1}\})W(\widehat{\mu}-\tilde{\mu}),
		\label{equation:completion-wwidehatsigma}
	\end{equation}
	and such that $\mathcal A(\widehat \sigma_B)\subseteq V$.
	\end{claim}
	
	\begin{proofclaim}
		Recalling the definitions, we have $\widehat{\mu}-\tilde{\mu} \leq \widehat{\mu} = \mu - \mu' \leq \mu$.
		Since $\mu$ is a matching running between $U$ and $V$ by assumption, the same is true for $\widehat{\mu}-\tilde{\mu}$.
		Hence, the claim follows by applying \Cref{lem:new-extending-out} with $\widehat{\mu}-\tilde{\mu}, \gamma_B, V, U$ playing the roles of $\mu, \gamma, U, V$.
	\end{proofclaim}
	
	Let $\sigma_B:=\widehat \sigma_B+\sigma_B'$.
	Since $\widehat \sigma_B \trianglelefteq \hat \mu - \tilde \mu \leq \mu - \mu'$ and $\sigma'_B \trianglelefteq \mu'$, it quickly follows that $\sigma_B$ is a $\gamma_B$-skew-matching.
	Recall that $\sigma_A$ was defined at the beginning of the proof.
	We shall prove that $\sigma_A, \sigma_B$ satisfy the required \ref{item:completion-out1}--\ref{item:completion-out5}.
	
	To see \ref{item:completion-out1}, note that
	\[ \sigma_B+\sigma_A = \widehat{\sigma}_B + \sigma_B' + \sigma_A' + \widehat{\sigma}_A \trianglelefteq (\widehat{\mu} - \widetilde{\mu}) + \mu' + \widetilde{\mu} = \mu,\]
	where we used the definition of $\sigma_B$ and \eqref{eq:widehat_sigma_A} in the first equality,
	and then we used the choice of $\widehat{\sigma}_B$, \eqref{equation:completion-sigmaABtrianglemu} and \eqref{equation:completion-hatsigmaAvshatmu}.
	Point \ref{item:completion-out2} follows from \eqref{eq:lem-compl-WeightOfSigma_A}.
	
	To see \ref{item:completion-out3}, we shall estimate $W(\sigma_B)$ using all our work so far.
	We have
	\begin{align*}
		W(\sigma_B)
		& = W(\widehat \sigma_B) + W(\sigma'_B) \\
		& \overset{\eqref{equation:completion-wwidehatsigma}}{\geq} W(\widehat{\mu}-\tilde{\mu}) + W(\sigma'_B) \\
		& \overset{\eqref{equation:completion:wsigmab'}}{=} 2W(\mu')-W(\sigma_A')+W(\widehat{\mu}-\tilde{\mu}) \\
		& \overset{\eqref{eq:widthtildemu}}{\ge} 2W(\mu')-W(\sigma_A')+W(\widehat\mu)-W(\widehat\sigma_A)\\
		& \overset{\eqref{eq:deg_u_mu}}{\ge} \deg_w(u, \mu')-W(\sigma_A')+W(\widehat\mu)-W(\widehat\sigma_A)\\
		& \overset{\eqref{eq:widehat_sigma_A}}{=} \deg_w(u, \mu')+W(\widehat \mu)-W(\sigma_A)\\
		& = \sum_{x \in V(H)} \min \{ w(\oriented{ux}), \mu'(x) \} +W(\widehat \mu)-W(\sigma_A) \\
		& = \sum_{x \in U} \min \{ w(\oriented{ux}), \mu'(x) \} + \sum_{x \in V} \min \{ w(\oriented{ux}), \mu'(x) \} +W(\widehat \mu)-W(\sigma_A) \\
		& \overset{\eqref{eq:mu'_y_leq_w_uy}}{=} \sum_{x \in U} \min \{ w(\oriented{ux}), \mu'(x) \} + \sum_{x \in V} \mu'(x) +W(\widehat \mu)-W(\sigma_A) \\
		& = \sum_{x \in U} \min \{ w(\oriented{ux}), \mu'(x) \} + \sum_{x \in V} \mu'(x) + \sum_{x \in V} \hat \mu(x) -W(\sigma_A) \\
		& \overset{\eqref{eq:widehat_mu}}{=} \sum_{x \in U} \min \{ w(\oriented{ux}), \mu'(x) \} + \sum_{x \in V}
		\mu(x) -W(\sigma_A) \\
		& \overset{\eqref{eq:for_mu'}}{=} \sum_{x \in U} \min \{ w(\oriented{ux}), \mu(x) \} + \sum_{x \in V}
		\mu(x) -W(\sigma_A) \\
		& \overset{\eqref{eq:mu'_y_leq_w_uy}}{=} \sum_{x\in U\cup V}\min\{w(\oriented{ux}), \mu(x)\}-W(\sigma_A) \\
		& = \deg_{w}(u, \mu) - W(\sigma_A),
	\end{align*}
	which, together with the obvious $W(\sigma_B) \geq 0$, gives \ref{item:completion-out3}.
	
	Point \ref{item:completion-out4} follows from \eqref{equation:completion-anchorofsAisinU}.
	Finally, to verify \ref{item:completion-out5}, we need to check that $\sigma^1_B(x) \leq w(\oriented{ux})$ for all $x \in V(H)$.
	It suffices to check it for $x \in U \cup V$.
	Suppose first that $x \in U$.
	Note that $\mathcal A(\widehat\sigma_B)\subseteq V$ follows from our choice in \Cref{cl:I_need_a_reference}, so we have $\widehat\sigma^1_B(x) = 0$, and $\mathcal A(\sigma_B')$ fits in the $w$-neighbourhood of $u$ by our choice in \Cref{cl:upTomu'}.
	Hence we have $\sigma^1_B(x) = \sigma'_B(x) \leq w(\oriented{ux})$, as required.
	Now, suppose that $x \in V$.
	Here, we have $\sigma_B^1(x) \leq \mu(x)\leq w(\oriented{ux})$, where in the first inequality we used \ref{item:completion-out1}, and the second inequality follows from \eqref{eq:mu'_y_leq_w_uy}.
	Thus \ref{item:completion-out5} holds.

We have thus proven \ref{item:completion-1}--\ref{item:completion-out5} and it only remains to prove \ref{item:completion-moreover}. So we may assume that $\gamma_B\le 1$, as otherwise there is nothing to show.

	Let us explain our strategy here.
	We will need to divide the proof in two cases, depending if $\gamma_A \geq 1$, or not. In any case, we shall keep the definition of $\sigma_A, \sigma_A'$, and $\widehat\sigma_A$ as in~\eqref{eq:new-completion-def-sigma_A},~\eqref{eq:def-simga_A'}, and~\eqref{eq:def-widehatsimga_A},  as well as we shall keep the definition of $\sigma_B'$   from \Cref{cl:upTomu'} and $\bar \sigma_B$ from \Cref{cl:new-balancing}. We keep $\widehat\sigma_B$ from \Cref{cl:I_need_a_reference} as well, but we shall use it to determine $\sigma_B$ only in some cases.
	This is the part where the value of $\gamma_A$ becomes relevant.
	Indeed, while verifying \ref{item:completion-moreover} is natural in the case when $\gamma_A\ge 1$ and $\gamma_B\le 1$, one has to define things differently when $\gamma_A<1$.  In this case, we shall not use $\widehat\sigma_B$ in the definition of $\sigma_B$ as above. The main problem is that $\widehat \sigma_A$  does not saturate $\tilde \mu$ in $V$. Therefore, it has to be suitably completed by some $\gamma_B$-skew-matching. To achieve this, we shall proceed similarly as in \Cref{cl:new-balancing}. This $\gamma_B$-skew-matching is then completed with a $\gamma_B$-skew-matching in a standard way to fill the left-over of $\widehat \mu$ in $V$. This is possible, as $\gamma_B\le 1$.
	 After that, we glue the respective $\gamma_A$ and $\gamma_B$-skew-matchings together and check that indeed fulfill the required properties. 
	Now we turn to the details. \medskip

	\noindent \emph{Step 5: Saturating $V$, case $\gamma_A \geq 1$.}
	We consider first the easiest case when $\gamma_A \geq 1$.
	Observe that $\mathcal{A}(\widehat{\sigma}_B) \subseteq V$ implies that $\widehat{\sigma}_B(\oriented{xy}) = 0$ for all $x\in U, y\in V$.
	From $\widehat{\sigma}_B \trianglelefteq \widehat \mu - \widetilde \mu$ we have $\widehat{\sigma}_B(y) \leq \widehat \mu(y) - \widetilde \mu (y)$ for all $y \in V$.
	From $\gamma_B\le 1$ and \eqref{equation:completion-wwidehatsigma} we get $W(\widehat{\sigma}_B) = (1 + \gamma_B) W(\widehat{\mu} - \widetilde{\mu})$, and therefore
	\[ \sum_{y \in V} \widehat{\sigma}_B(y) = \frac{W(\widehat{\sigma}_B)}{1+\gamma_B}
	= W(\widehat{\mu} - \widetilde{\mu}) = \sum_{y \in V} (\widehat{\mu}(y) - \widetilde{\mu}(y)),
	 \]
	 which together with the above, implies that for all $y\in V$, we have
	 \[\widehat\sigma_B(y)= \widehat\mu(y)-\tilde \mu(y).\] 
	 Now note that from $\widehat{\sigma}_A \trianglelefteq\tilde \mu$ we have $\widehat{\sigma}_A(y) \leq \tilde{\mu}(y)$ for all $y \in V$.
	 Note that from $\gamma_A \geq 1$ and \eqref{eq:widthtildemu}, we obtain $\gamma_A W(\widehat \sigma_A) = (1 + \gamma_A) W(\tilde \mu)$.
	 Together with $\mathcal{A}(\widehat{\sigma}_A) \subseteq \mathcal{A}(\sigma_A) \subseteq U$, a similar argument as above gives, for all $y \in V$,
	 \[\widehat\sigma_A(y)=\tilde \mu(y).\]
	 Next, by \Cref{cl:upTomu'} 
	we have $\sigma_B'+\sigma_A'\trianglelefteq \mu'$ and $W(\sigma_B')+W(\sigma_A' )=2W(\mu')$.
	By a similar reasoning as above, we have, for all $y\in V$,
	\begin{equation} \label{eq:moreover_sigma_mu_bar}
		\sigma_B'(y)+\sigma_A'(y)=\mu'(y).
	\end{equation}
	Using the three equalities above, together with $\sigma_B:=\widehat \sigma_B+\sigma_B'$ and \eqref{eq:widehat_sigma_A}, \eqref{eq:widehat_mu}, we obtain
	\begin{equation*}
		\sigma_A(y)+\sigma_B(y) = 
		\widehat\sigma_B(y)+\sigma_B'(y)+\widehat \sigma_A(y)+\sigma_A'(y) = \widehat \mu(y) + \mu'(y) = \mu(y)
	\end{equation*} for all $y\in V$, thus giving \ref{item:completion-moreover} in the case $\gamma_A \geq 1$. \medskip
%

	\noindent \emph{Step 6: Saturating $V$, case $\gamma_A < 1$.}
	From now on, we can assume, in addition to $\gamma_B \leq 1$, that $\gamma_A < 1$.
	By~\eqref{eq:max-min}, in fact we have that $\gamma_B<1$.
	Observe now that $(1-\gamma_A)(1-\gamma_B)>0$, as in the second case of \Cref{cl:new-balancing}.
	We shall define a new $\gamma_B$-skew-matching in place of $\sigma_B$ to conclude.
	
	Let $\sigma$ be the maximal $\gamma_B$-skew-matching with $\sigma+\widehat\sigma_A\trianglelefteq \widehat\mu$ such that for every $x\in U$ and every $y\in V$ we have 
	\begin{align}
		\label{eq:nas-stary-b'')}
		\frac{1-\gamma_A}{1+\gamma_A}\widehat \sigma_A(\oriented{xy})\ge \frac{1- \gamma_B}{1+\gamma_B} \sigma(\oriented{yx}),
	\end{align} and $\sigma=0$ on all other edges.
	Analogously as in the second case of \Cref{cl:new-balancing}, the definition implies that we obtain equality in~\eqref{eq:nas-stary-b'')}.
	The calculations go verbatim, replacing $\bar\sigma_B$ by $\sigma$, $\sigma_A'$ by $\widehat\sigma_A$, and  $\mu'$ by $\widehat\mu$.

	Next, let $\tilde \mu'\le \widehat \mu$ be the fractional matching such that
	\begin{align} \label{eq:mu_tilde_ap}
		\tilde \mu'(\fmat{xy}):= \frac{\widehat\sigma_A(\oriented{xy})}{1+\gamma_A}+\frac{\gamma_B\sigma(\oriented{yx})}{1+\gamma_B}=\frac{\sigma(\oriented{yx})}{1+\gamma_B}+\frac{\gamma_A\widehat\sigma_A(\oriented{xy})}{1+\gamma_A},
	\end{align}
	where the equality follows from the fact that we have equality in~\eqref{eq:nas-stary-b'')} for each $\oriented{xy}$. 
	From \eqref{eq:mu_tilde_ap}, we have 
	\begin{align} \label{eq:W_mu_tilde_ap}
		W(\widehat\sigma_A)+W(\sigma)=2W(\tilde \mu')
	\end{align} 
	and, similarly as before, we also obtain 
	\begin{align} \label{eq:mu_tilde_ap_y}
	\sigma(y)+\widehat\sigma_A(y)=\sum_{x\in U}\frac{\sigma(\oriented{yx})}{1+\gamma_B}+\frac{\gamma_A\widehat\sigma_A(\oriented{xy})}{1+\gamma_A}=\tilde \mu'(y)
	\end{align} 
	for all $y\in V$.
	From \eqref{eq:mu_tilde_ap} we also have that
	\begin{equation}
		\label{equation:completion-sigawidehattildemuprime}
		\sigma + \widehat{\sigma}_A \trianglelefteq \tilde \mu'.
	\end{equation}
	
	Recall that $\gamma_B < 1$.
	Apply \Cref{lem:new-extending-out} with $V, U, \widehat\mu-\tilde \mu', \gamma_B$ playing the roles of $U, V, \mu, \gamma$, respectively.
	We obtain a $\gamma_B$-skew-matching $\sigma'\trianglelefteq\widehat\mu-\tilde \mu'$ in $H^\leftrightarrow$ of weight 
\begin{align} \label{eq:W_sigma_ap}
		W(\sigma')=(1+\gamma_B)W(\widehat\mu-\tilde \mu'),
\end{align}
 with $\mathcal A(\sigma')\subseteq V$.
 In particular, this means that $\sigma'^2(y) = 0$ for all $y \in V$.
 Using $\gamma_B<1$ and \eqref{eq:W_sigma_ap} we obtain, for all $y\in V$, that
 \begin{align}\label{eq:sigma_ap_y}
 	\sigma'(y)=\widehat\mu(y)-\tilde \mu'(y).
 \end{align}
 
	Now we shall `glue' the skew-matchings together and check we indeed have the required properties.
	Set $\sigma''_B:= \sigma_B'+\sigma+\sigma'$.
	We now claim that $\sigma_A, \sigma''_B$ are the required skew-matchings, i.e. they satisfy \ref{item:completion-out1}--\ref{item:completion-moreover}.
	
	To see \ref{item:completion-out1}, from \eqref{equation:completion-sigmaABtrianglemu} we have $\sigma_A'+\sigma_B'\trianglelefteq  \mu'$, from \eqref{equation:completion-sigawidehattildemuprime} we have $\sigma + \widehat{\sigma}_A \trianglelefteq \tilde \mu'$, and from the definition of $\sigma'$ we have $\sigma' \trianglelefteq \widehat \mu - \tilde \mu'$.
	Therefore, by \eqref{eq:widehat_mu}, we have
	\[\sigma_A+\sigma''_B = (\sigma_A' + \widehat \sigma_A) + (\sigma_B' + \sigma + \sigma') 
	\trianglelefteq \mu' + \tilde \mu' + (\widehat{\mu} - \tilde \mu') = \mu\]
	so \ref{item:completion-out1} holds.
	Item \ref{item:completion-out2} follows from \eqref{eq:lem-compl-WeightOfSigma_A}.
	
	To see \ref{item:completion-out3}, we estimate
	\begin{align*}
		W(\sigma''_B)
		& = W(\sigma_B')+W(\sigma)+W(\sigma')\\
		& \overset{\eqref{equation:completion:wsigmab'}}{=} 2W(\mu') - W(\sigma_A')+W(\sigma)+W(\sigma') \\
		& \overset{\eqref{eq:W_mu_tilde_ap}}{=} 2W(\mu') - W(\sigma_A')+2 W(\tilde \mu') - W(\widehat \sigma_A)+W(\sigma') \\
		& \overset{\eqref{eq:W_sigma_ap}}{=}
		 2W(\mu') - W(\sigma_A')+2 W(\tilde \mu') - W(\widehat \sigma_A)+(1+\gamma_B)W(\hat \mu - \tilde \mu') \\
		& \overset{\eqref{eq:widehat_sigma_A}}{\ge} 2W(\mu')+W(\widehat\mu)-W(\sigma_A)\\
		& = \sum_{x\in U}\mu'(x)+\sum_{y\in V}(\mu'(y)+\widehat\mu(y))-W(\sigma_A)\\
		&		\overset{\eqref{eq:for_mu'}}{=} \sum_{x\in U}\min\{w(\oriented{ux}), \mu(x)\}+\sum_{y\in V}(\mu'(y)+\widehat\mu(y))-W(\sigma_A)\\
		& \overset{\eqref{eq:widehat_mu}}{=} \sum_{x\in U}\min\{w(\oriented{ux}), \mu(x)\} + \sum_{y\in V}\mu(y)-W(\sigma_A) \\
		& \overset{\eqref{eq:mu'_y_leq_w_uy}}{=} \sum_{x\in U}\min\{w(\oriented{ux}), \mu(x)\} + \sum_{y\in V}\min\{ w(\oriented{uy}), \mu(y)\}-W(\sigma_A) \\
		& = \deg_w(u, \mu)-W(\sigma_A),
	\end{align*}	
	which together with the trivial $W(\sigma''_B) \geq 0$ gives \ref{item:completion-out3}.
	Property \ref{item:completion-out4} follows from \eqref{equation:completion-anchorofsAisinU}.
	
	Now we verify \ref{item:completion-out5}.
	Our definitions of $\sigma$ and $\sigma'$ imply that $\mathcal{A}(\sigma), \mathcal{A}(\sigma') \subseteq V$.
	Hence, for $x \in U$, we have that $(\sigma''_B)^1(x) = (\sigma'_B)^1(x) \leq w(\oriented{ux})$, where the inequality follows because $\mathcal{A}(\sigma'_B)$ fits in the $w$-neighbourhood of $u$.
	For $y \in V$, from \ref{item:completion-out1} and \eqref{eq:mu'_y_leq_w_uy} we have $(\sigma''_B)^1(y) \leq \mu(y) \leq w(\oriented{uy})$.
	Thus \ref{item:completion-out5} holds.
	
	Finally, we verify \ref{item:completion-moreover}.
	For any $y\in V$, we have 
	\begin{align*}
		\sigma''_B(y)+\sigma_A(y)&=\sigma_A'(y)+\sigma_B'(y)+\widehat\sigma_A(y)+\sigma(y)+\sigma'(y)\\
		& \overset{\eqref{eq:moreover_sigma_mu_bar},\eqref{eq:mu_tilde_ap_y},\eqref{eq:sigma_ap_y}}{=} \mu'(y)+\tilde \mu'(y)+\widehat\mu(y)-\tilde \mu'(y)
		 = \mu'(y)+ \widehat\mu(y)\overset{\eqref{eq:widehat_mu}}{=}\mu(y),
	\end{align*}
	which implies \ref{item:completion-moreover}.
	This concludes the proof of \Cref{lem:new-completion}.
\end{proof}

\subsection{Proof of the Greedy Lemmas}
Now we shall prove the `Greedy Lemmas', i.e. Lemmas \ref{prop:weighted-greedy-1}--\ref{prop:weighted-greedy-3}.

\begin{proof}[Proof of \Cref{prop:weighted-greedy-1}]
	We shall define the required skew-matching $\sigma'_A$ by a greedy iterative process.
	We will define a sequence $\sigma'_0, \dotsc, \sigma'_t$ of $\gamma_A$-skew-matchings, where $\sigma'_i$ will be obtained from $\sigma'_{i-1}$ by increasing the value in exactly one edge.
	Then $\sigma'_A$ will be the final skew-matching obtained at the end of this procedure.
	
	The following claim ensures we can carry one iterative step of the process.
	
	\begin{claim}
		Let $\sigma'$ be a $\gamma_A$-skew-matching supported in $H[V,U]$ with $\mathcal{A}(\sigma') \subseteq V$.
		If $W(\sigma') < (1 + \gamma_A)\kappa$, then there exists $x \in N_w(u) \cap V$ such that $w(\oriented{ux})>\sigma_A(x)+\sigma_B(x)+\sigma'(x)$, and $y\in N_w(x)\cap U$ that is not covered by $\sigma_A+\sigma_B+\sigma_A'$, i.e. that $\sigma_A(y)+\sigma_B(y)+\sigma_A'(y)<1$.
	\end{claim}
	
	\begin{proofclaim}
		Suppose the desired $x$ does not exist.
		Then, for all $x \in V$ we have
		$w(\oriented{ux}) \leq \sigma_A(x) + \sigma_B(x) + \sigma'(x)$.
		Taking the sum over all $x \in V$ and using \ref{item:greedy-1-condition-1} we get $\kappa \leq \sum_{x \in V} \sigma'(x) = W(\sigma')/(1 + \gamma_A)$, a contradiction.
		Now suppose the desired $y$ does not exist.
		Then for all $y \in N_w(x) \cap U$ we have $\sigma_A(y) + \sigma_B(y) + \sigma'(y) \geq 1$.
		Taking the sum over all $y \in N_w(x) \cap U$, we obtain $|N_w(x) \cap U| \leq \sum_{y \in N_w(x) \cap U } ( \sigma_A(y) + \sigma_B(y) + \sigma'(y) ) \leq \sum_{y \in U } ( \sigma_A(y) + \sigma_B(y) + \sigma'(y) )$.
		Together with \ref{item:greedy-1-condition-2} we obtain $\kappa \gamma_A \leq \sum_{y \in U } \sigma'(y) = \gamma_A W(\sigma') / (1 + \gamma_A)$, a contradiction.
	\end{proofclaim}

	Initially, let $\sigma'_0$ be the identically zero $\gamma_A$-skew-matching.
	Next, given $i \geq 0$, suppose that we are given $\sigma'_i$, which is supported on $H[V,U]$, with $\mathcal{A}(\sigma'_i) \subseteq V$, which is disjoint from $\sigma_A+\sigma_B$, and such that $\sigma_A(x) + \sigma_B(x) + \sigma'_{i}(x) \leq w(\oriented{ux})$ for all $x \in V$.
	We run the following algorithm:
	\begin{enumerate}
		\item If $W(\sigma'_i) \geq (1 + \gamma_A) \kappa$, then we set $\sigma'_A := \sigma'_i$ and we finalise the construction.
		\item Otherwise, we have $W(\sigma'_i) < (1 + \gamma_A) \kappa$.
		By the claim, there exists $x \in V$ and $y \in N_w(x) \cap U$ such that $w(\oriented{ux}) > \sigma_A(x) + \sigma_B(x) + \sigma'_i(x)$ and $\sigma_A(y) + \sigma_B(y) + \sigma'_i(y) < 1$.
		We define $\sigma'_{i+1}$ from $\sigma'_i$ by increasing the weight maximally on $\oriented{xy}$ such that $\sigma'_{i+1}$ is still a $\gamma_A$-skew-matching disjoint from $\sigma_A+\sigma_B$, and such that $\sigma_A(x) + \sigma_B(x) + \sigma'_{i+1}(x) \leq w(\oriented{ux})$.
	\end{enumerate}
	Observe that this process strictly increases the weight of the current matching in each iteration, and further, it chooses each edge $\oriented{xy}$, for $x \in V$ and $y \in U$, at most once.
	Hence the process must stop in at most $|U||V|$ steps.
	
	The $\gamma_A$-skew-matching $\sigma'_A$, by construction, satisfies \ref{item:greedy-1-out-1}--\ref{item:greedy-1-out-3}.
	It remains to verify \ref{item:greedy-1-out-4}, for which we need to verify that \ref{itdef:disjoint}--\ref{itdef:anchorpartition} hold for $(\sigma_A + \sigma'_A, \sigma_B)$ and $\oriented{uv}$.
	Property \ref{itdef:disjoint} follows directly from our construction, since $\sigma'_A$ is disjoint from $\sigma_A + \sigma_B$.
	Property \ref{itdef:anchord} is true by assumption, because $(\sigma_A, \sigma_B)$ is a $(\gamma_A, \gamma_B)$-skew-matching pair anchored in $\oriented{uv}$.
	To see \ref{itdef:anchorc}, we need to prove that $\sigma_A+\sigma'_A$ fits in the $w$-neighbourhood of $u$.
	Note that for every $x \in V(H)$ for which $\sigma'^1_A(x) > 0$, by construction we have
	\[ w(\oriented{ux}) \geq \sigma'_A(x)+\sigma_A(x)+\sigma_B(x). \]
	Together with the fact that $(\sigma_A, \sigma_B)$ is a $(\gamma_A, \gamma_B)$-skew-matching pair anchored in $\oriented{uv}$, this gives \ref{itdef:anchorc}.
	Property \ref{itdef:anchorpartition} also follows from this inequality and the fact that $(\sigma_A, \sigma_B)$ is a $(\gamma_A, \gamma_B)$-skew-matching pair anchored in $\oriented{uv}$.
	This finishes the proof.
\end{proof}

\begin{proof}[Proof of \Cref{prop:weighted-greedy-2}]
	We first shall reserve some space for the anchor of $\sigma_A'$.
	We do this by following a greedy process that defines a sequence of functions $A_0, A_1, \dotsc$, each from $V(H)$ to $[0,1]$, and supported on $V$ (i.e. we will have $A_i(x) = 0$ for all $x \notin V$).
	The next claim will guarantee that we can carry one step of the process.
	
	\begin{claim}
		Let $A: V(H) \rightarrow [0,1]$ be a function supported on $V$.
		Suppose that $\sum_{x \in V} A(x) < \kappa$.
		Then there exists $x \in V$ such that $w(\oriented{ux}) > \sigma_A(x)+\sigma_B(x)+A(x)$.
	\end{claim}
	
	\begin{proofclaim}
		Suppose otherwise.
		Then we have $w(\oriented{ux}) \leq \sigma_A(x)+\sigma_B(x)+A(x)$ for all $x \in V$.
		Summing over all $x \in V$, we obtain $\deg_w(u, V) \leq \sum_{x \in V}(\sigma_A(x)+\sigma_B(x)) + \sum_{x \in V} A(x)$, which together with $\sum_{x \in V} A(x) < \kappa$ contradicts \ref{item:greedy-2-condition-1}.
	\end{proofclaim}
	
	Let $A_0: V(H) \rightarrow [0,1]$ be the identically-zero function.
	Next, given $i \geq 0$ and a function $A_i: V(H) \rightarrow [0,1]$ supported on $V$, we execute the following process:
	\begin{enumerate}
		\item If $\sum_{x \in V} A_i(x) = \kappa$; then set $A := A_i$, and halt the process.
		\item Otherwise, by the claim, there exists $x \in V$ such that $w(\oriented{ux}) > \sigma_A(x)+\sigma_B(x)+A(x)$.
		Define $A_{i+1}: V(H) \rightarrow [0,1]$ from $A_i$ by increasing the value maximally on $x$, subject to $\sigma_A(x)+\sigma_B(x)+A_{i+1}(x) \leq w(\oriented{ux})$ and $\sum_{x \in V}A_{i+1}(x) \leq \kappa$.
	\end{enumerate}
	Since each $x \in V$ can be chosen at most once during the process, we must reach $\sum_{x \in V}A_i(x) = \kappa$ after at most $|V|$ steps.
	Thus the process ends up defining $A: V(H) \rightarrow [0,1]$ supported on $V$, such that $\sum_{x \in V} A(x) = \kappa$.
	
	Now we shall use $A$ to define our desired skew-matching $\sigma'_A$.
	To do this, we will also use an iterative process, constructing skew-matchings $\sigma'_0, \sigma'_1, \dotsc$ and functions $A'_0, A'_1, \dotsc$, and we will maintain the invariants $\sum_{x \in V}(A'_i(x) + \sigma'^1_i(x)) = \kappa$ and $\mathcal{A}(\sigma'_i) \subseteq V$ during the process.
	Again, we state a claim ensuring we can carry with one iteration of the process.

	\begin{claim}
		Suppose $A: V(H) \rightarrow [0,1]$ is supported on $V$ and such that $\sum_{z \in V} A(z) > 0$.
		Suppose $\sigma'$ is a $\gamma_A$-skew-matching such that $\sum_{z \in V} (A(z)+\sigma'^1(z)) = \kappa$ and $\mathcal{A}(\sigma') \subseteq V$.
		Then there exists $x \in V$ and $y \in N_w(x) \cap (U \cup V)$ with $A(x) > 0$ and $\sigma_A(y) + \sigma_B(y) + A(y) + \sigma'(y) < 1$.
	\end{claim}
	
	\begin{proofclaim}
		Select any $x \in V$ with $A(x) > 0$.
		Aiming at a contradiction, suppose for all $y \in N_w(x) \cap (U \cup V)$ we have $\sigma_A(y) + \sigma_B(y) + A(y) + \sigma'(y) \geq 1$.
		Taking the sum of $\sigma_A(y) + \sigma_B(y) + A(y) + \sigma'(y)$ over all $y \in N_w(x) \cap (U \cup V)$, we get
		\begin{equation}
			\sum_{y \in N_w(x) \cap (U \cup V)} (\sigma_A(y) + \sigma_B(y) + A(y) + \sigma'(y))
			\geq |N_w(x) \cap (U \cup V)|.
			\label{equation:greedy2-iteration2}
		\end{equation}
		On the other hand, from \ref{item:greedy-2-condition-2} we get
		\begin{align*}
			|N_w(x) \cap (U \cup V)|
			& \geq (1 + \gamma_A) \kappa + \sum_{y \in U \cup V} (\sigma_A(y) + \sigma_B(y)) \\
			& = \gamma_A \kappa + \sum_{z \in V} (A(z)+\sigma'(z)) + \sum_{z \in U \cup V} (\sigma_A(z) + \sigma_B(z)) \\
			& = \gamma_A \kappa + \sum_{z \in V}\sigma'^1(z) + \sum_{z \in U \cup V} (\sigma_A(z) + \sigma_B(z) + A(z)) \\
			& > (1 + \gamma_A) \sum_{z \in V}\sigma'^1(z) + \sum_{z \in U \cup V} (\sigma_A(z) + \sigma_B(z) + A(z)),
		\end{align*}
		where in the last step we used $\sum_{z \in V}\sigma'^1(z) < \kappa$.
		From $\mathcal{A}(\sigma') \subseteq V$, we get that 
		\[ (1 + \gamma_A) \sum_{z \in V}\sigma'^1(z) = \sum_{z \in V} \sum_{x \in N(z)} \sigma'(\oriented{xz}) = \sum_{z \in U \cup V} \sigma'(z). \]
		Combining the last two inequalities we obtain the desired contradiction with \eqref{equation:greedy2-iteration2}, which proves the claim.
	\end{proofclaim}
	
	The process is described as follows.
	Initially, let $\sigma'_0$ be the identically-zero skew-matching, and let $A'_0 := A$.
	This clearly satisfies $\sum_{x \in V}(A'_0(x) + \sigma'^1_0(x)) = \sum_{x \in V} A(x) = \kappa$ and $\mathcal{A}(\sigma'_0) \subseteq V$, as required.
	
	Next, given $i \geq 0$ and $A'_i, \sigma'_i$ with $\sum_{x \in V}(A'_i(x) + \sigma'^1_i(x)) = \kappa$ and $\mathcal{A}(\sigma'_i) \subseteq V$, we do the following.
	\begin{enumerate}
		\item If $\sum_{x \in V} A'_i(x) = 0$, then we let $A' := A'_i$ and $\sigma'_A := \sigma'_i$, and halt the process.
		\item Otherwise, by the previous claim there exists $x \in V$ and $y \in N_w(x) \cap (U \cup V)$ with $A_i(x) > 0$ and $\sigma_A(y) + \sigma_B(y) + A_i(y) + \sigma'_i(y) < 1$.
		
		Let $\delta$ be maximum such that $\delta \leq A_i(x)$ and $\gamma_A \delta \leq 1 - (\sigma_A(y) + \sigma_B(y) + A_i(y) + \sigma'_i(y))$.
		We define $A_{i+1}$ from $A_i$ by decreasing the value of $A_i(x)$ in $\delta$, and we define $\sigma'_{i+1}$ from $\sigma'_i$ by adding $(1 + \gamma_A) \delta$ to the weight of $\oriented{xy}$.
		This ensures that $\sigma'^1_i(x)$ increases by $\delta$ and $\sigma'^2_i(y)$ increases by $\gamma_A \delta$; so $\sum_{x \in V}(A'_{i+1}(x) + \sigma'^1_{i+1}(x)) = \kappa$ holds.
	\end{enumerate}
	
	Note that each $\oriented{xy}$ with $x \in V, y \in U \cup V$ is chosen at most once during the iterations, so the process finishes after at most $|V||V \cup U|$ steps.
	So the process halts, and at the end of the process we have obtained a $\gamma_A$-skew-matching $\sigma'_A$, with $\mathcal{A}(\sigma'_A) \subseteq V$ and $\sum_{x \in V} \sigma'^1_A(x) = \kappa$.
	
	We verify the required properties \ref{item:greedy-2-out-1}--\ref{item:greedy-2-out-4}.
	Properties \ref{item:greedy-2-out-1} and \ref{item:greedy-2-out-3} follow straight from the construction; and property \ref{item:greedy-2-out-2} follows from  $\mathcal{A}(\sigma'_A) \subseteq V$ together with $\sum_{x \in V} \sigma'^1_A(x) = \kappa$.
	To see \ref{item:greedy-2-out-4}, we need to verify that \ref{itdef:disjoint}--\ref{itdef:anchorpartition} hold for $(\sigma_A+ \sigma'_A, \sigma_B)$ and $\oriented{uv}$.
	The construction of $\sigma'_A$ implies that $\sigma_A+\sigma'_A$ is a $\gamma_A$-skew-matching and for all $x \in V(H)$ with $\sigma'_A(x) > 0$, we have $\sigma^1_A(x)+\sigma_B(x)+\sigma'^1_A(x) \leq w(\oriented{ux})$.
	This, together with the fact that $(\sigma_A, \sigma_B)$ is a $(\gamma_A, \gamma_B)$-skew-matching anchored in $\oriented{uv}$, readily gives \ref{itdef:disjoint}--\ref{itdef:anchorpartition}, and finishes the proof.
\end{proof}

\begin{proof}[Proof of \Cref{prop:weighted-greedy-3}]	
	The proof proceeds analogously to the proof of~\Cref{prop:weighted-greedy-1}, 
	with the difference that as long as $W(\sigma_A')<(1+\gamma_A)\kappa$,
	 we first find a vertex $y\in U$ such that $\sigma_A(y)+\sigma_B(y)+\sigma_A'(y)<1$ and a vertex $x\in N_w(y)\cap V$ with $w(\oriented{ux})>\sigma_A(x)+\sigma_B(x)+\sigma_A'(x)$. We increment $\sigma_A'$ by putting a maximal weight on $\oriented{xy}$, so that $\sigma_A'$ stays a $\gamma_A$-skew-matching, $\sigma_A'$ is disjoint from $\sigma_A+\sigma_B$ and $\sigma_A'(x)+\sigma_A(x)+\sigma_B(x)\le w(\oriented{ux})$.
	
	 Hence, $\mathcal A(\sigma_A')\subseteq V$, the support edges of $\sigma_A'$ run from $V$ to $U$, the anchor of $\sigma'_A$ fits in the $w$-neighbourhood of $u$, we have $W(\sigma_A')=(1+\gamma_A)\kappa$, and $(\sigma_A+\sigma_A', \sigma_B)$ is a $(\gamma_A, \gamma_B)$-skew-matching anchored in $\oriented{uv}$.	
\end{proof}

\subsection{Proof of the $(k, k/2)$-Lemma}

Now we prove the $(k,k/2)$-Lemma, \Cref{lem:new-k-k/2}.
We shall need the following auxiliary lemma. 
In some easy cases, the fractional matching can be partitioned into two disjoint matchings that will host each one of the skew-matchings from our desired skew-matching pair; the next lemma finds the required skew-matchings in this situation.

\begin{lemma}[Filling disjoint matchings]\label{lem:degreetodisjointmatching}
	Let $(H,w)$ be a weighted graph and $\fmat{uv}\in E(H)$. 
	Let $\mu_1, \mu_2, \mu$ be fractional matchings satisfying $\mu_1+\mu_2\le \mu$, and let $(H,\bar w)$ be the $\mu_1$-truncated weighted graph obtained from $(H,w)$. 
	Then, for any $\gamma_1, \gamma_2>0$, there is a $(\gamma_1,\gamma_2)$-skew-matching pair $(\sigma_1,\sigma_2)$ in $(H,w)$ anchored in $\oriented{uv}$ 
	with
	\stepcounter{propcounter}
	\begin{enumerate}[\upshape{(\Alph{propcounter}\arabic*)},topsep=0.7em, itemsep=0.5em]
		\item \label{item:filling-out1} $W(\sigma_1)= \deg_w(u, \mu_1)$,
		\item \label{item:filling-out2} $W(\sigma_2)=\deg_{\bar w}(v, \mu_2)$, and
		\item \label{item:filling-out3} $\sigma_1+\sigma_2\trianglelefteq \mu$. 
	\end{enumerate}
\end{lemma}

\begin{proof}
	First, we construct $\sigma_1$.
	To do so, we use \Cref{lem:new-combination} (Combination) with
	\begin{center}
		\begin{tabular}{c|c|c|c|c}
			object & $(H, w)$  & $u$ & $\mu_1$ & $\gamma_1$  \\
			\hline
			in place of & $(H, w)$ & $v$ & $\mu$ & $\gamma$ 
		\end{tabular}
	\end{center}
	We obtain a $\gamma_1$-skew-matching $\sigma_1$ in $H^{\leftrightarrow}$ with $\sigma \trianglelefteq \mu_1$ and such that the anchor $\mathcal A(\sigma_1)$ fits in the $w$-neighbourhood of $u$.
	By decreasing the weight if necessary, we can also assume that $W(\sigma_1)=\deg_w(u, \mu_1)$.
	
	Secondly, we build $\sigma_2$.
	We use \Cref{lem:new-combination} (Combination) a second time, now with
		\begin{center}
		\begin{tabular}{c|c|c|c|c}
			object & $(H, \bar w)$  & $v$ & $\mu_2$ & $\gamma_2$  \\
			\hline
			in place of & $(H, w)$ & $v$ & $\mu$ & $\gamma$ 
		\end{tabular}
	\end{center}
	We obtain a $\gamma_2$-skew-matching $\sigma_2$ in $H^{\leftrightarrow}$ such that $\sigma_2\trianglelefteq\mu_2$ and
	$\mathcal A(\sigma_2)$ fits in the $\bar w$-neighbourhood of $v$. 
	Again, decreasing the weight if necessary we can assume that $W(\sigma_2)=\deg_{\bar w}(v, \mu_2)$.
	
	Now we combine both skew-matchings to form a skew-matching pair.
	Let $\sigma_\emptyset$ be the empty skew-matching.
	Applying \Cref{prop:adding-skew-matchings} (recall \Cref{crem:adding-skew-matchings}) with
	\begin{center}
	\begin{tabular}{c|c|c|c|c|c|c|c|c|c}
		object & $(H, w)$ & $\fmat{uv}$ & $\mu_1$ & $\mu_2$ & $\bar w$ & $(\sigma_1, \sigma_\emptyset)$ & $(\sigma_\emptyset, \sigma_2)$ & $\gamma_1$ & $\gamma_2$ \\
		\hline
		in place of & $(G, w)$ & $\fmat{uv}$ & $\mu$ & $\bar \mu$ & $\bar w$ & $(\sigma_A, \sigma_B)$ & $(\bar\sigma_A, \bar\sigma_B)$ & $\gamma_A$ & $\gamma_B$ 
	\end{tabular}
	\end{center}
	 we obtain that $(\sigma_A, \sigma_B)$ is a $(\gamma_A, \gamma_B)$-skew-matching pair in $(H, w)$ anchored in $\oriented {uv}$ with $\sigma_1+\sigma_2\trianglelefteq\mu$. 
	 By construction, it satisfies \ref{item:filling-out1}--\ref{item:filling-out3}.
\end{proof}

\begin{proof}[Proof of \Cref{lem:new-k-k/2}]
Without loss of generality, we may assume that 
\begin{equation}\label{eq:deg_w_c_d}
	\deg_w(c, \mu)=k \quad \quad \text{ and } \quad \quad \deg_w(d, \mu)=k/2,
\end{equation} as we can decrease the weight function $w$ so that it  is satisfied.
A second assumption we may do without loss of generality is 
\begin{equation}\label{eq:assumption-alpha_isnammler}
\alpha_1+\alpha_2\le \frac k2,
\end{equation}
as otherwise we can switch the roles of $\alpha_1, \alpha_2$ and $\beta_1, \beta_2$.

The proof consists of five phases.
In the first phase, we partition the fractional matching into meaningful parts, where we will apply some of the matching lemmas from this section. You can refer to \Cref{fig:k-k2} for illustration of how the fractional matchings are defined.
In the second phase, we can construct the required matchings quickly in some `easy' cases, mostly relying on applications of \Cref{lem:degreetodisjointmatching}.
After this is done, we get some extra assumptions about the weights of our auxiliary matchings.
Now it will not be immediately clear whether the $(\gamma_A, \gamma_B)$-skew-matching pair will be anchored at $\oriented{cd}$ or $\oriented{dc}$; and this decision is made in the third phase.
In the fourth phase, we construct partial skew-matchings inside the meaningful parts.
Finally, in the fifth phase, we combine these partial skew-matchings to form the desired skew-matching pair. \medskip

\noindent \emph{Step 1: Partitioning the fractional matching into meaningful bits.} 
Each part of the fractional matching is connected to vertices $c$ and $d$ in different ways.
To handle them systematically, we divide $\mu$ into `homogeneous bits' --groups with similar properties--, so that they can be treated in the same way.
The fractional matching $\mu_d$ will represent the portion that is strongly connected to vertex $d$ (from both sides) and can be efficiently packed using a skew-matching anchored at $d$.
On the other hand, the fractional matching $\bar{\mu}_d$ will be connected to $d$ from only one side and is further split into $\mu'_d$, which will also be connected to vertex $c$, and $\mu^*$, which will be exclusive or `private' to $d$.
Finally, $\mu_c$ will be the part of the fractional matching that is private to $c$, and it will include sections that are connected to $c$ either from both sides or from only one side.

Let $\mu_d\leq \mu$ be a fractional matching of maximal weight such that $\mu_d(x)\leq w(\oriented{dx})$ for all $x\in V(H)$.
Note that this implies that if $\mu_d(xy) > 0$, then $x, y \in N_w(d)$.
We have 
\begin{align}\label{eq:sizemu_d}
\deg_w(d, \mu_d)=2W(\mu_d)\ge \deg_w(c, \mu_d).
\end{align}
Next, let $(H,w')$ be the $\mu_d$-truncated weighted graph obtained from $(H,w)$ (\Cref{def:truncated-weighted-graph}). 
Set $U:= N_{w'}(d)$ and $V:= V(H)\setminus U$.
We have
\begin{equation} \label{eq:V_mu_d}
	\mu_d(x) = w(\oriented{dx})
\end{equation}
for all $x \in V$.
Also observe that for all $x \in U$, from the definition of $w'$ we have
\begin{equation} \label{equation:kk2-wprimevsmud}
	w'(\oriented{dx}) + \mu_d(x) = w(\oriented{dx}). 
\end{equation}

Now, let $\bar \mu_d\le \mu-\mu_d$ be a maximal fractional matching supported only on edges intersecting the set $U$, and such that
\begin{equation}
	\bar \mu_d(x)\le w'(\oriented{dx})
	\label{equation:barmudw'}
\end{equation}
for all $x\in U$. 
We claim that $\bar \mu_d$ is supported in $H[U, V]$.
Indeed, otherwise there is an edge $xy$, completely contained in $U$, with $\bar \mu_d(\fmat{xy})>0$. 
This implies that $w'(\oriented{dx}) > 0$ and $w'(\oriented{dy}) > 0$; and therefore $w(\oriented{dx}) > \mu_d(x)$ and $w(\oriented{dy}) > \mu_d(y)$; but this contradicts the maximality of $\mu_d$.
Now, the fact that $\bar \mu_d$ runs between $U$ and $V$ implies that
\begin{equation} \label{eq:deg_w_prime_d_bar_mu_d}
	\deg_{w'}(d,\bar\mu_d) = W(\bar\mu_d).
\end{equation}
It also gives
\begin{equation*}
	\bar \mu_d(x) + \mu_d(x) = \min\{w(\oriented{dx}), \mu(x)\}
\end{equation*}
for all $x \in U$, and
\begin{equation}
	\nonumber
	w(\oriented{dx})\overset{\eqref{eq:V_mu_d}}{=}\mu_d(x)\le \bar \mu_d(x) + \mu_d(x) \le  \mu(x)
\end{equation}
for all $x \in V$. This yields
\begin{equation} \label{eq.deg_w_mu_mu_d_bar_mu_d}
	\deg_w(d,\mu_d+\bar\mu_d) = \deg_w(d,\mu) = \frac{k}{2}.
\end{equation}
Moreover, observe that by construction we have $\mu_d(x) = \min\{ \mu_d(x), w(\oriented{dx}) \}$ for all $x \in V(H)$.
Hence, we have
\begin{align}
	\nonumber
	2W(\mu_d)+W(\bar \mu_d)
	&\overset{\eqref{eq:sizemu_d},\eqref{eq:deg_w_prime_d_bar_mu_d}}{=}\deg_w(d, \mu_d)+\deg_{w'}(d, \bar\mu_d)\\ \nonumber
	&= \sum_{x\in V(H)}\mu_d(x)+\sum_{x\in U}\min\{w'(\oriented{dx}),\bar\mu_d(x)\} \\ \nonumber 
	& = \sum_{x \in V} \mu_d(x) + \sum_{x\in U} \left( \mu_d(x) + \min\{w'(\oriented{dx}),\bar\mu_d(x)\} \right) \\ \nonumber 
	& = \sum_{x \in V} \mu_d(x) + \sum_{x\in U} \min\{w'(\oriented{dx})+\mu_d(x),\bar\mu_d(x)+\mu_d(x) \}  \\ \nonumber 
	& \overset{\eqref{equation:kk2-wprimevsmud}}{=} \sum_{x \in V} \mu_d(x) + \sum_{x\in U}  \min\{w(\oriented{dx}),\bar\mu_d(x)+\mu_d(x) \}\\ \nonumber
	& \overset{\eqref{eq:V_mu_d}}{=} \sum_{x \in V} \min\{w(\oriented{dx}), \mu_d(x)+\bar\mu_d(x)\} + \sum_{x\in U}  \min\{w(\oriented{dx}),\bar\mu_d(x)+\mu_d(x) \} \\
	&=\deg_w(d, \mu_d+\bar \mu_d). \label{eq:size_mu_d_bar_mu_d}
\end{align}

Let $\mu_d'\le \bar \mu_d$ be a fractional matching such that, for all $y \in V$,
\begin{equation}\label{eq:lessthanw'}
\mu_d'(y)=\min\{ w'(\oriented{cy}), \bar\mu_d(y)\}.
\end{equation}
Let $\mu^*:= \bar\mu_d-\mu_d'$. 
Since $\mu^\ast, \mu'_d \leq \bar \mu_d$ and $\bar \mu_d$ is supported in $H[U,V]$,
we have that $\mu_d'$ and $\mu^*$ run between $U$ and $V$ as well.
Straight from the definition, we have
\begin{equation}\label{eq:W_mu_star_bar_d_d_prime}
	W(\mu^*) = W(\bar \mu_d) - W(\mu_d').
\end{equation}
Moreover, from $\bar\mu_d = \mu^* + \mu_d'$ we have
\begin{align}\label{eq:degreefromc}
\deg_w(c, \mu_d+\bar \mu_d)\le  \deg_w(c, \mu_d+\mu_d')+W(\mu^*).
\end{align}
We finalise the construction of objects by setting $\mu_c:= \mu-(\mu_d+\bar \mu_d)$. 
Let $(H, \bar w)$ be the $\bar \mu_d$-truncated graph obtained from $(H, w')$.
Observe that $(H, \bar w)$ also corresponds to the $(\mu_d+\bar\mu_d)$-truncated graph obtained from $(H, w)$.

Observe that 
\begin{align}
\nonumber \deg_{\bar w}(c, \mu_c)&=\deg_w(c, \mu)-\deg_w(c, \mu_d+\bar\mu_d)\overset{\eqref{eq:deg_w_c_d}, \eqref{eq:degreefromc}}{\ge}k-\deg_w(c, \mu_d+\mu_d')-W(\mu^*)\\
\nonumber &\overset{\eqref{eq:deg_w_c_d}}{\ge } 2\deg_w(d, \mu)-2\Big(W(\mu_d)+W(\mu_d')\Big)-W(\mu^*)\\
\nonumber & = 2\Big(\deg_w(d, \mu_d)+\deg_{w'}(d, \bar\mu_d)\Big)-2W(\mu_d)-2W(\mu_d')-W(\mu^*)\\
\nonumber &\overset{\eqref{eq:sizemu_d},\eqref{eq:deg_w_prime_d_bar_mu_d}}{=}4W(\mu_d)+2W(\bar\mu_d) -2W(\mu_d)-2W( \bar \mu_d)+W(\mu^*)\\
&=2W(\mu_d)+W(\mu^*).\label{eq:deg(c,mu_c)}
\end{align}

\begin{figure}[ht]
\centering
\begin{tikzpicture}
\definecolor{myblue}{RGB}{31,119,180}
\definecolor{mypurple}{RGB}{148,103,189}
\definecolor{myred}{RGB}{214,39,40}
\definecolor{myteal}{RGB}{23,190,207}
\definecolor{mybrown}{RGB}{140,86,75}

  \def\R{1.5}   
  \def\G{5.0}   
  \def\eps{0.05}  

  \coordinate (Cdown) at (0,0);
  \coordinate (Cup)   at (0,\G);

  \draw (Cdown) circle (\R);
  \draw (Cup)   circle (\R);

\pgfmathsetmacro{\xthree}{\R/2 }

\begin{scope}
  \clip (Cup) circle (\R);
  \clip ($(Cup)+(\xthree,-\R)$) rectangle ($(Cup)+(\R,\R)$);

  \foreach \t in {-4,-3.6,...,6}{
    \draw[green!75!black, thick] ($(Cup)+(\t,-2*\R)$) -- ($(Cup)+(\t+4,2*\R)$);
  }
\end{scope}

\pgfmathsetmacro{\xone}{-\R/2}

\begin{scope}
  \clip (Cdown) circle (\R);
  \clip ($(Cdown)+(\xone,-\R)$) rectangle ($(Cdown)+(\R,\R)$);

  \foreach \t in {-4,-3.6,...,6}{
    \draw[green!75!black, thick] ($(Cdown)+(\t,-2*\R)$) -- ($(Cdown)+(\t+4,2*\R)$);
  }
\end{scope}


\begin{scope}
  \clip (Cup) circle (\R);
  \clip ($(Cup)+(0,-\R)$) rectangle ($(Cup)+(\R,\R)$);

   \foreach \t in {-6,-5.5,...,6}{
    \draw[orange!75!black, thick]
      ($(Cup)+(\t,-3*\R)$) -- ++(135:{6*\R});
  }
\end{scope}

  \draw [black, thick](-\R,0) -- (-\R,\G);
  \draw [myred, thick] ( \R,0) -- ( \R,\G);

\pgfmathsetmacro{\axone}{-\R + (2*\R/4)}
\pgfmathsetmacro{\byone}{sqrt(\R*\R - \axone*\axone)}
\coordinate (Top1) at ($(Cup)+(\axone,\byone)$);
\coordinate (Bot1) at ($(Cdown)+(\axone,-\byone)$);
\draw[myteal, thick] (Top1) -- (Bot1);

\pgfmathsetmacro{\xblack}{\axone - \eps}
\pgfmathsetmacro{\yblack}{sqrt(\R*\R - \xblack*\xblack)}
\coordinate (Topblack) at ($(Cup)+(\xblack,\yblack)$);
\coordinate (Botblack) at ($(Cdown)+(\xblack,-\yblack)$);
\draw[black, thick] (Topblack) -- (Botblack);

\pgfmathsetmacro{\axtwo}{-\R + 2*(2*\R/4)}
\pgfmathsetmacro{\bytwo}{sqrt(\R*\R - \axtwo*\axtwo)}
\coordinate (Top2) at ($(Cup)+(\axtwo,\bytwo)$);
\coordinate (Bot2) at ($(Cdown)+(\axtwo,-\bytwo)$);
\draw [mypurple, thick](Top2) -- (Bot2);

\pgfmathsetmacro{\xbrown}{\axtwo - \eps}
\pgfmathsetmacro{\ybrown}{sqrt(\R*\R - \xbrown*\xbrown)}
\coordinate (Topbrown) at ($(Cup)+(\xbrown,\ybrown)$);
\coordinate (Botbrown) at ($(Cdown)+(\xbrown,-\ybrown)$);
\draw[myteal, thick] (Topbrown) -- (Botbrown);

\pgfmathsetmacro{\axthree}{-\R + 3*(2*\R/4)}
\pgfmathsetmacro{\bythree}{sqrt(\R*\R - \axthree*\axthree)}
\coordinate (Top3) at ($(Cup)+(\axthree,\bythree)$);
\coordinate (Bot3) at ($(Cdown)+(\axthree,-\bythree)$);
\draw[mypurple, thick] (Top3) -- (Bot3);

\pgfmathsetmacro{\xred}{-\R + 3*(2*\R/4) + \eps}  
\pgfmathsetmacro{\yred}{sqrt(\R*\R - \xred*\xred)}

\coordinate (Topred) at ($(Cup)+(\xred,\yred)$);
\coordinate (Botred) at ($(Cdown)+(\xred,-\yred)$);

\draw[myred, thick] (Topred) -- (Botred);

\node at ($(Cup)+(-\R-0.7,0)$) {$y\in V$};
\node at ($(Cdown)+(-\R-0.7,0)$) {$x\in U$};

\node [draw=none, text=myred] at ($(Cdown)+(1.1,2.5)$) {$\mu_d$};
\node [draw=none, text=mypurple] at ($(Cdown)+(0.4,2.5)$) {$\mu_d'$};
\node [draw=none, text=myteal] at ($(Cdown)+(-0.4,2.5)$) {$\mu^*$};
\node [draw=none, text=black] at ($(Cdown)+(-1.1,2.5)$) {$\mu_c$};

\begin{scope}[shift={(3*\R,5.5)}] 
  \begin{scope}
    \clip (-0.3,-0.3) rectangle (0.3,0.3);
    \foreach \t in {-1,-0.8,...,1}{
      \draw[green!50!black, thick] (\t,-0.5) -- ++(45:1.5);
    }
  \end{scope}
  \node[right=4pt] at (0.3,0) {$w(\oriented{dx}),\; w(\oriented{dy})$};

  \begin{scope}[yshift=-1cm]
    \clip (-0.3,-0.3) rectangle (0.3,0.3);
    \foreach \t in {-1,-0.8,...,1}{
      \draw[orange!75!black, thick] (\t,-0.5) -- ++(135:1.5);
    }
  \end{scope}
  \node[right=4pt] at (0.3,-1cm) {$w(\oriented{cy})$};
  \node at (0.85, -2) {$\bar\mu_d=\mu_d'+\mu^*$};
 \end{scope}

\pgfmathsetmacro{\angL}{180}                        
\pgfmathsetmacro{\angBlack}{atan2(\yblack,\xblack)} 
\pgfmathsetmacro{\angOne}{atan2(\byone,\axone)}     
\pgfmathsetmacro{\angBrown}{atan2(\ybrown,\xbrown)} 
\pgfmathsetmacro{\angTwo}{atan2(\bytwo,\axtwo)}     
\pgfmathsetmacro{\angThree}{atan2(\bythree,\axthree)} 
\pgfmathsetmacro{\angRed}{atan2(\yred,\xred)}       
\pgfmathsetmacro{\angR}{0}                           

\newcommand{\drawSeg}[3]{%
  \ifdim #1pt>#2pt
    \draw[#3, line width=1.2pt] ($(Cup)+(#2:\R)$) arc[start angle=#2, end angle=#1, radius=\R];
    \draw[#3, line width=1.2pt] ($(Cdown)+({-(#1)}:\R)$) arc[start angle={-(#1)}, end angle={-(#2)}, radius=\R];
  \else
    \draw[#3, line width=1.2pt] ($(Cup)+(#1:\R)$) arc[start angle=#1, end angle=#2, radius=\R];
    \draw[#3, line width=1.2pt] ($(Cdown)+({-(#2)}:\R)$) arc[start angle={-(#2)}, end angle={-(#1)}, radius=\R];
  \fi
}

\drawSeg{\angL}{\angBlack}{black}
\drawSeg{\angBlack}{\angOne}{myteal}
\drawSeg{\angOne}{\angBrown}{myteal}
\drawSeg{\angBrown}{\angTwo}{mypurple}
\drawSeg{\angTwo}{\angThree}{mypurple}
\drawSeg{\angThree}{\angRed}{myred}
\drawSeg{\angRed}{\angR}{myred}
\newcommand{\drawSegOpp}[3]{%
  \ifdim #1pt>#2pt
    \draw[#3, line width=1.2pt] ($(Cup)+({-(#2)}:\R)$) arc[start angle={-(#2)}, end angle={-(#1)}, radius=\R];
    \draw[#3, line width=1.2pt] ($(Cdown)+(#1:\R)$) arc[start angle=#1, end angle=#2, radius=\R];
  \else
    \draw[#3, line width=1.2pt] ($(Cup)+({-(#1)}:\R)$) arc[start angle={-(#1)}, end angle={-(#2)}, radius=\R];
    \draw[#3, line width=1.2pt] ($(Cdown)+(#1:\R)$) arc[start angle=#1, end angle=#2, radius=\R];
  \fi
}


\drawSegOpp{\angL}{\angBlack}{black}
\drawSegOpp{\angBlack}{\angOne}{myteal}
\drawSegOpp{\angOne}{\angBrown}{myteal}
\drawSegOpp{\angBrown}{\angTwo}{mypurple}
\drawSegOpp{\angTwo}{\angThree}{mypurple}
\drawSegOpp{\angThree}{\angRed}{myred}
\drawSegOpp{\angRed}{\angR}{myred}

\end{tikzpicture}
\caption{A representation of one edge $xy$ and how it is partitioned in different fractional matchings.
Depending on $w(\protect\vv{dx}), w(\protect\vv{dy}), w(\protect\vv{cy})$ the fractional matchings $\mu_d, \mu_d', \mu^*, \mu_c$ will look differently.
Note that $w(\protect\vv{cx})$ has not been represented, as it is irrelevant for the definition of the different fractional matchings. Also  $\mu_d(y)$ may not necessarily be ``covered'' by $N_w(c)$ (i.e. $\mu_d(y)<w(\protect\vv{cy})$). However, in this case $\mu_d'(xy)=0$. Also note that $\mu_c(y)$ may be partially ``covered'' by $N_w(c)$. However, if this is the case, than $\mu^*(xy)=0$.}
\label{fig:k-k2}
\end{figure}
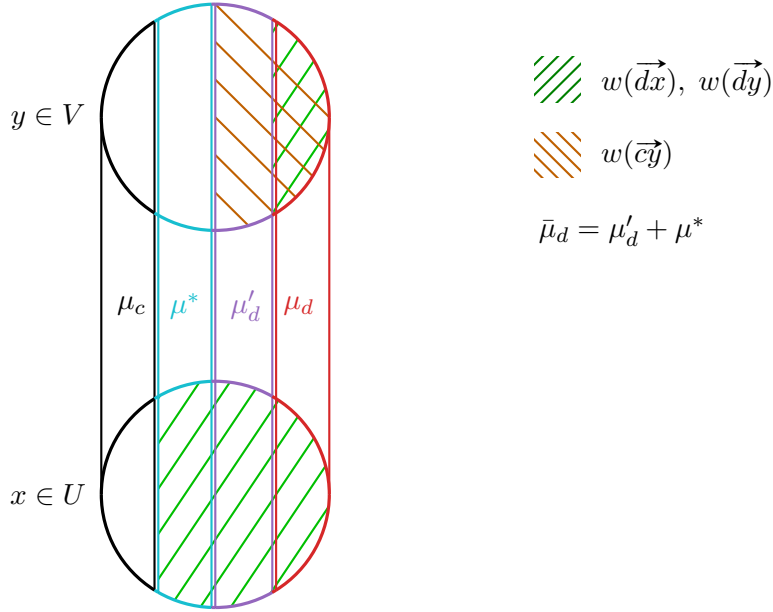

\noindent \emph{Step 2: The simple cases.}
The most delicate part of the proof involves handling the fractional matching $\mu_d'$, as its space might be shared by both skew-matchings: one anchored at $d$ and the other at $c$.
Before diving into this difficulty, let us first get rid of some easier cases, when the fractional matching $\mu_d'$ does not need to be shared by the skew-matchings.
The main outcome after this step is that we will get some numerical information about the weights of the auxiliary matchings, encapsulated in inequalities \eqref{eq:bpositive} and \eqref{eq:apositive}.

The first case we consider is when $2W(\mu_d)+W(\mu^*)\ge \beta_1+\beta_2$ holds.
If this holds, then by~\eqref{eq:deg(c,mu_c)} we obtain that $\deg_{\bar w}(c, \mu_c)\ge \beta_1+\beta_2$.
On the other hand,~\eqref{eq.deg_w_mu_mu_d_bar_mu_d} and~\eqref{eq:assumption-alpha_isnammler} together imply that $\deg_w(d, \mu_d+\bar\mu_d)\ge \alpha_1+\alpha_2$.
We recall that, by definition of $\mu_c$, we have $\mu_c + \mu_d+\bar \mu_d \trianglelefteq \mu$; and also that $(H, \bar w)$ is the $(\mu_d+\bar\mu_d)$-truncated graph obtained from $(H, w)$. 
Thus we can apply \Cref{lem:degreetodisjointmatching} (Filling disjoint matchings) with
\begin{center}
\begin{tabular}{c|c|c|c|c|c|c|c|c}
	object & $(H, w)$ & $\fmat{dc}$ & $\mu_d+\bar\mu_d$ & $\mu_c$ & $\mu$ & $\bar w$ & $\gamma_A$ & $\gamma_B$ \\
	\hline
	in place of & $(H,w)$ & $\fmat{uv}$ & $\mu_1$ & $\mu_2$ & $\mu$ & $\bar w$ & $\gamma_1$ & $\gamma_2$ 
\end{tabular}
\end{center}
to obtain the required $(\gamma_A, \gamma_B)$-skew-matching pair $(\sigma_A, \sigma_B)$ with $\sigma_A+\sigma_B\trianglelefteq \mu_c+\mu_d+\bar\mu_d = \mu$, and we are done in this case.
So, from now on, we may assume that 
\begin{equation}\label{eq:bpositive}
2W(\mu_d)+W(\mu^*) <\beta_1+\beta_2.
\end{equation}

The second simple case we consider is when $2W(\mu_d)\ge \alpha_1+\alpha_2$.
If this holds, then we recall \eqref{eq:sizemu_d} and we use it to select $\mu''_d\le \mu_d$ to be such that $\deg_w(d, \mu''_d)=2W(\mu''_d)=\alpha_1+\alpha_2$.
Also, let $\mu''_c=\mu-\mu''_d$.
Let $(H, w'')$ be the $\mu''_d$-truncated graph obtained from $(H,w)$.
By applying \Cref{prop:adding-degree-truncated} with $(H,w), \mu_d'', \mu, w''$ playing the role of $(G,w), \mu', \mu, w'$,  we see that
\[\deg_{w''}(c, \mu''_c)=\deg_w(c, \mu)-\deg_w(c, \mu''_d)\ge k-2W(\mu''_d)=\beta_1+\beta_2.\]
Hence, we can apply \Cref{lem:degreetodisjointmatching} with
\begin{center}
	\begin{tabular}{c|c|c|c|c|c|c|c|c}
		object & $(H, w)$ & $\fmat{dc}$ & $\mu''_d$ & $\mu''_c$ & $\mu$ & $w''$ & $\gamma_A$ & $\gamma_B$\\
		\hline
		in place of & $(H,w)$ & $\fmat{uv}$ & $\mu_1$ & $\mu_2$ & $\mu$ & $\bar w$ & $\gamma_1$ & $\gamma_2$
	\end{tabular}
\end{center}
to obtain the required $(\gamma_A, \gamma_B)$-skew-matching pair $(\sigma_A, \sigma_B)$ with  $\sigma_A+\sigma_B\trianglelefteq \mu''_c+\mu''_d = \mu$. 
Thus, from now on, we can assume that $2W(\mu_d)<\alpha_1+\alpha_2$ holds.

Finally, the last simple case we will consider is when $\alpha_1+\alpha_2\le 2W(\mu_d)+W(\mu^*)$.
Since $2W(\mu_d)<\alpha_1+\alpha_2$, we can select $\hat \mu_d\le \mu^*$ to be such that 
\begin{equation}
	2W(\mu_d)+W(\hat\mu_d)=\alpha_1+\alpha_2.
	\label{equation:2wmudwhatmud}
\end{equation}

We recall that $(H,w')$ is the $\mu_d$-truncated weighted graph obtained from $(H,w)$.
Then, we can apply \Cref{lem:degreetodisjointmatching} with
\begin{center}
	\begin{tabular}{c|c|c|c|c|c|c|c|c}
		object &  $(H, w)$ & $\fmat{dc}$ & $\mu_d$ & $\mu_d'$ & $\mu_d+\mu_d'$ & $w'$ & $\gamma_A$ & $\gamma_B$ \\
		\hline
		in place of & $(H,w)$ & $\fmat{uv}$ & $\mu_1$ & $\mu_2$ & $\mu$ & $\bar w$ & $\gamma_1$ & $\gamma_2$
	\end{tabular}
\end{center}
to obtain a $(\gamma_A, \gamma_B)$-skew-matching pair $(\sigma_A', \sigma_B')$ 
that is anchored in $\oriented{dc}$ with
\begin{equation}
	W(\sigma_A')=\deg_w(d, \mu_d)=2W(\mu_d) \label{equation:wsigmaa'}
\end{equation}
and 
\begin{equation}\label{equation:wsigmab'}
W(\sigma_B')=\deg_{w'}(c, \mu_d')\ge \deg_w(c, \mu_d+\mu_d')-2W(\mu_d)=\deg_w(c, \mu_d+\mu_d')-W(\sigma_A').
\end{equation}
Moreover, we have
\begin{equation}
	\sigma'_A + \sigma'_B \trianglelefteq \mu_d + \mu'_d. \label{equation:sigmaprimeaplussigmaprimebtrianglemudmuprimed}
\end{equation}

Next, let $(H, \hat w)$ be the $\mu_d'$-truncated weighted graph obtained from $(H,w')$.
We recall that $U = N_{w'}(d)$, and since $\hat w \leq w'$ we get $N_{\hat{w}}(d) \subseteq U$. 

We claim that
\begin{equation}\label{eq:cSeesMu*OnlyFromOneSide}
	V(\mu^\ast) \cap N_{\hat w}(c, V) = \emptyset.
\end{equation}
Indeed, from \eqref{eq:lessthanw'} we see that for any $y\in V$ with $\mu^*(y)=\bar\mu_d(y)-\mu_d'(y)>0$, we get $\hat w(\oriented{cy}) \leq w'(\oriented{cy})-\mu_d'(y)=\mu_d'(y)-\mu_d'(y)=0$.

We also recall that, by the definition of $\hat \mu_d$ and $\bar \mu_d$, we have $\hat \mu_d(x) \leq \bar \mu_d(x) - \mu'_d(x) \leq w'(\oriented{dx}) - \mu'_d(x) \leq \hat w(\oriented{dx})$ for all $x \in U$. 
Therefore, we have
\begin{equation}
	\deg_{\hat w}(d, \hat \mu_d) = \sum_{x \in U} \min\{ \hat w(\oriented{dx}), \hat \mu_d(x) \} = \sum_{x \in U} \hat \mu_d(x) = W(\hat\mu_d),
	\label{equation:deghatwdhatmud}
\end{equation}
where in the last equality we used that $\hat \mu_d \leq \bar \mu_d$, and that $\bar \mu_d$ is supported in $H[U,V]$.

We apply \Cref{lem:new-combination} (Combination) with
\begin{center}
	\begin{tabular}{c|c|c|c|c}
		object &  $(H, \hat w)$ & $d$ & $\hat \mu_d$ & $\gamma_A$ \\
		\hline
		in place of & $(H,w)$ & $v$ & $\mu$ & $\gamma$
	\end{tabular}
\end{center}
The outcome is an $\gamma_A$-skew-matching $\sigma_A''$ in $(H, \hat w)$, such that $\sigma''_A \trianglelefteq \hat\mu_d$, and with weight $W(\sigma_A'') \geq \deg_{\hat{w}}(d, \hat \mu_d)$.
Using \eqref{equation:deghatwdhatmud}, \eqref{equation:2wmudwhatmud} and \eqref{equation:wsigmaa'}, we obtain
\begin{equation} \label{eq:sigma_Ainhatmu_d}
W(\sigma_A'') \geq \deg_{\hat{w}}(d, \hat \mu_d) = W(\hat \mu_d) = \alpha_1 + \alpha_2 - 2 W(\mu_d) = \alpha_1 + \alpha_2 - W(\sigma'_A), 
\end{equation}
so we can assume, decreasing the weight of $\sigma''_A$ if necessary, that in fact $W(\sigma''_A) = \alpha_1 + \alpha_2 - W(\sigma'_A)$ holds.

Now we wish to apply \Cref{prop:adding-skew-matchings} (recall \Cref{crem:adding-skew-matchings}) with
\begin{center}
	\begin{tabular}{c|c|c|c|c|c|c|c|c|c|c}
		object &  $(H, w)$ & $(H, \hat w)$ & $\fmat{dc}$ & $\mu_d + \mu'_d$ & $\hat\mu_d$ & $\hat w$ & $(\sigma'_A, \sigma'_B)$ & $(\sigma''_A, \sigma_\emptyset)$ & $\gamma_A$ & $\gamma_B$ \\
		\hline
		in place of & $(G,w)$ & $(G, \bar w)$ & $\fmat{uv}$ & $\mu$ & $\bar \mu$ & $\bar w$ & $(\sigma_A, \sigma_B)$ & $(\bar \sigma_A, \bar \sigma_B)$ & $\gamma_A$ & $\gamma_B$
	\end{tabular}
\end{center}
Let us quickly verify that our choice of objects satisfies the required conditions.
We have that $\mu_d + \mu'_d$ and $\hat\mu_d$ are disjoint, because $\hat\mu_d \leq \mu^\ast = \bar \mu_d - \mu'_d$ and $\bar \mu_d \leq \mu - \mu_d$.
We observe that the choice $\hat w$ (as the $\mu'_d$-truncated weighted graph from $(H,w')$) also implies that $(H, \bar w)$ is the $(\mu_d + \mu'_d)$-truncated weighted graph  obtained from $(H,w)$, as required. 
Given this, then \ref{item:addingskew-in-1} follows from \eqref{equation:sigmaprimeaplussigmaprimebtrianglemudmuprimed}, and \ref{item:addingskew-in-2} follows from the choice of $\sigma''_A$.

From the application of \Cref{prop:adding-skew-matchings} we obtain a $(\gamma_A, \gamma_B)$-skew-matching pair $(\sigma_A, \sigma_B')$ in $(H, w)$ anchored in $\oriented{dc}$, where $\sigma_A:=\sigma_A'+\sigma_A''$.

Now, we observe that by applying \Cref{prop:adding-degree-truncated} twice (first with $(H, w), \mu_d', \mu_d+\mu_d'+\hat\mu_d$, and $w'$ playing the role of $(G,w), \mu', \mu,$, and $w'$, and then with $(H, w'), \mu_d', \mu_d'+\hat\mu_d$, and $\hat w$ playing the role of $(G,w), \mu', \mu$, and $w'$)  we obtain
\begin{align*}
\deg_w(c, \mu_d+\mu_d'+\hat\mu_d)&=\deg_{w'}(c, \mu_d'+\hat\mu_d)+\deg_w(c, \mu_d)\\
&= \deg_{\hat w}(c, \hat \mu_d) +\deg_{w'}(c, \mu_d')+\deg_w(c,\mu_d)\\
&\le \deg_{\hat w}(c, \hat\mu_d) +\deg_{w'}(c, \mu_d')+ 2W(\mu_d),
\end{align*}
where in the last inequality we used $\deg_w(c, \mu_d)\le 2W(\mu_d)$.
Next, from \eqref{eq:cSeesMu*OnlyFromOneSide} and $\hat \mu_d \leq \mu^\ast \leq \bar \mu_d \subseteq H[U,V]$, we get that for every edge $uv \in E(H)$ with $u \in U$, $v \in V$ and $\hat \mu_d(uv) > 0$, it must hold that $\hat{w}(\oriented{cv}) = 0$.
Hence, we obtain that $\deg_{\hat w}(c, \hat\mu_d) \leq W(\hat \mu_d)$.
Combining this with the last displayed inequality, we obtain
\begin{align*}
	\deg_w(c, \mu_d+\mu_d'+\hat\mu_d)
\leq W(\hat\mu_d)+\deg_{w'}(c, \mu_d')+2W(\mu_d) = W(\sigma_B')+\alpha_1+\alpha_2,
\end{align*}
where the last equality is now direct by using \eqref{eq:sigma_Ainhatmu_d}, \eqref{equation:wsigmab'} and \eqref{equation:wsigmaa'}.

Last, let $(H, w'')$ be the $(\mu_d+\mu_d'+\hat\mu_d)$-truncated graph obtained from $(H,w)$.
Using our upper bound for $\deg_w(c, \mu_d+\mu_d'+\hat\mu_d)$, we see that
\[\deg_{w''}(c, \mu-\hat\mu_d-\mu_d'-\mu_d)=\deg_w(c, \mu)-\deg_w(c, \mu_d+\mu_d'+\hat\mu_d)\ge k-W(\sigma_B')-\alpha_1-\alpha_2.\] 
We apply \Cref{lem:new-combination} (Combination) with
\begin{center}
	\begin{tabular}{c|c|c|c|c}
		object &  $(H,  w'')$ & $c$ & $\mu-(\mu_d+\mu_d'+\bar\mu_d)$ & $\gamma_B$ \\
		\hline
		in place of & $(H,w)$ & $v$ & $\mu$ & $\gamma$
	\end{tabular}
\end{center}
to obtain a $\gamma_B$-skew-matching $\sigma_B''$ satisfying $\sigma''_B \trianglelefteq \mu-(\mu_d+\mu_d'+\bar\mu_d)$, and that has weight $W(\sigma_B'')=k-W(\sigma_B')-\alpha_1-\alpha_2 = \beta_1+\beta_2-W(\sigma_B')$. 
Combined with \Cref{prop:adding-skew-matchings} with
\begin{center}
	\begin{tabular}{c|c|c|c|c|c|c|c}
		object &  $(H,w)$ & $\fmat{dc}$ & $\mu_d+\mu_d'+ \hat\mu_d$ & $\mu-(\mu_d+\mu_d'+ \hat\mu_d)$ & $w''$ & $(\sigma_A, \sigma_B')$ & $( \emptyset,\sigma_B')$ \\
		\hline
		in place of & $(G, w)$ & $\fmat{uv}$ & $\mu$ & $\bar\mu$ & $\bar w$ & $(\sigma_A, \sigma_B)$ & $(\bar \sigma_A, \bar\sigma_B)$
	\end{tabular}
\end{center}
we get the desired $(\gamma_A, \gamma_B)$-skew-matching pair $(\sigma_A, \sigma_B)$, with $\sigma_B:=\sigma_B'+\sigma_B''$ and 
$\sigma_A+\sigma_B\trianglelefteq \mu$.

Hence, from now on we may assume that 
\begin{equation}\label{eq:apositive}
2W(\mu_d)+W(\mu^*)<\alpha_1+\alpha_2.
\end{equation}

\noindent \emph{Step 3: Deciding how to anchor the skew-matching pair.}
To manage this effectively, we will use the sophisticated \Cref{lem:new-completion} (Completion) in the next step. At this point, we determine which configuration satisfies the conditions, which, in turn, dictates which skew-matching should be anchored to each vertex.
The fact that $\mu_d'$ governs how the skew-matching pair is anchored, while the other fractional matchings are irrelevant, follows from the property established below in~\eqref{eq:degmuc=degmud}, which states that we can fit the same amount of skew-matching in $\mu_c$ as in $\mu_d$.

Let
\begin{align*}
	a & := \alpha_1+\alpha_2-2W(\mu_d)-W(\mu^*), \\
	b & := \beta_1+\beta_2-2W(\mu_d)-W(\mu^*).
\end{align*} By~\eqref{eq:bpositive} and~\eqref{eq:apositive}, we obtain that $a>0$ and $b>0$.
Set
\[ a_1:= \frac{a}{1+\gamma_A}, \qquad a_2:=\gamma_A a_1, \qquad b_1:= \frac{b}{1+\gamma_B}, \qquad b_2:= \gamma_B b_1. \]
Clearly, $a_1, a_2, b_1, b_2 \geq 0$.
We have 
\begin{align}
\nonumber a_1+a_2+b_1+b_2
& = \alpha_1+\alpha_2+\beta_1+\beta_2-4W(\mu_d)-2W(\mu^*)\\\nonumber 
& = k-4W(\mu_d)-2W(\mu^*)\\ \nonumber
& \overset{\eqref{eq:deg_w_c_d}}{=} 2\deg_w(d,\mu) -4W(\mu_d)-2W(\mu^*)\\
\nonumber
& \overset{\eqref{eq.deg_w_mu_mu_d_bar_mu_d},\eqref{eq:size_mu_d_bar_mu_d}}{=}
2W(\bar\mu_d) - 2W(\mu^*)\\
& \overset{\eqref{eq:W_mu_star_bar_d_d_prime}}{=}
2W(\mu_d').\label{eq:sum}
\end{align}

Now, assume that 
\begin{align}\label{eq:assumption}
\max\{a_1,a_2\}+\min\{b_1,b_2\}\le W(\mu_d').
\end{align}
Then, we shall find below a suitable $(\gamma_A,\gamma_B)$-skew-matching anchored in $\oriented{dc}$.
In the case~\eqref{eq:assumption} does not  hold, we find a suitable $(\gamma_A,\gamma_B)$-skew-matching anchored in $\oriented{cd}$ 
using the same argumentation as below, just interchanging $\alpha, \alpha_i, \gamma_A, a, a_i$ with $\beta, \beta_i, \gamma_B, b, b_i$. \medskip

\noindent \emph{Step 4: Building the skew-matchings.}
In this phase, we apply the relevant matching lemmas from this section, based on the specific properties of the fractional matchings we are completing.

By \Cref{lem:new-balancing} (Balancing-out) (with $\mu_d, V(H)$ in place of $\mu, U$) there is a $\gamma_A$-skew-matching $\tilde \sigma_A\trianglelefteq\mu_d$ in $H^\leftrightarrow$ of weight 
$W(\tilde \sigma_A)=2W(\mu_d)$. 
By the definition of $\mu_d$ the anchor of $\tilde \sigma_A$ fits in the $w$-neighbourhood of $d$. 

Now, we wish to apply \Cref{lem:new-completion} (Completion) with
\begin{center}
\begin{tabular}{c|c|c|c|c|c|c|c|c|c}
	object & $(H, w')$ & $U$ & $V$ & $a_1$ & $a_2$ & $b_1$ & $b_2$ & $\mu'_d$ &$c$  \\
	\hline
	in place of & $(H, w)$ & $U$ & $V$ & $\alpha_1$ & $\alpha_2$ & $\beta_1$ & $\beta_2$ & $\mu$ & $u$
\end{tabular}
\end{center}
To be able to do so, we need to verify the corresponding \ref{item:completion-1}--\ref{eq:VybalancovaniMozne}.
Properties \ref{item:completion-1} and \ref{eq:neq-fit} follow from~\eqref{eq:lessthanw'} and \eqref{eq:assumption}, respectively; and \ref{eq:VybalancovaniMozne} follows by combining \eqref{eq:assumption} with~\eqref{eq:sum}. 
The application of the Completion Lemma gives a $\gamma_A$-skew-matching $\sigma_A'$ and a $\gamma_B$-skew-matching $\sigma_B'$ satisfying \ref{item:completion-out1}--\ref{item:completion-out5}.
In particular, we have $\sigma_A'+\sigma_B'\trianglelefteq \mu_d'$; we have $W(\sigma_A')=a_1+a_2=a$ and $W(\sigma_B')\ge \deg_{w'}(c, \mu_d')-W(\sigma_A')$.
Moreover, we have $\mathcal A(\sigma_A')\subseteq U$, and the anchor of $\sigma_B'$ fits in the $w'$-neighbourhood of $c$.
 
We claim that $(\sigma'_A, \sigma'_B)$ is a $(\gamma_A, \gamma_B)$-skew-matching pair anchored in $\oriented{dc}$, with respect to $w'$.
For this, we need to verify the required \ref{itdef:disjoint}--\ref{itdef:anchorpartition}.
Property \ref{itdef:disjoint} follows from $\sigma'_A + \sigma'_B \trianglelefteq \mu'_d$; and we already verified \ref{itdef:anchord}.
Using \eqref{equation:barmudw'}, we see that for every $x \in U$, 
$\sigma'^1_A(x) + \sigma'^1_B(x) \leq \mu'_d(x) \leq \bar \mu_d(x) \leq w'(\oriented{dx})$.
Since $\mathcal{A}(\sigma'_A) \subseteq U$ and $N_{w'}(d) = U$, this gives both \ref{itdef:anchorc} and \ref{itdef:anchorpartition}.

Since $(H, w')$ is the $\mu_d$-truncated weighted graph obtained from $(H,w)$, and $\tilde \sigma_A \trianglelefteq \mu_d$, we can apply \Cref{prop:adding-skew-matchings} (with $(H,w), \fmat{dc} \mu_d, \mu_d', w', (\tilde\sigma_A, \sigma_\emptyset)$, and $(\sigma_A', \sigma_B')$  playing the role of $(G,w), \mu, \bar\mu, \bar w, (\sigma_A, \sigma_B)$, and $(\bar\sigma_A, \bar \sigma_B)$) to see that the pair $(\tilde \sigma_A+\sigma_A', \sigma_B')$ is a $(\gamma_A, \gamma_B)$-skew-matching pair in $H^{\leftrightarrow}$ anchored in $\oriented{dc}$ with $\tilde \sigma_A+\sigma_A'+\sigma_B'\trianglelefteq \mu_d+\mu_d'$.

Let $(H, w^*)$ be the $(\mu_d+\mu_d')$-truncated graph  obtained from $(H, w)$.
We wish to apply \Cref{lem:new-extending-out} (Extending-out Lemma) with $U,V, H,\mu^*$ playing the roles of $U,V,H,\mu$; we can do this because $\mu^\ast$ is supported in $H[U,V]$.
From the lemma, we obtain a $\gamma_A$-skew-matching $\sigma^\ast_A$ such that 
$\sigma_A^*\trianglelefteq \mu^*$ and its anchor is contained in $U$.
We can also assume that
\[
W(\sigma_A^*)= W(\mu^*)
\]
(the lemma gives a skew-matching with larger weight, which we can scale down, and this scaling down does not break the property $\sigma_A^*\trianglelefteq \mu^*$).

Now we claim that the anchor of $\sigma_A^*$ fits in the $w^*$-neighbourhood of $d$.
Indeed, since $\mathcal{A}(\sigma^\ast_A) \subseteq A$ and $\sigma^\ast_A \trianglelefteq \mu^\ast$, it suffices to verify that $\mu^\ast(x) \leq w^\ast(x)$ holds for $x \in U$.
Since $\mu^\ast = \bar \mu_d - \mu'_d$, from \eqref{eq:lessthanw'} and \eqref{equation:kk2-wprimevsmud} we obtain $\mu^\ast(x) \leq w(\oriented{dx}) - (\mu_d(x) + \mu'_d(x))$, from which the desired inequality follows by the definition of $w^\ast$.
This allows us to apply \Cref{prop:adding-skew-matchings} (with $(H,w), \fmat{dc}, \mu_d+\mu_d', \mu^*, w^*, (\tilde\sigma_A+\sigma_A', \sigma_B')$, and $(\sigma_A^*, \sigma_\emptyset)$ playing the role of $(G,w), \mu, \bar\mu, \bar w, (\sigma_A, \sigma_B)$, and $(\bar\sigma_A, \bar \sigma_B)$) again, to obtain that the pair $(\tilde \sigma_A+\sigma_A'+\sigma_A^*, \sigma_B')$ is a $(\gamma_A, \gamma_B)$-skew-matching pair in $(H,w)$ anchored in $\oriented{dc}$ with $\tilde \sigma_A+\sigma_A'+\sigma_A^*+\sigma_B'\trianglelefteq \mu_d+\bar\mu_d$.

Now, let $(H, \bar w)$ be the $(\mu_d+\bar \mu_d)$-truncated weighted graph  obtained from $(H, w)$.
Recall that $\mu_c = \mu - (\mu_d + \bar \mu_d)$.
We have the following:
\begin{align}
\nonumber \deg_{\bar w}(c, \mu_c)
& = \sum_{x\in V(H)}\min\{\bar w(\oriented{cx}), \mu_c(x)\}\\ \nonumber
& = \sum_{x\in V(H)}\max\{0, \min\{w(\oriented{cx}), \mu(x)\}-(\mu_d(x)+\bar \mu_d(x))\}\}\\ \nonumber
& \ge \sum_{x\in V(H)}\min\{w(\oriented{cx}), \mu(x)\}-2W(\mu_d+\bar \mu_d)  \\ \nonumber
&=\deg_w(c, \mu)-2W(\mu_d+\bar \mu_d) \overset{\eqref{eq:deg_w_c_d}}{=} 2\deg_w(d, \mu) -2W(\mu_d+\bar \mu_d)\\
 &\overset{\eqref{eq:size_mu_d_bar_mu_d},\eqref{eq.deg_w_mu_mu_d_bar_mu_d}}{=} 2W(\mu_d) \overset{\eqref{eq:sizemu_d}}{=}\deg_w(d, \mu_d).\label{eq:degmuc=degmud}
\end{align}

By \Cref{lem:new-combination} (Combination) (with $ H, \bar w,\mu_c,c$ playing the roles of $H,w,\mu,v$) there exists a $\gamma_B$-skew-matching $\tilde \sigma_B\trianglelefteq\mu_c$ in $\bar H^\leftrightarrow$ with weight $W(\tilde \sigma_B)\ge \deg_{\bar w}(c, \mu_c)$, and its anchor $\mathcal A(\tilde \sigma_B)$ fits in the $\bar w$-neighbourhood of~$c$. \medskip

\noindent \emph{Step 5: Gluing skew-matchings together.}
The final step is to combine all the skew-matchings and verify that the resulting skew-matching pair satisfies the required properties.

We check a quick inequality before proceeding.
Note that, since $w'$ is the $\mu_d$-truncated weight  obtained from $w$, we have $w'(\oriented{cx}) \geq w(\oriented{cx}) - \mu_d(x)$ for each $x \in V(H)$, and therefore
\begin{equation} \label{equation:degfwcdfjkgfsfgfd}
	\deg_{w'}(c, \mu_d')+2W(\mu_d) \geq \sum_{x \in V(H)} \left( \min\{ w'(\oriented{cx}), \mu'_{d}(x) \} + \mu_d(x) \right) \geq \deg_{w}(c, \mu_d+\mu_d')
\end{equation}

By \Cref{prop:adding-skew-matchings} (with $(H,w), \fmat{dc}, \mu_d+\bar\mu_d, \mu_c, \bar w, (\tilde\sigma_A+\sigma_A'+\sigma_A^*, \sigma_B')$, and $(\sigma_\emptyset, \tilde\sigma_B)$ playing the role of $(G,w), \mu, \bar\mu, \bar w, (\sigma_A, \sigma_B)$, and $(\bar\sigma_A, \bar \sigma_B)$) the pair $( \tilde \sigma_A+\sigma_A'+\sigma_A^*, \tilde \sigma_B+\sigma_B')$ is a $( \gamma_A, \gamma_B)$-skew-matching pair in $H^\leftrightarrow$ anchored in $\oriented{dc}$ with $\tilde \sigma_A+\sigma_A'+\sigma_A^*+\tilde \sigma_B+\sigma_B'\trianglelefteq\mu$.  Set $\sigma_A:= \tilde \sigma_A+\sigma_A'+\sigma_A^*$ and $\sigma_B:= \tilde \sigma_B+\sigma_B'$. Then $\sigma_A+\sigma_B\trianglelefteq \mu$.
We have that
\begin{align*}
W(\sigma_A)=W(\tilde \sigma_A)+W(\sigma_A')+W(\sigma_A^*)=2W(\mu_d)+W(\mu^*)+a=\alpha_1+\alpha_2,
\end{align*}
and
\begin{align*}
W(\sigma_B)
& = W(\tilde \sigma_B)+W(\sigma_B')\ge \deg_{\bar w}(c, \mu_c)+\deg_{w'}(c, \mu_d')-W(\sigma_A')\\
& = \deg_{\bar w}(c, \mu_c)+\deg_{w'}(c, \mu_d')-\left(\alpha_1+\alpha_2-2W(\mu_d)-W(\mu^*)\right)\\
& \overset{\eqref{equation:degfwcdfjkgfsfgfd}}{\ge} \deg_{\bar w}(c, \mu_c)+\deg_{w}(c, \mu_d+\mu_d')+W(\mu^*)-(\alpha_1+\alpha_2) \\
& \overset{\eqref{eq:degreefromc}}{\ge} \deg_{\bar w}(c, \mu_c)+\deg_w(c, \mu_d + \bar \mu_d)-(\alpha_1+\alpha_2) \\
& \overset{\textrm{Prop. \ref{prop:adding-degree-truncated}}}{\ge} \deg_w(c, \mu)-(\alpha_1+\alpha_2)\overset{\eqref{eq:deg_w_c_d}}{\geq}
 k-(\alpha_1+\alpha_2)=\beta_1+\beta_2.
\end{align*}
If $W(\sigma_B)>\beta_1+\beta_2$, we can scale it down.
This finishes the proof of \Cref{lem:new-k-k/2}.	
\end{proof}

\section{The Structural Proposition}\label{sec:new-proof-structural}
The goal of this section is to prove the Structural Proposition (\Cref{prop:weighted-structural}), which we restate here for convenience.

\restatestructural*

The proof is split across this section and spans several subsections.
We begin by giving a sketch of the proof, highlighting the various cases that will arise during our analysis.
Then we proceed with the main proof in the remaining subsections.

\subsection{Sketch of the proof}
We aim to find appropriate neighbouring vertices $c,d$ in $V(H)$ and to construct a sufficiently large skew-matching pair anchored at those two vertices.
In the course of the proof, we will choose $c$ to have total degree at least $k$, and vertex $d$ to have total degree at least $k/2$.
The proof splits into several cases, considering different configurations of the graph $H$, as well as the different possible structures of the pair $(\sigma_A, \sigma_B)$ of skew-matchings.

\begin{enumerate}
	\item \emph{The fractional matching cover case.}
	In the first two claims (\Cref{cl:new-cinS,cl:new-cnotinS}), 
	we consider the situation where a fractional matching obtained from a Gallai-Edmonds triple is sufficient to build the entire skew-matching pair.
	This arises when the fractional matching covers well the neighbourhood of vertex~$c$.
	
	The case is divided into two claims, depending on the position of $c$. \Cref{cl:new-cinS} considers the configuration where $c$ lies in the separator $S$. Then finding a suitable vertex $d$, whose neighbourhood is also completely covered by the fractional matching, is easy.
	The skew-matching pair is obtained by applying the $(k,k/2)$-Lemma (\Cref{lem:new-k-k/2}).
	
	On the other hand, \Cref{cl:new-cnotinS} considers the configuration where $c$ does not lie in the separator $S$. Then, having chosen a neighbouring vertex $d$, it may be necessary to slightly alter the given fractional matching to ensure it covers well its neighbourhood, while still keeping a good coverage of the neighbourhood of~$c$. After building a suitable fractional matching, we again apply the $(k,k/2)$-Lemma (\Cref{lem:new-k-k/2}). Assuming that we are not in the fractional matching case provides important structural information about the graph, in particular on the existence of the set $\mathcal R$ and $ S_{\mathcal R}$ from \Cref{def:weighted-R}.
	
	\item \emph{The easy skew case}.
	In \Cref{cl:new-weighted-a_1b_1>k/2}, we consider the case where the skew-matching structure has favorable parameters, making it easy to build.
	Then most of the skew-matching pair is built using the fractional matching from the GE triple, and the left-over can be handled with a greedy argument.
	This allows us to make some basic assumptions  on the structure of the skew-matching pair we aim to embed for the rest of the proof. 
	
	\item \emph{The skew-matching cover case}.
	In the next two claims (\Cref{cl:new-weighted-largetildesigma} and \Cref{cl:new-weighted-coveringR}) we aim to replace part of the (unsatisfactory) fractional matching by ``blowing it'' into a $\gamma_B$-skew-matching in order to cover the neighbourhood of $c$ as much as possible.
	For this purpose, we will use a GE pair (\Cref{def:GE-weighted-pair}).
	Here, \Cref{cl:new-weighted-largetildesigma} will take care of the case where the ``blowed part'' can actually accommodate the whole skew-matching $\sigma_B$.
	The skew-matching $\sigma_A$ will be built within the remaining fractional matching.
	Secondly, \Cref{cl:new-weighted-coveringR} covers the complementary situation where the obtained skew-matching is perhaps not large enough to accommodate the whole one part of the pair, but similarly as in the \emph{fractional matching case}, it covers well the neighbourhood of~$c$. 
	
	After this step, we will have found an optimal GE pair, which in particular allow us to use the Separating Lemmas (\Cref{lemma:separatinglemma-1} and \Cref{lemma:separatinglemma-2}), to obtain two important ``Separating Claims'' in our situation (\Cref{cl:new-weighted-alternating-skew-1} and \Cref{cl:new-weighted-alternating-skew-2}).
	
	\item \emph{The balanced case}.
	Harnessing the assumptions we made on the configuration of the graph $H$, in \Cref{cl:new-weighted-smallb_2} we manage to build the required pair of skew-matchings, under the condition that it is reasonably balanced, i.e., none of the four parts $a_1,a_2,b_1,b_2$ exceeds half of the total weight $a_1+a_2+b_1+b_2$. This allows us  to assume from now on that the pair of skew-matchings has a huge $b_2$. 
	
	\item \emph{The large $ S_{\mathcal R}$-case}.
	This auxiliary case allows us to extract further assumptions on the configuration of the graph $H$.
	In \Cref{cl:new-weighted-smallintersection} we assume that the neighbourhoods of two elements of $\mathcal R$ do not intersect much, leading to the existence of a large $ S_{\mathcal R}$.
	The large size of $ S_{\mathcal R}$ allows us to build the required pair of skew-matching.
	This enables us to make a crucial assumption for the following case, i.e., that there is a `flabellum structure' emanating from $ S_{\mathcal R}$ and spreading to~$\mathcal R$.
	
	\item \emph{The flabellum case}.
	In \Cref{cl:new-weighted-evantail}, we assume the graph $H$ contains a large `flabellum structure'.
	Intuitively, a flabellum should have two parts: a smaller one, called the base, from which the structure `expands' to the second (larger) part.
	In our graph~$H$, this is represented by a bipartite graph: the smaller colour class fits in $ S_{\mathcal R}$ and forms the base of the flabellum; the larger colour class contains $\mathcal R$ and has the property that the neighbourhood of each of its vertices intersects the base substantially.
	We shall use the flabellum to  host the (very large) $\sigma_B$ and then build $\sigma_A$ somehow in the leftover of the graph.
	
	\item \emph{The avoiding case}.
	Finally, we treat the last possible configuration of the graph $H$. We assume that the `flabellum structure' emanating from $ S_{\mathcal R}$ is too small to accommodate the whole $\sigma_B$. We proceed as follows: We build a large  part of $\sigma_B$ within the fractional matching and skew-matching covering the base of the flabellum structure. Then we build the rest of $\sigma_B$ and the whole $\sigma_A$ working greedily from the part of the neighbourhood of $c$ that is not covered by the flabellum structure. We manage to do this by using the fact that we start from something that is not in the flabellum structure and thus we can avoid its base. We also heavily exploit the structural information obtained in the ``Separating Claims'' (\Cref{cl:new-weighted-alternating-skew-1} and \Cref{cl:new-weighted-alternating-skew-2}) to avoid the whole fractional matching and skew-matching covering the base of the flabellum.
\end{enumerate}

Now we begin our proof of \Cref{prop:weighted-structural}, which will take up the rest of this section. \medskip

\noindent \emph{Proof of \Cref{prop:weighted-structural}.}
For brevity, say that a \emph{good matching} is a $(\gamma_A, \gamma_B)$-skew-matching pair $(\sigma_A, \sigma_B)$ anchored in some edge $\oriented{xy}\in E(H^\leftrightarrow)$ with $w(\fmat{xy})>0$ such that $W(\sigma_A)=a_1+a_2$ and $W(\sigma_B)=b_1+b_2$.
Our goal is to show that $H^\leftrightarrow$ has a good matching.

Note that if $uv \in E(H)$ is such that $w(uv) = 0$, then we can remove $uv$ from $H$ and this does not affect any of our assumptions.
Thus we can assume that $w(uv) > 0$ for each $uv \in E(H)$.
In particular, for every $v \in V(H)$ we have $N(v) = N_H(v) = N_w(v)$.
We will use this during the whole proof.
 Let $(H^\leftrightarrow,w)$ is a weighted directed graph associated with $(H,w)$, with its weight function inherited from $(H,w)$, i.e., $w(\oriented{uv})=w(\oriented{vu})=w(\fmat{uv})$.

\subsection{Proof of \Cref{prop:weighted-structural}: The fractional matching cover case}
\label{ssection:strucprop-step1}
	We begin by applying the Gallai--Edmonds theorem (\Cref{thm:Gallai-Edmond}) on $H$, which provides us with a Gallai--Edmonds triple $(H,S,M_S)$.
	Recall that $\mathcal{K}^\ast_S$ is the set of vertices in non-singleton components in $H-S$,
	and $U_S$ is the set of vertices which correspond to singleton components in $H - S$.
	
	By assumption, there exists a vertex $c$ with $\deg_w(c) \geq k$. 
	Depending on the location of $c$ and its neighbours within the set $S$, there are two situations where we can quickly find good matchings in $(H,w)$, as the next two claims show.
	Recall the meaning of a fractional Gallai--Edmonds triple (\Cref{definition:gallaiedmonds-fractional}).

	\begin{claim}\label{cl:new-cinS}
		Suppose there exists $c \in S$ and a fractional Gallai--Edmonds triple $(H,S,\mu)$ such that $\deg_w(c,\mu)\ge k$.
		Then $H^\leftrightarrow$ has a good matching.
	\end{claim}
	
	\begin{proofclaim}[Proof of \Cref{cl:new-cinS}]
		Let $d\in N(c)\setminus S$ (this exists, because in the Gallai--Edmonds triple, $M_S$ matches $c$ with some vertex not in $S$).
		Let $K$ be the component of $H - S$ that contains $d$, so that $N_H(d) \subseteq S \cup (V(K) \setminus \{d\})$.
		The definition of the fractional Gallai--Edmonds triple implies that $\mu$ covers $N_H(d)$, so we have $\deg_w(d, \mu) = \deg_w(d) \geq \delta_w(H) \geq k/2$.
		We obtain the required good matching immediately from an application of \Cref{lem:new-k-k/2} (the $(k, k/2)$-Lemma).
	\end{proofclaim}

	\begin{claim}\label{cl:new-cnotinS}
		Suppose there exists $c \in V(H) \setminus S$ with $\deg_w(c)\ge k$.
		Then $H^\leftrightarrow$ has a good matching.
	\end{claim}

\begin{proofclaim}[Proof of \Cref{cl:new-cnotinS}]
	Suppose first that there exists $d \in N(c)\setminus S$.
	In this case, both $c$ and $d$ are not in $S$.
	Let $\mu$ be any fractional matching such that $(H, S, \mu)$ is a fractional Gallai--Edmonds triple (at least one must exist, e.g., by applying \Cref{proposition:weighted-GEmatching-c} with an arbitrary $c' \in S$ in place of $c$).
	By definition, $\mu$ covers all of $S$ and each non-singleton component of $H - S$, and therefore it must cover $N(c)\cup N(d)$.
	Hence, we have $\deg_w(c,\mu)\ge \deg_w(c) \geq k$ and $\deg_w(d, \mu) \ge \deg_w(d) \geq \delta_w(H) \geq \frac{k}{2}$.
	Then, an application of the $(k, k/2)$-Lemma (\Cref{lem:new-k-k/2}) yields the existence of the required good matching.
	From now on, we let $M = M_S$ for simplicity of notation.
	
	Hence, we may assume from now on that $N(c)\subseteq S$ holds.
	Pick $d\in N_H(c)\subseteq S$ arbitrarily.
	Our objective now is to construct a fractional matching $\mu_d$ that allows us to apply the $(k, k/2)$-Lemma once more with input edge $\fmat{cd}$.
	The construction of $\mu_d$ will need several steps, the first of which is to define an auxiliary weighted graph $(H,w')$ where $w'$ is obtained from $w$ by decreasing the weight in some of its edges, as follows.
	Recall that, by assumption, $\deg_w(d) \geq \frac{k}{2}$.
	We obtain $w'$ from $w$ by decreasing its value on the edges incident to $d$ so that $\deg_{w'}(d)=\frac{k}{2}$, and such that $N_{w}(d) \cap U_S = N_{w'}(d) \cap U_S$.
	The only purpose of $w'$ is to define a suitable $\mu_d$ as an input of the $(k, k/2)$-Lemma. 
	
	Our second step to find $\mu_d$ is to construct an initial fractional matching $\mu_1$ with the help of the First Greedy Lemma (\Cref{prop:weighted-greedy-1}), as follows.
	Let $V = U_S \cap N(d)$, i.e. it consists of all vertices forming singleton components in $H- S$ that are incident to $d$.

	Let $\sigma_\emptyset: E(H^\leftrightarrow) \rightarrow [0,1]$ 
	be the function which is identically zero everywhere, and note that it is (trivially) a $1$-skew oriented fractional matching.
	Also, define $\kappa = \deg_{w'}(d, V)$.
	We wish to apply \Cref{prop:weighted-greedy-1}
	with
	\begin{center}
	\begin{tabular}{c|c|c|c|c|c|c|c|c|c}
			object & $(H, w')$ & $(d,c)$ & $S$ & $V$ & $\sigma_\emptyset$ & $\sigma_\emptyset$ & $1$ & $1$ & $\kappa$ \\
			\hline
			in place of & $(H,w)$ & $(u,v)$ & $U$ & $V$ & $\sigma_A$ & $\sigma_B$ & $\gamma_A$ & $\gamma_B$ & $\kappa$
	\end{tabular}
	\end{center}
	We can indeed do so: condition \ref{item:greedy-1-condition-1} reduces to $\deg_{w'}(d, V) \geq \kappa$, which trivially holds.
	To check \ref{item:greedy-1-condition-2}, we note that for any $x\in V$, by our choice of $w'$, we have $N_{w'}(x) = N_w(x) = N_H(x)$.
	Since $N_H(x) \subseteq S$, we have
	\begin{align*}
		|N_{w'}(x) \cap S|
		& = |N_H(x)| \geq \deg_{w}(x) \ge \delta_w(H) \ge \frac{k}{2} \\
		& =\deg_{w'}(d)\ge \deg_{w'}(d,V) = \kappa, 
	\end{align*}
	as required.
	
	Hence, from \Cref{prop:weighted-greedy-1} we obtain a $1$-skew oriented fractional matching $\sigma$ of weight $W(\sigma) \geq 2\kappa$, whose support is contained in $H[V,S]$,	
	 and such that $(\sigma, \sigma_\emptyset)$ is anchored in $\oriented{dc}$ (with respect to $w'$),
	 and moreover the anchor of $\sigma$ is contained in $V$.
	 Then we have $\kappa = \sum_{u \in  V} w'(\oriented{du}) \geq \sum_{u \in  V} \sigma^1(u) = \frac{1}{2} W(\sigma) = \kappa$. 
	 This means that $\sigma$ saturates $N_{w'}(d)\cap V$.
	 Forgetting the orientation (formally, by \Cref{lemma:fractionalfrom1skew}) we obtain from~$\sigma$ a fractional matching $\mu_1$ in~$H$ of weight $\deg_{w'}(d,V)$ that saturates $N_{w'}(d)\cap V$. 
	 
	Our next step is to obtain a new fractional matching $\mu_2$ to cover $S$.
	Recall that $M_S \subseteq H$ is a matching such that $(H, S, M_S)$ is a Gallai--Edmonds triple.
 	Let $\mu_2$ be a fractional matching supported on $M_S$, disjoint from $\mu_1$, and of maximum possible weight. We have that $\mu_1+\mu_2$ fully covers $S$.
 	
 	Next, we obtain a new fractional matching $\mu_3$ to cover the non-singleton components in $H - S$.
 	For any component~$K$ of $H - S$ that is not a singleton, there is at most one vertex $v\in V(K)$ for which $\mu_2(v) > 0$, and every other vertex in $K$ receives zero weight under $\mu_2$. 
	As $K$ is factor-critical, there is a perfect matching $M_K$ in $K-v$. By \Cref{prop:factor-critical->fractional matching} there is a fractional matching $\mu_K'$ that completely covers $K$.
	We shall build a fractional matching $\mu_K$ by \[\mu_K:= \mu_2(v)\cdot \mathbf{1}_{M_K}+(1-\mu_2(v))\cdot \mathbf{1}_{\mu_K'},\]
	so $\mu_K$ is disjoint from $\mu_1 + \mu_2$ and supported completely in $E(K)$.  
	Let $\mu_3 = \sum_{K} \mu_K$, where the sum ranges over all non-singleton components of $H - S$.
	Then $\mu_3$ is also disjoint from $\mu_1 + \mu_2$.
	
	We can finally define our desired matching as $\mu_d := \mu_1 + \mu_2 + \mu_3$.
	By construction, $\mu_d$ covers $S$ and also covers all non-singleton components of $H-S$. Moreover $\mu_d$ saturates $N_{w'}(d)\cap V$ and thus $\deg_w(d, \mu_d)\ge \deg_{w'}(d, \mu_d)=\deg_{w'}(d)=  k/2$. As $\mu_d$ covers $S$, we have $\deg_w(c, \mu_d)=\deg_w(c)\ge k$. We can now use the $(k,k/2)$-Lemma (\Cref{lem:new-k-k/2}), with
		\begin{center}
		\begin{tabular}{c|c|c|c|c|c|c|c|c|c}
			object & $(H,w)$ & $c$ & $d$ & $\mu_d$ & $k$ & $a_1$ & $a_2$ & $b_1$ & $b_2$ \\
			\hline
			in place of & $(H,w)$ & $c$ & $d$ & $\mu$ & $k$ & $\alpha_1$ & $\alpha_2$ & $\beta_1$ & $\beta_2$
		\end{tabular}
	\end{center}
	 to obtain the required good matching and conclude the proof.
\end{proofclaim}


Now, we fix important objects (vertices and fractional matchings) for the remainder of the proof and record key properties of those objects.
Fix an arbitrary vertex $c$ with $\deg_w(c)\ge k$. 
If $c \notin S$, we are done by \Cref{cl:new-cnotinS}.
Hence, we can assume that
\stepcounter{propcounter}
\begin{enumerate}[\upshape (\Alph{propcounter}\arabic*), topsep=0.7em, itemsep=0.5em]
	\item $c \in S$.
\end{enumerate}
Recall that $U_S$ is the set of vertices that correspond to singleton components in $H - S$.
Now we apply \Cref{proposition:weighted-GEmatching-c} to $(H,w)$, $(H, S, M)$ and $c$.
We obtain a $c$-optimal fractional matching $\mu$, which means the following:
for each $u \in U_S$ such that $\mu(u) < w(\oriented{cu})$, for each $\mu$-alternating path $P_u$ starting at $u$, and for each $v \in V(P_u) \cap S$, the set $N_\mu(v) = \{ x \in N(v) : \mu(\fmat{vx}) > 0 \}$ satisfies
\begin{enumerate}[\upshape (\Alph{propcounter}\arabic*), topsep=0.7em, itemsep=0.5em, resume]
	\item \label{item:SP-fge1} $N_\mu(v) \subseteq U_S$;
	\item \label{item:SP-fge2} $\mu(v) = 1$; and
	\item \label{item:SP-fge3} for every $y \in N_\mu(v)$, $0 <\mu(y) \leq w(\oriented{cy})$.
\end{enumerate}
From \cref{remark:coptimal-neighbourhoodc}, we get
\begin{enumerate}[\upshape (\Alph{propcounter}\arabic*), topsep=0.7em, itemsep=0.5em, resume]
	\item \label{item:SP-fge4} $N_\mu(v) \subseteq N(c)$. 
\end{enumerate}
Moreover, since $(H, S, \mu)$ is a fractional Gallai--Edmonds triple (\Cref{definition:gallaiedmonds-fractional}) we also have
\begin{enumerate}[\upshape (\Alph{propcounter}\arabic*), topsep=0.7em, itemsep=0.5em, resume]
	\item \label{item:SP-musupport} $\mu$ is supported in the Gallai--Edmonds support $E_{S,M}$, and
	\item \label{item:SP-mucovered} $S \cup \bigcup_{K \in \mathcal{K}^\ast_S} V(K)$ is covered by $\mu$.
	In particular, $S \cup \bigcup_{K \in \mathcal{K}^\ast_S} V(K) \subseteq V(\mu)$.
\end{enumerate}

Let $\mathcal{R}$ be the set of reachable vertices (\Cref{def:weighted-R}) with respect to $(H,S,\mu)$, also let $S_{\mathcal{R}}$ as in that definition, i.e. the neighbourhood of the set of reachable vertices.
By Observations \ref{ob:R=singletons} and~\ref{ob:R=neighbours}, we have that
\begin{enumerate}[\upshape (\Alph{propcounter}\arabic*), topsep=0.7em, itemsep=0.5em, resume]
	\item \label{item:SP-reachableareisolated} $\mathcal{R} \subseteq U_S \cap N_w(c)$, and
	\item \label{item:SP-neighbourreachable} $S_{\mathcal{R}} \subseteq S$.
\end{enumerate}

If $\deg_w(c, \mu) \geq k$, 
 we would be done by \Cref{cl:new-cinS}.
Hence, we can assume that
\begin{enumerate}[\upshape (\Alph{propcounter}\arabic*), topsep=0.7em, itemsep=0.5em, resume]
	\item \label{item:SP-weightc} $\deg_w(c, \mu) < k$.
\end{enumerate}
Note that in particular, \ref{item:SP-weightc} implies that $\mu$ does not saturate $N_w(c)$, that is, there exists $u \in N_w(c)$ such that $\mu(u) < w(\fmat{cu})$.
By \ref{item:SP-mucovered}, we have that $u \in U_S$, i.e., $\{u\}$ is a singleton component of $H-S$.
This implies that $cu$ is an $(E(H), \mu)$-alternating path, so $u \in \mathcal{R}$.
In summary,
\begin{enumerate}[\upshape (\Alph{propcounter}\arabic*), topsep=0.7em, itemsep=0.5em, resume]
	\item \label{item:SP-reachable-neighbourhood} $N_w(c) \cap \mathcal R \cap \{ u : \mu(u) < w(\oriented{cu}) \} \neq \emptyset$.
\end{enumerate}
	

Since $a_1 + a_2 + b_1 + b_2 = k$, without loss of generality we can suppose that 
\begin{equation}\label{assumption:new-WLOG}
a_2+b_1\le \frac{k}{2},
\end{equation}
as otherwise we can just swap the roles of $a_1, b_1$ with $a_2, b_2$.
Set $\gamma_A :=  \frac{a_2}{a_1}$ and $\gamma_B := \frac{b_2}{b_1}$. 

\subsection{Proof of \Cref{prop:weighted-structural}: The easy skew case}
\label{ssection:strucprop-step2}
Now we can find a good matching if the skew of the tree satisfies a favourable condition.

\begin{claim}\label{cl:new-weighted-a_1b_1>k/2}
	If $a_1+b_1\ge k/2$, then $H^\leftrightarrow$ has a good matching.
\end{claim}
\begin{proofclaim}[Proof of \Cref{cl:new-weighted-a_1b_1>k/2}]
	The proof has three parts.
	First we will define auxiliary objects: a weighted graph $(H, w_d)$ and a fractional matching $\mu_d$, obtained from $(H,w)$ and $\mu$, respectively; and then gather properties of those objects.
	The second part of the proof is to use the auxiliary objects as an input for an application of the Completion Lemma, which gives a skew-matching as an output.
	In the third part we apply the Combination and Greedy matching lemmas to obtain two further skew-matchings.
	The combination of all the skew-matchings we have constructed yields the result.
	
	\medskip \noindent \emph{Step 1: Setting the auxiliary objects.}
	Let $d\in N(c) \cap \mathcal R$ be an arbitrary vertex with $\mu(d) < w(\oriented{cd})$, which exists by \ref{item:SP-reachable-neighbourhood}.
	Since $d \in \mathcal{R}$, we have in fact that $N(d) \subseteq S$.	
%
	
	We obtain a new weighted graph $(H,w_d)$, where $w_d$ is obtained from $w$ by decreasing the weight function on edges incident to $d$ (other than $cd$) in such a way that $\deg_{w_d}(d, \mu)=k/2$ holds. 
	The goal of $w_d$ is to define a fractional matching $\mu_d$ and filling it completely with some skew-matching $\sigma_A'$ and $\sigma_B$. 
	
	Next, we define a new fractional matching $\mu_d \leq \mu$ using the following procedure.
	We recall that by \ref{item:SP-musupport} we have that $\mu$ is not supported in any edge with two endpoints in~$S$; and that $N_{w_d}(d)\subseteq N_w(d) \subseteq S$.
	Initially, we obtain $\mu_d$ from $\mu$ by setting its value on every edge which does not touch $N_{w_d}(d)$ to $0$.
	At this point, we have that $\mu_d$ is supported only in edges with one endpoint in $N_{w_d}(d)\subseteq S$ and the other in $V(H) \setminus S$.
	Next, we process every $x \in N_{w_d}(d)$ in turn.
	If $\mu_d(x) \le  w_d(\oriented{dx})$ we do nothing, otherwise we decrease the value of $\mu_d$ on the edges incident to $x$ so that $\mu_d(x) = w_d(\oriented{dx})$ holds.
	This finishes the construction of $\mu_d$.
	
	Note that by our construction, $\min\{ \mu_d(x), w_d(\oriented{dx}) \} = \mu_d(x)$ if $x \in N_{w_d}(d)$, and $\mu_d(x) = 0$ for $x \in S \setminus N_{w_d}(d)$.
	Also, if $\mu_d(x) < \mu(x)$, then $\mu_d(x) = w_d(\oriented{dx})$.
	This implies that
	\begin{equation} \text{$\min\{ \mu(x), w_d(\oriented{dx}) \} = \min\{ \mu_d(x), w_d(\oriented{dx}) \} = \mu_d(x)$ for every $x \in N_{w_d}(d)$.}
	\label{item:strucprop-wdmuclaim51}
	\end{equation} 
	Using all of this, we get
	\begin{align}
		W(\mu_d) & = \sum_{xy \in E(H)} \mu_d(\fmat{xy}) = \sum_{x \in N_{w_d}(d)} \mu_d(x) = \sum_{x \in N_{w_d}(d)} \min\{ \mu(x), w_d(\oriented{dx}) \} \nonumber \\
		& =\deg_{w_d}(d,\mu)=\frac{k}{2}. \label{eq:degwd(d)mu_d}
	\end{align}	
	
	Before continuing, we gather some inequalities.
	From \eqref{assumption:new-WLOG} and the assumption $a_1 + b_1 \geq k/2$, we deduce
	\begin{equation}\label{eqneq-neq:gamma_a<1-cl:weighted-a_1b_1>k/2}
	\gamma_A=\frac{a_2}{a_1}\le 1.
	\end{equation}
	Since $a_1 + b_1 \geq k/2$, we get $a_2 + b_2 \leq k/2$;
	similarly from~\eqref{assumption:new-WLOG} we also have $a_1 + b_2 \geq k/2$.
	Using all of this and \eqref{eq:degwd(d)mu_d}, we get
	\begin{equation}\label{eq-new:a_ib_i-cl:weighted-a_1b_1>k/2}
	\max\{b_1,b_2\}+a_2\le \frac{k}{2} = \deg_{w_d}(d, \mu)=W(\mu_d)\le \min\{b_1,b_2\}+a_1.
	\end{equation}
	
	Now we want to use \ref{item:SP-fge1}--\ref{item:SP-fge4} with $d$ playing the role of $u$, for which we verify the necessary conditions.
	Since $d\in {\mathcal{R}}$, we have $d \in U_{S}$ by \ref{item:SP-reachableareisolated}; and by the choice of $d$ we also have $\mu(d) < w(\oriented{cd})=w_d(\oriented{cd})$, so indeed $d$ can play the role of $u$ in that property.
	For any $x \in N_{w_d}(d)$, by \cref{remark:coptimal-neighbourhoodofR} we have that
	it can play the role of $v$ in \ref{item:SP-fge1}--\ref{item:SP-fge4}.
	We deduce that
	for any $x \in N_{w_d}(d)$ and any $y \in N_{\mu}(x)$, by \ref{item:SP-fge3} we have
	\begin{equation}
		\mu_d(y)\le \mu(y)\le w(\oriented{cy})=w_d(\oriented{cy}).
		\label{equation:claim51-conditionmudwd}
	\end{equation}

	\noindent \emph{Step 2: Applying the Completion Lemma.}
	Recall that $V(\mu_d)$ denotes the set of vertices $v$ such that $\mu_d(v) > 0$.
	Now we want to apply \Cref{lem:new-completion} (Completion) with
	\begin{center}
	\begin{tabular}{c|c|c|c|c|c|c|c|c|c}
		object & $(H, w_d)$ & $N_{w_d}(d)$ & $V(\mu_d) \setminus S$ & $b_1$ & $b_2$ & $a_1$ & $a_2$ & $\mu_d$ & $c$  \\
		\hline
		in place of & $(H, w)$ & $U$ & $V$ & $\alpha_1$ & $\alpha_2$ & $\beta_1$ & $\beta_2$ & $\mu$ & $u$
	\end{tabular}
	\end{center}
	
	Let us verify the required conditions \ref{item:completion-1}--\ref{eq:VybalancovaniMozne}.
	Since $N_{w_d}(d) \subseteq S$, we indeed have that $N_{w_d}(d) \cap (V(\mu_d) \setminus S) = \emptyset$.
	By construction, $\mu_d$ is supported only in edges with one endpoint in $N_{w_d}(d)$ and the other in $V(\mu_d) \setminus S$, as required.
	If $y \in V(\mu_d) \setminus S$, there exists $x \in N_{w_d}(d)$ such that $y \in N_{\mu}(x)$, so from \eqref{equation:claim51-conditionmudwd} we get that
	\begin{equation}
		\label{item:strucprop-mudvswprimeprime} \text{$\mu_d(y) \leq w_d(\oriented{cy})$ for all $y \in V(\mu_d) \setminus S$,}
	\end{equation} 
	which gives~\ref{item:completion-1}.
	Inequality~\eqref{eq-new:a_ib_i-cl:weighted-a_1b_1>k/2} implies~\ref{eq:neq-fit} and~\ref{eq:VybalancovaniMozne}.
	Finally,~\eqref{eqneq-neq:gamma_a<1-cl:weighted-a_1b_1>k/2} also allows us to use part~\ref{item:completion-moreover}.
	
	Hence, $H^\leftrightarrow$ admits a $\gamma_A$-skew-matching $\sigma'_A$ and a $\gamma_B$-skew-matching $\sigma_B$ satisfying \ref{item:completion-out1}--\ref{item:completion-moreover}.
	In particular, we have that
	$\sigma_A'+\sigma_B \trianglelefteq \mu_d$ and $W(\sigma_B)=b_1 + b_2$.
	We also have that
	\begin{equation}
		W(\sigma'_A) \geq \deg_{w_d}(c, \mu_d) - W(\sigma_B). \label{equation:strucprop-easyskew-weightsigmaaprime}
	\end{equation}
	Moreover, we have that the anchor $\mathcal{A}(\sigma_B)$ is contained in $N_{w_d}(d)$, the anchor $\mathcal{A}(\sigma'_A)$ fits in the $w_d$-neighbourhood of $c$, and
	\begin{equation} \label{item:strucprop-mudysigmaprimera} \text{for all $y \in V(\mu_d) \setminus S$, we have $\sigma'_A(y) + \sigma_B(y) = \mu_d(y)$.}
	\end{equation} 
	
	Note that for each $y \in V(\mu_d) \setminus S$, we have that there exists $x \in S \cap N_H(d)$ with $\mu(\fmat{xy}) \geq \mu_d(\fmat{xy}) > 0$.
	We have $d \in U_S$ and $\mu_d(d) \leq \mu(d) < w(\oriented{cd})$.
	Hence $d x y$ is a $\mu$-alternating path and $y$ is reachable, so $V(\mu_d) \setminus S \subseteq \mathcal{R}$.
	Therefore, together with \eqref{item:strucprop-mudysigmaprimera} and \eqref{eq:degwd(d)mu_d}, we have
	 \begin{equation}\label{eqnew:fills-up}
	 \sum_{y\in \mathcal R} \left( \sigma_A'(y)+\sigma_B(y) \right)
	  = \sum_{y\in V(\mu_d) \setminus S} \mu_d(y) = W(\mu_d)=\frac{k}{2}.
	 \end{equation} 
	 Using $\sigma_A'+\sigma_B\trianglelefteq \mu_d$ and the properties of $\sigma'_A, \sigma_B$, we have that $(\sigma_A', \sigma_B)$ is a $(\gamma_A, \gamma_B)$-skew pair in $(H, w_d)$, anchored in $\oriented{cd}$ (we use \eqref{item:strucprop-wdmuclaim51} to check \ref{itdef:anchord} and \ref{itdef:anchorpartition}).
	 
	 \medskip \noindent \emph{Step 3: Two more skew-matchings.}
	Let $(H, w')$ be the $\mu_d$-truncated weighted graph obtained from $(H,w_d)$. 
	Let $\mu_c= \mu-\mu_d$.
	
	By \Cref{lem:new-combination} (Combination) there is a $\gamma_A$-skew-matching $\sigma_A^*\trianglelefteq \mu_c$ with its anchor $\mathcal A(\sigma_A^*)$ fitting in the $w'$-neighbourhood of $c$ and of weight
	\begin{equation}
		W(\sigma_A^*)\ge  \deg_{w'}(c,\mu_c). \label{equation:strucprop-easyskewwsigmastar}
	\end{equation}
	By \Cref{prop:adding-skew-matchings} (with $(H,w_d), \fmat{cd}, \mu_d, \mu_c, w', \sigma_A, \sigma_B, \sigma_A^*$, and $\sigma_\emptyset$ playing the role of $(G,w), \fmat{uv}, \mu, \bar\mu, \bar w, \sigma_A, \sigma_B, \bar\sigma_A$, and $\bar\sigma_B$, respectively) we have that $(\sigma_A'+\sigma_A^*, \sigma_B)$ is a $(\gamma_A, \gamma_B)$-skew pair anchored in 
	$\oriented{cd}$ with $\sigma_A'+\sigma_A^*+\sigma_B\trianglelefteq \mu$.
	Moreover, we have (with explanations to follow)
	\begin{align*}
		W(\sigma_A^*+\sigma_A'+\sigma_B)
		& \ge W(\sigma_A^*) + W(\sigma'_A) + W(\sigma_B)
		 \overset{\eqref{equation:strucprop-easyskew-weightsigmaaprime},\eqref{equation:strucprop-easyskewwsigmastar}}{\ge } \deg_{w'}(c, \mu_c) + \deg_{w_d}(c, \mu_d) \\
		& 		= \deg_{w_d}(c,\mu) =\deg_w(c, \mu).
	\end{align*} 
	Here, the last equality follows from the construction of $w_d$, from which we infer $w_d(\oriented{cx})=w(\oriented{cx})$ for all $x\in V(H)$. 
	The penultimate inequality follows from \Cref{prop:adding-degree-truncated} (with $(H, w_d), \mu_d, \mu$, and $w'$ playing the role of $(G,w), \mu', \mu$, and $w'$, respectively).
	In summary, we obtain
	\begin{equation}
		W(\sigma_A^*+\sigma_A'+\sigma_B) \geq \deg_w(c,\mu).
		\label{equation:claimeasyskew-weightsigmaastarplussigmaaprime}
	\end{equation}	

	Let $\kappa' = (k-W(\sigma_B+\sigma_A'+\sigma_A^*))(1+\gamma_A)^{-1}$.
	We want to use the First Greedy Lemma (\Cref{prop:weighted-greedy-1}) with
	\begin{center}
	\begin{tabular}{c|c|c|c|c|c|c}
		object & $(H, w)$ & $\oriented{cd}$ & $(\sigma_A'+\sigma_A^*, \sigma_B)$ & $V(H)\setminus \mathcal R$ & $\mathcal R$ & $\kappa'$  \\
		\hline
			in place of & $(H, w)$ & $\oriented{uv}$ & $(\sigma_A, \sigma_B)$ & $U$ & $V$ & $\kappa$ 
	\end{tabular}
	\end{center}
	To verify we can do so, we check the required inequalities \ref{item:greedy-1-condition-1} and \ref{item:greedy-1-condition-2}.
	We note first that, all support edges of $\mu_d$ intersect $N_{w_d}(d)$, by construction.
	Hence, if $y \in \mathcal{R}$ and $\mu_d(y) > 0$, by \eqref{equation:claim51-conditionmudwd} we have $\mu(y) \leq w(\oriented{cy})$.
	Using this, together with $\sigma_A^*\trianglelefteq \mu_c$ and $\sigma_A'+\sigma_B\trianglelefteq \mu_d$, we obtain
	\begin{align}
\nonumber \sum_{x\in \mathcal R}\min\{w(\oriented{cx}), \mu(x)\}
& = \sum_{x\in \mathcal R} \mu(x) = \sum_{x\in \mathcal R} \mu_c(x) + \sum_{x\in \mathcal R} \mu_d(x)\\
&\ge \sum_{x\in \mathcal R} \left( \sigma_A^*(x)+ \sigma'_A(x)+\sigma_B(x) \right).\label{eq:intermediatecalculations}
	\end{align}
	Also note that if $y \in V(H)$ is such that $0 < \mu(y) < w(\oriented{cy})$,
	then $y \notin S$ (by \ref{item:SP-mucovered}) and hence $y \in \mathcal{R}$ by definition.
	Hence,
	\begin{equation}
		\sum_{x\in V(H)\setminus \mathcal R} w(\oriented{cx}) \leq \sum_{x\in V(\mu)\setminus \mathcal R}\min\{w(\oriented{cx}), \mu(x)\}.
		\label{eq:intermediatecalculations2}
	\end{equation}
	Using this, we can verify \ref{item:greedy-1-condition-1}. Indeed,

	\begin{align*}
	\deg_w(c,\mathcal R)
	& \ge \deg_w(c)-\sum_{x\in V(H)\setminus \mathcal R} w(\oriented{cx})\\
	& \overset{\eqref{eq:intermediatecalculations2}}{\ge } \deg_w(c)-\sum_{x\in V(\mu)\setminus \mathcal R}\min\{w(\oriented{cx}), \mu(x)\}\\
	& \ge k -\sum_{x\in V(\mu)\setminus \mathcal R}\min\{w(\oriented{cx}), \mu(x)\}\\
	&\ge k -\deg_w(c, \mu)+\sum_{x\in \mathcal R}\min\{w(\oriented{cx}), \mu(x)\} \\
&\overset{\eqref{equation:claimeasyskew-weightsigmaastarplussigmaaprime}}{\geq} k-W(\sigma_A'+\sigma_A^*+\sigma_B)+\sum_{x\in \mathcal R}\min\{w(\oriented{cx}), \mu(x)\} \\
	& \overset{\eqref{eq:intermediatecalculations}}{\ge} k-W(\sigma_A'+\sigma_A^*+\sigma_B) +\sum_{x\in \mathcal{R}}(\sigma_A'(x)+\sigma_A^*(x)+\sigma_B(x) ) \\
	& \ge \kappa' +\sum_{x\in \mathcal{R}}(\sigma_A'(x)+\sigma_A^*(x)+\sigma_B(x) ),
\end{align*}
	Next, for any $x\in\mathcal R$, we have $N_w(x) \subseteq S \setminus \mathcal{R}$, so
	\begin{align*}
		|N_w(x)\cap (V(H)\setminus \mathcal R)|
		& \geq \deg_w(x) \geq \frac{k}{2}
	\overset{\eqref{eqnew:fills-up}}{\geq} k-\sum_{x\in \mathcal R}(\sigma_B(x)+\sigma_A'(x)+\sigma_A^*(x))\\
	& \ge k -W(\sigma_B+\sigma_A'+\sigma_A^*)+\sum_{y\in V(H)\setminus \mathcal R}(\sigma_B(y)+\sigma_A'(y)+\sigma_A^*(y)) \\
	& \geq (1+\gamma_A) \kappa'+\sum_{y\in V(H)\setminus \mathcal R}(\sigma_B(y)+\sigma_A'(y)+\sigma_A^*(y)),
	\end{align*}
	which is \ref{item:greedy-1-condition-2}, as required.
	
	The outcome of \Cref{prop:weighted-greedy-1} is a $\gamma_A$-skew-matching $\tilde \sigma_A$
	of weight $W(\tilde \sigma_A)\ge (1+\gamma_A)\kappa' = k -W(\sigma_B+\sigma_A'+\sigma_A^*)$ such that $\sigma_A:= \sigma_A'+\sigma^*_A+\tilde \sigma_A$ is such that $(\sigma_A,\sigma_B)$ is a $(\gamma_A,\gamma_B)$-skew-matching anchored in $\oriented{cd}$.
	We have that $W(\sigma_B) = b_1+b_2$ and
	$W(\sigma_A) \geq k - (b_1+b_2) = a_1+a_2$, so we have found the required good matching,
	proving the claim.
	\end{proofclaim}

In the rest of the proof we may assume that 
\begin{equation}\label{assumption-new:a_1+b_1-small}
\max\{a_1, a_2\}+b_1< \frac{k}{2}.
\end{equation}
This inequality, together with $a_1 + a_2 + b_1 + b_2 = k$, implies $b_2 + \min\{a_1, a_2\} > k/2 > b_1 + \max\{a_1, a_2\} \geq b_1 + \min\{a_1, a_2\}$, which in turn implies $b_2 > b_1$.
Hence, we have
\begin{enumerate}[\upshape (\Alph{propcounter}\arabic*), topsep=0.7em, itemsep=0.5em, resume]
	\item \label{eqnew:gamma_B>1} $\gamma_B= b_2/b_1 >1$.
\end{enumerate}

\subsection{Proof of \Cref{prop:weighted-structural}: The skew-matching cover case}\label{ssec:skew-matching-case}
In the remainder of the proof we will need the concept of a \emph{GE pair} (\Cref{def:GE-weighted-pair}).
Recall that we have already defined the fractional matching $\mu$, which is $c$-optimal.
We will let $(\tilde \sigma, \tilde \mu)$ be a $(H, w, S, M, c, \gamma_B)$-GE pair, which means that

\begin{enumerate}[\upshape (\Alph{propcounter}\arabic*), topsep=0.7em, itemsep=0.5em, resume]
	\item \label{item:strucprop-GEskewmatching} $\tilde \sigma$ is a $\gamma_B$-skew-matching,
	\item \label{item:strucprop-GEpairdisjoint} $\tilde \mu$ is a fractional matching disjoint from $\tilde \sigma$,
	\item \label{item:strucprop-GEpairanchorS}  the anchor $\mathcal A(\tilde \sigma)$ of $\tilde \sigma$ is contained in $S_{\mathcal R}$,
	\item \label{item:strucprop-GEpairanchorfits} $\mathcal A(\tilde \sigma)$ fits in the $w$-neighbourhood of $c$,
	\item \label{item:strucprop-GEskewinR} $V(\tilde \sigma) \setminus \mathcal A(\tilde \sigma) \subseteq \mathcal R$,
	\item \label{item:structuprop-GEpairsaturation} for every $y\in \mathcal R$ we have $\tilde \mu(y)+\tilde \sigma(y)\le w(\oriented{cy})$,
	\item \label{item:strucprop-GEpaircovers-w} 
	 for every $y\in V(H)\setminus \mathcal R$ we have $\tilde \mu(y)+\tilde \sigma(y)\ge w(\oriented{cy})$,
	\item \label{item:strucprop-GEpaircovers} $\tilde \mu+\tilde \sigma$ covers  $S$, and
	\item \label{item:strucprop-GEpairsupport} $\tilde \mu$ equals $\mu$ when restricted to the graph $H-(\mathcal R\cup S_{\mathcal R})$ and any supporting edge of $\tilde \mu$ intersecting the set $\mathcal R\cup S_{\mathcal R}$ lies in the bipartite graph $H[\mathcal R, S_\mathcal R]$.
\end{enumerate}
Such objects exist by \Cref{lemma:gepairs-exist}.
\label{pageref:GE-pair}Over all possible $(H,w,S,M,c,\gamma_B)$-GE pairs $(\tilde \sigma, \tilde \mu)$, we can assume we choose $(\tilde \sigma, \tilde \mu)$ to be optimal, meaning that, over the possible choices,
\begin{enumerate}[\upshape (\Alph{propcounter}\arabic*), topsep=0.7em, itemsep=0.5em, resume]
	\item \label{item:strucprop-GEpairmaxsat} $\tilde \mu+\tilde \sigma$ maximises the saturation $\deg_w(c, \tilde \mu+\tilde \sigma)$ of $N(c)$.
\end{enumerate}

As we shall see now, we can conclude if $\tilde \sigma$ has enough weight.

\begin{claim}\label{cl:new-weighted-largetildesigma}
	Suppose that $W(\tilde \sigma)\ge b_1+b_2$.
	Then $H^\leftrightarrow$ has a good matching.
\end{claim}

\begin{proofclaim}[Proof of \Cref{cl:new-weighted-largetildesigma}]
	Let $d\in {\mathcal{R}}$ be arbitrary.
	By \ref{item:SP-reachableareisolated} we have that $d \in N(c)$.
	Next, let $\sigma_B$ be a $\gamma_B$-skew-matching obtained by scaling down $\tilde \sigma$, so that $\sigma_B \le \tilde \sigma$ and $W(\sigma_B)=b_1+b_2$.
	Let $\sigma'_B := \tilde \sigma - \sigma_B$, so $\sigma'_B$ is a $\gamma_B$-skew-matching of weight $W(\tilde \sigma) - (b_1 + b_2)$.
	
	We wish to apply \Cref{prop:weighted-extending-skew} (Extending-out skew-matching) with
	\begin{center}
	\begin{tabular}{c|c|c|c|c|c}
		object & $(H, w)$ & $c$ & $\sigma'_B$ & $1$ & $\gamma_B$ \\
		\hline
		in place of & $(H, w)$ & $u$ & $\sigma_B$ & $\gamma_A$ & $\gamma_B$ 
	\end{tabular}
	\end{center}
	Indeed we can do so: we have $\gamma_B \geq 1$ by \ref{eqnew:gamma_B>1}; and as $\sigma'_B\le \tilde \sigma$, 
	by \ref{item:strucprop-GEpairanchorfits}, we have that the anchor $\mathcal A(\sigma'_B)$ fits in the $w$-neighbourhood of $c$.
	From this application we obtain a $1$-skew-matching $\sigma'_A$ of weight $2 W(\sigma'_B)/(1 + \gamma_B)$, and such that $\sigma'_A \leq \sigma'_B$.
	Using \Cref{lemma:fractionalfrom1skew}, we obtain from $\sigma'_A$ a fractional matching $\mu_{\tilde \sigma}$ such that
	\begin{equation}
		W(\mu_{\tilde \sigma}) = \frac{1}{2} W(\sigma'_A) = \frac{W(\sigma'_B)}{1 + \gamma_B} = \frac{W(\tilde \sigma)}{1+\gamma_B}-b_1,
	\end{equation}
	and, for each $x \in V(H)$,
	\begin{equation}
		\mu_{\tilde \sigma}(x) = \sigma'_A(x),
		\label{equation:largetildesigma-mutildesigmavssigmaprimeA}
	\end{equation}
	and
	\begin{equation}
		\mu_{\tilde \sigma} \preceq \tilde \sigma-\sigma_B.
	\end{equation}
	Note that, since $\sigma'_A \leq \sigma'_B \leq \tilde \sigma$, from \ref{item:strucprop-GEpairanchorS}  we have that the anchor of $\sigma'_A$ is in $S_{\mathcal R}$. Then \ref{item:strucprop-GEskewinR} implies that $W(\sigma'_A) = \sum_{x \in S} \sum_{y \in N(x)} \sigma'_A(xy)$, and the analogous equation for $W(\sigma'_B)$ also holds.
	Using $\sigma'_A \leq \sigma'_B$, we get that
	\[ W(\sigma'_B) = \sum_{x \in S} \sum_{y \in N(x)} \sigma'_B(xy) \geq \frac{1 + \gamma_B}{2} \sum_{x \in S} \sum_{y \in N(x)} \sigma'_A(xy) = \frac{1 + \gamma_B}{2} W(\sigma'_A) = W(\sigma'_B), \]
	so $\sigma'_B(xy) = (1+\gamma_B)\sigma'_A(xy)/2$ holds for every edge.
	This implies that, for each $x \in S$,
	\begin{equation}
		\sigma'_A(x) = \sigma'_B(x).
		\label{equation:largetildesigma-sigmaprimeAsigmaprimeB}
	\end{equation}
	Similarly, the same argument we used above shows that for each $x \in S$ we have $\sigma_B(x) = \sigma^1_B(x)$. 
	Moreover, since $\sigma_B$ is $\gamma_B$-skew anchored in $S$, we have
	$(1 + \gamma_B)\sum_{x\in S}\sigma_B^1(x)= W(\sigma_B)$.
	Hence (since $N(d) \subseteq S$), we have
	\begin{equation}
	\deg_w(d, \sigma_B)\le \sum_{x\in S} \sigma_B^1(x) = \frac{1}{1 + \gamma_B} W(\sigma_B) = b_1.
	\label{equation:largetildesigma-degwdsigmaB}
	\end{equation}
	Now we claim that for each $x \in S$, we have
	\begin{equation}
		\mu_{\tilde \sigma}(x)+\tilde \mu(x) \geq w(\oriented{dx})-\sigma_B(x).
		\label{equation:largetildesigma-maximal}
	\end{equation}
	Indeed, let $x \in S$.
	We have
	\begin{align*}
		\mu_{\tilde \sigma}(x)+\tilde \mu(x)
		& \overset{\eqref{equation:largetildesigma-mutildesigmavssigmaprimeA}}{=}
		\sigma'_{A}(x)+\tilde \mu(x)
		\overset{\eqref{equation:largetildesigma-sigmaprimeAsigmaprimeB}}{=}
		\sigma'_{B}(x)+\tilde \mu(x)
		 = (\sigma'_{B}(x)+\tilde \mu(x)+\sigma_B(x)) - \sigma_B(x) \\
		& = 1 - \sigma_B(x) \geq w(\oriented{dx})-\sigma_B(x),
	\end{align*}
	in the last line we used that $\tilde \sigma+\tilde \mu=\sigma_B'+\sigma_B+\tilde \mu$ covers $S$ by \ref{item:strucprop-GEpaircovers}.
	This proves~\eqref{equation:largetildesigma-maximal}.
	
	Next, consider $\mu_d\le \mu_{\tilde \sigma}+\tilde \mu$ to be a maximal fractional matching such that for each $x\in N(d)$ we have $\mu_d(x) \leq \max\{0, w(\oriented{dx})-\sigma_B(x)\}$, and $\mu_d(\fmat{xy})=0$ if $\fmat{xy}$ does not intersect $N_w(d)$. 
	Since $\mu_{\tilde \sigma} \preceq \sigma'_B$ and $\sigma'_B \leq \tilde \sigma$, from \ref{item:strucprop-GEpairanchorS} and further since $N(d) \subseteq S$, we get that $\mu_d$ is only supported in edges with exactly one endpoint in $S_\mathcal R$.
	
	We claim that, in fact, for each $x \in N(d)$, we have
	\begin{equation}
		\mu_d(x) = \max\{0, w(\oriented{dx})-\sigma_B(x)\}, \label{eq:mu_dfitsInNw(d)}
	\end{equation}
	By the maximal choice of $\mu_d$, it is enough to check that for all such $x$, it holds that $\mu_{\tilde \sigma}(x)+\tilde \mu(x) \geq w(\oriented{dx})-\sigma_B(x)$, and this is true by \eqref{equation:largetildesigma-maximal}.
	This gives \eqref{eq:mu_dfitsInNw(d)}.
	
	We thus have 
	\begin{align*}
		W(\mu_d)
		& = \sum_{x \in N_w(d)} \mu_d(x) \overset{\eqref{eq:mu_dfitsInNw(d)}}{=} \sum_{x \in N_w(d)} \max\{0, w(\oriented{dx})-\sigma_B(x) \} \\
		& = \sum_{x \in N_w(d)} \left( w(\oriented{dx}) - \min \{ w(\oriented{dx}), \sigma_B(x) \} \right) =  \deg_w(d)-\deg_w(d,\sigma_B) \\
		& \geq \frac{k}{2} -\deg_w(d,\sigma_B) \overset{\eqref{equation:largetildesigma-degwdsigmaB}}{\ge} \frac{k}{2}-b_1\overset{\eqref{assumption-new:a_1+b_1-small}}{\ge }\max\{a_1,a_2\}.
	\end{align*}
	
We apply \Cref{lem:new-extending-out} (Extending-out) with \begin{center}
	\begin{tabular}{c|c|c|c|c|c}
		object & $H$ & $N_w(d)$ & $V(H) \setminus S$ & $\mu_d$ & $\gamma_A$ \\
		\hline
		in place of & $H$ & $U$ & $V$ & $\mu$ & $\gamma$
	\end{tabular}
\end{center} to obtain a $\gamma_A$-skew-matching $\sigma_A$ such that $\sigma_A \trianglelefteq \mu_d$, which has weight
	\[W(\sigma_A)\ge (1+\min\{\gamma_A,\gamma_A^{-1}\})W(\mu_d)\ge a_1+a_2,\]
and its anchor $\mathcal A(\sigma_A)$ is contained in $N_w(d)$.	
	
	We claim that $(\sigma_A,\sigma_B)$ is the desired $(\gamma_A,\gamma_B)$-skew-matching, anchored in $\oriented{dc}$.
	To do so, let us verify \ref{itdef:disjoint}--\ref{itdef:anchorpartition}.
	Indeed, by construction, $\sigma_A$ is disjoint from $\sigma_B$.
	We have that the anchor $\mathcal A(\sigma_A)$ is contained in $N_w(d)$.
	By \ref{item:strucprop-GEpairanchorfits} and $\sigma_B \leq \tilde \sigma$, we get that $\mathcal{A}(\sigma_B)$ fits in the $w$-neighbourhood of $c$.
	Also, for each $x\in N_w(d)\cap N_w(c)$ for which $\mu_d(x)>0$ we have that 
	\begin{align*}
	\sigma_A(x)\le \mu_d(x)\le w(\oriented{dx})-\sigma_B(x).
	\end{align*}
	In particular, for this we can deduce that the anchor $\mathcal A(\sigma_A)$ fits in the $w$-neighbourhood of~$d$; and moreover, that \ref{itdef:anchorpartition} holds.
	This proves the claim.
\end{proofclaim}

Owing to \Cref{cl:new-weighted-largetildesigma}, we may assume from now on that 
\begin{equation}\label{assumption:new-small-skew-weight}
W(\tilde \sigma)<b_1+b_2.
\end{equation}
Now we show that we are also done if $\tilde \sigma+\tilde \mu$ has sufficient weight in $N_w(c)$, i.e., if $\deg_w(c, \tilde \sigma+\tilde \mu) \geq k$ holds.

\begin{claim}\label{cl:new-weighted-coveringR}
	If $\deg_w(c, \tilde \sigma+\tilde \mu) \geq k$, then $H^\leftrightarrow$ has a good matching.
\end{claim}

\begin{proofclaim}[Proof of \Cref{cl:new-weighted-coveringR}]
	Let $d\in {\mathcal{R}}$ be arbitrary.
	We shall build a $(\gamma_A,\gamma_B)$-skew-matching pair anchored in $\oriented{dc}$.
	We summarise the proof strategy now.
	The neighbourhood of $c$ receives enough weight by $\tilde \mu+\tilde \sigma$, so we can build the whole pair within those objects.
	We shall use the full $\tilde \sigma$ to start building a $\gamma_B$-skew-matching with a fraction of the total required weight.
	In order to build the left-over of the skew pair in $\tilde \mu$, we shall first define a weight $w'$, somewhat similarly to a ``truncated weight'' (\Cref{def:truncated-weighted-graph}), but with the skew-matching $\tilde \sigma$ playing the role of the fractional matching $\mu$ in that definition.
	Next, we shall partition $\tilde \mu$ into three auxiliary fractional matchings $\mu_d, \bar{\mu}, \mu_c$, and $\mu'$, so that they are homogeneous with respect to the weights from $c$ and $d$. The fractional matching $\mu_d$ represents the part of $\tilde \mu$ that is anchored in $d$ and where we can place the $\gamma_A$-skew matching, and $\bar\mu$ its complement, where only the $\gamma_B$-skew matching can be placed. We further divide $\mu_d$ into two matchings $\mu_c$ and $\mu'$. The matching $\mu_c$ represents the part of $\mu_d$ that is not only accessible from $d$, but as well accessible from $c$ on both its side, allowing a lot of flexibility to build the $\gamma_B$-skew matching (using \Cref{lem:new-completion}). Before doing that however, we aim to build as much of the $\gamma_A$-skew matching in what remains of $\mu_d$ and $\mu'$ represents how much of $\mu_d-\mu_c$ we need for that. That is, if we manage to find the whole $\gamma_A$-skew matching in $\mu_d-\mu_c$, then $\mu'$ might be slightly smaller that $\mu_d-\mu_c$. This is a special (easier case) we treat separately in Step~4. If we do not manage to find the whole $\gamma_A$-skew matching in $\mu_d-\mu_c$, then the rest will have to share the space with the $\gamma_B$-skew matching in $\mu_c$ (and $\mu'=\mu_d-\mu_c$). The latter situation is the core of the proof and is in Step~5 further  divided in two cases (Cases~1 and Case~2) according to how much of the $\gamma_B$-skew matching we manage to build in $\mu_c$ beside the leftover of the $\gamma_A$-skew matching.

	Some readers may be slightly confused during the proof with the order in which the skew matchings and in which order the truncated weights are defined. The philosophy behind it is the following. The fractional matching $\mu_c$ is defined first, w.r.t. $w'$. It is then cut away (reserved for later) and $\mu_d-\mu_c$ is filled using the truncated weight~$\tilde w$, obtained from~$w'$. Next, the reserved fractional matching $\mu_c$ is used up using the weight $w'$. Last, if needed, we cut away the whole $\mu_d$, which is already completely used up, and fit whatever we can in $\bar\mu$ using the truncated weight $\bar w$. At the end we glue the different bits together.
	 \medskip
	
	\noindent \emph{Step 1: Defining the ``truncated weight''.}
	By \ref{item:SP-reachableareisolated} we have that $d \in N_w(c)$.
	We define an auxiliary weight $w'$ by setting 
	\begin{align*}
		w'(\oriented{ux}) & = \begin{cases}
			\max\{0,w(\oriented{ux})-\tilde \sigma(x)\} & \text{if $u \in \{c,d\}$ and $x \in V(H)\setminus \{c,d\}$,} \\
			w(\oriented{ux}) & \text{otherwise.}
		\end{cases}
	\end{align*}

	As mentioned before, $w'$ is analogous to a truncated weight (Definition~\ref{def:truncated-weighted-graph}), but taken instead with respect to a skew-matching $\tilde \sigma$ and somewhat simplified.
	Observe that 
	\begin{align}
	\nonumber \deg_{w'}(c, \tilde \mu)&= \sum_{v\in V(H)}\left(\min\{w'(\oriented{cx})+\tilde \sigma(x), \tilde \mu(x)+\tilde \sigma(x)\}-\tilde \sigma(x)\right)\\
	&\ge \sum_{x\in V(H)}\min\{w(\oriented{cx}), \tilde \mu(x)+\tilde \sigma(x)\}-W(\tilde \sigma) \nonumber \\
	&=\deg_w(c, \tilde \mu+\tilde \sigma)-W(\tilde \sigma).\label{eq:takingouttildesigma}
	\end{align}
 
	By \ref{item:strucprop-GEpaircovers} we have, for each $x \in S$, that $\tilde \mu(x) = 1 - \tilde \sigma(x) \geq w(\oriented{dx}) - \tilde \sigma(x)$; hence $w'(\oriented{dx}) \leq \tilde \mu(x)$ holds for each $x \in S$.
	This implies that
	\begin{align}
		\deg_{w'}(d, \tilde \mu)
		& = \sum_{x \in S} \min\{ w'(\oriented{dx}), \tilde \mu(x) \} = \sum_{x \in S} w'(\oriented{dx}) \geq \sum_{x \in S} (w(\oriented{dx}) - \tilde \sigma(x)) \nonumber \\
		& = \deg_{w}(d) - \sum_{x \in S_{\mathcal R}} \tilde \sigma(x) \geq \frac{k}{2} - \frac{W(\tilde \sigma)}{1 + \gamma_B}, \label{equation:skewmatchingcover2-degwprimedtildemu}
	\end{align}
	where in the last inequality we used \ref{item:strucprop-GEpairanchorS},  \ref{item:strucprop-GEskewinR} and $\deg_{w}(d) \geq k/2$.\medskip
	
	\noindent
	\emph{Step 2: Defining the auxiliary fractional matchings.}
	Let $\mu_d\le \tilde \mu$ be a maximal fractional matching such that for each $x\in N_{w'}(d)$, we have $\mu_d(x)\le w'(\oriented{dx})$ and its support edges intersect $N_{w'}(d)$.
	We have $\mu_d(x) = \min\{ w'(\oriented{dx}), \tilde \mu(x) \}$ for each $x \in S$, and therefore
	\begin{equation}\label{eq:mu_dAvoidingc,d}
	\deg_{w'}(d, \mu_d)= \deg_{w'}(d, \tilde \mu).
	\end{equation}
	Now, let $\mu_c\le \mu_d$ be a maximal fractional matching so that $\mu_c(x)\le w'(\oriented{cx})$ for all $x\in N_{w'}(d)$.
	Similarly, by construction we have that $\mu_c(x) = \min\{ w'(\oriented{cx}), \mu_d(x) \} = \min\{ w'(\oriented{cx}), w'(\oriented{dx}), \tilde \mu(x) \}$ for all $x \in S$. Let $\mu'\le \mu_d-\mu_c$ be maximal such that
	\begin{equation}\label{eq:struct-tildesigmaANotTooBig}
	W(\mu')\le \frac{a_1+a_2}{1+\min\{\gamma_A, \gamma_A^{-1}\}}.
	\end{equation}
	Last, let $( H, \tilde w)$ be the $\mu_c$-truncated graph obtained from $(H,w')$, let $(H, \bar w)$ be the $\mu_d$-truncated weighted graph obtained from $(H, w')$, and set $\bar \mu= \tilde \mu-\mu_d$. 
	
	We now claim the following property, that we will need repeatedly later:
	\begin{equation}
		\text{for all $v \notin S$, $\mu_c(v) \leq w'(\oriented{cv})$.}
		\label{equation:vnotSmucwprime}
	\end{equation}
	Indeed, if $\mu_c(v) = 0$ there is nothing to check, so assume $\mu_c(v) > 0$.
	We have $0 < \mu_c(v) \leq \mu_d(v) \leq \tilde \mu(v)$.
	In particular, there exists $u \in S$ such that $\mu_d(uv) > 0$.
	By definition of $\mu_d$, we must have $u \in N_{w'}(d)\subseteq S_{\mathcal R}$ and also $0 < \mu_d(uv) \leq \tilde \mu(uv)$. By \ref{item:strucprop-GEpairsupport}, we get
 $v\in \mathcal R$. Hence by \ref{item:structuprop-GEpairsaturation}  we have $\tilde \mu(v) +\tilde \sigma(v) \leq w(\oriented{cv})$, and thus $\mu_c(v)\le \tilde \mu (v)\le w'(\oriented{cv})$. This shows that \eqref{equation:vnotSmucwprime} holds.\medskip

\noindent
	\emph{Step 3: Starting building the $\gamma_A$-skew-matching in $\mu'$}.
	We apply \Cref{lem:new-extending-out} (Extending-out) with
	\begin{center}	
	\begin{tabular}{c|c|c|c|c|c}
		object & $H$ & $S$ & $V(H) \setminus S$ & $\mu'$ & $\gamma_A$ \\
		\hline
		in place of & $H$ & $U$ & $V$ & $\mu$ & $\gamma$
	\end{tabular}
	\end{center}
	By doing so, we obtain a $\gamma_A$-skew-matching $\tilde \sigma_A$ in $ H^\leftrightarrow$ such that
	\begin{equation}
		\tilde \sigma_A \trianglelefteq\mu'\le \mu_d-\mu_c,
		\label{equation:skewmatchingcover2-tildesigmatriangle}
	\end{equation}
	of weight 
	\begin{equation}
		W(\tilde \sigma_A)=(1+\min\{\gamma_A, \gamma_A^{-1}\})W(\mu')\overset{\eqref{eq:struct-tildesigmaANotTooBig}}{\le} a_1+a_2,
		\label{equation:skewmatchingcover2-weigthtildesigmaA}
	\end{equation}
	and whose anchor $\mathcal A(\tilde \sigma_A)$ is contained in $S$. 
	We claim that
	\begin{equation}
		\text{the anchor $\mathcal A(\tilde\sigma_A)$ fits in the $\tilde w$-neighbourhood of $d$.}
		\label{equation:skewmatchingcover2-anchortildesigmafits}
	\end{equation}
	Indeed, for all $x\in N_{w'}(d)\subseteq S$, using \eqref{equation:skewmatchingcover2-tildesigmatriangle} we have 
	\begin{align*}
	\tilde \sigma_A^1(x)&\le \mu'(x) \leq \mu_d(x)-\mu_c(x)\overset{\eqref{equation:vnotSmucwprime}}{=} \mu_d(x)-\min\{w'(\oriented{cx}),\mu_d(x)\}\\
	&=\max\{0,\mu_d(x)-w'(\oriented{cx})\} \le \max\{0,w'(\oriented{dx})-\mu_c(x)\} = \tilde w(\oriented{dx}),
	\end{align*}
	where the last inequality follows from $\mu_d(x) \leq w'(\oriented{dx})$ and $\mu_c(x) \leq w'(\oriented{cx})$.

	Now we claim that 
\begin{equation}\label{eq:weightTidleSigmaA}
W(\tilde \sigma_A)\ge \deg_{\tilde w}(c, \mu').
\end{equation}
Indeed, $\deg_{\tilde w}(c, \mu')$ is the sum of the terms $\min\{\tilde w(\oriented{cx}), \mu'(x)\}$, and those are all zero for $x \in S$.
Hence the sum incorporates the weight $\mu'$ of each edge at most once, and thus is at most $\deg_{\tilde w}(c, \mu') \leq W(\mu')$.
Then we get the desired inequality by \eqref{equation:skewmatchingcover2-weigthtildesigmaA}.\medskip

 \noindent \emph{Step 4: The special (easy) case when $W(\mu_d - \mu_c)$ is large.}
 Before continuing, we will show that we can conclude in the case where the inequality $W(\mu_d-\mu_c)>W(\mu')$ holds.
 If the fractional matching $\mu'$ does not spend the whole weight of $\mu_d-\mu_c$, it means it is big enough to accommodate a whole $\gamma_A$-skew matching of weight $a_1+a_2$.
 Indeed, we have in this case equality in~\eqref{eq:struct-tildesigmaANotTooBig}, leading to equality in~\eqref{equation:skewmatchingcover2-weigthtildesigmaA}, i.e. we have that $W(\tilde \sigma_A) = a_1 + a_2$.
 This means that all is left to do is to build the rest of the $\gamma_B$-skew fractional matching and glue everything together.
 
 First, we apply \Cref{lem:new-balancing} (Balancing-Out) with
 \begin{center}	
 \begin{tabular}{c|c|c|c|c}
 	object & $H$ & $V(\mu_c)$ & $\mu_c$ & $\gamma_B$ \\
 	\hline
 	in place of & $H$ & $U$ & $\mu$ & $\gamma$
 \end{tabular}
 \end{center}
 to obtain a $\gamma_B$-skew matching $\widehat \sigma_B$ such that $\widehat \sigma_B \trianglelefteq \mu_c$ and $W( \widehat \sigma_B) = 2 W(\mu_c)$.
 Since $\widehat \sigma_B \trianglelefteq \mu_c$, together with the fact that $\mu_c(x) \leq w'(\oriented{cx})$ for $x \in S$ and \eqref{equation:vnotSmucwprime}, implies that in fact $\widehat \sigma_B$ fits in the $w'$-neighbourhood of $c$.
 
 Secondly, we let $(H, \hat w)$ be the $\mu'$-truncated graph obtained from $(H, \tilde w)$.
 We claim that
 \begin{equation}
 	\deg_{\widehat w}(c, \tilde \mu-\mu'-\mu_c) \ge b_1+b_2-W(\widehat \sigma_B)-W(\tilde \sigma).
 	\label{equation:deghatwctildemumuprimemuc}
 \end{equation}
 Indeed, using \Cref{prop:adding-degree-truncated} twice in the first two equalities 
\begin{center}	
 \begin{tabular}{c|c|c|c|c}
 	object & $(H, \tilde w)$ & $\mu'$ & $\tilde \mu-\mu_c$ & $\hat w$ \\
 	\hline
		object & $(H, w')$ & $\mu_c$ & $\tilde \mu$ & $\tilde w$ \\
 	\hline
 	in place of & $(G,w)$ & $\mu'$ & $\mu$ & $w'$
 \end{tabular}
 \end{center}

 we have
 \begin{align*}
 	\deg_{\widehat w}(c, \tilde \mu-\mu'-\mu_c)
 	& = \deg_{\tilde w}(c, \tilde \mu - \mu_c) - \deg_{\tilde w}(c,\mu') \\
 	& = \deg_{w'}(c, \tilde \mu) - \deg_{w'}(c, \mu_c) - \deg_{\tilde w}(c,\mu') \\
 	& \geq \deg_{w'}(c, \tilde \mu) - 2 W(\mu_c) - \deg_{\tilde w}(c,\mu') \\
 	& = \deg_{w'}(c, \tilde \mu) - W(\hat \sigma_B) - \deg_{\tilde w}(c,\mu') \\
 	& \overset{\eqref{eq:takingouttildesigma}}{\ge} \deg_{w'}(c, \tilde \mu + \tilde \sigma) - W(\tilde \sigma) - W(\hat \sigma_B) - \deg_{\tilde w}(c,\mu') \\
 	& \overset{\eqref{eq:weightTidleSigmaA}}{\ge} \deg_{w'}(c, \tilde \mu + \tilde \sigma) - W(\tilde \sigma) - W(\hat \sigma_B) - W(\tilde \sigma_A) \\
 	& \overset{\eqref{equation:skewmatchingcover2-weigthtildesigmaA}}{\ge} \deg_{w'}(c, \tilde \mu + \tilde \sigma) - W(\tilde \sigma) - W(\hat \sigma_B) - (a_1 + a_2) \\
 	& = \deg_{w'}(c, \tilde \mu + \tilde \sigma) - k - W(\tilde \sigma) - W(\hat \sigma_B) + (b_1 + b_2),
 \end{align*}
 which gives \eqref{equation:deghatwctildemumuprimemuc} by recalling that $\deg_{w'}(c, \tilde \mu + \tilde \sigma) \geq k$ holds by the claim assumption.
 Thanks to \eqref{equation:deghatwctildemumuprimemuc}, we can use \Cref{lem:new-combination} (Combination) with
  \begin{center}	
 	\begin{tabular}{c|c|c|c|c}
 		object & $(H,\hat w)$ & $c$ & $\tilde \mu - \mu_c - \mu'$ & $\gamma_B$ \\
 		\hline
 		in place of & $(H,w)$ & $v$ & $\mu$ & $\gamma$
 	\end{tabular}
 \end{center}
 to obtain a $\gamma_B$-skew matching $\bar \sigma_B$ such that $\bar \sigma_B \trianglelefteq \tilde \mu - \mu_c - \mu'$, $W(\bar \sigma_B) \geq b_1 + b_2 - W(\widehat \sigma_B) - W(\tilde \sigma)$, and such that $\mathcal{A}(\bar \sigma_B)$ fits in the $\hat w$-neighbourhood of $c$.
 
 We have finalised the construction of the skew matchings, and we need to `glue' them; namely, we need to argue that $\tilde \sigma + \widehat{\sigma}_B + \bar \sigma_B$ forms a skew-matching pair together with $\tilde \sigma_A$.
 Recall that $\sigma_\emptyset$ is the empty skew-matching.
 We apply \Cref{prop:adding-skew-matchings} with
 \begin{center}	
 	\begin{tabular}{c|c|c|c|c|c|c|c|c|c}
 		object & $(H,\tilde{w})$ & $(d,c)$ & $\mu'$ & $\tilde{\mu} - \mu_c - \mu'$ & $(H, \hat{w})$ & $\tilde \sigma_A$ & $\sigma_\emptyset$ & $\sigma_\emptyset$ & $\bar \sigma_B$ \\
 		\hline
 		in place of & $(H,w)$ & $(u,v)$ & $\mu$ & $\bar \mu$ & $(H, \bar w)$ & $\sigma_A$ & $\sigma_B$ & $\bar \sigma_A$ & $\bar \sigma_B$
 	\end{tabular}
 \end{center}
 to obtain that $(\tilde \sigma_A, \bar \sigma_B)$ is a $(\gamma_A, \gamma_B)$-skew matching pair in $(H, \tilde{w})$, anchored in $\oriented{dc}$, with $\tilde \sigma_A + \bar \sigma_B \trianglelefteq \tilde{\mu} - \mu_c$.
 Next, we apply \Cref{prop:adding-skew-matchings} again, now with
 \begin{center}	
 	\begin{tabular}{c|c|c|c|c|c|c|c|c|c}
 		object & $(H,w')$ & $(d,c)$ & $\mu_c$ & $\tilde{\mu} - \mu_c$ & $(H, \tilde{w})$ & $\sigma_\emptyset$ & $\widehat \sigma_B$ & $\tilde \sigma_A$ & $\bar \sigma_B$ \\
 		\hline
 		in place of & $(H,w)$ & $(u,v)$ & $\mu$ & $\bar \mu$ & $(H, \bar w)$ & $\sigma_A$ & $\sigma_B$ & $\bar \sigma_A$ & $\bar \sigma_B$
 	\end{tabular}
 \end{center}
 to obtain that $(\tilde \sigma_A, \widehat{\sigma}_B + \bar \sigma_B)$
  is a $(\gamma_A, \gamma_B)$-skew matching pair in $(H, w')$, anchored in $\oriented{dc}$, with $\tilde \sigma_A + \bar \sigma_B + \widehat{\sigma}_B \trianglelefteq \tilde{\mu}$.
 
 We would like to apply \Cref{prop:adding-skew-matchings} again to incorporate $\tilde \sigma$, but formally this does not work because $w'$ is not a truncated weight obtained from $w$.
 Nevertheless, the argument is morally the same and we sketch it for completeness.
 We need to show that $(\bar \sigma_A, \widehat{\sigma}_B + \bar \sigma_B + \tilde \sigma)$ is a $(\gamma_A, \gamma_B)$-skew matching pair in $(H, w)$, anchored in $\oriented{dc}$.
 We shall verify \ref{itdef:disjoint}--\ref{itdef:anchorpartition}.
 Property \ref{itdef:disjoint} follows because $\tilde \sigma_A + \bar \sigma_B + \widehat{\sigma}_B \trianglelefteq \tilde{\mu}$ and \ref{item:strucprop-GEpairdisjoint}.
 We already have \ref{itdef:anchorc}.
 To see \ref{itdef:anchord}, let $x \in N_w(c)$.
 Using property \ref{itdef:anchord} for $(\bar \sigma_A, \widehat{\sigma}_B + \bar \sigma_B)$ and $w'$, we see that $\widehat{\sigma}^2_B(x) + \bar \sigma^1_B(x) + \tilde \sigma^1(x) \leq w'(\oriented{cx}) + \tilde \sigma(x) = w(\oriented{cx})$, where the last inequality is because the anchor of $\tilde \sigma$ fits in the $w$-neighbourhood of $c$.
 Lastly, \ref{itdef:anchorpartition} follows as in the proof of \Cref{prop:adding-skew-matchings}, considering three cases depending if $x \in N_{w'}(c)$ and $x \in N_{w'}(d)$, or not; we omit further details.
 
 To finish, we argued at the beginning of this step that the assumption $W(\mu_d-\mu_c)>W(\mu')$ implies that $W(\tilde \sigma_A) = a_1 + a_2$, and by construction we have $W(\tilde \sigma) + W(\widehat \sigma_B) + W(\bar \sigma_B) = b_1 + b_2$.
 Hence, we are done with the construction in this case. \medskip
 
  \noindent \emph{Step 5: The main case distinction.}
  Hence, from now on, we can assume that $W(\mu_d-\mu_c) \leq W(\mu')$.
  Since $\mu' \leq \mu_d - \mu_c$, this implies that
  \begin{equation}\label{eq:mu'=mu_d-mu_c2}
  	\mu'=\mu_d-\mu_c.
  \end{equation}
   Now, we set the auxiliary parameters
  \begin{align*}
  	\alpha_1 & := a_1 - \frac{W(\tilde\sigma_A)}{1+\gamma_A}, &
  	\alpha_2 & := \gamma_A \alpha_1, \\
  	\beta_1 & := b_1-\frac{W(\tilde \sigma)}{1+\gamma_B}, &
  	\beta_2 & := \gamma_B \beta_1.
  \end{align*}
  Note that $\alpha_1, \alpha_2 \geq 0$ thanks to \eqref{equation:skewmatchingcover2-weigthtildesigmaA}, and $\beta_1, \beta_2 \geq 0$ thanks to \eqref{assumption:new-small-skew-weight}.
  We will use the parameters $\alpha_1, \alpha_2, \beta_1, \beta_2$ to build auxiliary skew-matchings in the cases that follow.
 
 We do some preliminary calculations that will be useful.
 Using that $\mu_c(x) \leq \mu_d(x) \leq w'(\oriented{dx})$ for all $x \in S$, we get that $W(\mu_c) = \deg_{w'}(d,\mu_c)$ and that $W(\mu_d) = \deg_{w'}(d, \mu_d)$.
From this, we can observe that 
	\begin{align}\nonumber
	W(\mu_c)
	& = \deg_{w'}(d, \tilde \mu) + W(\mu_c) - \deg_{w'}(d, \tilde \mu) \\ \nonumber
	& \overset{\eqref{eq:mu_dAvoidingc,d}}{= } \deg_{w'}(d, \tilde \mu) + W(\mu_c) - \deg_{w'}(d, \mu_d)\\ \nonumber
	& = \deg_{w'}(d, \tilde \mu) - W(\mu_d - \mu_c) \\
	\nonumber 
	& \overset{\eqref{eq:mu'=mu_d-mu_c2}}{=} \deg_{w'}(d, \tilde \mu) - W(\mu') \\
	\nonumber 
	& \overset{\eqref{equation:skewmatchingcover2-degwprimedtildemu}}{\geq} \frac{k}{2}-\frac{W(\tilde \sigma)}{1+\gamma_B} - W(\mu') \\
	\nonumber 
	& \overset{\eqref{equation:skewmatchingcover2-weigthtildesigmaA}}{=} \frac{k}{2}-\frac{W(\tilde \sigma)}{1+\gamma_B}-\max\{1,\gamma_A\}\frac{W(\tilde \sigma_A)}{1+\gamma_A}\\
	& \nonumber \overset{\eqref{assumption-new:a_1+b_1-small}}{> }
	\max\{a_1,a_2\}+b_1-\frac{W(\tilde \sigma)}{1+\gamma_B}-\max\{1,\gamma_A\}\frac{W(\tilde\sigma_A)}{1+\gamma_A}\\
	& =
	\max\{1,\gamma_A\}\left( a_1 - \frac{W(\tilde\sigma_A)}{1+\gamma_A}\right) + \beta_1 \nonumber \\
	& = \max\{\alpha_1,\alpha_2\} + \beta_1  \nonumber \\
	& \geq \max\{\alpha_1,\alpha_2\}+\min\{\beta_1,\beta_2\}.\label{eqnew:deg->mu_d}
	\end{align}

From now on, we separate the proof of the claim into two cases, depending if $W(\mu_c)$ is sufficiently large with respect to $\alpha_1, \alpha_2, \beta_1, \beta_2$ or not. \medskip

\noindent \emph{Case 1: $W(\mu_c)$ is large.}
In this case, we will assume that
\begin{align}\label{eq:minmax<mu_c}
	\min\{\alpha_1,\alpha_2\}+\max\{\beta_1,\beta_2\}<W(\mu_c).
\end{align}

Since $\mu_c$ is large enough, we can accommodate the left-over of the skew-matching pair in it and will not need to use $\bar \mu$ at all. However, the setting is such that we cannot guarantee the use of \Cref{lem:new-completion} (Completion), and therefore, we shall have to build the skew-matching pair in $\mu_c$ ``by hand''. \medskip

\noindent \emph{Case 1, Step I: Building a skew-matching pair in $\mu_c$.}
We will define two auxiliary skew-matchings, $\widehat\sigma_A$ and $\widehat\sigma_B$, which are $\gamma_A$-skew and $\gamma_B$-skew, respectively.
For every $xy \in E(H)$ with $x \in S$, $y \notin S$, we let
\begin{align*}
\widehat\sigma_A(\oriented{xy}):= \frac{\alpha_1+ \alpha_2}{W(\mu_c)}\mu_c(\fmat{xy
}),
\end{align*}
and $\widehat\sigma_A$ takes the value $0$ in any other case.
Note that
\begin{align}
	W(\widehat\sigma_A) = \sum_{x \in S} \sum_{y \notin S} \frac{\alpha_1+ \alpha_2}{W(\mu_c)}\mu_c(\fmat{xy
	}) = \alpha_1 + \alpha_2.
	\label{equation:claimsaturated-weigthwidehatsigmaA}
\end{align}

The definition of $\widehat\sigma_B$ depends if $\gamma_A < 1$ or not.
For brevity, we use the Iverson bracket notation: for a logical statement $P$ the value $[P]$ equals $1$ if $P$ is true, and $0$ otherwise.
Now consider any $xy \in E(H)$ with $x \in S$ and $y \notin S$.
We define
\begin{align*}
\widehat\sigma_B(\oriented{xy}):= [\gamma_A < 1]\frac{\beta_1+\beta_2}{W(\mu_c)}\mu_c(\fmat{xy}),
\end{align*}
as well as
\begin{align*}
\widehat\sigma_B(\oriented{yx}):= [\gamma_A \geq 1]\frac{\beta_1+\beta_2}{W(\mu_c)}\mu_c(\fmat{xy});
\end{align*}
and $\widehat\sigma_B$ takes the value $0$ in every other ordered edge.
Similarly as before, we can observe that
\begin{align}
	W(\widehat\sigma_B) = \beta_1 + \beta_2.
	\label{equation:claimsaturated-weigthwidehatsigmaB}
\end{align}

Now, we claim that
\[\widehat\sigma_A+\widehat\sigma_B\trianglelefteq\mu_c.\]
To do this, let $xy \in E(H)$ be arbitrary such that $x \in S$ and $y \notin S$ (as otherwise there is nothing to check).
Observe that 
\begin{align*}
\frac{\widehat\sigma_A(\oriented{xy})+\gamma_A\widehat\sigma_A(\oriented{yx})}{1+\gamma_A}
&+ \frac{\widehat\sigma_B(\oriented{xy})+\gamma_B\widehat\sigma_B(\oriented{yx})}{1+\gamma_B}\\
& = \left( \frac{\alpha_1 + \alpha_2}{1+\gamma_A} + ([\gamma_A < 1] + \gamma_B [\gamma_A \geq 1]) \frac{\beta_1 + \beta_2}{1+\gamma_B} \right) \frac{\mu_c(\fmat{xy})}{W(\mu_c)}\\
& = \left( \alpha_1 + \beta_1([\gamma_A < 1] + \gamma_B [\gamma_A \geq 1]) \right) \frac{\mu_c(\fmat{xy})}{W(\mu_c)}.
\end{align*} 
By \ref{eqnew:gamma_B>1}, we have that $\beta_2 = \gamma_B \beta_1 \geq \beta_1$.
A brief moment of reflection (considering the cases $\gamma_A \geq 1$ and $\gamma_A < 1$ separately) then implies that $\alpha_1 + \beta_1([\gamma_A < 1] + \gamma_B [\gamma_A \geq 1])$ is equal either to $\max\{\alpha_1, \alpha_2\} + \min\{\beta_1, \beta_2\}$ or $\min\{\alpha_1, \alpha_2\} + \max\{\beta_1, \beta_2\}$.
In any case, either by \eqref{eqnew:deg->mu_d} or \eqref{eq:minmax<mu_c}, we can conclude this term is at most $W(\mu_c)$.
This implies that we have
\[ \frac{\widehat\sigma_A(\oriented{xy})+\gamma_A\widehat\sigma_A(\oriented{yx})}{1+\gamma_A}
+ \frac{\widehat\sigma_B(\oriented{xy})+\gamma_B\widehat\sigma_B(\oriented{yx})}{1+\gamma_B} \leq \mu_c(\fmat{xy}), \]
as desired.
An analogous calculation implies that
\begin{align*}
\frac{\widehat\sigma_A(\oriented{yx})+\gamma_A\widehat\sigma_A(\oriented{xy})}{1+\gamma_A} + \frac{\widehat{\sigma}_B(\oriented{yx}) + \gamma_B \widehat{\sigma}_B(\oriented{xy})}{1 + \gamma_B} \leq \mu_c(\fmat{xy}),
\end{align*}
thus indeed $\widehat\sigma_A+\widehat\sigma_B\trianglelefteq\mu_c$ holds.\medskip

\noindent \emph{Case 1, Step II: Gluing the skew-matchings together.}
Now, we claim that $(\widehat\sigma_A, \widehat\sigma_B)$ is a $(\gamma_A, \gamma_B)$-skew-matching pair in $(H, w')$ anchored in $\oriented{dc}$, with respect to $w'$.
We check the required properties \ref{itdef:disjoint}--\ref{itdef:anchorpartition}.
The disjointness of $\widehat\sigma_A$ and $\widehat\sigma_B$ follows from $\widehat\sigma_A+\widehat\sigma_B\trianglelefteq\mu_c$, so \ref{itdef:disjoint} holds.
We check that $\widehat \sigma_A$ fits in the $w'$-neighbourhood of $d$.
Since $\widehat\sigma_A$ is supported only in directed edges $xy$ with $x \in S$ and $y \notin S$, we have $\widehat \sigma^1_A(v) = 0$ for each $v \notin S$.
For each $v \in S$, again $\widehat\sigma_A+\widehat\sigma_B\trianglelefteq\mu_c$ implies that $\widehat \sigma^1_A(v) \leq \mu_c(v) \leq \mu_d(v) \leq w'(\oriented{dv})$ (where we used $\mu_c \leq \mu_d$, the definition of $\mu_d$, and $v \in S$ in the last inequalities); this implies that $\widehat\sigma_A$ fits in the $w'$-neighbourhood of $d$.
Now, we check that $\widehat \sigma_B$ fits in the $w'$-neighbourhood of $c$.
If $v \in S$, a similar argument as before (using the definition of $\mu_c$) gives that $\widehat \sigma^1_B(v) \leq \mu_c(v) \leq w'(\oriented{cv})$.
Hence, we can suppose that $v \notin S$.
We have $0 < \widehat \sigma^1_B(v) \leq \mu_c(v) \leq w'(\oriented{cv})$, where the last inequality is from \eqref{equation:vnotSmucwprime}.
Thus indeed $\widehat \sigma_B$ fits in the $w'$-neighbourhood of $c$; hence \ref{itdef:anchorc}--\ref{itdef:anchord} hold.
Finally, for each $v \in N_{w'}(c) \cap N_{w'}(d) \subseteq S$ we have $\widehat\sigma_A^1(v)+\widehat\sigma^1_B(v) \leq \mu_c(v) \leq w'(\oriented{cv})$ which gives \ref{itdef:anchorpartition}.
We have then checked that $(\widehat\sigma_A, \widehat\sigma_B)$ is a $(\gamma_A, \gamma_B)$-skew-matching pair in $(H, w')$ anchored in $\oriented{dc}$.

Set $\sigma_A:= \widehat\sigma_A+\tilde \sigma_A$, and let $\sigma_\emptyset$ be the identically-zero skew-matching.
Recalling~\eqref{equation:skewmatchingcover2-tildesigmatriangle} and~\eqref{equation:skewmatchingcover2-anchortildesigmafits} reveals that we can use \Cref{prop:adding-skew-matchings} with
\begin{center}	
	\begin{tabular}{c|c|c|c|c|c|c|c|c|c}
		object & $(H,w')$ & $(d,c)$ & $\mu_c$ & $\mu_d - \mu_c$ & $(H, \tilde w)$ & $\widehat \sigma_A$ & $\widehat \sigma_B$ & $\tilde \sigma_A$ & $\sigma_\emptyset$ \\
		\hline
		in place of & $(H,w)$ & $(u,v)$ & $\mu$ & $\bar \mu$ & $(H, \bar w)$ & $\sigma_A$ & $\sigma_B$ & $\bar \sigma_A$ & $\bar \sigma_B$
	\end{tabular}
\end{center}
and we get that $(\sigma_A, \widehat\sigma_B)$ is a $(\gamma_A, \gamma_B)$-skew-matching in $H^\leftrightarrow$ anchored in $\oriented{dc}$ (with respect to $w'$) with $\sigma_A+\widehat\sigma_B\trianglelefteq\mu_d$.

We set $\sigma_B:= \tilde \sigma+\widehat\sigma_B$.
Using \eqref{equation:claimsaturated-weigthwidehatsigmaA} and the definition of $\alpha_1, \alpha_2$, we have
\begin{align*}
	W(\sigma_A)& = W(\tilde \sigma_A)+\alpha_1+\alpha_2 = W(\tilde \sigma_A) + (1 + \gamma_A)\left(a_1 - \frac{W(\tilde\sigma_A)}{1+\gamma_A}\right) = a_1+a_2,
\end{align*}
and similarly, using \eqref{equation:claimsaturated-weigthwidehatsigmaB} and the definition of $\beta_1, \beta_2$, we have
\begin{align*}
	W(\sigma_B)& = W(\tilde \sigma)+\beta_1+\beta_2=b_1+b_2.
\end{align*}

To conclude in this case, we just need to check that $(\sigma_A, \sigma_B)$ is a $(\gamma_A, \gamma_B)$-skew-matching in $H^\leftrightarrow$ anchored in $\oriented{dc}$.
Recall that $\mu_d \leq \tilde \mu$.
Since $\tilde \mu$ and $\tilde \sigma$ are disjoint, so are $\sigma_A$ and $\sigma_B$, which gives \ref{itdef:disjoint}.
Point \ref{itdef:anchorc} follows from the fact (that we have already checked) that $\sigma_A$ is a $\gamma_A$-skew-matching which fits in the $w'$-neighbourhood of $d$; and $w' \leq w$.
To see \ref{itdef:anchord} we need to check that $\sigma_B$ is a $\gamma_B$-skew-matching, whose anchor fits in the $w$-neighbourhood of $c$.
We already know that $\widehat\sigma_B$ is a $\gamma_B$-skew-matching whose anchor fits in the $w'$-neighbourhood of $c$.
Because of \ref{item:strucprop-GEpairanchorS} and \ref{item:strucprop-GEskewinR}, we just need to check the property for $x \in S_{\mathcal R}$.
In such a case, we have
\begin{align}\label{eq:anchor-sigma_B}
	\sigma^1_B(x)
	& = \widehat \sigma^1_B(x) + \tilde \sigma^1(x)
	= \widehat \sigma^1_B(x) + \tilde \sigma(x) \leq w'(\oriented{cx}) + \tilde \sigma(x).
\end{align}
We check two cases depending on the definition of $w'$.
If $w'(\oriented{cx}) = 0$, then \eqref{eq:anchor-sigma_B} gives $\sigma^1_B(x) \leq \tilde \sigma(x) \leq w(\oriented{cx})$, as desired.
Otherwise, we have $w'(\oriented{cx}) = w(\oriented{cx}) - \tilde \sigma(x)$, so $\sigma^1_B(x) \leq w(\oriented{cx})$, in which case again we are done.
This gives \ref{itdef:anchord}.
Finally, to see \ref{itdef:anchorpartition} let $x \in N(c) \cap N(d) \subseteq S$.
Note that $\sigma^1_A(x) + \sigma^1_B(x) \leq \mu_d(x) + \tilde \sigma(x) \leq w'(\oriented{dx}) + \tilde \sigma(x) \leq \max\{ w(\oriented{dx}), \tilde \sigma(x) \} \leq \max\{w(\oriented{dx}), w(\oriented{cx})\}$, where in the second to last inequality we used the definition of $w'$, and in the last inequality we used \ref{item:strucprop-GEpairanchorfits}.
This finishes the proof in Case~1. \medskip

\noindent \emph{Case 2: $W(\mu_c)$ is small.}
We can assume that Case 1 does not hold. Hence \eqref{eq:minmax<mu_c} fails to hold, and we can assume that
\begin{align} \label{eq:enoughmaterial}
	\min\{\alpha_1,\alpha_2\}+\max\{\beta_1,\beta_2\}\ge W(\mu_c).
\end{align}
In this case, since $W(\mu_c)$ is small, we shall need possibly to complement whatever we build in $\mu_d$ by a $\gamma_B$-skew-matching in $\bar\mu$.
On the other hand, the setting is such that for building a skew-matching pair in $\mu_c$, we may use \Cref{lem:new-completion} (Completion). \medskip

\noindent \emph{Case 2, Step I: Building a skew-matching pair in $\mu_c$.}
We wish to apply \Cref{lem:new-completion} (Completion) with
\begin{center}	
\begin{tabular}{c|c|c|c|c|c|c|c|c|c}
	object & $(H,w')$ & $N_{w'}(d)$ & $V(H)\setminus S$ & $c$ & $\mu_c$ & $\alpha_1$ & $\alpha_2$ & $\beta_1$ & $\beta_2$ \\
		\hline
	in place of & $(H,w)$ & $U$ & $V$ & $u$ & $\mu$ & $\alpha_1$ & $\alpha_2$ & $\beta_1$ & $\beta_2$
\end{tabular}
\end{center}
We check the required hypothesis.
To check \ref{item:completion-1}, we need that for all $y \notin S$, $\mu_c(y) \leq w'(\oriented{cy})$, and this follows from \eqref{equation:vnotSmucwprime}.
The other required inequalities~\ref{eq:neq-fit} and~\ref{eq:VybalancovaniMozne} follow from \eqref{eqnew:deg->mu_d} and~\eqref{eq:enoughmaterial}, respectively.

The application of \Cref{lem:new-completion} gives us a $\gamma_A$-skew-matching $ \sigma_A'$ and a $\gamma_B$-skew-matching $\sigma_B'$ with $\sigma_A'+\sigma_B'\trianglelefteq \mu_c$ such that
\begin{equation}
	W(\sigma_A')=a_1+a_2-W(\tilde \sigma_A)
	\label{equation:skewmatchingcover2-weightsigmaAprime}
\end{equation}
and
\begin{equation}
	W(\sigma_B')\ge \deg_{w'}(c, \mu_c)-W(\sigma_A'),
	\label{equation:skewmatchingcover2-weightsigmaBprime}
\end{equation} the anchor $\mathcal A(\sigma_B')$ fits in the $w'$-neighbourhood of~$c$, and $\mathcal A(\sigma_A')$ is contained in $N_{w'}(d)$.
 
We claim that $(\sigma_A', \sigma_B')$ forms a $(\gamma_A, \gamma_B)$-skew pair in $H^\leftrightarrow$ anchored in $\oriented{dc}$, with respect to $w'$.
Indeed, \ref{itdef:disjoint} follows from $\sigma_A'+\sigma_B'\trianglelefteq \mu_c$;  we have that $\mathcal A(\sigma_B')$ fits in the $w'$-neighbourhood of~$c$, implying \ref{itdef:anchord}.
To see \ref{itdef:anchorc}, we need to check that $\sigma'_A$ is a $\gamma_A$-skew-matching whose anchor fits in the $w'$-neighbourhood of $d$.
This follows because of $\sigma_A'+\sigma_B'\trianglelefteq \mu_c$ and the fact that $\mu_c(x) \leq \mu_d(x) \leq w'(\oriented{xd})$ holds for all $x \in S$, by the construction of $\mu_c, \mu_d$.
Finally, to see \ref{itdef:anchorpartition} we can argue similarly using $\sigma_A'+\sigma_B'\trianglelefteq \mu_c$. \medskip

\noindent \emph{Case 2, Step II: Completing the skew-matching pair in $\bar\mu$.}
By \Cref{lem:new-combination} with 
\begin{center}
\begin{tabular}{c|c|c|c|c}
	object & $(H,\bar w)$ & $c$ & $\bar\mu$ & $\gamma_B$ \\
		\hline
	in place of & $(H,w)$ & $v$ & $\mu$ & $\gamma$ 
\end{tabular}
\end{center}
there is a $\gamma_B$-skew-matching $\sigma_B^*\trianglelefteq\bar \mu$ of weight
\begin{equation}
	W(\sigma_B^*)\ge \deg_{\bar w}(c, \bar \mu)
	\label{equation:skewmatchingcover2-weightsigmaBstar}
\end{equation}
with its anchor $\mathcal A(\sigma_B^*)$ fitting in the $\bar w$-neighbourhood of~$c$. \medskip

\noindent \emph{Case 2, Step III: Gluing the skew-matching pairs together.}
Recall that $\sigma_\emptyset$ is the zero skew-matching.
Together with~\eqref{equation:skewmatchingcover2-tildesigmatriangle} and~\eqref{equation:skewmatchingcover2-anchortildesigmafits}, we can use \Cref{prop:adding-skew-matchings} with
\begin{center}	
	\begin{tabular}{c|c|c|c|c|c|c|c|c|c}
		object & $(H,w')$ & $(d,c)$ & $\mu_c$ & $\mu'$ & $(H, \tilde w)$ & $\sigma'_A$ & $\sigma'_B$ & $\tilde \sigma_A$ & $\sigma_\emptyset$ \\
		\hline
		in place of & $(H,w)$ & $(u,v)$ & $\mu$ & $\bar \mu$ & $(H, \bar w)$ & $\sigma_A$ & $\sigma_B$ & $\bar \sigma_A$ & $\bar \sigma_B$
	\end{tabular}
\end{center}
and obtain that the pair $(\sigma_A'+\tilde \sigma_A, \sigma_B')$ forms a $(\gamma_A, \gamma_B)$-skew-matching in $H$ anchored in $\oriented{dc}$ (w.r.t. $w'$) with $\sigma_A'+\tilde \sigma_A+\sigma_B'\trianglelefteq\mu_d$.

We apply \Cref{prop:adding-skew-matchings} once again, this time with
\begin{center}	
	\begin{tabular}{c|c|c|c|c|c|c|c|c|c}
		object & $(H,w')$ & $(d,c)$ & $\mu_d$ & $\bar \mu$ & $(H, \bar w)$ & $\sigma_A'+\tilde \sigma_A$ & $\sigma_B'$ & $\sigma_\emptyset$ & $\sigma_B^\ast$ \\
		\hline
		in place of & $(H,w)$ & $(u,v)$ & $\mu$ & $\bar \mu$ & $(H, \bar w)$ & $\sigma_A$ & $\sigma_B$ & $\bar \sigma_A$ & $\bar \sigma_B$
	\end{tabular}
\end{center}
and by doing so we obtain that $(\sigma_A'+\tilde \sigma_A, \sigma_B'+\sigma_B^*)$ is a $(\gamma_A, \gamma_B)$-skew pair in $H^\leftrightarrow$ anchored in $\oriented{dc}$ (w.r.t. $w'$) with $\sigma_A'+\tilde \sigma_A+\sigma_B'+\sigma_B^*\le \tilde \mu$.
Note that we have 
\begin{equation}
W(\sigma_A'+\tilde \sigma_A) \overset{\eqref{equation:skewmatchingcover2-weightsigmaAprime}}{=} a_1+a_2-W(\tilde \sigma_A)+W(\tilde \sigma_A)=a_1+a_2.\label{eq:W(sigmaA'+tidlesigmaA)}
\end{equation}
Recalling that $( H, \tilde w)$ is the $\mu_c$-truncated graph obtained from $(H,w')$, and that $(H, \bar w)$ is the $\mu_d$-truncated weighted graph obtained from $(H, w')$; we also have
\begin{align*}
W(\sigma_B'+\sigma_B^*)
&\overset{\eqref{equation:skewmatchingcover2-weightsigmaBprime}, \eqref{equation:skewmatchingcover2-weightsigmaBstar}}{\ge} \deg_{\bar{w}}(c, \bar \mu)+\deg_{w'}(c, \mu_c)-W(\sigma_A')\\
&\overset{\eqref{equation:skewmatchingcover2-weightsigmaAprime}}{=}\deg_{\bar w}(c, \tilde \mu-\mu_d)+\deg_{w'}(c,\mu_c)-(a_1+a_2)+W(\tilde \sigma_A)\\
&\overset{\eqref{eq:weightTidleSigmaA}}{\ge}\deg_{\bar w}(c, \tilde \mu-\mu_d)+\deg_{w'}(c,\mu_c)
+\deg_{\tilde w}(c, \mu_d-\mu_c)-(a_1+a_2)
\end{align*}
Now we apply \Cref{prop:adding-degree-truncated} twice, to deduce
\begin{align}
\nonumber
W(\sigma_B'+\sigma_B^*)
& \ge \deg_{\bar w}(c, \tilde \mu-\mu_d)+\deg_{w'}(c, \mu_d)-(a_1+a_2). \\
\nonumber&\ge  \deg_{w'}(c, \tilde \mu)-(a_1+a_2)
\overset{\eqref{eq:takingouttildesigma}}{\ge} \deg_w(c, \tilde \mu+\tilde\sigma)-W(\tilde \sigma)-(a_1+a_2)\\
&\geq k-W(\tilde \sigma)-(a_1+a_2),
\label{eq:W(sigmaB'+simgaB*)}
\end{align}
where the last inequality holds by our assumption in the statement of the claim.

We can finally define our final skew-matchings that will allow us to conclude.
Let  $\sigma_B:= \tilde \sigma+\sigma_B'+\sigma_B^*$ and $\sigma_A=\sigma_A'+\tilde \sigma_A$.
By~\eqref{eq:W(sigmaA'+tidlesigmaA)} and~\eqref{eq:W(sigmaB'+simgaB*)}, we have $W(\sigma_A) \geq a_1 + a_2$ and $W(\sigma_B)\ge k-W(\tilde \sigma)-(a_1+a_2)+W(\tilde \sigma)= b_1+b_2$.
We claim that $(\sigma_A, \sigma_B)$ is a $(\gamma_A, \gamma_B)$-skew-matching anchored in $\oriented{dc}$; this gives the existence of the desired good matching in $H^{\leftrightarrow}$.

We check the required properties \ref{itdef:disjoint}--\ref{itdef:anchorpartition}.
To see \ref{itdef:disjoint}, note that since $\tilde \mu$ and $\tilde \sigma$ are disjoint and $\sigma_A+\sigma_B'+\sigma_B^*\trianglelefteq \tilde\mu$, we have that $\sigma_A$ and $\sigma_B$ are disjoint.
As $(\sigma_A'+\tilde\sigma_A, \sigma_B'+\sigma_B^*)$ is a $(\gamma_A, \gamma_B)$-skew pair anchored in $\oriented{dc}$ w.r.t. $w'$, the anchor of $\sigma_A = \sigma_A'+\tilde \sigma_A$ fits in the $w$-neighbourhood of $d$, so \ref{itdef:anchorc} holds.

Now we check that \ref{itdef:anchord} holds.
We have already checked that the anchor of $\sigma_B'+\sigma^\ast_B$ fits in the $w'$-neighbourhood of $c$.
This means that we only need to check that for any $x \in V(H)$ with $\tilde \sigma^1(x)>0$, we have that 
\begin{equation}
	w(\oriented{cx})\ge \sigma_B^1(x).
	\label{equation:prop59-sigmaBwfitsc}
\end{equation}
To see this, let $x$ be a vertex with $\tilde \sigma^1(x)>0$.
By \ref{item:strucprop-GEpairanchorS}-\ref{item:strucprop-GEpairanchorfits}, we have that $x \in N_w(c) \cap S_{\mathcal R}$.
Suppose first that $w'(\oriented{cx})= 0$.
In this case, $\mu_c(x)=0$ and thus $\sigma'_B(x)=0$.
Also we have then that $\bar w(\oriented{cx})=0$, and therefore $\sigma_B^{*1}(x)=0$.
Hence, we have that $\sigma^1_B(x) = \tilde \sigma^1_B(x)$.
By \ref{item:strucprop-GEpairanchorfits}, we have that $\tilde \sigma^1(x) \le w(\oriented{cx})$, and therefore we obtain
\begin{align*}
	w(\oriented{cx}) \ge \tilde \sigma(x)=\tilde \sigma^1(x)+\sigma_B^{*1}(x)+\sigma_B'^1(x)=\sigma_B^1(x),
\end{align*}
so \eqref{equation:prop59-sigmaBwfitsc} holds in this case.
Hence, we can assume that $w'(\oriented{cx})>0$.
In this case (by the definition of $w'$) we have $w'(\oriented{cx}) = w(\oriented{cx}) - \tilde \sigma(x)$.
Hence, using that $\sigma_B'+\sigma^\ast_B$ fits in the $w'$-neighbourhood of $c$, we get
\begin{equation*}
	\sigma_B^1(x) = \tilde \sigma^1(x)+\sigma_B'^1(x)+\sigma_B^{*1}(x)
	\leq \tilde \sigma(x)+w'(\oriented{cx})= w(\oriented{cx}).
\end{equation*}
Hence \eqref{equation:prop59-sigmaBwfitsc} holds in all cases, and therefore, \ref{itdef:anchord} holds.

It remains to check that \ref{itdef:anchorpartition} holds.
Let $x \in N_w(c) \cap N_w(d)$ be arbitrary.
Recall that we know already that $(\sigma_A'+\tilde \sigma_A, \sigma_B'+\sigma_B^*)$ is a $(\gamma_A, \gamma_B)$-skew pair with respect to $w'$.
This implies that
\[ \sigma_A'^1(x)+\tilde \sigma^1_A(x)+ \sigma_B'^1(x)+\sigma_B^{*1}(x) \leq \max\{ w'(\oriented{cx}), w'(\oriented{dx})\}. \]
Hence, we have
\begin{align*}
	\sigma_A^1(x)+\sigma_B^1(x)
	& = \sigma_A'^1(x)+\tilde \sigma^1_A(x)+ \sigma_B'^1(x)+\sigma_B^{*1}(x) + \tilde \sigma^1(x)  \\
	& \leq \max\{ w'(\oriented{cx}), w'(\oriented{dx})\} + \tilde \sigma(x) \\
	& = \max\{ w(\oriented{cx})-\tilde \sigma(x), w(\oriented{dx})-\tilde \sigma(x)\} + \tilde \sigma(x) \\
	& = \max\{ w(\oriented{cx}), w(\oriented{dx}) \},
\end{align*}
where in the second-to-last step we have used the definition of $w'$, and the fact that $\tilde \sigma$ fits in the $w$-neighbourhood of $c$ (so $w(\oriented{cx}) \geq \tilde \sigma^1(x) = \tilde \sigma(x)$ holds).
Since $x$ was arbitrary, this implies \ref{itdef:anchorpartition} holds.
This finishes the proof of the claim that $(\sigma_A, \sigma_B)$ is a good matching; which in turn finishes the proof in Case~2, and the proof of \Cref{cl:new-weighted-coveringR}.
\end{proofclaim}

Let \[ \tilde {\mathcal{R}}:=\{d\in V(H)\::\: \tilde \sigma(d)+\tilde \mu(d)<w(\oriented{cd})\}.\]
If $\tilde {\mathcal{R}}$ is empty, then we would have $\deg_w(c, \tilde \sigma + \tilde \mu) = \sum_{x \in V} w(\oriented{cx}) = \deg_w(c) \geq k$, so we would be done by \Cref{cl:new-weighted-coveringR}.
Hence, 
from now on, we can and will assume that
\begin{enumerate}[\upshape (\Alph{propcounter}\arabic*), topsep=0.7em, itemsep=0.5em, resume]
	\item \label{assumption:N(c)notsaturated}
$\tilde {\mathcal{R}}\neq \emptyset$,
\end{enumerate}
Observe that by \ref{item:strucprop-GEpaircovers-w} we have $\tilde {\mathcal R} \subseteq \mathcal{R}$.
Since $(\tilde \sigma, \tilde \mu)$ is an optimal $(H, w, S, M, c, \mu, \gamma_B)$-GE pair by \ref{item:strucprop-GEpairmaxsat}, we can apply the Separating Lemmas (\Cref{lemma:separatinglemma-1} and \Cref{lemma:separatinglemma-2}) to it with any choice of $d \in \tilde{\mathcal{R}}$.
From those applications we directly obtain the following claims.

\begin{claim}[First Separating Claim]\label{cl:new-weighted-alternating-skew-1}
	For any $d\in \tilde {\mathcal{R}}$  we have  $\tilde \sigma(x)= w(\oriented{cx})$  for all $x\in N_w(c)\cap N_w(d)$.
\end{claim}

\begin{claim}[Second Separating Claim]\label{cl:new-weighted-alternating-skew-2}
	Let $d\in \tilde {\mathcal{R}}$, and let  $\sigma_d\le \tilde \sigma$ and $\mu_d\le \tilde \mu$ be the $\gamma_B$-skew-matching and fractional matching, respectively, obtained from $\tilde \sigma$ and $\tilde \mu$, respectively, by considering only the edges of their support that intersect $N_w(d)$.
	Then there is no edge between the sets \[\{x\in N_w(c)\cap S\::\: \tilde \sigma(x)<w(\oriented{cx})\}\]
	and 
	\[\{y\in N_w(c)\setminus S\::\: \sigma_d(y)+\mu_d(y)>0\}.\]
\end{claim}

\subsection{Proof of \Cref{prop:weighted-structural}: The balanced case}
Now we have the tools to conclude in the case where $b_2 \leq k/2$.
We discuss a bit about our strategy here.
This case is still somewhat manageable, in the sense that the $\gamma_B$-skew-matching would fit in the fractional matching built on top of the neighbourhood of $d$, for some $d\in \mathcal R$. However, we do not have the condition that $\min\{b_1,b_2\}+\max\{a_1,a_2\}\le k/2$ in order to use directly the Completion Lemma (\Cref{lem:new-completion}).
To overcome this problem and to ensure that we can actually embed in the mentioned fractional matching over $N_w(d)$ at least as much as the degree of $c$ into it, we shall leverage the fact we can fill things more efficiently by using the skew-matching $\tilde \sigma$.

\begin{claim}\label{cl:new-weighted-smallb_2}
	Suppose that  $b_2\le k/2$. Then 
	$H^\leftrightarrow$ admits a good matching. 
\end{claim}

\begin{proofclaim}[Proof of \Cref{cl:new-weighted-smallb_2}] 
	The proof has three steps.
	We begin by picking any $d\in \tilde {\mathcal{R}}$ (this exists, by \ref{assumption:N(c)notsaturated}).
	\medskip
	
	\noindent \emph{Step 1: Building the $\gamma_B$-skew-matching $\sigma_B$.}
	Let $\sigma_B'\le \tilde \sigma$ be a maximal $\gamma_B$-skew-matching such that for every $x\in N_w(d)\cap N_w(c)$ we have $\sigma_B'(x)=\min\{\tilde \sigma(x),  w(\oriented{dx})\}$.
	We recall that by \ref{item:strucprop-GEpairanchorfits}, we have $\sigma'_B(x) \leq \tilde \sigma(x) = 0$ for $x \notin N_w(c)$.
	Hence, by construction, the anchor $\mathcal A(\sigma_B')$ fits in the $ w$-neighbourhood of~$d$.
	
	Let $w'$ be the weight function defined as
	\[
	w'(\oriented{uv})=\begin{cases}
		w(\oriented{dx})-\sigma'_B(x) & \text{if } \{u,v\}=\{d,x\}\text{ with }x\in S,\\
		\max\{0,\,w(\oriented{cx})-\sigma'_B(x)\} & \text{if } \{u,v\}=\{c,x\}\text{ with }x\in V(H)\setminus\{c\},\\
		w(\oriented{uv}) & \text{otherwise.}
	\end{cases}
	\] 
	The definition of $w'$ and $\sigma'_B$, and a quick case analysis, shows that for every $x \in S$ we have $w'(\oriented{dx}) = w(\oriented{dx}) - \sigma'_B(x) \leq 1 - \tilde \sigma(x)$.
	For every $x \in S$ we have $1 - \tilde \sigma(x) = \tilde \mu (x)$, by \ref{item:strucprop-GEpaircovers}, which combined with the previous bound gives, for all $x \in S$,
	\begin{equation}
		w'(\oriented{dx}) \leq \tilde \mu(x).
		\label{equation:strucprop-balancedcase-wprimeversusmu}
	\end{equation}
	
	Now, let $\tilde \mu_d\le \tilde  \mu$  be a maximal fractional matching whose support edges intersect $N_{w'}(d)$, and such that for all $x\in N_{w'}(d)\subseteq S$, we have $\tilde \mu_d(x)\le  w'(\oriented{dx})$ and 
	\begin{align}\label{eqdef:tildemu_d}
	W(\tilde \mu_d) & \le \gamma_B\left(b_1-\frac{W(\sigma_B')}{1+\gamma_B}\right).
	\end{align}
	By assumption, we have $b_2 \leq k/2$, and therefore
	\begin{align*}
	\gamma_B\left(b_1-\frac{W(\sigma_B')}{1+\gamma_B}\right)
	& = b_2 -\frac{\gamma_B}{1+\gamma_B}W(\sigma_B')
	\leq \frac{k}{2} -\frac{\gamma_B}{1+\gamma_B}W(\sigma_B')
	 \overset{\text{\ref{eqnew:gamma_B>1}}}{\le} \frac{k}{2}-\frac{W(\sigma_B')}{1+\gamma_B} \\
	& \leq \sum_{x \in S} \left( w(\oriented{dx}) - \sigma'_B(x) \right) = \sum_{x \in S} w'(\oriented{dx}) \overset{\eqref{equation:strucprop-balancedcase-wprimeversusmu}}{\le} \deg_{w'}(d, \tilde \mu).
	\end{align*}
	Hence, by construction, in fact we attain equality in~\eqref{eqdef:tildemu_d}.
	We thus have
	\begin{align}
		W(\tilde \mu_d) = \gamma_B\left(b_1-\frac{W(\sigma_B')}{1+\gamma_B}\right).
		\label{eqdef:tildemu_dequal}
	\end{align}
	
	We apply \Cref{lem:new-extending-out} (Extending-out), with
	\begin{center}	
		\begin{tabular}{c|c|c|c|c}
			object & $H$ & $N_{w'}(d)$ & $V(H) \setminus S$ & $\tilde \mu_d$ \\
			\hline
			in place of & $H$ & $U$ & $V$ & $\mu$ 
		\end{tabular}
	\end{center}
	The outcome of this application is an $\gamma_B$-skew-matching $\tilde \sigma_B$ in $H^\leftrightarrow$ with
	\begin{equation}
		\tilde \sigma_B \trianglelefteq\tilde \mu_d,
		\label{equation:strucprop-balanced-stildebdom}
	\end{equation}
	and of weight
	\begin{align}
	W(\tilde \sigma_B) & =(1+\min\{\gamma_B, \gamma_B^{-1}\})W(\tilde \mu_d) \overset{\text{\ref{eqnew:gamma_B>1}}}{=} (1+\gamma_B^{-1})W(\tilde \mu_d) \nonumber \\
	& \overset{\eqref{eqdef:tildemu_dequal}}{=}b_1+b_2-W(\sigma_B'), \label{eq:weighttildesigmaB}
	\end{align}
	and whose anchor satisfies $\mathcal A(\tilde \sigma_B)\subseteq N_{w'}(d)\subseteq S$.
	By the way $w'$ and $\tilde \mu_d$ were defined and \eqref{equation:strucprop-balanced-stildebdom} we have that
	\begin{equation}\text{the anchor $\mathcal A(\tilde \sigma_B)$ fits in the $w'$-neighbourhood of $d$.}
		\label{equation:c62-sigmaBtildeanchor}
	\end{equation}
	
	We finish the construction of $\sigma_B$ by setting
	\begin{equation}
		\sigma_B:= \sigma_B'+\tilde \sigma_B.
		\label{equation:c62-defsigmaBtilda}
	\end{equation}
	By \eqref{eq:weighttildesigmaB}, we have 
	\begin{align}
		W(\sigma_B) = W(\sigma_B')+W(\tilde \sigma_B) = b_1+b_2.
		\label{equation:c62-weightsigmab}
	\end{align}
	
	\noindent \emph{Step 2: Starting to build $\sigma_A$.}  In this step, we try to build a $\gamma_A$-skew-matching, which is as large as possible, and is contained in the leftover of $\tilde \mu$ and in the leftover of $\tilde \sigma$.
	
	By construction, the edges supporting $\tilde \mu_d$ intersect $N_{w'}(d)$. By \Cref{cl:new-weighted-alternating-skew-1}
 we have that $\tilde \sigma(x)=w(\oriented{cx})$ for all $x\in N_w(c) \cap N_{w'}(d)$,
 and by \ref{item:strucprop-GEpairanchorfits} we have $\tilde \sigma(x) = w(\oriented{cx}) = 0$ for all $x \notin N_w(c)$.
 Then (in particular) the equality $\tilde \sigma(x)=w(\oriented{cx})$ holds for all $x \in N_{w'}(d) \cap V(\tilde \mu_d)$. 
	For $x\in S\cap N_w(c)$, by the definition of $\tilde \mu_d$ we have that if $\tilde \mu_d(x)>0$, then  $w'(\oriented{dx}) > 0$.
	This, by definition of $w'$, means that $w(\oriented{dx}) > \sigma'_B(x)$ for such $x$.
	By the definition of $\sigma'_B$, then we have $\sigma_B'(x)=\tilde\sigma(x)=w(\oriented{cx})$, which then means that $w'(\oriented{cx})=0$.
	Hence, we have that
	\begin{equation}
		\deg_{w'}(c, V(\tilde \mu_d)\cap S)=0.
		\label{equation:structprop-balanced-degwprimecmud}
	\end{equation}
	
	Now, let $y \in V(\tilde \mu_d) \setminus S$.
	For such $y$, there exists $v \in S \cap N_w(d)$ such that $0 < \tilde \mu_d(vy)$. Because $\tilde \mu_d \leq \tilde \mu$, by \ref{item:strucprop-GEpairsupport} $\tilde \mu_d$ is supported in the edges between $S_\mathcal R$ and $\mathcal R$, and thus 
 $y\in \mathcal R$ and thus \ref{item:structuprop-GEpairsaturation} implies $\tilde \mu_d(y)+\tilde \sigma(y)\le \tilde \mu(y)+\tilde \sigma(y)\le  w(\oriented{cy})$, implying $\tilde \mu_d(y)\le w(\oriented{cy})-\tilde \sigma(y)\le w'(\oriented{cy})$. Thus we have
	 
	\begin{align}
	\deg_{w'}(c, \tilde \mu_d)
	& = \sum_{x \in N_{w'}(c)} \min\{ w'(\oriented{cx}), \tilde \mu_d(x) \} \overset{\eqref{equation:structprop-balanced-degwprimecmud}}{=} \sum_{x \in N_{w'}(c) \setminus S} \min\{ w'(\oriented{cx}), \tilde \mu_d(x) \} \nonumber \\
	& = \sum_{x \in \mathcal R} \tilde \mu_d(x) {\leq} W(\tilde \mu_d)
	\overset{\eqref{eq:weighttildesigmaB}}{\le} W(\tilde \sigma_B),
	\label{eq:saturationtildemu_d}
	\end{align}
using that  $\tilde \mu_d$ is supported in the edges between $S_\mathcal R$ and $\mathcal R$.

	Let $( H, \bar w)$ be the $\tilde \mu_d$-truncated weighted graph obtained from $(H, w')$.
	We want to apply \Cref{prop:weighted-extending-skew} (Extending-out skew-matching) with
	\begin{center}	
		\begin{tabular}{c|c|c|c|c|c}
			object & $(H, \bar w)$ & $1$ & $\gamma_B$ & $c$ & $\tilde \sigma - \sigma'_B$  \\
			\hline
			in place of & $(H,w)$ & $\gamma_A$ & $\gamma_B$ & $u$ & $\sigma_B$
		\end{tabular}
	\end{center}
	To do so, we recall that $\gamma_B > 1$, by \ref{eqnew:gamma_B>1}.
	We also need to check that the anchor of $\tilde \sigma - \sigma'_B$ fits in the $\bar w$-neighbourhood of $c$.
	Indeed, we have $\mathcal{A}(\tilde \sigma - \sigma'_B) \subseteq S_{\mathcal R}$ by \ref{item:strucprop-GEpairanchorS} and $\sigma'_B \leq \tilde \sigma$, by definition.
	For any $x \in S$, if $\tilde \sigma(x) \leq w(\oriented{dx})$, we have $\sigma'_B(x) = \tilde \sigma(x)$, so there is nothing to check.
	So we can assume $\tilde \sigma(x) > w(\oriented{dx})$, so $\sigma'_B(x) = w(\oriented{dx})$.
	Then we have
	\begin{align*}
		\tilde \sigma^1(x) - (\sigma'_B)^1(x)
		& = \tilde \sigma(x) - \sigma'_B(x) = \tilde \sigma(x) - w(\oriented{dx})  \leq w(\oriented{ cx }) - w(\oriented{dx}) \\
		& = w(\oriented{ cx }) - \sigma'_B(x) + \sigma'_B(x) - w(\oriented{dx}) \\
		& \leq w'(\oriented{cx}) + \sigma'_B(x) - w(\oriented{dx})
		= w'(\oriented{ cx }) - w'(\oriented{dx})  \\
		& \leq w'(\oriented{ cx }) - \tilde \mu_d(x) = \bar w ( \oriented{cx}),
	\end{align*}
	as required.
	The application of \Cref{prop:weighted-extending-skew} then gives a $1$-skew-matching $\sigma_c$ of weight $W(\sigma_c) = 2W(\tilde \sigma - \sigma'_B)/(1 + \gamma_B)$ whose anchor fits in the $\bar w$-neighbourhood of $c$, and $\sigma_c \leq \tilde \sigma - \sigma'_B$.
	By \Cref{lemma:fractionalfrom1skew} with 
	
	\begin{center}	
		\begin{tabular}{c|c|c|c|c}
			object & $H$ & $\sigma_c$ & $\tilde \sigma-\sigma_B'$ & $\gamma_B$   \\
			\hline
			in place of & $G$ & $\sigma$ & $\sigma'$ & $\gamma$ 
		\end{tabular}
	\end{center}
	we get a fractional matching $\mu_c$ such that
	\begin{equation} \label{equation:strucprop-balanced-mucdom}
		\mu_c\preceq\tilde \sigma-\sigma_B',
	\end{equation}
	and
	\begin{equation*}
		W(\mu_c) = \frac{W(\tilde \sigma - \sigma'_B)}{1 + \gamma_B}.
	\end{equation*}
	
	The last equality and \eqref{equation:strucprop-balanced-mucdom} imply that, in fact, for each $x \in S$, we have
	\begin{equation*}
		\mu_c(x) = \tilde \sigma(x) - \sigma'_B(x),
	\end{equation*}
	and therefore, for each $x \in S$,
	\begin{equation*}
		\mu_c(x) + \tilde \mu(x) + \sigma'_B(x) = \tilde \sigma(x) + \tilde \mu(x) = 1 \geq w(\oriented{cx}).
	\end{equation*}
	Observe that by \ref{item:strucprop-GEskewinR} (and the fact that $\mu_c \preceq \tilde \sigma - \sigma'_B$ and $\sigma'_B \leq \tilde \sigma$), we have that $\mu_c(x)=\sigma_B'(x)=0$ for all $x\in V(H)\setminus (S\cup \mathcal R)$. 
	Together with \ref{item:strucprop-GEpaircovers-w}, we get that for each $x \in V(H)\setminus (S_{\mathcal R} \cup \mathcal R)$, we have $\mu_c(x) + \tilde \mu(x) + \sigma'_B(x) \ge w(\oriented{cx})$.
	Together with the last displayed inequality, we in fact have deduced that, for all $x \in V(H) \setminus \mathcal{R}$, we have
	\begin{equation}
		\mu_c(x) + \tilde \mu(x) + \sigma'_B(x) \ge w(\oriented{cx}).
		\label{equation:strucprop-balanced-mucx-pointwise}
	\end{equation}

	As $\mu_c \preceq \tilde \sigma-\sigma_B'$, using \ref{item:structuprop-GEpairsaturation} we can deduce that for each $x \in \mathcal{R}$,
	\begin{equation}
		\mu_c(x) + \tilde \mu(x) + \sigma'_B(x) \le \tilde \sigma(x) + \tilde \mu(x) \leq w(\oriented{cx}).
		\label{equation:strucprop-balanced-mucx-pointwise-R}
	\end{equation}
	
	By $\tilde \mu_d \leq \tilde \mu$, \eqref{equation:strucprop-balanced-mucdom} and \ref{item:strucprop-GEpairdisjoint}, we have that $\mu_c$ and $\tilde \mu - \tilde \mu_d$ are disjoint and thus $\mu_c + \tilde \mu - \tilde \mu_d$ is a fractional matching.
	We apply \Cref{lem:new-combination} (Combination)
	 \begin{center}	
 	\begin{tabular}{c|c|c|c|c}
 		object & $(H, \bar w)$ & $c$ & $\mu_c + \tilde \mu - \tilde \mu_d$ & $\gamma_A$ \\
 		\hline
 		in place of & $(H,w)$ & $v$ & $\mu$ & $\gamma$
 	\end{tabular}
 \end{center}
	to obtain a $\gamma_A$-skew-matching $\sigma_A'$ such that
	\begin{equation}
		\sigma'_A \trianglelefteq \mu_c+ \tilde \mu-\tilde \mu_d,
		\label{equation:strucprop-balanced-sprimeadom}
	\end{equation}
	with weight
	\begin{equation}
		W(\sigma'_A) \ge \deg_{\bar w}(c, \mu_c+ \tilde \mu-\tilde \mu_d),
		\label{equation:strucprop-balanced-wsigmaprimea}
	\end{equation}
	and \begin{equation}
		\text{the anchor $\mathcal A(\sigma_A')$ fits in the $\bar w$-neighbourhood of $c$.}
		\label{equation:strucprop-balanced-wsigmaprimea-anchor}
	\end{equation}

	Recall that $\sigma_\emptyset$ is the identically-zero matching.
	Owing to \eqref{equation:strucprop-balanced-stildebdom}, \eqref{equation:c62-sigmaBtildeanchor}, \eqref{equation:strucprop-balanced-sprimeadom}, \eqref{equation:strucprop-balanced-wsigmaprimea-anchor}, and the fact that $\bar w$ is the $\tilde \mu_d$-truncated weight obtained from $w'$, we can apply \Cref{prop:adding-skew-matchings} with
	\begin{center}	
	\begin{tabular}{c|c|c|c|c|c|c|c|c|c}
		object & $(H,w')$ & $\fmat{cd}$ & $\tilde \mu_d$ & $\mu_c + \tilde \mu - \tilde \mu_d$ & $\bar w$ & $\sigma_\emptyset$ & $\tilde \sigma_B$ & $\sigma'_A$ & $\sigma_\emptyset$ \\
			\hline
		in place of & $(G,w)$ & $\fmat{uv}$ & $\mu$ & $\bar \mu$ & $\bar w$ & $\sigma_A$ & $\sigma_B$ & $\bar \sigma_A$ & $\bar \sigma_B$
	\end{tabular}
	\end{center}
	to obtain that $(\sigma_A', \tilde \sigma_B)$ is a $(\gamma_A, \gamma_B)$-skew-matching pair in $H^\leftrightarrow$, anchored in $\oriented{cd}$ (with respect to $w'$) with $\sigma_A'+\tilde \sigma_B\trianglelefteq\mu_c+ \tilde \mu$.
	
	We now give further estimates on $W(\sigma'_A)$.
	We have
 \begin{align}
\nonumber 
 W(\sigma_A')& \overset{\eqref{equation:strucprop-balanced-wsigmaprimea}}{\geq} \deg_{\bar w}(c, \mu_c+ \tilde \mu-\tilde \mu_d)
  = \sum_{x \in N_{\bar w}(c)} \min \{ \bar w(\oriented{cx}), \mu_c(x)+ \tilde \mu(x)-\tilde \mu_d(x) \} \\
 \nonumber & = \sum_{x \in N_{w'}(c)} \min \{ \bar w(\oriented{cx}), \mu_c(x)+ \tilde \mu(x)-\tilde \mu_d(x) \} \\
\nonumber & \geq \sum_{x \in N_{w'}(c)} (\min \{ w'(\oriented{cx}), \mu_c(x)+ \tilde \mu(x)\}-\tilde \mu_d(x) ) \\
\nonumber & = \deg_{w'}(c,\mu_c+\tilde \mu) - \sum_{x \in N_{w'}(c)} \tilde \mu_d(x) 
 \overset{\eqref{equation:structprop-balanced-degwprimecmud}}{\ge} \deg_{w'}(c,\mu_c+\tilde \mu)-W(\tilde \mu_d)\\
& \overset{\eqref{eq:saturationtildemu_d}}{ \ge} \deg_{w'}(c,\mu_c+\tilde \mu)-W(\tilde \sigma_B),
 \end{align}
 here in the second line we used that $N_{\overline{w}}(c) \subseteq N_{w'}(c)$, and then that $\overline{w}(\oriented{cx}) \geq w'(\oriented{cx}) - \tilde{\mu}_d(x)$, both things holding by definition of $\overline{w}$.
	
Now we want to show that $(\sigma_A', \sigma_B)$ is a $(\gamma_A, \gamma_B)$-skew-matching pair in $H^\leftrightarrow$ anchored in $\oriented{cd}$ (with respect to $w$). 
For this, we need to verify properties \ref{itdef:disjoint}--\ref{itdef:anchorpartition}.
By \eqref{equation:strucprop-balanced-sprimeadom}, \eqref{equation:strucprop-balanced-mucdom},
\eqref{equation:strucprop-balanced-stildebdom}, we can appeal to the definition of $\mu_c$ and $\sigma_B'$ to deduce that $\sigma_A'$ is disjoint from $\sigma_B'+\tilde \sigma_B=\sigma_B$, which gives \ref{itdef:disjoint}.
As $(\sigma_A', \tilde \sigma_B)$ is a $(\gamma_A, \gamma_B)$-skew pair in $H^\leftrightarrow$ anchored in $\oriented{cd}$ (with respect to $w'$), we automatically have that the anchor $\mathcal A(\sigma_A')$ fits in the  $w$-neighbourhood of $c$, which gives \ref{itdef:anchorc}.
To see~\ref{itdef:anchord}, we use $\sigma_B'+\tilde \sigma_B=\sigma_B$ together with \eqref{equation:c62-sigmaBtildeanchor} to get $\sigma^1_B(x) = (\sigma'_B)^1(x) + \tilde \sigma^1_B(x) \leq (\sigma'_B)^1(x) + w'(\oriented{dx})$ for all $x \in V$.
From the definition of $w'$ we get that the last term is at most $w(\oriented{dx})$, which gives~\ref{itdef:anchord}. 

It remains to check \ref{itdef:anchorpartition}.
Let $x\in N_w(d)\cap N_w(c)$ be arbitrary.
We note first that
\begin{align*}
	\sigma_A'^1(x)+\tilde \sigma_B^1(x)\le \bar w(\oriented{cx})+\tilde \mu_d(x) \leq \max\{ w'(\oriented{cx}), \tilde \mu_d(x) \} \leq \max\{ w'(\oriented{cx}), w'(\oriented{dx}) \},
\end{align*}
where in the first inequality we used that the anchor of $\sigma'_A$ fits in the $\bar w$-neighbourhood of $c$ together with  \eqref{equation:strucprop-balanced-stildebdom}, 
in the second we used that $\bar w$ is the $\bar \mu_d$-truncated weight obtained from $w'$, and finally, we used the definition of $\tilde \mu_d$.
Using this, we get
\begin{align*}
	\sigma_A'^1(x)+\sigma_B^1(x) \leq \max\{ w'(\oriented{cx}), w'(\oriented{dx}) \} + \sigma'_B(x) \leq \max\{w(\oriented{cx}), w(\oriented{dx}), \sigma'_B(x)\},
\end{align*}
where we used $\sigma_B'+\tilde \sigma_B=\sigma_B$ and the definition of $w'$.
The last term is at most $\max\{w(\oriented{cx}), w(\oriented{dx})\}$, because $\sigma'_B(x) \leq \tilde \sigma(x) \leq w(\oriented{cx})$ by \ref{item:strucprop-GEpairanchorS}--\ref{item:strucprop-GEpairanchorfits}.
This implies \ref{itdef:anchorpartition}. \medskip
%

\noindent \emph{Step 3: Complementing $\sigma_A$.}
Summarising our efforts so far, we have obtained that $(\sigma'_A, \sigma_B)$ is a $(\gamma_A, \gamma_B)$-skew-matching anchored in $\oriented{cd}$, with respect to $w$.
By \eqref{equation:c62-weightsigmab}, we only need to focus on $\sigma'_A$.
In this step, we use a greedy argument to complement the already-obtained $\gamma_A$-skew-matching $\sigma'_A$.
By doing so, we will obtain a $\gamma_A$-skew-matching $\widehat{\sigma_A}$, so we can conclude our construction by setting $\sigma_A = \sigma'_A + \widehat{\sigma_A}$.

We proceed as follows.
If $W(\sigma_A')\ge a_1+a_2$, then we are actually done (since then $(\sigma'_A, \sigma_B)$ is the desired good matching), so we can finish by setting $\widehat\sigma_A\equiv 0$.
Hence, we can assume that $W(\sigma_A') < a_1+a_2$.

We also note that from \eqref{equation:strucprop-balanced-stildebdom}, \eqref{equation:c62-defsigmaBtilda}, and \eqref{equation:strucprop-balanced-sprimeadom}, we get
\begin{equation}
	\sigma'_A + \sigma_B -\sigma_B'\trianglelefteq \tilde \mu  + \mu_c.
	\label{equation:strucprop-balanced-sprimeAsigmaBdom}
\end{equation}

Let
\[ \kappa' := a_1 - \frac{W(\sigma'_A)}{1 + \gamma_A}. \]
We wish to apply \Cref{prop:weighted-greedy-1} (First Greedy Lemma), with
	\begin{center}	
	\begin{tabular}{c|c|c|c|c|c|c|c|c}
		object & $(H,w)$ & $\sigma'_A$ & $\sigma_B$ & $c$ & $d$ & $\mathcal{R}$ & $S$ & $\kappa'$ \\
		\hline
		in place of & $(H,w)$ & $\sigma_A$ & $\sigma_B$ & $u$ & $v$ & $V$ & $U$ & $\kappa$ 
	\end{tabular}
\end{center}
To do so, we verify the required properties \ref{item:greedy-1-condition-1}--\ref{item:greedy-1-condition-2}.
Indeed, we have (using $\sigma'_B \leq \tilde \sigma$ in the fourth line)
\begin{align*}
	\deg_w(c, \mathcal R)
	& = \deg_w(c)-\deg_w(c, V(H)\setminus \mathcal{R})
	 = \deg_w(c)- \sum_{y \notin \mathcal{R}} w(\oriented{cy}) \\
	& \overset{\eqref{equation:strucprop-balanced-mucx-pointwise}}{=}
	\deg_w(c)- \sum_{y \notin \mathcal{R}} \min \{ w(\oriented{cy}), \tilde \mu(y)+\sigma_B'(y)+\mu_c(y)\} \\
	& =
	\deg_w(c)- \deg_w(c, \tilde \mu+\sigma'_B+\mu_c)+\sum_{y\in \mathcal R}\min\{w(\oriented{cy}), \tilde \mu(y)+\sigma_B'(y)+\mu_c(y)\} \\
	& \geq k - \deg_w(c, \tilde \mu+\tilde \sigma+\mu_c)+\sum_{y\in \mathcal R}\min\{w(\oriented{cy}), \tilde \mu(y)+\sigma_B'(y)+\mu_c(y)\} \\
	& \overset{\eqref{equation:strucprop-balanced-mucx-pointwise-R}}{=} k - \deg_w(c, \tilde \mu+\tilde \sigma+\mu_c)+\sum_{y\in \mathcal R} (\tilde \mu(y)+\sigma_B'(y)+\mu_c(y)) \\
	&\overset{\eqref{equation:strucprop-balanced-wsigmaprimea}}{\geq} (a_1+a_2)-W(\sigma_A')+\sum_{y\in \mathcal R}(\tilde\mu(y)+\sigma_B'(y)+\mu_c(y))\\
	& \overset{\eqref{equation:strucprop-balanced-sprimeAsigmaBdom}}{\geq} a_1-\frac{W(\sigma_A')}{1+\gamma_A}+\sum_{y\in \mathcal R}\left(\sigma_A'(y)+\sigma_B(y)\right),
\end{align*}
which verifies \ref{item:greedy-1-condition-1}.
To see  \ref{item:greedy-1-condition-2}, we need to verify, for any $x\in \mathcal R$, that we have
\begin{equation}
\deg_w(x, S) \geq \sum_{y\in S}\left(\sigma_B(y)+\sigma_A'(y)\right)+\left(a_2-\frac{\gamma_A W(\sigma_A')}{1+\gamma_A}\right).
\label{equation:strucprop-c64-finalverificationgreedy}
\end{equation}
Indeed, let $x \in \mathcal{R}$ be arbitrary.
We have
\begin{align*}
\deg_w(x, S)
& = \deg_w(x)\ge \frac{k}{2} \overset{\eqref{assumption-new:a_1+b_1-small}}{>} b_1+\max\{a_1,a_2\}
\overset{\eqref{equation:c62-weightsigmab}}{=} \frac{1}{1 + \gamma_B} W(\sigma_B) + \max\{a_1,a_2\}\\
& \geq \sum_{y\in S} \sigma_B(y) + \max\{a_1,a_2\} = \sum_{y\in S} \sigma_B(y) + (\max\{a_1,a_2\} - a_2) + a_2.
\end{align*}
We also note that
\begin{align*}
	\sum_{y \in S} \sigma_A'(y)
	& \leq \frac{\gamma_A}{1 + \gamma_A} \sum_{(u,y), y \in S} \sigma'_A(\oriented{uy}) +  \frac{1}{1 + \gamma_A} \sum_{(y,u), y \in S} \sigma'_A(\oriented{yu}) \\
	& \leq \frac{\gamma_A}{1 + \gamma_A} W(\sigma'_A) +  \frac{1 - \gamma_A}{1 + \gamma_A} \sum_{(y,u), y \in S} \sigma'_A(\oriented{yu}) \leq \frac{1}{1 + \gamma_A} W(\sigma'_A),
\end{align*}
so to reach \eqref{equation:strucprop-c64-finalverificationgreedy} it suffices to show that
\begin{equation}
	\frac{1}{1 + \gamma_A} W(\sigma'_A) \leq \frac{\gamma_A}{1 + \gamma_A} W(\sigma'_A) + (\max\{a_1,a_2\} - a_2)
	\label{equation:strucprop-c64-finalverificationgreedy2}
\end{equation}
We do this as follows.
Suppose first that $a_2 \geq a_1$, so $\gamma_A \geq 1$ and $\max\{a_1,a_2\} - a_2 = 0$.
In this case, \eqref{equation:strucprop-c64-finalverificationgreedy2} follows immediately from $\gamma_A \geq 1$.
Thus we can assume that $a_2 < a_1$, so $\gamma_A < 1$, and $\max\{a_1,a_2\} - a_2 = a_1 - a_2$.
Recall that we assume $W(\sigma'_A) < a_1 + a_2$, so we have
\begin{align*}
	\frac{1}{1 + \gamma_A} W(\sigma'_A)
	& \leq \frac{\gamma_A}{1 + \gamma_A} W(\sigma'_A) + \frac{1 - \gamma_A}{1 + \gamma_A} (a_1 + a_2)
	= \frac{\gamma_A}{1 + \gamma_A} W(\sigma'_A) + (a_2 - a_1),
\end{align*}
which again gives \eqref{equation:strucprop-c64-finalverificationgreedy2}.
Thus indeed \eqref{equation:strucprop-c64-finalverificationgreedy} holds, and we can apply \Cref{prop:weighted-greedy-1} with the above-mentioned desired parameters.

The outcome of \Cref{prop:weighted-greedy-1} is a $\gamma_A$-skew-matching $\widehat \sigma_A$ in $H^\leftrightarrow$ of weight $W(\widehat\sigma_A) \geq a_1+a_2-W(\sigma_A')$ such that $(\sigma_A'+\widehat\sigma_A, \sigma_B)$ is a $(\gamma_A, \gamma_B)$-skew-matching anchored in $\oriented{cd}$.
Set $\sigma_A:= \sigma_A'+\widehat\sigma_A$.
We have 
\begin{align*}
W(\sigma_A) = W(\sigma_A')+W(\widehat\sigma_A) \geq W(\sigma_A')+a_1+a_2-W(\sigma_A') = a_1+a_2,
\end{align*}
which, together with \eqref{equation:c62-weightsigmab}, implies we have found our desired good matching.
\end{proofclaim}

In the following, we may assume that 
\begin{equation}\label{assumption:b_2large}
	b_2>k/2,
\end{equation}and thus
\begin{equation}\label{eq:assumption:a_i+b_1<k/2}
	a_1+a_2+b_1< k/2.
\end{equation}

\subsection{Proof of \Cref{prop:weighted-structural}: The large $S_\mathcal{R}$ case}
In this case we will argue that if the neighbourhoods of two elements of $\mathcal R$ do not intersect much, the set $ S_{\mathcal R}$ will be unusually large, allowing us to conclude.

\begin{claim}\label{cl:new-weighted-smallintersection}
	If there are $d_1,d_2\in {\mathcal{R}}$ such that 
	\begin{equation}\label{eq:smallintersection}
	\sum_{u\in N_w(d_2)}  w(\oriented{d_1 u})\le b_1,
	\end{equation}
	then $H^\leftrightarrow$ has a good matching. 
\end{claim}

In the proof of \Cref{cl:new-weighted-smallintersection} we will not use the GE pair $(\tilde \mu, \tilde \sigma)$, but only the $c$-optimal fractional matching $\mu$. 

\begin{proofclaim}[Proof of \Cref{cl:new-weighted-smallintersection}]
	Suppose $d_1, d_2$ are as in the statement.
	The strategy here is to first build $\sigma_A$, which is anchored in $S$, and its anchor $w$-fits in the neighbourhood of $d_1 \in \mathcal{R}$.
	Then we will set up two auxiliary matchings $\mu_A$ and $\mu'$.
	The third step is to use the condition \eqref{eq:smallintersection} to deduce that there is enough space in $S_\mathcal{R}$ to build $\sigma_B$, anchored in $\mathcal{R}$, that also $w$-fits in the neighbourhood of $c$.
	Both matchings will be found by invoking \Cref{lem:new-extending-out} (Extending-out).
	\medskip
	
	\noindent \emph{Step 1: Building $\sigma_A$}.
	Let $d_1, d_2 \in \mathcal{R}$ be such that they satisfy \eqref{eq:smallintersection}.
	By \Cref{ob:R=singletons}, we have $N_w(d_1)\subseteq S_{\mathcal{R}}\subseteq S$.
	By \ref{item:SP-mucovered}, each $x \in N_w(d_i)$ is covered by $\mu$.
	In particular, $1 = \mu(x) \geq w(\oriented{d_i x})$ holds.
	Now, let $\mu_{d_1}\le \mu$	be a maximal fractional matching such that for any $u\in N_w(d_1)$, we have $\mu_{d_1}(u)\le w(\oriented{d_1u})$.
	This implies that for any such $u$, in fact we have $\mu_{d_1}(u) = w(\oriented{d_1 u})$, so it follows that \begin{equation}\label{eq:W(d_1)}
	W(\mu_{d_1})=\deg_w(d_1,\mu_{d_1})=\deg_w(d_1)\ge \frac k2.
	\end{equation}
	
	Now we apply \Cref{lem:new-extending-out} (Extending-out), with
	\begin{center}	
	\begin{tabular}{c|c|c|c|c|c}
		object & $H$ & $N_w(d_1)$ & $V(H) \setminus S$ & $\mu_{d_1}$ & $\gamma_A$ \\
		\hline
		in place of & $H$ & $U$ & $V$ & $\mu$ & $\gamma$
	\end{tabular}
	\end{center}
	 and deduce the existence of a $\gamma_A$-skew-matching $\sigma_A\trianglelefteq \mu_{d_1}$ in $H^\leftrightarrow$ of weight
	 \begin{align*}
	W(\sigma_A) = (1+\min\{\gamma_A, \gamma_A^{-1}\})W(\mu_{d_1}) \geq W(\mu_{d_1}) \overset{\eqref{eq:W(d_1)}}{\ge} \frac k2 \overset{\eqref{eq:assumption:a_i+b_1<k/2}}{\ge} a_1+a_2,
	\end{align*}
	and such that its anchor $\mathcal A(\sigma_A)$ is contained in $N_w(d_1) \subseteq S$.
	By decreasing the weight in the edges of $\sigma_A$ if necessary, we can assume that in fact
	\begin{equation}
		W(\sigma_A) =a_1+a_2.
		\label{eq:degdomud_1}
	\end{equation}
	By the construction of $\mu_{d_1}$  and the fact that $\sigma_A\trianglelefteq \mu_{d_1}$, the anchor $\mathcal A(\sigma_A)$ fits in the $w$-neighbourhood of~$d_1$. \medskip
	
	\noindent \emph{Step 2: Constructing the auxiliary fractional matchings.}
	Now we will set up several auxiliary fractional matchings.
	Define $\mu_A$ as a fractional matching such that $\mu_A\le \mu_{d_1}$ and minimum such that $\sigma_A\trianglelefteq \mu_A$.
	Note that, by the fact that the anchor of $\sigma_A$ is contained in $S$, then the only directed edges with non-zero weight in $\sigma_A$ are of the form $(x,y)$ with $x \in S$ and $y \notin S$.
	Then the minimality of $\mu_A$ implies that for such a pair $xy \in E(H)$,
	\begin{align*}
		\mu_A(\fmat{xy})
		& = \frac{1}{1 + \gamma_A} \max \{ \sigma_A( \oriented{xy}) + \gamma_A \sigma_A(\oriented{yx}), \sigma_A( \oriented{yx}) + \gamma_A \sigma_A(\oriented{xy}) \} \\
		& =  \frac{1}{1 + \gamma_A} \max \{ \sigma_A( \oriented{xy}) , \gamma_A \sigma_A(\oriented{xy}) \} \\
		& = \sigma_A(\oriented{xy}) \frac{\max\{1, \gamma_A\}}{1 + \gamma_A}
		= \sigma_A(\oriented{xy}) \frac{\max\{a_1, a_2\}}{a_1 + a_2},
	\end{align*}
	which implies that
	\[ W(\mu_A) = W(\sigma_A) \frac{\max\{a_1, a_2\}}{a_1 + a_2} = \max\{a_1, a_2\}.\]
	Therefore, we conclude that
	\begin{equation}\label{eq:degd1muA}
		\deg_w(d_1, \mu_A)\le W(\mu_A)=\max\{a_1,a_2\}\le a_1+a_2.
	\end{equation}
	
	Let $(H, w')$ be the $\mu_A$-truncated weighted graph obtained from $(H, w)$.
	Now we will define another fractional matching $\mu'$, by
	\[ \mu'(xy) = \begin{cases}
		\mu(xy) - \mu_A(xy) & xy \in H[\mathcal{R}, S_\mathcal{R}] \\
		0 & \text{otherwise.}
	\end{cases} \]
	In words, $\mu'$ is equal to $\mu - \mu_A$ whenever the corresponding edge touches $\mathcal{R}$, and zero otherwise.
	Note that $\mu' \leq \mu - \mu_A \leq \mu$.
	
	To understand $\mu'$ we recall a few consequences of the definition of reachable vertices (\Cref{def:weighted-R}).
	By \Cref{observation:reachable-wversusmu}, together with the definition of $w'$, we get that
	\begin{equation}\label{eq:wprimecyR}
		w'(\oriented{cy}) \geq w(\oriented{cy}) - \mu_A(y) \geq \mu'(y) \text{ for all $y \in \mathcal{R}$,}
	\end{equation}
	and recall from \Cref{observation:reachable-weightsR-SR} that we have
	$\sum_{x \in S_\mathcal{R}} \mu(x) = \sum_{x \in \mathcal{R}} \mu(x)$.
	
	Using this, we obtain
	\begin{align*}
	\deg_{w'}(c,\mu')
	& \geq \sum_{y\in \mathcal R}\min\{w'(\oriented{cy}), \mu'(y)\}
	\overset{\eqref{eq:wprimecyR}}{\geq} \sum_{y\in \mathcal R} \mu'(y) \\
	& = \sum_{y\in \mathcal R}\left(\mu(y)-\mu_A(y)\right)
	\\
	& \geq \sum_{y\in \mathcal R} \mu(y) - W(\mu_A)
	 = \sum_{x\in S_{\mathcal R}}\mu(x)-W(\mu_A) \\
	& \overset{\text{\ref{item:SP-fge2}}}{=}|S_{\mathcal{R}}|-W(\mu_A)
	\overset{\eqref{eq:degd1muA}}{\ge} |S_{\mathcal{R}}| - (a_1 + a_2) \\
	& \geq |N_w(d_1)\cup N_w(d_2)|-(a_1+a_2)\\
	& = |N_w(d_2)| + |N_w(d_1) \setminus N_w(d_2)| - (a_1 + a_2) \\
	& \geq \deg_w(d_2) + \deg_w(d_1, S_\mathcal{R} \setminus N_w(d_2)) - (a_1 + a_2) \\
 	& = \deg_w(d_2) + \deg_w(d_1) - \deg_w(d_1, N_w(d_2)) - (a_1 + a_2) \\
 	& \overset{\eqref{eq:smallintersection}}{\ge} \deg_w(d_1)+\deg_w(d_2)-(a_1+a_2)-b_1\\
 	& \geq \frac{k}{2}+\frac{k}{2}-(a_1+a_2)-b_1= b_2= \max\{b_1,b_2\}. 
	\end{align*}
	
	\noindent \emph{Step 3: Building $\sigma_B$.}
	Now we will use \Cref{lem:new-extending-out} (Extending out), with
	\begin{center}	
		\begin{tabular}{c|c|c|c|c|c}
			object & $H$ & $\mathcal{R}$ & $S_\mathcal{R}$ & $\mu'$ & $\gamma_B$ \\
			\hline
			in place of & $H$ & $U$ & $V$ & $\mu$ & $\gamma$
		\end{tabular}
	\end{center}
	which implies the existence of a $\gamma_B$-skew-matching $\sigma_B$ such that $\sigma_B \trianglelefteq \mu'$, of weight \begin{align*}
		W(\sigma_B)=b_1+b_2\le  (1+\min\{\gamma_B,\gamma_B^{-1}\})\deg_{w'}(c, \mu'),\end{align*} 
	and with its anchor $\mathcal A(\sigma_B)$ contained in $\mathcal R$.
	Thanks to this, together with the fact that $\sigma_B\trianglelefteq \mu'$ and the definition of $(H, w')$, we can infer that the anchor $\mathcal A(\sigma_B)$ fits in the $w'$-neighbourhood of~$c$.
	By \Cref{prop:adding-skew-matchings} 
	\begin{center}	
	\begin{tabular}{c|c|c|c|c|c|c|c|c|c}
		object & $(H,w)$ & $\fmat{d_1c}$ & $\mu_A$ & $\mu'$ & $w'$ & $\sigma_A$ & $\sigma_\emptyset$ & $\sigma_\emptyset$ & $\sigma_B$ \\
			\hline
		in place of & $(G,w)$ & $\fmat{uv}$ & $\mu$ & $\bar \mu$ & $\bar w$ & $\sigma_A$ & $\sigma_B$ & $\bar \sigma_A$ & $\bar \sigma_B$
	\end{tabular}
	\end{center}
	 we obtain that $(\sigma_A, \sigma_B)$ is a $(\gamma_A, \gamma_B)$-skew pair anchored in $\oriented{d_1 c}$. 
\end{proofclaim}

In the rest of the proof, we may assume that for any $d, d'\in {\mathcal{R}}$ we have  that
\begin{equation}\label{assumption:bigintersection}
\sum_{u\in N_w(d')}  w(\fmat{du})>b_1,
\end{equation}
as otherwise we would be done by \cref{cl:new-weighted-smallintersection}.

\subsection{Proof of \Cref{prop:weighted-structural}: The flabellum case}
Before proceeding with the proof, let us summarise our strategy here.
As $b_2$ is so large (by inequality \eqref{assumption:b_2large}),
the biggest challenge is to find enough space for a $\gamma_B$-skew-matching $\sigma_B$ of weight $b_1+b_2$. For any $d\in \mathcal R$, we can easily accommodate part (at least half) of such a $\gamma_B$-skew-matching using edges between $\mathcal R$ and $N_H(d)$.
This is because inequality~\eqref{assumption:bigintersection} makes it actually very easy to find space for its anchor in $S_\mathcal{R} \subseteq S$.

Unfortunately, $\mathcal R$ might be too small to host the whole \emph{tail} (i.e.,\ the non-anchor part) of $\sigma_B$.
To solve this problem, we shall try to extend the space where to place the \emph{tail} of $\sigma_B$.
To this aim, we define below an auxiliary set $X_d$.
Here we specifically use the fact that we aim to use it only for at most half of the \emph{tail} of $\sigma_B$.
\Cref{cl:new-weighted-evantail} then says that if $X_d\cup \mathcal R$ is large enough to fit the whole tail of $\sigma_B$, then we can find the required $(\gamma_A, \gamma_B)$-skew-matching in $H$.

Referencing the `flabellum structure' described informally at the beginning of the section when we sketched the proof;
the small `base' of the flabellum will be contained in $S_{\mathcal{R}}$, and the large part will be $X_d \cup \mathcal R$.
These two structures will be used to host $\sigma_B$; and $\sigma_A$ will be built after that, in the leftover.

Here is the key definition for this step.
For any $d\in {\mathcal{R}}$, set
\[X_d:=\left\{u\in N_w(c)\setminus ({\mathcal{R}}\cup N_w(d))\::\: \sum_{x\in N_w(u)}  w(\oriented{dx})\ge \frac{b_1}{2}\right\}.\] 

\begin{claim}\label{cl:new-weighted-evantail}
	Suppose there is  $d\in \tilde {\mathcal{R}}$ such  that $|{\mathcal{R}}\cup X_d|\ge b_2$.
	Then $H^\leftrightarrow$ admits a good matching.
\end{claim}

\begin{proofclaim}[Proof of \Cref{cl:new-weighted-evantail}]
The proof has three steps.
Choosing $d$ as in the statement, in the first step, we build a $\gamma_B$-skew-matching $\sigma^\ast_B$ completely outside $\mathcal{R}$, with anchor in $N_w(d)$ and tail in $X_d$.
This matching will contain only a fraction of the weight of the required $\gamma_B$-skew-matching.
In the second step, we build a $\gamma_B$-skew-matching $\sigma'_B$ which completes the required weight, and $\sigma'_B$ will be such that its tail is in $\mathcal{R}$.
In the last step, we build the $\gamma_A$-skew-matching~$\sigma_A$.

Recall that $(\tilde \mu, \tilde \sigma)$ is the  $(H,w,S,M,c,\mu,\gamma_B)$-GE pair,
chosen on 
\Cref{ssec:skew-matching-case} (`The skew-matching cover case'),
and it satisfies \ref{item:strucprop-GEpairanchorS}--\ref{item:strucprop-GEpairmaxsat}.
Let $\mu_d\le \tilde \mu$ and $\sigma_d\le \tilde \sigma$ be defined as in \Cref{cl:new-weighted-alternating-skew-2}, i.e. discarding from $\tilde \mu$ and $\tilde \sigma$ all support edges not intersecting $N_w(d)$.
We define $\mathcal{R}_{\textrm{cov}} = \mathcal{R} \cap V(\mu_d+\sigma_d)$. 
Observe that 
\begin{equation}
	V(\mu_d + \sigma_d) = N_w(d) \cup \mathcal{R}_{\textrm{cov}}
	\label{equation:flabellum-coveredbymudsigmad}
\end{equation}
Indeed, by \ref{item:strucprop-GEpaircovers} we clearly have $V(\mu_d + \sigma_d) \cap S = N_w(d)$.
Since $\sigma_d \leq \tilde \sigma$, by \ref{item:strucprop-GEskewinR}, we have $V(\sigma_d) \setminus \mathcal A(\sigma_d) \subseteq \mathcal{R}$;
and since $\mu_d \leq \tilde \mu$ by \Cref{observation:reachable-weightsR-SR} we have $V(\mu_d) \setminus S \subseteq \mathcal{R}$, which implies the equality above. \medskip

\noindent \emph{Step 1: Starting building $\gamma_B$-skew-matching using $X_d \cup ( \mathcal{R} \setminus \mathcal{R}_{\mathrm{cov}})$ for the tail.}
Let $d$ be as in the statement.
The goal of this step is to build a $\gamma_B$-skew-matching $\sigma_B^*$.
If $b_2<|\mathcal{R}_{\textrm{cov}}|$, we simply set $\sigma_B^*\equiv \sigma_\emptyset$ to be the empty skew-matching, and finalise the construction.
Otherwise, we have that $b_2\ge |\mathcal{R}_{\textrm{cov}}|$.
Let \[\kappa' := \frac{b_2 - |\mathcal{R}_{\textrm{cov}}|}{\gamma_B}.\]
We wish to apply \Cref{prop:weighted-greedy-3} (Third Greedy Lemma) with
\begin{center}	
\begin{tabular}{c|c|c|c|c|c|c|c|c}
	object & $X_d \cup ( \mathcal{R} \setminus \mathcal{R}_{\mathrm{cov}})$ & $ N_w(d)$ & $d$ & $c$ & $\gamma_B$ & $\gamma_A$ & $(\sigma_\emptyset, \sigma_\emptyset)$ & $\kappa'$ \\
		\hline
	in place of & $U$ & $V$ & $u$ & $v$ & $\gamma_A$ & $\gamma_B$ & $(\sigma_A, \sigma_B)$ & $\kappa$
\end{tabular}
\end{center}
Observe that $N_w(d) \subseteq S$, so $\mathcal{R}$ and $N_w(d)$ are disjoint; also $X_d$ is disjoint from $N_w(d)$ by definition.
Hence $X_d \cup ( \mathcal{R} \setminus \mathcal{R}_{\mathrm{cov}})$ is disjoint from $N_w(d)$.
We need to check \ref{item:greedy-3-condition-1}--\ref{item:greedy-3-condition-2}.
With our choice of parameters, \ref{item:greedy-3-condition-1} translates to
\[ |X_d \cup ( \mathcal{R} \setminus \mathcal{R}_{\mathrm{cov}})| = |X_d \cup \mathcal{R}| - |\mathcal{R}_{\mathrm{cov}}| \geq b_2 - |\mathcal{R}_{\textrm{cov}}|, \]
which follows from \eqref{equation:flabellum-coveredbymudsigmad}, the fact that $X_d$ is disjoint from $\mathcal{R} \cup N_w(d)$, and our main assumption $| \mathcal{R} \cup X_d| \geq b_2$.

To see \ref{item:greedy-3-condition-2}, we gather one extra inequality.
From \eqref{equation:flabellum-coveredbymudsigmad} we can deduce that $\sigma_d$ is supported between $N_w(d)$ and $\mathcal{R}$, with the anchor in $N_w(d)$.
Since $\gamma_B > 1$ by \ref{eqnew:gamma_B>1}, we must have that $\sum_{x \in N(d)} \sigma_d(x) \leq \sum_{x \in \mathcal{R}} \sigma_d(x)$.
Similarly, using $\mu_d \leq \tilde \mu$ and \Cref{observation:reachable-weightsR-SR} we have $\sum_{x \in N(d)} \mu_d(x) \leq \sum_{x \in \mathcal{R}} \mu_d(x)$.
Combining these two bounds, we obtain
\begin{equation}
	\frac{k}{2} \leq |N(d)| \leq \sum_{x \in N(d)} (\sigma_d(x) + \mu_d(x)) \leq \sum_{x \in \mathcal{R}} (\sigma_d(x) + \mu_d(x)) \leq |\mathcal{R}_{\textrm{cov}}|
	\label{equation:flabellum-sigmamu}
\end{equation}
Next, note that for any $y\in X_d$ we have (using the definition of $X_d$ in the first inequality),
\begin{align*}
	\deg_w(d, N_w(y)\cap N_w(d))
	&= \sum_{x\in N_w(y) \cap N_w(d)} w(\oriented{dx}) = \sum_{x\in N_w(y)}w(\oriented{dx})\ge \frac{b_1}{2}=\frac{b_2}{2\gamma_B}\\
	& \ge \frac{b_2-k/2}{\gamma_B} 
	\overset{\eqref{equation:flabellum-sigmamu}}{\ge} \frac{b_2-|\mathcal{R}_{\textrm{cov}}|}{\gamma_B},
\end{align*}
giving the required inequality in this case.
Finally, for $y \in \mathcal{R}$, we have
\begin{align*}
	\deg_w(d, N_w(y)\cap N_w(d))
	&= \sum_{x\in N_w(y) \cap N_w(d)} w(\oriented{dx}) = \sum_{x\in N_w(y)}w(\oriented{dx}) \overset{\eqref{assumption:bigintersection}}{\ge} b_1,
\end{align*}
which allows us to conclude as before.
This gives the desired inequality for all $y \in X_d \cup ( \mathcal{R} \setminus \mathcal{R}_{\mathrm{cov}})$, so we obtain \ref{item:greedy-3-condition-2}.

As a consequence of the application of \Cref{prop:weighted-greedy-3}, we deduce the existence of a $\gamma_B$-skew-matching $\sigma_B^*$, of weight \[W(\sigma_B^*)=(1+\gamma_B)\frac{b_2-|\mathcal{R}_{\textrm{cov}}|}{\gamma_B}\le  (1+\gamma_B)\frac{b_1}{2},
\]
with its anchor $\mathcal A(\sigma_B^*)$ fitting in the $w$-neighbourhood of $d$ and the support edges of $\sigma_B^*$ lying in $H[N_w(d), X_d \cup ( \mathcal{R} \setminus \mathcal{R}_{\mathrm{cov}})]$.
Therefore, only the anchor $\mathcal A(\sigma_B^*)$ intersects $N_w(d)$, and $\sigma_B^*$ does not intersect $\mathcal{R}_{\mathrm{cov}}$.

Note that, no matter if $b_2 < |\mathcal{R}_{\textrm{cov}}|$ or not, in both cases our construction of $\sigma^\ast_B$ ensures that we have
\begin{equation}
	W(\sigma_B^*)=(1+\gamma_B)\frac{\max\{0, b_2-|\mathcal{R}_{\textrm{cov}}|\}}{\gamma_B},
	\label{equation:flabellum-weightsigmaBstar}
\end{equation}
and ensures that all vertices from $\mathcal{R}_{\textrm{cov}}$ do not receive any weight from $\sigma^\ast_B$. \medskip
	
	\noindent \emph{Step 2: Completing the $\gamma_B$-skew-matching using $\mathcal{R}_{\mathrm{cov}}$ for the tail.}
Let $\mathcal R_d\subseteq \mathcal{R}_{\textrm{cov}} $ be of size $|\mathcal R_d|=\min\{|\mathcal{R}_{\textrm{cov}}|, b_2\}$.
We wish to apply \Cref{prop:weighted-greedy-3} (Third Greedy Lemma) again, this time with
\begin{center}	
	\begin{tabular}{c|c|c|c|c|c|c|c|c}
		object & $\mathcal R_d$ & $N_w(d)$ & $d$ & $c$ & $\gamma_B$ & $\gamma_A$ & $(\sigma^\ast_B, \sigma_\emptyset)$ & $|\mathcal R_d|/\gamma_B$ \\
		\hline
		in place of & $U$ & $V$ & $u$ & $v$ & $\gamma_A$ & $\gamma_B$ & $(\sigma_A, \sigma_B)$ & $\kappa$
	\end{tabular}
\end{center}
Again, we need to check \ref{item:greedy-3-condition-1}--\ref{item:greedy-3-condition-2}.
With our choice of parameters, \ref{item:greedy-3-condition-1} translates to
\[ |\mathcal R_d| \geq |\mathcal R_d| + \sum_{y \in \mathcal R_d} \sigma^\ast_B(y), \]
which indeed holds, since in any case our construction ensures that $\sigma^\ast_B$ has zero weight in all vertices of $\mathcal R_d \subseteq \mathcal{R}_{\textrm{cov}}$, so the second sum is identically zero.

To verify \ref{item:greedy-3-condition-2}, we note that for every $y \in \mathcal R_d$ we have
\begin{align*}
\deg_w(d, N_w(y)\cap N_w(d))
&=\sum_{x\in N_w(y)}w(\oriented{dx})
\overset{\eqref{assumption:bigintersection}}{>}b_1
 = \frac{b_2}{\gamma_B} = \frac{|\mathcal R_d|+(b_2-|\mathcal R_d|)}{\gamma_B} \\
& =  \frac{|\mathcal R_d|}{\gamma_B} +\frac{\max\{ 0, b_2 - |\mathcal{R}_{\textrm{cov}}|\}}{\gamma_B}\\
&\overset{\eqref{equation:flabellum-weightsigmaBstar}}{=}\frac{|\mathcal R_d|}{\gamma_B} + \frac{W(\sigma^\ast_B)}{1 + \gamma_B} = \frac{|\mathcal R_d|}{\gamma_B}+\sum_{x\in N_w(d)}\sigma_B^*(x),
\end{align*}
where the last equality follows from the fact the anchor of $\sigma^\ast_B$ is contained in $N_w(d)$. 

Thanks to \Cref{prop:weighted-greedy-3}, we obtain a $\gamma_B$-skew-matching $\sigma_B'$ of weight
\begin{equation}\label{eq:weightinR}
W(\sigma_B')=(1+\gamma_B)\frac{|\mathcal R_d|}{\gamma_B} \leq \frac{1 + \gamma_B}{\gamma_B} |\mathcal{R}_{\textrm{cov}}|
\end{equation}
such that $\sigma_B:=\sigma_B^*+\sigma_B'$ is a $\gamma_B$-skew-matching in $H^\leftrightarrow$, with its anchor $\mathcal A(\sigma_B)$ fitting in the $w$-neighbourhood of $d$. 
Note that we have
\begin{equation}
	W(\sigma_B)=W(\sigma_B^*)+W(\sigma_B')=(1+\gamma_B)\frac{b_2}{\gamma_B}=b_1+b_2.
	\label{eq:flabellum-weightsigmaB}
\end{equation}
	
\noindent \emph{Step 3: Building the $\gamma_A$-skew-matching.}
In order to build the $\gamma_A$-skew-matching, we shall heavily leverage the `separation' properties obtained in \Cref{cl:new-weighted-alternating-skew-1} and \Cref{cl:new-weighted-alternating-skew-2}.
Indeed, we shall use the fact that there is a ``separation'' between the (fractional, and skew) matching covering the neighbourhood  $N_w(d)$ and the matching built on top of the rest of $S_{\mathcal R}$.
In this way we can guarantee that the $\gamma_A$-skew-matching we build greedily from its anchor can avoid the already-built $\gamma_B$-skew-matching.

Let $S_{\tilde \sigma}:= \{x\in S\::\: \tilde \sigma(x)=w(\oriented{cx})\}$.
Observe that $N_w(d)\subseteq S_{\tilde \sigma}$, thanks to \Cref{cl:new-weighted-alternating-skew-1} and $d \in \tilde{\mathcal{R}}$.

Define
\begin{align*}
	    V' & = N_w(c)\setminus (\mathcal R_d\cup S_{\tilde \sigma}), \\
	U' & = V(H) \setminus ( \mathcal{R}_{\textrm{cov}} \cup V').
\end{align*}
Clearly, $U'$ and $V'$ are disjoint.
Since $\mathcal R_d \subseteq \mathcal{R}_{\textrm{cov}} \subseteq \mathcal{R} \subseteq N_w(c)$, we have
\begin{equation}
	U' \cup V' = V(H) \setminus \mathcal{R}_d.
	\label{equation:flabellum-uprimevprime}
\end{equation}

Our goal now is to apply \Cref{prop:weighted-greedy-2} (Second Greedy Lemma) with
\begin{center}	
	\begin{tabular}{c|c|c|c|c|c}
		object & $V'$ & $U'$ & $\oriented{cd}$ & $(\sigma_\emptyset, \sigma_B)$ & $a_1$ \\
		\hline
		in place of & $V$ & $U$ & $\oriented{uv}$ & $(\sigma_A, \sigma_B)$ & $\kappa$
	\end{tabular}
\end{center}

Let us verify the required \ref{item:greedy-2-condition-1}--\ref{item:greedy-2-condition-2} conditions hold.
By definition, 
we have for all $x\in S_{\tilde \sigma}$ that $w(\oriented{cx})= \tilde \sigma(x)$.
Using \ref{item:strucprop-GEpairanchorS} we also recall that the anchor of $\tilde \sigma$ is in $S_{\mathcal R}$, which contains $S_{\tilde \sigma}$. Therefore,
\begin{align}
	    \deg_w(c, S_{\tilde \sigma})
    & = \sum_{x \in S_{\tilde \sigma}} w(\oriented{cx})
    = \sum_{x \in S_{\tilde \sigma}} \tilde \sigma(x)
	\leq \frac{W(\tilde \sigma)}{1+\gamma_B}\overset{\eqref{assumption:new-small-skew-weight}}{<}b_1.
	\label{eq:smalldeginN_H(d)}
\end{align}
Using this, we have
\begin{align*}
	    \deg_w(c, N_w(c)\setminus (\mathcal R_d\cup S_{\tilde \sigma}))
    & \ge \deg_w(c)-\deg_w(c, S_{\tilde \sigma})-|\mathcal R_d|\\
	& \overset{\eqref{eq:smalldeginN_H(d)}}{>}
	\deg_w(c)-b_1-|\mathcal R_d|\\
	& \overset{\eqref{eq:weightinR}}{\geq} k-b_1-\frac{\gamma_B W(\sigma_B')}{1+\gamma_B}\\
	& \geq a_1 + b_2 -\frac{\gamma_B W(\sigma_B')}{1+\gamma_B}\\
	& \overset{\eqref{eq:flabellum-weightsigmaB}}{=} a_1 + b_2 -\frac{\gamma_B }{1+\gamma_B} \left( b_1 + b_2 - W(\sigma_B^\ast) \right)\\
	& = a_1 + \frac{\gamma_B}{1 + \gamma_B} W(\sigma^\ast_B) \\
	& = a_1 +\sum_{x\in X_d \cup (\mathcal{R} \setminus \mathcal{R}_{\textrm{cov}})}\sigma_B^*(x) \\
	    & = a_1 +\sum_{x\in N_w(c)\setminus (\mathcal R_d\cup S_{\tilde \sigma})}\sigma_B^*(x) \\
    & = a_1 +\sum_{x\in N_w(c)\setminus (\mathcal R_d\cup S_{\tilde \sigma})}\sigma_B(x)
	\end{align*}
where the third to last line follows since, by construction, $\sigma^\ast_B$ only puts weight from the tail in $X_d$, i.e. $\sigma^\ast_B(x) = (\sigma^\ast_B)^2(x)$ for $x \in X_d$;
the second to last line follows from the fact that $V(\sigma^\ast_B) \cap \mathcal{R}_d = \emptyset$ and $X_d \cup (\mathcal{R} \setminus \mathcal{R}_{\textrm{cov}}) \subseteq N_w(c) \setminus S_{\tilde \sigma}$;
and the last line follows from the fact that $\sigma_B - \sigma^\ast_B = \sigma'_B$ does not have any weight outside $\mathcal{R}_d \cup S_{\tilde \sigma}$.
This gives~\ref{item:greedy-2-condition-1}.

To see that \ref{item:greedy-2-condition-2} holds, we proceed as follows.
First, we have 
\begin{align}
	\sum_{y\in \mathcal{R}_{\textrm{cov}}}\sigma_B(y)
	& =\sum_{y\in \mathcal R}\sigma_B'(y)
	=\frac{\gamma_B}{1+\gamma_B}W(\sigma_B') \overset{\eqref{eq:weightinR}}{=}|\mathcal R_d|
	=\min\{|\mathcal{R}_{\textrm{cov}}|, b_2\}
	\ge \frac{k}{2},
	\label{eq:ObsahInRcapmu_dsigma_d}
\end{align}
where the last inequality holds by \eqref{equation:flabellum-sigmamu} and \eqref{assumption:b_2large}.
Next, we note that, since $\mathcal{R}$ consists of singletons and $\mathcal{R}_{d} \subseteq \mathcal{R}$, every neighbour of $\mathcal{R}_d$ must be in $S$.
Note that \Cref{cl:new-weighted-alternating-skew-2} forbids edges from $\mathcal{R}_{\textrm{cov}}$ to $S \setminus S_{\tilde \sigma}$, and $\mathcal R_d\subseteq\mathcal{R}_{\textrm{cov}}$.
We conclude that there is no edge between the set $N_w(c)\setminus (\mathcal R_d\cup S_{\tilde \sigma})$ and the set $\mathcal{R}_{d}$.
Therefore, for any $x\in N_w(c)\setminus (\mathcal R_d\cup S_{\tilde \sigma})$, we have 
\begin{align*}
	|N_w(x)\cap (V(H)\setminus \mathcal{R}_{d})|
	& \ge \deg_w(x, V(H)\setminus \mathcal{R}_{d})\\
	&= \deg_w(x)\ge \frac{k}{2}\overset{\eqref{eq:ObsahInRcapmu_dsigma_d}}{>}k-\sum_{y\in \mathcal{R}_{\mathrm{cov}}}\sigma_B(y)\\
	&= (a_1+a_2)+W(\sigma_B)-\sum_{y\in \mathcal{R}_{\mathrm{cov}}}\sigma_B(y)\\
	&=(a_1+a_2)+\sum_{y\in V(H)\setminus \mathcal{R}_{\mathrm{cov}}}\sigma_B(y) \\
	&=(a_1+a_2)+\sum_{y\in V(H)\setminus \mathcal{R}_{d}}\sigma_B(y),
\end{align*}
where the last line follows from the fact that $V(\sigma_B) \cap \mathcal{R}_{\mathrm{cov}} \subseteq \mathcal{R}_d$.
Together with \eqref{equation:flabellum-uprimevprime}, this gives \ref{item:greedy-2-condition-2}.

As a consequence of the application of \Cref{prop:weighted-greedy-2}, we obtain a $\gamma_A$-skew-matching $\sigma_A$ of weight $W(\sigma_A)=a_1+a_2$ such that the pair $(\sigma_A, \sigma_B)$ is a $(\gamma_A, \gamma_B)$-skew-matching pair in $H^\leftrightarrow$ anchored in $\oriented{cd}$.
Together with \eqref{eq:flabellum-weightsigmaB}, this finishes the proof of \Cref{cl:new-weighted-evantail}. 
\end{proofclaim}

\subsection{Proof of \Cref{prop:weighted-structural}: The avoiding case} \label{ssection:strucprop-last}

For the rest of the proof, we may assume that for every $d\in \tilde {\mathcal{R}}$, we have
\begin{equation}\label{new:assumption:smallX_d}
|{\mathcal{R}}\cup X_d|<b_2,
\end{equation}
as otherwise we would be done by \cref{cl:new-weighted-evantail}.

Now we will finalise the proof of \Cref{prop:weighted-structural} by constructing a good matching.
The construction will have four main steps.
First, we select a suitable region in $S_{\mathcal R}\cup \mathcal R$ that is ``large enough''. After picking this special region, we pick the anchor $d$ for the $\gamma_A$-skew-matching to be a vertex ``avoiding'' this region. In the second step, we shall construct a partial $\gamma_B$-skew-matching using this region and estimate its size w.r.t. the degree of $d$ into this region. In the third step, we shall construct the $\gamma_A$-skew-matching by leveraging the Separation  Claims (\Cref{cl:new-weighted-alternating-skew-1} and \Cref{cl:new-weighted-alternating-skew-2}) to avoid the large, already built, part of the $\gamma_B$-skew-matching.
In the last step, we complement the already built $\gamma_B$-skew-matching to obtain the full pair. This last step relies on the same avoiding strategy as in Step~3.
 \medskip

\noindent \emph{Step 1: Defining a separating region and finding $d$.}
Slightly counter-intuitively, we first pick an auxiliary vertex $d'$.
This vertex and, in particular, its large neighbourhood in $S_{\mathcal R}$, will allow us to build a pair $(\mu_{d'}, \sigma_{d'})$ which we shall fill by a $\gamma_B$-skew-matching in the second step.

Pick any $d'\in \tilde {\mathcal{R}}$.
Let $\mu_{d'}\le \tilde \mu$ and $\sigma_{d'}\le \tilde \sigma$ be as in \Cref{cl:new-weighted-alternating-skew-2}, i.e. they are obtained from $\mu$ and $\sigma$ by setting the weight to zero in all edges not touching $N_w(d')$, and keeping every other edge intact.

We recall, from \ref{item:strucprop-GEpairanchorS}--\ref{item:strucprop-GEskewinR}, that $\tilde \sigma$ is supported only in directed edges $(u,v)$ with $u \in N_w(c) \cap S$ and $v \in \mathcal{R}$.
Together with \eqref{assumption:new-small-skew-weight}, we have
\begin{equation}
	b_1 > \frac{W(\tilde \sigma)}{1+\gamma_B} =\sum_{x\in N_w(c) \cap S}\tilde \sigma(x).
	 \label{equation:avoiding-wtildesigma}
\end{equation}
Also, observe that by \Cref{cl:new-weighted-alternating-skew-1}, and by the definition of $\sigma_{d'}$, we have 
\begin{equation}
	\deg_w(c, N_w(d')) = \sum_{x \in N_w(d')} w(\oriented{cx}) = \sum_{x \in N_w(d') \cap N_w(c)} \tilde \sigma(x) = \frac{W( \sigma_{d'})}{1+\gamma_B}. \label{equation:avoiding-degwcnwdprime}
\end{equation}
Thus, 
\begin{align*}\nonumber
	\sum_{x\in N_w(c)\setminus (\mathcal R\cup X_{d'}\cup N_w(d'))}w(\oriented{cx})
	& = \deg_w(c, V(H)\setminus (\mathcal R\cup X_{d'}\cup N_w(d'))) \\
	& \geq \deg_w(c) - |\mathcal{R} \cup X_{d'}| - \deg_w(c, N_w(d')) \\
	& \geq k - |\mathcal{R} \cup X_{d'}| - \deg_w(c, N_w(d')) \\
	& \overset{\eqref{new:assumption:smallX_d}}{>} k - b_2 - \deg_w(c, N_w(d')) \\
	& \overset{\eqref{equation:avoiding-degwcnwdprime}}{=} k-b_2-\sum_{x\in N_w(d') \cap N_w(c)}\tilde \sigma(x) \\
	& \geq b_1-\sum_{x\in N_w(d') \cap N_w(c)}\tilde \sigma(x) \\
	& \overset{\eqref{equation:avoiding-wtildesigma}}{>} \sum_{x \in N_w(c) \cap S} \tilde \sigma(x) - \sum_{x \in N_w(d') \cap N_w(c)} \tilde \sigma(x) \\
	& = \sum_{x \in (N_w(c) \setminus N_w(d')) \cap S} \tilde \sigma(x) \\
	& \geq \sum_{x \in  N_w(c) \setminus (\mathcal{R} \cup X_{d'} \cup N_w(d'))} \tilde \sigma(x),
\end{align*}
where the last inequality follows because every $x$ counted in the last sum and not in the previous sum must belong to $N_w(c) \setminus (S \cup \mathcal{R})$, and for those vertices we have $\tilde \sigma(x) = 0$.
This chain of inequalities implies that there is at least one vertex $d\in  N_w(c)\setminus (\mathcal R\cup X_{d'}\cup N_w(d'))$ such that $w(\oriented{cd})>\tilde \sigma(d)$.
Pick such a $d$ such that $\deg_w(d, N_w(d'))$ is maximum.

Now, since $d \in N_w(c) \setminus (\mathcal{R} \cup N_w(d')\cup X_{d'})$, from the definition of $X_{d'}$ we have
\[\deg_w(d',N_w(d))= \sum_{x \in N_w(d)} w(\oriented{d'x}) < \frac{b_1}{2}. \]

Then,
	\begin{align}
W(\mu_{d'})+\frac{W(\sigma_d')}{1+\gamma_B}
& = \sum_{x \in N_w(d')} ( \mu_{d'}(x) + \sigma_{d'}(x) )
= \sum_{x \in N_w(d')} ( \tilde \mu(x) + \tilde \sigma(x) ) \nonumber \\
& \overset{\text{\ref{item:strucprop-GEpaircovers}}}{=}  |N_w(d')|\label{eq:lemavoiding-mezMud'}\\
\nonumber&\ge \deg_w(d, N(d'))+\deg_w(d')-\deg_w(d', N(d)) \\
&> \deg_w(d, N(d'))+\frac k2-\frac{b_1}{2}.\label{eq:lemavoiding-deg(d,N_H(d'))}
	\end{align}

\noindent \emph{Step 2: Starting to build the $\gamma_B$-skew-matching.}
As explained earlier, we build a $\gamma_B$-skew-matching $\sigma_B^*$ in $V(\mu_{d'})$, which is then completed by $\sigma_{d'}$.  

We note first that
\begin{equation}
	V(\mu_{d'}) \setminus S \subseteq \mathcal{R},
	\label{equation:avoiding-supportmudprime}
\end{equation}
that is, $\mu_{d'}$ is supported between $\mathcal{R}$ and $S_\mathcal{R}$.
Indeed, by construction, $V(\mu_{d'}) \cap S \subseteq N_w(d') \subseteq S_\mathcal{R}$, since $d' \in \mathcal{R}$.
Since $\mu_{d'} \leq \tilde \mu$, by \ref{item:strucprop-GEpairsupport} for each $v \in V(\mu_{d'}) \cap S$, each edge with non-zero weight in $\mu_{d'}$ joined to $v$ must have its other endpoint in $\mathcal{R}$.

This allows us to apply \Cref{lem:new-extending-out} (Extending out) with
\begin{center}	
\begin{tabular}{c|c|c|c|c|c}
	object & $H$ & $\mathcal{R}$ & $S$ & $\mu_{d'}$ & $\gamma_B$ \\
	\hline
	in place of & $H$ & $U$ & $V$ & $\mu$ & $\gamma$
\end{tabular}
\end{center}
This application yields a $\gamma_B$-skew-matching $\sigma_B^*\trianglelefteq\mu_{d'}$ of weight \begin{align*}
W(\sigma_B^*)=(1+\gamma_B^{-1})W(\mu_{d'})
\end{align*} with its anchor $\mathcal A(\sigma_B^*)$ contained in $\mathcal R$.

Let $x \in \mathcal{A}(\sigma_B^\ast) \subseteq \mathcal{R}$ be arbitrary, and suppose $\sigma^\ast_B(x) > 0$.
Since $\sigma^\ast_B \trianglelefteq \mu_{d'} \leq \tilde \mu $, by \ref{item:structuprop-GEpairsaturation} we obtain that $(\sigma^\ast_B)^1(x) \leq \tilde \mu(x) \leq w(\oriented{cx})$.
We deduce that the anchor $\mathcal A(\sigma_B^*)$ fits in the $w$-neighbourhood of $c$. 
	
We have 
	\begin{align*}
\nonumber 
	W(\sigma_B^*+\sigma_{d'})
	& = (1+\gamma_B^{-1})W(\mu_{d'})+(1+\gamma_B) \frac{W( \sigma_{d'})}{1+\gamma_B}\\
\nonumber
&=W(\mu_d')+\frac{W(\sigma_{d'}) }{1+\gamma_B} +\gamma_B^{-1}W(\mu_{d'})+\frac{\gamma_B}{1+\gamma_B}W(\sigma_{d'})\\
&\overset{\eqref{eq:lemavoiding-deg(d,N_H(d'))}}>
\nonumber \deg_w(d, N_w(d'))+\frac{k-b_1}{2}+
\gamma_B^{-1}
\left(W(\mu_d')+\frac{W(\sigma_{d'})}{1+\gamma_B}\right)\\
\nonumber &\overset{\eqref{eq:lemavoiding-mezMud'}}\ge \deg_w(d, N_w(d'))+\frac{k-b_1}{2}+\gamma_B^{-1}|N_w(d')|\\
\nonumber &\ge \deg_w(d, N_w(d'))+\frac{k-b_1+\gamma_B^{-1}b_2}{2}\\
&=\deg_w(d, N_w(d'))+\frac k2.
\end{align*}
Hence, by the maximal choice of $d$, we have that
\begin{equation}
W(\sigma_B^*+\sigma_{d'}) \ge \deg_w(x, N_w(d'))+\frac k2
\label{eq:clavoid-W(sigma_d')}
\end{equation}
holds for any $x\in N_w(c)\setminus (\mathcal R\cup X_{d'}\cup N_w(d'))$ that also satisfies $w(\oriented{cx})>\tilde \sigma(x)$. \medskip

\noindent \emph{Step 3: Building the $\gamma_A$-skew-matching $\sigma_A$.}
We find a fractional matching $\widehat \mu$, completely avoiding the already build $\gamma_B$-skew-matching, and we exploit the unique property of the anchor vertex $d$ to  deduce that $d$ has large degree into this matching, allowing us to build a sufficiently large $\gamma_A$-skew-matching $\sigma_A$.

Using that $\sigma_B^\ast \trianglelefteq \mu_{d'}$, we check that $\frac{\gamma_B}{1 + \gamma_B} W(\sigma_B^\ast) \leq W(\mu_{d'})$.
Observe that
\begin{align*}
\frac{\gamma_B}{1+\gamma_B}W(\sigma_B^*+ \tilde \sigma)
& = \frac{\gamma_B}{1+\gamma_B}W(\sigma_B^\ast) + \frac{\gamma_B}{1+\gamma_B} W(\tilde \sigma) 
& \leq W(\mu_{d'})+\frac{\gamma_B}{1+\gamma_B}W(\tilde \sigma) \\
& \le |\mathcal R|<b_2,
\end{align*}
where the last line follows because $\mu_{d'}$ and $\tilde \sigma$ are disjoint, and supported between $\mathcal{R}$ and $S_{\mathcal R}$, by \eqref{equation:avoiding-supportmudprime} and \ref{item:strucprop-GEskewinR}; and the last inequality follows from \eqref{new:assumption:smallX_d}.
Hence, we deduce 
\[ W(\tilde \sigma+\sigma_B^*) < b_1+b_2.\]

Now, note that we have $V(\mu_{d'} + \sigma_{d'}) \cap S = N_w(d')$ and $V(\mu_{d'} + \sigma_{d'}) \setminus S \subseteq \mathcal{R}$.
Since $\mathcal{R}$ consists of singletons, all of its neighbours are in $S$.
Hence, since $w(\oriented{cd})>\tilde \sigma(d)$, by \Cref{cl:new-weighted-alternating-skew-2}, we have that $N_w(d) \cap V(\mu_{d'} + \sigma_{d'}) = N_w(d')$.
Therefore, we get
\begin{align}
\nonumber \deg_w(d, V(H)\setminus V(\mu_{d'}+\sigma_{d'}))&= \deg_w(d)-\deg_w(d, N_w(d'))\\
\nonumber  & \overset{\eqref{eq:clavoid-W(sigma_d')}}{>}\frac k2-\left(W(\sigma_B^*+\sigma_{d'})-\frac{k}{2}\right)\\
\nonumber  &=k-W(\sigma_B^*+\sigma_{d'})\\
\nonumber  &=a_1 + a_2 + b_1 + b_2 - W(\sigma_B^* + \tilde \sigma) + W(\tilde \sigma - \sigma_{d'})\\
&> a_1+a_2+W(\tilde \sigma-\sigma_{d'}).\label{eq:avoiding-deg}
\end{align}

Let $\widehat\mu=  \tilde \mu-\mu_{d'}$. Observe that $V(\widehat \mu)$ does not intersect $V(\mu_{d'}+\sigma_{d'}) \cap S$. Indeed, if $0<\sigma_{d'}(x)<1$, then $x\in N_w(d')$, and thus $\tilde \mu(x)=\mu_{d'}(x)$. Hence, $\widehat \mu(x)=0$.
Define a weight function $\widehat w$ by
\[ \widehat w(\oriented{uv}) = \begin{cases}
	\max\{ 0, w(\oriented{ux})-\tilde \sigma(x)-\mu_{d'}(x) \} & 
	\text{for }u \in \{c,d\},  v=x \in V(H), \\
	w(\oriented{uv})
	& \text{otherwise}.
\end{cases} \]
We claim that we have
\begin{equation}
	\deg_{\widehat w}(d, \widehat\mu) > a_1+a_2.
	\label{equation:avoiding-degwhatdhatmu}
\end{equation}
Indeed,
\begin{align*}
	\deg_{\widehat w}(d, \widehat\mu)
	&= \sum_{x \in V(H)} \min\{ \widehat w(\oriented{dx}), \widehat \mu(x) \} \\
	&\ge \sum_{x \in V(H)\setminus V(\mu_{d'}+\sigma_{d'})} \min\{ \widehat w(\oriented{dx}), \widehat \mu(x) \} \\
	&\ge \sum_{x \in V(H)\setminus V(\mu_{d'}+\sigma_{d'})} \min\{w(\oriented{dx}) - \tilde \sigma(x) - \mu_{d'}(x),\ \tilde \mu(x) - \mu_{d'}(x)\} \\
	&= \sum_{x \in V(H)\setminus V(\mu_{d'}+\sigma_{d'})} \min\{ w(\oriented{dx}) - \tilde \sigma(x),\ \tilde \mu(x)\} \\
	&= \sum_{x \in V(H)\setminus V(\mu_{d'}+\sigma_{d'})} \bigl( w(\oriented{dx}) - \tilde \sigma(x) \bigr) \\
	&= \sum_{x \in V(H)\setminus V(\mu_{d'}+\sigma_{d'})} \bigl( w(\oriented{dx}) - \tilde \sigma(x) + \sigma_{d'}(x) \bigr) \\
	&\ge \deg_w\!\bigl(d, V(H)\setminus V(\mu_{d'}+\sigma_{d'})\bigr)
	- \sum_{x \in V(H)\setminus V(\mu_{d'}+\sigma_{d'})} \bigl( \tilde \sigma(x) - \sigma_{d'}(x) \bigr) \\
	&\overset{\eqref{eq:avoiding-deg}}{>} a_1 + a_2 + W(\tilde \sigma - \sigma_{d'}) 
	- \sum_{x \in V(H)\setminus V(\mu_{d'}+\sigma_{d'})} \bigl( \tilde \sigma(x) - \sigma_{d'}(x) \bigr) \\
	&\ge a_1 + a_2.
\end{align*}
where in the fifth line we used $N_w(d) \subseteq S_{\mathcal R}$ and \ref{item:strucprop-GEpaircovers-w}. In the sixth line we used that $x \notin V(\mu_{d'}+\sigma_{d'})$.

We apply \Cref{lem:new-combination} (Combination) with
\begin{center}	
	\begin{tabular}{c|c|c|c|c}
		object & $(H, \widehat{w})$ & $d$ & $\widehat{\mu}$ & $\gamma_A$ \\
		\hline
		in place of & $(H,w)$ & $v$ & $\mu$ & $\gamma$
	\end{tabular}
\end{center}
and we obtain as a consequence a $\gamma_A$-skew-matching $\sigma_A$, which satisfies $\sigma_A \trianglelefteq \widehat \mu$, and whose weight (scaling down if necessary, using \eqref{equation:avoiding-degwhatdhatmu}) is 
\[W(\sigma_A) = a_1+a_2.\]
Moreover, its anchor $\mathcal A(\sigma_A)$ fits in the $\widehat w$-neighbourhood of $d$. 
 
 We shall argue now that $(\sigma_A,\tilde\sigma+\sigma_B^*)$ is a $(\gamma_A, \gamma_B)$-skew-matching pair anchored in $\oriented{dc}$.
 We need to check properties \ref{itdef:disjoint}--\ref{itdef:anchorpartition}.
 The disjointedness of $\sigma_A$ with $\tilde \sigma+\sigma_B^*$ follows from the fact that $\sigma_A \trianglelefteq \widehat \mu =  \tilde \mu - \mu_{d'}$ and $\sigma^\ast_B \trianglelefteq \mu_{d'}$,
together with the fact that $\tilde \mu$ is disjoint from $\tilde \sigma$.
This gives \ref{itdef:disjoint}.
 We have already stated above that $\mathcal A(\sigma_A)$ fits in the $\widehat{w}$-neighbourhood of $d$, which gives \ref{itdef:anchorc} of \Cref{def:weighted-anchored-pair}.
 We have also said that the anchor $\mathcal A(\sigma_B^*)$ fits in the $w$-neighbourhood of $c$, and also we have $\mathcal A(\sigma_B^*) \subseteq \mathcal R$.
 By \ref{item:strucprop-GEpairanchorS}-\ref{item:strucprop-GEpairanchorfits}, we have that the anchor $\mathcal A(\tilde\sigma)$ is contained in $S$ and fits in the $w$-neighbourhood of $c$.
 As $\mathcal R$ and $S_{\mathcal R}$ are disjoint, we have that $\mathcal A(\tilde \sigma+\sigma_B^*)$ fits in the $w$-neighbourhood of $c$, so we have \ref{itdef:anchord}.
 
 Finally, to see \ref{itdef:anchorpartition}, we note first that $N(c) \cap N(d) \subseteq N(d) \subseteq S_\mathcal{R}$.
 Since $\sigma_B^\ast$ is anchored in $\mathcal{R}$, which is disjoint from $S_\mathcal{R}$, for each $x \in N(c) \cap N(d)$ we have $(\sigma_B^*)^{1}(x) = 0$.
 Hence, to see \ref{itdef:anchorpartition} it suffices to show, for every $x \in N(c) \cap N(d)$, that
 \begin{equation*}
 	\max \{ w(\oriented{cx}), w(\oriented{dx}) \} \geq \sigma_A^1(x)+\tilde \sigma^1(x)
 \end{equation*}
 holds.
 Suppose first that $\sigma^1_A(x) = 0$.
 Then to get the desired inequality it suffices to check that $w(\oriented{cx}) \geq \tilde \sigma^1(x)$ holds, and this is true because the anchor of $\tilde \sigma$ fits in the $w$-neighbourhood of $c$.
 Hence, we can assume that $\sigma^1_A(x) > 0$.
 Since $\sigma_A \trianglelefteq \widehat \mu = \tilde \mu - \mu_{d'}$, we have that $\tilde \mu(x) \neq \mu_{d'}(x)$, meaning that $x \notin N(d')$, and therefore, $\mu_{d'}(x) = 0$.
 On the other hand, since the anchor of $\sigma_A$ fits in the $\widehat w$-neighbourhood of $d$, from $\sigma^1_A(x) > 0$ we also have that $\widehat w (\oriented{dx}) \geq \sigma^1_A(x) > 0$.
 From the definition of $\widehat w$, this implies that $\widehat w (\oriented{dx}) = w(\oriented{dx}) - \tilde \sigma(x)$ (since $\mu_{d'}(x) = 0$).
 Putting all together, we have
 \[ w(\oriented{dx}) = \widehat w (\oriented{dx}) + \tilde \sigma(x) \geq \sigma^1_A(x) + \tilde \sigma^1(x), \]
 as desired.
 This finishes the verification of \ref{itdef:anchorpartition}.
\medskip

\noindent \emph{Step 4: Completing the $\gamma_B$-skew-matching.}
Having found $\sigma_A$, we need to extend $\tilde \sigma+\sigma_B^*$ in order to obtain $\sigma_B$.
We use a similar trick as in the previous step, but instead of picking an anchor $d$ that avoids $\sigma_B^*+\sigma_{d'}$,
we shall select  ``avoiding vertices'' to place the anchor of $\widehat \sigma_B$,
the missing part of the $\gamma_B$-skew-matching.
Using their avoiding properties, we can greedily build then the tail of $\widehat \sigma_B$.
 
Define \begin{align*}
	V &:= N_{\widehat w}(c)\setminus (\mathcal R\cup X_{d'}\cup N_w(d')), \text{ and} \\
	U &:=V(H)\setminus ( V(\sigma_{d'}+\mu_{d'})\cup V).
\end{align*}
We gather first some useful observations.
We note first that, since $V(\sigma_{d'} + \mu_{d'}) \subseteq N_w(d') \cup \mathcal{R}$, we have that $V \cap V(\sigma_{d'} + \mu_{d'}) = \emptyset$.
This implies that
\begin{equation}
	U \cup V = V(H) \setminus V(\sigma_{d'} + \mu_{d'}).
	\label{equation:avoiding-uunionv}
\end{equation}

Suppose that $x \notin \mathcal{R} \cup X_{d'}$.
The claim is that $x \in N(d') \cup V$, or $w(\oriented{cx}) = 0$.
Indeed, suppose that $x \notin \mathcal{R} \cup X_{d'} \cup N(d') \cup V$.
In particular we have $x \notin N_{\widehat w}(c)$, so $\widehat w(\oriented{cx}) = 0$.
If $x \in S$, we have $\mu_{d'}(x) = 0$ (because $x \notin N(d')$) and $w(\oriented{cx}) \geq \tilde \sigma(x)$ (since $\tilde \sigma$ fits in the $w$-neigbourhood of $c$), then by definition of $\hat w$ we have $0 = \widehat w(\oriented{cx}) = w(\oriented{cx}) - \tilde \sigma(x) \geq 0$, hence $w(\oriented{cx}) = 0$, as desired.
On the other hand, if $x \notin S$, then we have that $\mu_{d'}(x) = \tilde \sigma(x) = 0$ (by \ref{item:strucprop-GEskewinR}, \eqref{equation:avoiding-supportmudprime} and $x \notin \mathcal{R}$).
Then again by definition of $\hat w$ we have $0 = \widehat w(\oriented{cx}) = w(\oriented{cx}) = 0$, as desired.
We conclude that $V(H) \setminus (\mathcal{R} \cup X_{d'}) \subseteq N(d') \cup V \cup \{ x : w(\oriented{cx}) = 0 \}$, so in particular
\begin{equation}
	\sum_{x\in V(H)\setminus  (\mathcal R\cup X_{d'})}w(\oriented{cx}) \leq \sum_{x\in N_w(d')}w(\oriented{cx}) + \sum_{x \in V} w(\oriented{cx}).
	\label{equation:avoinding-inequalitysumwcV}
\end{equation}

Similarly, again using \ref{item:strucprop-GEskewinR} we deduce that for each $x \in (V \cup N_w(d)) \setminus S$, we have $\tilde \sigma(x) = 0$.
In particular, since $V$ and $N_w(d')$ are disjoint, we have
\begin{equation}
	\sum_{x\in V} \tilde \sigma(x) + \sum_{x \in N(d')} \tilde \sigma(x) \leq \sum_{x \in S} \tilde \sigma(x).
	\label{equation:avoinding-inequalitysumsigmaS}
\end{equation}

Finally, we observe the following.
Suppose $x \in N_{\widehat{w}}(c) \setminus ( \mathcal{R} \cup N(d') )$; and consider $Y = \mathcal{R} \cap (V(\mu_{d'} + \sigma_{d'}))$.
We claim that
\begin{equation}
	N(x) \cap Y = \emptyset.
	\label{equation:usign-separationg}
\end{equation}
Indeed, since $N(Y) \subseteq N(\mathcal{R}) = S_\mathcal{R} \subseteq S$, it suffices to analyse the case where $x \in S$.
Since $x \in  N_{\widehat{w}}(c)$, we have $\widehat{w}(\oriented{cx}) > 0$.
By definition of $\widehat{w}$, this implies that $w(\oriented{cx}) > \tilde \sigma(x) + \mu_{d'}(x) \geq \tilde \sigma(x)$.
Hence, by \Cref{cl:new-weighted-alternating-skew-2}, we have that $x$ has no neighbours in $Y$, as desired.

We wish to apply \Cref{prop:weighted-greedy-2} (Second Greedy Lemma) with
\begin{center}	
	\begin{tabular}{c|c|c|c|c|c|c|c}
		object & $V$ & $U$ & $\oriented{cd}$ & $(\tilde \sigma + \sigma^\ast_B, \sigma_A)$ & $\gamma_B$ & $\gamma_A$ & $((b_1+b_2)-W(\sigma_B^*+\tilde \sigma)) / (1 + \gamma_B)$ \\
		\hline
		in place of & $V$ & $U$ & $\oriented{uv}$ & $(\sigma_A, \sigma_B)$ & $\gamma_A$ & $\gamma_B$ & $\kappa$
	\end{tabular}
\end{center}
Note that $U$ and $V$ are disjoint.
We verify the required conditions \ref{item:greedy-2-condition-1}--\ref{item:greedy-2-condition-2}.
We have
 \begin{align*}
 \deg_{w}(c, V)& = \sum_{x \in V} w(\oriented{cx})
 \overset{\eqref{equation:avoinding-inequalitysumwcV}}{\ge}
 \sum_{x\in V(H)\setminus  (\mathcal R\cup X_{d'})}w(\oriented{cx})-\sum_{x\in N_w(d')}w(\oriented{cx})\\
 & \geq \sum_{x \in V(H)} w(\oriented{cx}) - |\mathcal{R} \cup X_{d'}| -\sum_{x\in N_w(d')}w(\oriented{cx})\\
 & \overset{\text{Claim~\ref{cl:new-weighted-alternating-skew-1}}}{=} \sum_{x\in V(H)}w(\oriented {cx})-|\mathcal R\cup X_{d'}|-\sum_{x\in N_w(d')}\tilde \sigma(x)\\
 &  = \deg_w(c) - |\mathcal R\cup X_{d'}| - \sum_{x\in N_w(d')} \tilde \sigma(x) \\
& \overset{\eqref{equation:avoinding-inequalitysumsigmaS}}{\ge} \deg_w(c)-|\mathcal R\cup X_{d'}|-\sum_{x\in S}\tilde \sigma (x)+\sum_{x\in V}\tilde \sigma(x)\\
& \overset{\eqref{new:assumption:smallX_d}}{\ge} k - b_2 -\sum_{x\in S_{\mathcal R}}\tilde \sigma (x)+\sum_{x\in V}\tilde \sigma(x)\\
 &\overset{\text{\ref{item:strucprop-GEpairanchorS}\&\ref{item:strucprop-GEskewinR}}}{\geq} k -b_2 - \frac{W(\tilde \sigma)}{1+\gamma_B}+\sum_{x\in V}\tilde \sigma(x)
 \\
 & = b_1 - \frac{W(\tilde \sigma)}{1+\gamma_B} + (a_1 + a_2) + \sum_{x\in V}\tilde \sigma(x) \\
 & \geq b_1 - \frac{W(\tilde \sigma)}{1+\gamma_B} + \sum_{x \in V} \sigma_A(x) + \sum_{x\in V}\tilde \sigma(x) \\
& \geq b_1-\frac{W(\tilde \sigma+\sigma_B^*)}{1+\gamma_B}+\sum_{x\in V}(\sigma_A(x)+\tilde \sigma(x)+\sigma_B^*(x)),
 \end{align*}
 where in the last line we used that $\sigma^\ast_B \trianglelefteq \mu_{d'}$, and for every $x \in V$ we have $\mu_{d'}(x) = 0$, since $\mu_{d'}$ is supported between $N_w(d')$ and $\mathcal{R}$.
 This implies that $\sigma^\ast_B(x) = 0$ for all $x \in V$; and so indeed we have \ref{item:greedy-2-condition-1}.
 
Now we check \ref{item:greedy-2-condition-2}.
Recall from \eqref{equation:avoiding-uunionv} that $U \cup V = V(H) \setminus V(\sigma_{d'} + \mu_{d'})$, and that $\mu_{d'}$ and $\sigma_{d'}$ are supported in $N(d') \cup \mathcal{R}$. Also notice that if $x\in V$, then $\widehat w(\oriented{cx})>0$, i.e., $w(\oriented{cx})>\tilde \sigma(x)$ and thus~\eqref{eq:clavoid-W(sigma_d')} applies to all $x\in V$. 
Using this, we have, for all $x \in V$,
\begin{align*}
	|N(x) \cap (U \cup V)|
	& \geq \deg_w(x, V\cup U) \\
	& \geq \deg_w(x, V(H)\setminus V(\sigma_{d'}+\mu_{d'}))\\
    & =\deg_w(x)-\deg_w(x, \mathcal R\cap V(\sigma_{d'}+\mu_{d'}))-\deg_w(x, N_w(d')) \\
    & \geq \frac{k}{2}-\deg_w(x, \mathcal R\cap V(\sigma_{d'}+\mu_{d'}))-\deg_w(x, N_w(d')) \\
    & \overset{\eqref{equation:usign-separationg}}{=} \frac{k}{2}-\deg_w(x, N_w(d')) \\
	& \overset{\eqref{eq:clavoid-W(sigma_d')}}{\ge}
	\frac{k}{2}-\left(W(\sigma_B^*+\sigma_{d'})-\frac k2\right)
	= k-W(\sigma_B^*+\sigma_{d'})\\
	& = b_1 + b_2 -W(\sigma_B^*+\tilde \sigma) + W(\tilde \sigma - \sigma_{d'}) + (a_1 + a_2) \\
	& \geq b_1 + b_2 -W(\sigma_B^*+\tilde \sigma) + W(\tilde \sigma - \sigma_{d'}) + \sum_{x\in V\cup U}\sigma_A(x)\\
	& \geq b_1 + b_2 -W(\sigma_B^*+\tilde \sigma) +  \sum_{x\in V\cup U} (\sigma_A(x) + \tilde \sigma(x) - \sigma_{d'}(x)) \\
	& = (b_1+b_2-W(\sigma_B^*+\tilde \sigma))+\sum_{x\in V\cup U}\left(\sigma_A(x)+\tilde \sigma(x)+\sigma_B^*(x)\right),
\end{align*}
where in the last inequality we use the fact that $\sigma_{d'}(x) = \sigma_B^*(x)=0$ for $x\in U\cup V$ (this follows from $\sigma_B^\ast \trianglelefteq \mu_{d'}$ together with \eqref{equation:avoiding-uunionv}).

The outcome of the application of \Cref{prop:weighted-greedy-2} is a $\gamma_B$-skew-matching $\widehat\sigma_B$ in $H^\leftrightarrow$ of weight $W(\widehat\sigma_B)=(b_1+b_2)-W(\sigma_B^*+\tilde \sigma)$ such that $(\widehat \sigma_B+\sigma_B^*+\tilde \sigma, \sigma_A)$ is a $(\gamma_B, \gamma_A)$-skew-matching anchored in $\oriented{cd}$.
We set $\sigma_B:= \widehat\sigma_B+\sigma_B^*+\tilde \sigma$ and have 
\[W(\sigma_B)=b_1+b_2.\]
Hence, $(\sigma_A, \sigma_B)$ is a good matching.
This finishes the proof of \Cref{prop:weighted-structural}. \hfill \qedsymbol

\section{Embedding the Tree}\label{sec:embed}

In this section, we prove the Tree Embedding Lemma (\Cref{lemma:treeembedding}), which we restate for convenience.

\restatetreembedding*

The section is organised as follows.
We gather useful results about probability and the regularity method in~\cref{ssec:prob}.
In \Cref{ssec:shrub-embedding}, we prove a lemma that embeds a `shrub' in a regular pair and will be used repeatedly during the proof.
We give a sketch of the proof in \Cref{ssec:proof-embed-sketch}.
Then we proceed with the main proof, which is split in several steps, and takes the remainder of the section.

\subsection{Preliminaries}\label{ssec:prob}

We will need a bounded-differences inequality~\cite{McDiarmid1989}.

\begin{lemma} \label{lemma:mcdiarmid}
	Let $X_1, \dotsc, X_N$ be independent random variables, with $X_i \in \Lambda_i$.
	Let $f: \prod_{i=1}^N \Lambda_i \rightarrow \mathbb{R}$ 
	be a function such that for any $z, z' \in \prod_{i=1}^N \Lambda_i$ 
	which differ only in the $k$th coordinate, we have $|f(z) - f(z')| \leq c_k$.
	Then, the random variable $X = f(X_1, \dotsc, X_N)$ satisfies, for any $t \geq 0$,
	\[ \probability[X \geq \expectation[X] + t] \leq \exp \left( - \frac{2 t^2}{\sum_{i=1}^N c_i^2} \right). \]
\end{lemma}

We now collect a few results about graph regularity and regular pairs.
Let $(X,Y)$ be an $\varepsilon$-regular pair with density $\dred$, and $Y' \subseteq Y$.
Say that $x \in X$ is \emph{typical to $Y'$} if $x$ has at least $(\dred - \varepsilon)|Y'|$ neighbours in $Y'$.
The following fact is well-known (cf. \cite[Lemma 7.5.1]{Diestel}).

\begin{lemma} \label{lemma:typicalvertex}
	Let $(X,Y)$ be an $\varepsilon$-regular pair, and let $Y' \subseteq Y$ with $|Y'| \geq \varepsilon |Y|$.
	Then all but at most $\varepsilon |X|$ vertices in $X$ are typical to $Y'$.
\end{lemma}

We generalise this to deduce that many vertices are typical to many sets.

\begin{lemma} \label{lemma:typicalvertex2}
	Let $X, Y_1, \dotsc, Y_N$ be such that $(X,Y_i)$ is an $\varepsilon$-regular pair for each $1 \leq i \leq N$,
	and for each $i$ let $Y'_i \subseteq Y_i$ with $|Y'_i| \geq \varepsilon |Y_i|$.
	Then all but at most $\sqrt{\varepsilon} |X|$ vertices in $X$ are typical to all but at most $\sqrt{\varepsilon}N$ of the sets $Y'_1, \dotsc, Y'_N$.
\end{lemma}

\begin{proof}
	Let $B$ be an auxiliary bipartite graph with classes $X$ and $[N] = \{1, \dotsc, N\}$, where $x \in X$ and $i \in [N]$ are joined by an edge if $x$ is \emph{not} typical to $Y'_i$.
	By \Cref{lemma:typicalvertex}, each $i$ is adjacent to at most $\eps |X|$ edges in $B$,
	so $B$ has at most $\varepsilon |X| N$ edges.
	A double-counting argument shows that the number of vertices of $X$ which have degree at least $\sqrt{\varepsilon}N$ in $B$ is at most $\sqrt{\varepsilon} |X|$, as required.
\end{proof}

\subsection{Embedding shrubs}\label{ssec:shrub-embedding}
Now, we state and prove an auxiliary lemma, which concerns the embedding of a `shrub' in a regular pair.
The proof follows from very standard arguments about embedding in regular pairs,
so the reader familiar with graph regularity may safely skip its proof.

\begin{definition}
	We say that $(F,r,x)$ is a \emph{rooted shrub} if $F$ is a tree,
	 $r \in V(F)$, and either $x = \emptyset$ or $x \in V(F)$ and the distance of $x$ to $r$ is even and at least $4$.
	We say that $r$ is the \emph{root} of $F$,
	and if $x \neq \emptyset$ we call $x$ the \emph{adventitious root} of $F$ (and say that $F$ has an adventitious root).
\end{definition}

\begin{lemma}[Embedding a shrub]\label{lem:embed-embedding the shrubs}
	Suppose that $2 \eps \leq \dred \leq 1/3 $ and $\varepsilon <\frac{\dred^2\tilde \eta}{8}$.
	Let $(F, r, x)$ be a rooted shrub.
	Let $(X,Y)$ be an $\eps$-regular pair in a graph $G$ with $d(X,Y) \geq \dred$ and $|X| = |Y|$.
	Let $U\subseteq X \cup Y$ and $P\subset X$ be disjoint sets, and $v\in X \setminus P$, and suppose
	\stepcounter{propcounter}
	\begin{enumerate}[\upshape{(\Alph{propcounter}\arabic*)},topsep=0.7em, itemsep=0.5em]
		\item \label{item:shrubembed-XY} $|X\cap U|, |Y\cap U| \ge \tilde \eta|X|$,
		\item \label{item:shrubembed-P} $|P|\ge 2\varepsilon|X|$,
		\item \label{item:shrubembed-degv} $\deg_G(v, Y\cap U)\ge 2\varepsilon |Y|$,
		\item \label{item:shrubembed-Tsmol} $F$ has at most $\varepsilon |X|$ vertices.
	\end{enumerate}
	Then there is an embedding $\varphi$ of $F$ in $G$ with $\varphi(V(F)\setminus \{x, r\})\subseteq U$, $\varphi(r) = v$, and $\varphi(x)\in P$.
\end{lemma}

\begin{proof}
	Let $m = |X| = |Y|$.
	We suppose first that $x \neq \emptyset$, i.e. that $F$ contains an adventitious root $x$, and we define $\varphi(x)$ first.
	Since $|Y \cap U| \geq \tilde{\eta} m \geq \eps m$, then \Cref{lemma:typicalvertex} implies that all but at most $\eps m$ vertices in $X$ are typical to $Y \cap U$.
	In particular, we may select one such typical vertex $w$ which lies in $P$, and we set $\varphi(x) = w \in P$.
	By the choice of $w$, we have $\deg(w, Y \cap U) \ge (\dred - \eps) |Y \cap U| \geq \frac{\dred\tilde \eta }{2}m$.
	Set $N_w:=N(w)\cap Y\cap U$.
	If there is no adventitious root in $F$, we simply set $N_w:= Y\cap U$.
	In any case, $N_w \subseteq Y \cap U$ is defined and it holds that $|N_w| \geq \frac{\dred\tilde \eta }{2}m$.
	
	Let $X'$ be the set of vertices in $X$ that are typical to $N_w$, and let $Y'$ be the set of vertices in $Y$ that are typical to $X\cap U$.
	Again, \Cref{lemma:typicalvertex} implies that $|X\setminus X'|\le \varepsilon m$ and $|Y\setminus Y'|\le \varepsilon m$. 
	Observe that for every $u\in X'$, we have 
	\[\deg(u, N_w\cap Y')\ge (\dred-\varepsilon)|N_w|-|Y\setminus Y'|\ge \frac{\dred^2\tilde \eta}{4}m-\varepsilon m \ge \eps m \ge |V(F)|;\]
	for every $u\in Y'$, we have 
	\[\deg(u, X'\cap U)\ge (\dred-\varepsilon)|X\cap U|-|X\setminus X'|\ge  \frac{\dred\tilde \eta}{2}m -\varepsilon m \ge \varepsilon m \ge |V(F)|, \]
	and 
	\[\deg(v, Y'\cap U)\ge \deg(v, Y\cap U)-|Y\setminus Y'|\ge 2\varepsilon m -\varepsilon m \ge |V(F)|.\]
	Therefore we can greedily embed $V(F)\setminus \{x\}$ in $U\cup \{v\}$ so that $r$ is mapped to $v$,
	the neighbours of $r$ are mapped to $Y'$,
	the vertices of even distance from $r$ in $X'$
	and the vertices of odd distance at least $3$ from $r$ to $N_w\cap Y'$. 
\end{proof}

\subsection{Sketch of the proof}\label{ssec:proof-embed-sketch}
We will find an embedding of a tree $T$  in a graph $G$, provided $G$ has a regular partition whose weighted $\dred$-reduced graph has a suitable skew-matching pair.
The proof will consist of seven steps.

\begin{enumerate}
	\item \emph{Setting the stage.} First, we will set the stage and constants.
	By assumption, the tree $T$ to be embedded has a corresponding fine partition $(W_A, W_B, \mathcal{F}_A, \mathcal{F}_B)$ which decomposes $T$ into `seeds' and `shrubs'.
	The reduced host graph has a skew-matching pair $(\sigma_A, \sigma_B)$ which is anchored to an edge $(V_c, V_d)$.
	We select buffer and reservoir sets $U$ and $Q$ in $V(G)\setminus (V_c\cup V_d)$.

	\item \emph{Embedding the seeds.} Next, we will embed the `seeds' of $T$, meaning we embed $W_A \cup W_B$ into $V_c \cup V_d$.
	
	\item \emph{Allocating the shrubs.} Then we will allocate (but not embed yet) the shrubs from	
	$\mathcal{F}_A \cup \mathcal{F}_B$ in a suitable way, according to the weights indicated by $(\sigma_A, \sigma_B)$.
	This is done using a randomised procedure.
	At the end of this step, for each shrub $F$, we will have selected a cluster $V_i$ where one colour class of the shrub will be embedded.
	
	\item \emph{Allocating the roots.} We reserve, for each shrub $F$, enough space in the selected cluster $V_i$.
	We will select pairwise disjoint sets $R_F \subseteq V_i$, 
	one for each $F \in \mathcal{F}_A \cup \mathcal{F}_B$.
	This $R_F$ will be a `private' set for the shrub $F$, with enough space for the embedding of $F$. 
	
	\item \emph{Finding suitable clusters.} Having chosen $V_i$ and $R_F$ for all shrubs $F$, we will argue that there is always (even after having embedded some shrubs) a way to choose another `suitable' cluster $V_\ell$ for each $F$, in such a way that there is enough space in $G[V_i, V_\ell]$ to embed $F$.
	
	\item \emph{First Embedding Phase.} In the last two steps, we will carry the actual embedding of the shrubs.
	Here, we will try to embed the shrubs in their respective chosen bipartite graphs $G[V_i,V_j]$, avoiding (nearly completely) the reservoir set $Q$. 	However, this embedding can fail for some (few) shrubs.
	
	\item \emph{Second Embedding Phase.}
	In this phase, the remaining few shrubs, which failed to be embedded in the previous phase, are finally embedded in $Q$.
	For this purpose, we use \Cref{lem:embed-embedding the shrubs}.
\end{enumerate}

The rest of this section is dedicated to the proof of \Cref{lemma:treeembedding}, and spans several subsections.

\subsection{Proof of \Cref{lemma:treeembedding}: Setting the stage}
Suppose $\dred, \eta, q >0$ and $t\in \mathbb N$ are given.
Set 
\[
\eta':= {q\eta^3}.
\]
We assume $0<\dred \ll \eta' \ll \eta \ll 1$; otherwise we decrease their values to satisfy this hierarchy.
Set 
\[
\varepsilon:= \frac{\dred^2\cdot (\eta')^2}{100}, \qquad 
\rho:= \min\left\lbrace \frac{(\eta')^4\eta^2}{10^3t^3},\frac{\dred^2\eta}{100}, \frac{\varepsilon}{2t} \right\rbrace, \qquad
n_0:=\frac{4\cdot 10^3\cdot  t}{\rho^2}.
\]
  The parameters satisfy 
 \begin{equation}\label{eq:param-embedding}
 0<1/n_0 \ll \rho\ll \varepsilon\ll \dred\ll \eta'\ll \eta\ll 1. \end{equation}

Let $G$ be a given graph on $n$ vertices, and let $\Gamma := \Gamma_{\dred, \varepsilon}$ be the reduced graph with  $|V(\Gamma_{\dred, \varepsilon})|=N\le t$,
and let $\mathcal{P} = \{V_0, V_1, \dotsc, V_{N}\}$ be the $\eps$-regular equitable partition which yields $\Gamma$, so $i \in V(\Gamma)$ corresponds to the cluster $V_i \subseteq V(G)$. Let $w:E(\Gamma)\rightarrow \{0\} \cup [\dred,1]$ be the weight function defined by $w(\fmat{ij}):= d(V_i,V_j)$, if $d(V_i,V_j)\ge\dred$ and $w(\fmat{ij})=0$, otherwise. 
Let $\oriented{cd}$ be the edge in $E(\Gamma^\leftrightarrow)$ where $(\sigma_A, \sigma_B)$ is anchored.

Let $m := |V_1|$ be the common size of all clusters $V_1, \dotsc, V_N$.
Note that since $\mathcal{P}$ is an $\varepsilon$-regular equitable partition, we have
\begin{align}
	(1 - \varepsilon) \frac{n}{N} \leq m \leq \frac{n}{N}. \label{equation:clustersize}
\end{align}

Let $T \in \mathcal T_{a_1, a_2, b_1, b_2}^{\rho}$ be an arbitrary tree, which we need to embed into $G$.
Let $(W_A,W_B, \mathcal F_{A},\mathcal F_{B})$ be the $(\rho|V(T)|)$-fine partition of $T$ that witnesses $T\in \mathcal T_{a_1, a_2, b_1, b_2}^{\rho}$.
By assumption, the values $a_1,a_2, b_1, b_2, k$ must satisfy
\begin{align}
	a_1 + a_2 + b_1 + b_2 & = k \geq qn \label{equation:klarge} \\
	a_1, a_2, b_1, b_2 & \geq \eta k. \label{equation:aibilarge}
\end{align}
Let $\gamma_A := a_2/a_1$ and $\gamma_B := b_2/b_1$.
From \eqref{equation:aibilarge}, we easily get
\begin{align}
	\gamma_A, \gamma_B & \geq \eta. \label{equation:gammalarge}
\end{align}

We recall that the definition of $\sigma^1$ (\Cref{def:sigma}) implies that, for every $i \in V(\Gamma)$,
\begin{equation}
	\sum_{j \in N_\Gamma(i)} \sigma_A(\oriented{ij}) = (1 + \gamma_A) \sigma^1_A(i),
	\label{equation:sigma1Ai}
\end{equation}
and similarly,
\begin{equation}
	\sum_{j \in N_\Gamma(i)} \sigma_B(\oriented{ij}) = (1 + \gamma_B) \sigma^1_B(i).
	\label{equation:sigma1Bi}
\end{equation}

By assumption, we have that $(\sigma_A, \sigma_B)$ is a $(\gamma_A, \gamma_B)$-skew-matching pair anchored in $\oriented{cd}$ with weights $W(\sigma_A), W(\sigma_B)$.
In particular, from \ref{itdef:anchorc} we have that if $\sigma^1_A(i) > 0$, then $i \in N_\Gamma(c)$, and therefore $\sigma_A(\oriented{ij}) > 0$ implies that $i \in N_\Gamma(c)$ and $j \in N_\Gamma(i)$.
This implies that \[W(\sigma_A) = \sum_{e \in E(\Gamma^\leftrightarrow)} \sigma_A(e) = \sum_{i \in N_\Gamma(c)} \sum_{j \in N_\Gamma(i)} \sigma_A(\oriented{ij}) = (1 + \gamma_A) \sum_{i \in N_\Gamma(c)} \sigma^1_A(i),\]
 where we used \eqref{equation:sigma1Ai} in the last step.
We also have as assumption that $W(\sigma_A) n \geq (1 + \eta)(a_1 + a_2) N = (1 + \eta)a_1(1 + \gamma_A) N$.
Combined with the previous calculations, we can write
\begin{equation}
	\sum_{i \in N_\Gamma(c)} \sigma^1_A(i) \geq (1 + \eta) a_1 \frac{N}{n}. \label{equation:sumweightsigma1A}
\end{equation}
The analogous reasoning also yields
\begin{equation}
	\sum_{i \in N_\Gamma(d)} \sigma^1_B(i) \geq (1 + \eta) b_1 \frac{N}{n}. \label{equation:sumweightsigma1B}
\end{equation}

Let $1 \leq i \leq N$ be arbitrary.
By \ref{item:fp-sizeWAWB}, we have 
\begin{equation}\label{eq:W_AcupW_B}
|W_A\cup W_B| \leq \frac{672 |V(T)|}{ \rho |V(T)|}<\frac{10^3}{\rho}\le  \frac{\rho n}{4N} 
< \frac{\rho |V_i|}{2},
\end{equation}
where the second to last inequality follows from the choice of $n_0$ and the inequalities $n \geq n_0$ and $N \leq t$;
and the last inequality follows from \eqref{equation:clustersize}.

For any $\mathcal F\subseteq \mathcal F_A$ denote by $V_1(\mathcal F)$ and $V_2(\mathcal F)$ the set of vertices in $\mathcal F$ at an odd distance, and respectively even distance, from $W_A$.
Similarly, for $\mathcal F\subseteq \mathcal F_B$, we define $V_1(\mathcal F)$ and $V_2(\mathcal F)$ to be the vertices of $\mathcal F$ at odd, and respectively even, distance from~$W_B$.

Choose $Q\subseteq V(G)\setminus (V_c\cup V_d)$ arbitrarily under the condition that $|Q\cap V_i|=\frac{\eta}{4}|V_i|$ for every $V_i \in \mathcal{P}\setminus \{V_c,V_d\}$.
We shall use this set in a second phase of the embedding to deal with a small proportion of the shrubs we fail to embed in a first attempt.
Let $U\subseteq V(G)\setminus (Q\cup V_c\cup V_d)$ be an arbitrary set such that $|U\cap V_i|=\frac{\eta}{4}|V_i|$ for every $V_i \in \mathcal{P}\setminus \{V_c,V_d\}$.
When we embed the small shrubs of $\mathcal F_A\cup \mathcal F_B$, we shall use the set $U$ as a buffer that will guarantee us we have enough free neighbours to choose a typical vertex from.
Thus, by construction, we have, for each $i \in [N]\setminus\{c,d\}$,
\begin{equation}
	|V_i \cap U| = |V_i \cap Q| = \frac{\eta}{4}|V_i| = \frac{\eta}{4}m. \label{equation:viUviR}
\end{equation}

\subsection{Proof of \Cref{lemma:treeembedding}: Embedding the seeds}
In the following, we find a suitable embedding of the seeds of $T$.
More precisely, we shall find an embedding $\varphi_0$ of $T[W_A\cup W_B]$ in $G[V_c,V_d]$,
and subsets $I_x \subseteq V(\Gamma)$, one for each $x \in W_A\cup W_B$,
that satisfy

\stepcounter{propcounter}
\begin{enumerate}[(\Alph{propcounter}\arabic*),topsep=0.7em, itemsep=0.5em]
	\item \label{it:embed-locationcd} $\varphi_0(W_A) \subseteq V_c$, and $\varphi_0(W_B) \subseteq V_d$;
	\item \label{it:embed-locationix} for each $x \in W_A$, $I_x \subseteq N_\Gamma(c)$,
	 and for each $x \in W_B$, $I_x \subseteq N_\Gamma(d)$;
	\item \label{it:embed-typical-A}$|N({\varphi_0}(x))\cap V_i \setminus (Q\cup U)|\ge (w(\fmat{ci})-\varepsilon)|V_i \setminus (Q\cup U)|$ and\\
	$|N({\varphi_0}(x))\cap V_i \cap Q|\ge 7\varepsilon |V_i|$ for each $x\in W_A$ and for each $i \in I_x$, 
	\item \label{it:embed-typical-B}$|N({\varphi_0}(x))\cap V_i \setminus (Q\cup U)|\ge (w(\fmat{di})-\varepsilon)|V_i \setminus (Q\cup U)|$ and\\
	$|N({\varphi_0}(x))\cap V_i \cap Q|\ge 7\varepsilon |V_i|$ for each $x\in W_B$ and for each $i \in I_x$, 
	\item \label{it:embed-degree-A}$\sum_{i \not \in I_x\cap I_{x'}}\sigma_A^1(i) |V_i| \le 5\sqrt{\varepsilon} n$, for every $x, x'\in W_A$, 
	\item \label{it:embed-degree-B} $\sum_{i \not\in I_x\cap I_{x'}}\sigma_B^1(i) |V_i| \le 5\sqrt{\varepsilon}n$, for every $x, x'\in W_B$.
\end{enumerate}

Let us digest the meaning of some of these properties.
Intuitively, \ref{it:embed-degree-A} says that for every $x,x' \in W_A$, both of the sets $I_x$ and $I_x'$ (and therefore also their intersection) contains most of $N_\Gamma(c)$;
and \ref{it:embed-typical-A} says that every $x \in W_A$ will be embedded in a vertex $\varphi_0(x) \in V_c$ such that for each $i \in I_x$, $\varphi_0(x)$ has large degree both to $V_i \setminus (Q \cup U)$ and $V_i \cap Q$.
\ref{it:embed-degree-B} and \ref{it:embed-typical-B} state the analogous properties for $W_B$.

To find the required embedding, we will embed the seeds $W_A\cup W_B$ in $G[V_c, V_d]$ one after the other.
Recall that by \Cref{definition:reducedgraph}, for any $j\in [N]$,  $N_\Gamma(j) \subseteq V(\Gamma)$ denotes the set of indices $i$ such that the cluster $V_i$ forms an $\eps$-regular pair together with $V_j$ of density at least $\dred$.

Let $V_c'\subseteq V_c$ be the set of vertices that are typical to all but at most $\sqrt{\varepsilon}|N_\Gamma(c)|$ of the sets $V_i \setminus (Q\cup U)$ with $i \in N_\Gamma(c)\setminus \{d\}$, typical
to all but at most $\sqrt{\varepsilon}|N_\Gamma(c)|$ sets $Q\cap V_i$ with $i \in N_\Gamma(c)\setminus \{d\}$,
and typical to $V_d$.
Define $V_d'\subseteq V_d$ analogously.

We estimate the size of $V_c'$.
By \Cref{lemma:typicalvertex2}, we see that at most $\sqrt{\eps} |V_c|$ vertices in $V_c$ are not included in $V'_c$ because of the first requirement; and similarly at most $\sqrt{\eps} |V_c|$ vertices fail to be included in $V'_c$ because of the second requirement.
By \Cref{lemma:typicalvertex}, at most $\eps |V_c| \leq \sqrt{\eps} |V_c|$ are lost in $V'_c$ because of the third requirement.
The analysis for $V'_d$ is the same, and thus, we get that
\begin{align*}
|V_c'|, |V_d'|\ge (1-3\sqrt{\varepsilon})|V_c| = (1-3\sqrt{\varepsilon})|V_d|.
\end{align*} 
From the definition of $V'_c$ and the previous bound, we have that every $v\in V_c'$ satisfies
\begin{align*}
\deg_G(v, V_d')\ge (\dred-\varepsilon)|V_d|-|V_d\setminus V_d'| \ge \frac{\dred}{2}|V_d| \overset{\eqref{eq:W_AcupW_B}}{> }  |W_A\cup W_B|,
\end{align*}
where in the last step, we also used $\dred \geq \rho$.
Similarly, for any vertex $v\in V_d'$, we have $\deg(v, V_c')\ge |W_A\cup W_B|$.
We can thus define an embedding $\varphi_0$ of $T[W_A \cup W_B]$ in $G[V_c', V_d']$ vertex by vertex, in a greedy fashion,
in such a way that $\varphi_0(W_A) \subseteq V'_c$ and $\varphi_0(W_B) \subseteq V'_d$.
Thus, \ref{it:embed-locationcd} holds.

For each $x\in W_A\cup W_B$, let $j_x\in \{c,d\}$ be such that $\varphi_0(x)\in V_{j_x}$.
We define $I_x\subseteq N_\Gamma(j_x)\setminus \{c,d\}$ as the set of indices $i$ of the clusters $V_i$ for which $\varphi_0(x)$ is typical to $V_i \setminus (Q\cup U)$ and typical to $V_i \cap Q$.
By construction, \ref{it:embed-locationix} holds.

Now, we check that the sets $I_x$ satisfy \ref{it:embed-typical-A} and \ref{it:embed-typical-B}.
Indeed, for each $i \in I_x$, the vertex $\varphi_0(x)$ is typical to the set $V_i\setminus (Q\cup U)$, and thus
\begin{align*}\deg_G(\varphi_0(x), V_i\setminus (Q\cup U))&\ge (d(V_{j_x}, V_i)-\varepsilon)|V_i\setminus (Q\cup U)| = (w(\fmat{j_x i})-\varepsilon)|V_i\setminus (Q\cup U)|.
\end{align*}
Similarly, $\varphi_0(x)$ is typical to $V_i \cap Q$, and thus
\begin{align*}
	\deg_G(\varphi_0(x), V_i\cap Q)
	& \ge (d(V_{j_x}, V_i)-\varepsilon)|V_i\cap Q|
  = 
 (w(\fmat{j_x i})-\varepsilon)|V_i\cap Q| \\
 & = (w(\fmat{j_x i})-\varepsilon)\frac{\eta}{4} |V_i| 
  \geq 
 \frac{\dred}{2}\frac{\eta}{4}|V_i|\overset{\eqref{eq:param-embedding}}{>}  7\varepsilon |V_i|,
\end{align*}
where in the second to last inequality, we used that $i\in N_\Gamma(j_x)$, and therefore $w(\fmat{j_x i})=d(V_{j_x}, V_i)\ge \dred$; together with $\dred \gg \eps$.
Thus we have obtained \ref{it:embed-typical-A} and \ref{it:embed-typical-B}.

Now, we show that the sets $I_x$ satisfy \ref{it:embed-degree-A} and \ref{it:embed-degree-B}.
We focus on showing \ref{it:embed-degree-A}, as the proof of \ref{it:embed-degree-B} follows mutatis mutandis.
Let $x, x' \in W_A$ be arbitrary.
By the definition of $V_c'$, there are at most $\sqrt{\varepsilon}|N_\Gamma(c)| \leq \sqrt{\eps} N$ indices $i\in N_\Gamma(c)\setminus \{d\}$ such that $\varphi_0(x)$ is not typical to $V_i\setminus (Q\cup U)$,
and at most $\sqrt{\varepsilon}|N_\Gamma(c)| \leq \sqrt{\varepsilon}N$ indices $i\in N_\Gamma(c)\setminus \{d\}$ for which $\varphi_0(x)$ is not typical with respect to $V_i\cap Q$.
The same can be said of $\varphi_0(x')$.

Therefore, $|N_\Gamma(c) \setminus (I_x\cup \{d\})| \leq 2 \sqrt{\eps} N$ and $|N_\Gamma(c) \setminus (I_{x'}\cup \{d\})| \leq 2 \sqrt{\eps} N$ hold,
from which we can deduce that
\[
 |N_\Gamma(c) \setminus (I_x \cap I_{x'})| \leq 4 \sqrt{\eps} N+1. 
 \]
Next, note that since $(\sigma_A, \sigma_B)$ is a $(\gamma_A, \gamma_B)$-skew-matching pair anchored in $\oriented{cd}$,
then $\sigma_A$ is a $\gamma_A$-skew-matching anchored in $N_\Gamma(c)$.
In particular, we have that $\sigma^1_A(i)\le w(\fmat{ci}) \leq 1$ for all $i\in N_\Gamma(c)$,
and $\sigma_A^1(i)=0$ for every $i \not \in N_\Gamma(c)$.
Therefore, we have 
\begin{equation*}\label{eq:I_x}
\sum_{i \not\in I_x\cap I_{x'}}\sigma^1_A(i) \leq \sum_{i \in  (N_\Gamma(c) \setminus I_x) \cup (N_\Gamma(c) \setminus I_{x'})\cup\{d\}}\sigma^1_A(i) \le |N_\Gamma(c) \setminus (I_x \cap I_{x'})|\leq 5 \sqrt{\varepsilon}N.
\end{equation*}
The bounds in \eqref{equation:clustersize} imply that $|V_i| = m \leq n/N$ holds for each $i \in [N]$, so the previous inequality yields \ref{it:embed-degree-A}.

\subsection{Proof of \Cref{lemma:treeembedding}: Allocating the shrubs}
To decide where we shall embed each shrub $F$ of $\mathcal F_A\cup \mathcal F_B$, we first decide where the colour class containing the root of $F$ should go.
We will use a probabilistic argument to assign a cluster $i_{F}$ for each shrub $F \in \mathcal F_A\cup \mathcal F_B$.
This will be done in such a way that the neighbours of $F$ (in $T$), which must be in $W_A \cup W_B$ and are already embedded via $\varphi_0$, can be appropriately joined whenever we embed the shrub $F$ in $V_{i_F}$.

Now, we turn to the details.
Recall that $\eta':= q\eta^3$. 
We will construct partitions $\mathcal F_A= \mathcal{F}_A^1\cup\dotsb\cup \mathcal{F}_A^N$
and $\mathcal{F}_B = \mathcal F_B^1\cup \dotsb \cup \mathcal{F}_B^N$ such that for each $i \in [N]$, and each ${\AB} \in \{A, B\}$, we have 
\stepcounter{propcounter}
\begin{enumerate}[(\Alph{propcounter}\arabic*),topsep=0.7em, itemsep=0.5em]
	\item \label{it:grouping-A}$|V_1(\mathcal{F}_{\AB}^i)|\le (1-\eta')\sigma_{\AB}^1(i)|V_i\setminus (Q\cup U)|$,
	\item \label{it:grouping-v2} $|V_2(\mathcal{F}_{\AB}^i)|\le (1-\eta')\sum_{j \in [N]}\frac{\gamma_{\AB}}{1+\gamma_{\AB}}\sigma_{\AB}(\oriented{ij})|V_j\setminus (Q\cup U)|$,
	\item \label{item:grouping-etaprime} if $\mathcal{F}_{\AB}^i\neq\emptyset$, then $\sigma^1_{\AB}(i)\ge \eta'$, 
		\item \label{item:grouping-neighbours} for each $F \in \mathcal{F}_A^i\cup \mathcal{F}_B^i$, we have that $i \in \bigcap_{x \in X_F} I_x$, where we set $X_F := \varphi_0( N_{T}(V(F))\cap (W_A\cup W_B) )$. 
\end{enumerate}

We show the existence of a suitable partition of $\mathcal F_A$, since the argument for $\mathcal F_B$ can be done in the same way.
Let
\[\Xi_c:=\{ i \in N_{\Gamma}(c) : \sigma^1_A(i)\ge \eta'\} \subseteq V(\Gamma).\]
In words, $\Xi_c$ are the indices of the clusters which we can use to allocate trees from $\mathcal{F}_A$ while complying with \ref{item:grouping-etaprime}.
For any $i \notin \Xi_c$, we will set $\mathcal{F}^i_A := \emptyset$.

Let $F \in \mathcal{F}_A$, and recall the definition of $X_F$ from \cref{item:grouping-neighbours}.
By \ref{item:fp-locationroots}--\ref{item:fp-seedsinWAWB}, we have 
$N_{T}(V(F)) \cap W_B=\emptyset$,
and therefore $X_F = \varphi_0( N_{T}(V(F))\cap (W_A\cup W_B) ) \subseteq \varphi_0(W_A) \subseteq V_c$.
Observe that by \ref{item:fp-twoseeds}, we have $|X_F| \leq 2$.
For each $F \in \mathcal{F}_A$, define
\begin{equation}\label{eq:I_T}
I_F := \Xi_c \cap \bigcap_{x\in X_F} I_x,
\end{equation}
so $I_F \subseteq \Xi_c$ are the indices corresponding to the clusters which are permissible for the allocation of $F$.
We will ultimately allocate 
$F$ to $\mathcal{F}^i_A$ while ensuring that $i \in I_F$ holds,
and this choice will ensure that \ref{item:grouping-neighbours} holds.

Before continuing with the proof, we need the following crucial estimate.

\begin{claim}
	For any $F \in \mathcal{F}_A$, we have
	\begin{align}
		\tilde{a}_1 := \left( 1 + \frac{\eta}{3} \right) a_1 \leq \sum_{i \in I_F} \sigma^1_A(i)|V_i \setminus (Q\cup U)|
		\label{equation:a1tilde}
	\end{align}
\end{claim}

\begin{proofclaim}
Let $S_F$ be the term in the right-hand side of \eqref{equation:a1tilde}.
Using \eqref{equation:viUviR},
\begin{align*}
	S_F 
	& \geq \left(1-\frac{\eta}{2} \right) m \sum_{i \in I_F} \sigma^1_A(i) \\
	& \overset{\eqref{eq:I_T}}{\ge} \left(1-\frac{\eta}{2}\right)m \left(\sum_{i\in N_\Gamma(c)}\sigma^1_A(i)-\sum_{i\not \in \Xi_c}\sigma_A^1(i)-\sum_{i\not\in \bigcap_{x\in X_F}I_x}\sigma_A^1(i)\right).
\end{align*}
We estimate the last two sums in the last expression.
From the definition of $\Xi_c$ and \eqref{equation:clustersize}, we have that $\sum_{i\not \in \Xi_c}\sigma_A^1(i) m \leq \eta' m (N - |\Xi_c| ) \leq \eta' Nm \leq \eta' n$;
and from \ref{it:embed-degree-A}, we get that $m \sum_{i\not\in \bigcap_{x\in X_F}I_x}\sigma_A^1(i)$ is at most $5 \sqrt{\eps} n$.
Using this, we obtain
\begin{align*}
		S_F \geq \left(1-\frac{\eta}{2} \right) \left( m \sum_{i \in N_\Gamma(c)}\sigma_A^1(i)-\eta' n-5\sqrt{\varepsilon} n \right)
	 \geq \left(1-\frac{\eta}{2} \right) \left( m \sum_{i \in N_\Gamma(c)}\sigma_A^1(i)-2\eta' n \right),
\end{align*}
where in the last inequality, we used the definition of $\varepsilon$ to deduce $5 \sqrt{\eps} \leq  \eta'$.
Next, using \eqref{equation:sumweightsigma1A} and \eqref{equation:clustersize}, we get
\begin{align*}
	S_F \geq \left(1-\frac{\eta}{2} \right) \left( m (1 + \eta) a_1 \frac{N}{n} -2\eta' n \right)
	 \geq \left(1-\frac{\eta}{2} \right) \left( (1 - \eps)(1 + \eta) a_1 -2\eta' n \right).
\end{align*}
Using \eqref{equation:klarge} and \eqref{equation:aibilarge}, we get that $a_1 \geq \eta k \geq \eta q n$,
and together with the choice of $\eta' = q \eta^3$, we deduce that
$2 \eta' n \leq 2 \eta^2 a_1$.
This gives
\begin{align*}
	S_F \geq \left(1-\frac{\eta}{2} \right) \left( (1 - \eps)(1 + \eta) -2\eta^2 \right) a_1
	\geq \left( 1 + \frac{\eta}{3} \right) a_1 = \tilde{a}_1,
\end{align*}
where in the last line, we used that $\eps \ll \eta \ll 1$, according to \eqref{eq:param-embedding}.
This proves the claim.
\end{proofclaim}

We will assign each shrub $F \in \mathcal{F}_A$ to $\mathcal{F}_A^i$, by choosing some $i \in I_F$, independently at random, with probability 
\begin{align*}
p^i_F := \probability[F \in \mathcal{F}_A^i]=\frac{\sigma^1_A(i)}{\sum_{\ell \in  I_F} \sigma^1_A(\ell)}.
\end{align*}
We note that, in particular, \eqref{equation:a1tilde} implies $I_F\neq\emptyset$, so the probabilities $p_F^i$ are well-defined.
Note that by construction, any such assignment will satisfy \ref{item:grouping-etaprime}--\ref{item:grouping-neighbours}.
We shall show that with non-zero probability also \ref{it:grouping-A}--\ref{it:grouping-v2} are satisfied.
This suffices to prove the existence of the desired partition.

We will estimate the expectation of $|V_1(\mathcal{F}_A^i)|$ and $|V_2(\mathcal{F}_A^i)|$ for each $i \in \Xi_c$,
then we will show that those values are concentrated around its expectation.
Note that $|V_1(\mathcal{F}_A^i)|$ can be written as the sum of the independent random variables $\{ X^i_F : F \in \mathcal{F}_A \}$, each of which takes the value $|V_1(F)|$ with probability $p^i_F$, and takes the value $0$ otherwise.
Thus, using the linearity of expectation, we have that
\begin{align*}
	\expectation[|V_1(\mathcal{F}_A^i)|]
	 & = \sum_{F \in \mathcal{F}_A} p^i_F |V_1(F)|
	 = \sum_{F \in \mathcal{F}_A} \frac{\sigma^1_A(i)}{\sum_{\ell \in I_F} \sigma^1_A(\ell)} |V_1(F)| \\
	 & \leq \sum_{F \in \mathcal{F}_A} \frac{\sigma^1_A(i)|V_i \setminus (Q\cup U)|}{\tilde a_1}|V_1(F)|,
\end{align*}
where in the last inequality, we used \eqref{equation:a1tilde}, and also the fact that $|V_i \setminus (Q \cup U)|$ is equal for each $i$ (which follows from \eqref{equation:viUviR} and from the fact that $Q\cap U = \emptyset$).
We have that $\sum_{F \in \mathcal{F}_A} |V_1(F)| = a_1$, and using this, we get
\begin{align*}
\expectation[|V_1(\mathcal{F}_A^i)|]
& \leq \frac{a_1}{\tilde a_1}
\sigma^1_A(i)|V_i \setminus (Q\cup U)|
= \frac{\sigma^1_A(i)|V_i \setminus(Q\cup U)|}{1+\frac{\eta}{3}} \\
& \le \left(1-\frac{\eta}{4}\right)\sigma^1_A(i)|V_i \setminus(Q\cup U)|.
\end{align*}
Using these estimates, we have
\begin{align*}
(1 - \eta') \sigma^1_A(i)|V_i \setminus (Q \cup U)| - \expectation[|V_1(\mathcal{F}^i_A)|]
& \geq \left( \frac{\eta}{4} - \eta' \right) \sigma^1_A(i)|V_i \setminus (Q \cup U)| \geq \frac{\eta\eta'}{16 N} n,
\end{align*}
where we used $\sigma^1_A(i) \geq \eta'$ (as $i \in \Xi_c$) and $|V_i \setminus (Q \cup U)| \geq |V_i|/2$ (which follows from \eqref{equation:viUviR}) in the last inequality.

Now, note that changing the allocation of a single shrub $F$ changes the value of $|V_1(\mathcal{F}_A^i)|$ by at most $c_F := |V_1(F)|$.
Since the partition is a $(\rho |V(T)|)$-fine-partition, we have by \ref{item:fp-smallshrubs} that $c_F \leq \rho |V(T)| = \rho k$ for every $F \in \mathcal{F}_A$.
Moreover, we have $\sum_{F \in \mathcal{F}_A} c_F \leq k$.
Therefore,
\[ \sum_{F \in \mathcal{F}_A} c_F^2 \leq (\rho k) \sum_{F \in \mathcal{F}_A} c_F \leq \rho k^2. \]
Putting all together, using \Cref{lemma:mcdiarmid} and recalling $\eta' = q \eta^3$ and $k \leq n$, we have
\begin{align*}
\probability\left[ |V_1(\mathcal{F}^i_A)| > (1 - \eta') \right. & \left. \!\! \sigma^1_{A}(i) |V_i \setminus (Q \cup U)| \right]
\leq \probability\left[ |V_1(\mathcal{F}^i_A)| > \expectation[|V_1(\mathcal{F}^i_A)|] + \frac{\eta \eta'}{16 N} n \right] \\
& \leq \exp \left( - \frac{2 \left( \frac{\eta \eta' n}{16 N} \right)^2 }{\sum_F c^2_F} \right)
=  \exp \left( - \frac{2 (\frac{q\eta^4}{16 N} n)^2 }{\sum_F c^2_F} \right) \\
& \leq \exp \left( - \frac{2 (\frac{q\eta^4}{16 N} n)^2 }{\rho k^2} \right) \leq \exp \left( - \frac{q^2\eta^8}{128 N^2 \rho} \right) < \frac{1}{2N},
\end{align*}
where in the last inequality, we used $\rho < q^2\eta^8/(128 N^2 \ln(2N))$, which indeed is valid by the choice of $\rho$.

Similarly, to control $|V_2(\mathcal{F}_A^i)|$, we need an upper bound on its expectation.
The main point here is that by \eqref{equation:sigma1Ai}, we have that $\sum_{j \in [N]} \frac{\gamma_A}{1 + \gamma_A} \sigma_A(\oriented{ij}) = \gamma_A \sigma^1_A(i)$,
and we can use that $|V_j \setminus (Q \cup U)|$ have the same size for each choice of $j \in [N]$ again.
Thus, \ref{it:grouping-v2} will hold if we ensure that $|V_2(\mathcal{F}^i_A)| \leq (1-\eta')\gamma_A \sigma^1_A(i) |V_i \setminus (Q \cup U)|$ holds for every $i \in [N]$.
To get this, we argue in the same way as in the calculation of $\expectation[|V_1(\mathcal{F}_A^i)|]$, but now, we use that $a_2 = \sum_{F \in \mathcal{F}_A} |V_2(F)|$ and also $a_2 = \gamma_A a_1$. We have that
\begin{align*}
\expectation[|V_2(\mathcal{F}_A^i)|]
& = \sum_{F \in \mathcal{F}_A} p^i_F |V_2(F)| \leq \frac{a_2}{\tilde{a}_1} \sigma_A^1(i) |V_i \setminus (Q \cup U)| \\
& \leq \left( 1 - \frac{\eta}{4} \right) \gamma_A \sigma^1_A(i) |V_i \setminus (Q \cup U)|.
\end{align*}
From this, essentially the same argument (using \Cref{lemma:mcdiarmid}, and our choice of $\rho$) as before gives that
\begin{align*}
\probability\left[ |V_2(\mathcal{F}_A^i)| > (1-\eta')\gamma_A \sigma^1_A(i) |V_i \setminus (Q \cup U)| \right] < \frac{1}{2N}.
\end{align*}

Thus, using a union bound over the at most $N$ possible values of $i$,
we see that the random allocation satisfies \ref{it:grouping-A}--\ref{it:grouping-v2} simultaneously for every $i$ with positive probability: this implies that the desired partition exists.

\subsection{Proof of \Cref{lemma:treeembedding}: Allocating the roots}
Now, we will determine suitable sets to `root' the shrubs.
Given a shrub $F \in \mathcal{F}_A \cup \mathcal{F}_B$,
recall that by \ref{item:fp-seedsinWAWB}, $F$ has a \emph{seed} $s_F \in W_A \cup W_B$.
In $T$, there must be a unique neighbour $r_F$ of $s_F$ inside $V(F)$; we will say that $r_F$ is the \emph{root} of $F$.
Alternatively, $r_F$ is the unique vertex in $V(F)$, which is closest to the root of $T$. 
Recall that the seed $s_F$ is already embedded in $V_c \cup V_d$ via $\varphi_0$.

This step aims to define sets $R_F$, one for each $F \in \mathcal{F}_A \cup \mathcal{F}_B$, that satisfy the following properties.
\stepcounter{propcounter}
\begin{enumerate}[(\Alph{propcounter}\arabic*),topsep=0.7em, itemsep=0.5em]
	\item \label{item:embedanchors-rightcluster} if $i \in [N]$ is such that $F \in \mathcal{F}^i_A \cup \mathcal{F}^i_B$, then $R_F \subseteq V_i \setminus (Q \cup U)$,
	\item \label{item:sizeofrJ} $|R_F| = |V_1(F)|$,
	\item \label{it:neighbourhood} $R_F \subseteq N_G(\varphi_0(s_F))$, where $s_F \in W_A \cup W_B$ is the seed of $F$,
	\item \label{item:embedanchors-disjoint} all the sets $R_F$ are pairwise-disjoint.
\end{enumerate}

We proceed as follows.
Let $i \in [N]$ for which $\mathcal{F}_A^i \cup \mathcal{F}_B^i$ is non-empty.
Suppose first that $i \in N_{\Gamma}(c)\setminus N_{\Gamma}(d)$.
In this case, we have $\sigma^1_B(i)=0<\eta'$ (the anchor of $\sigma_B$ fits in the $w$-neighbourhood of $d$) and thus $\mathcal{F}_B^i=\emptyset$ by~\ref{item:grouping-etaprime}.
Hence, in this case, we have $\mathcal{F}_A^i  \cup \mathcal{F}_B^i \subseteq \mathcal{F}_A$.
By~\ref{item:grouping-etaprime} again, we have that $\sigma_A^1(i) \ge  \eta'$, and therefore $\sigma_A^1(i) - \eps \geq (1 - \eta')\sigma_A^1(i)$ (where we used $\eps \leq (\eta')^2$).
Let $F \in \mathcal{F}^i_A$,
and let $s_F \in W_A$ be the seed of $F$.
As the anchor of $\sigma_A$ fits in the $w$-neighbourhood of $c$ by~\ref{itdef:anchorc},
we have that $\sigma_A^1(i)\le w(\fmat{ci})= d(V_c, V_i)$. 
By~\ref{item:grouping-neighbours}, \ref{it:embed-typical-A} and~\ref{it:grouping-A}, we get
\begin{align*}
|N_G(\varphi_0(s_F))\cap V_i \setminus (Q\cup U)|
& \ge (w(\fmat{ci}) - \varepsilon) |V_i \setminus (Q\cup U)|
\ge (\sigma^1_A(i) - \varepsilon) |V_i \setminus (Q\cup U)|\\
& \ge (1 - \eta')\sigma^1_A(i)|V_i \setminus (Q\cup U)|
 \ge |V_1(\mathcal{F}_A^i)|.
\end{align*}
This means that we can find a set $R_F \subseteq \bigl(N_G(\varphi_0(s_F))\cap V_i\bigr)\setminus (Q\cup U)$ of size precisely $|R_F|=|V_1(F)|$.
Moreover, we can do it in such a way that all the sets $R_{F'} \subseteq V_i$, for $F' \in \mathcal{F}_A^i$, are pairwise-disjoint.
(Indeed, for this fixed $i$ we have  $|N_G(\varphi_0(s_F))\cap V_i\setminus(Q\cup U)| \ge |V_1(\mathcal{F}_A^i)|$,
so choosing the $R_F$ one by one uses at most the remaining total demand and preserves enough capacity in $V_i\setminus(Q\cup U)$ for the shrubs still to be assigned.)
The argument is analogous if $i \in N_{\Gamma}(d) \setminus N_{\Gamma}(c)$, 
and in that case we can also find pairwise-disjoint sets $R_F \subseteq V_i$, for all shrubs $F \in \mathcal{F}_B^i$. 

It is only left to consider the case where $i \in N_{\Gamma}(c)\cap  N_{\Gamma}(d)$.
By \ref{itdef:anchorpartition}, we have that one of $w(\fmat{c i})$ or $w(\fmat{d i})$ is at least $\sigma^1_A(i) + \sigma^1_B(i)$.
Without loss of generality, we can assume that the former inequality holds, i.e. $w(\fmat{c i}) \ge \sigma_A^1(i)+\sigma_B^1(i)$
(the proof in the complementary case is analogous).
In this case, we first define $R_F$ for all shrubs $F \in \mathcal{F}_B^i$ as above, which we can do since $\sigma_B$ is anchored in $N_\Gamma(d)$.
Next, for the remaining shrubs $F \in \mathcal{F}_A^i$ with seeds $s_F$, we note that 
\begin{align*}
|N_G(\varphi_0(s_F))\cap V_i \setminus (Q\cup U) |
& \ge (w(\fmat{c i}) - \eps)|V_i \setminus (Q\cup U)| \\
& \geq (\sigma_A^1(i)+\sigma_B^1(i) - \eps )|V_i \setminus (Q\cup U)|
 \geq |V_1(\mathcal{F}_A^i)|+|V_1(\mathcal{F}_B^i)|.
\end{align*}
So we can again find a set $R_F \subseteq N(\varphi_0(s_F))\cap V_i \setminus (Q\cup U)$ of size $|R_F|=|V_1(F)|$.
Moreover, we can find such set that is also disjoint from every previously chosen $R_{F'} \subseteq V_i$.
(Here $(w(\fmat{ci})-\varepsilon)|V_i\setminus(Q\cup U)| \ge |V_1(\mathcal{F}_A^i)|+|V_1(\mathcal{F}_B^i)|$,	 so the same greedy capacity vs demand argument ensures pairwise disjointness.)
We have thus ensured \ref{item:embedanchors-rightcluster}--\ref{item:embedanchors-disjoint} hold.

\subsection{Proof of \Cref{lemma:treeembedding}: Finding suitable clusters}
Summarising what we have done so far: we have already embedded $W_A \cup W_B$ in $G[V_c \cup V_d]$ via $\varphi_0$,
we have already allocated each shrub $F \in \mathcal{F}_A \cup \mathcal{F}_B$ to a cluster $i \in V(\Gamma)$ where $V_1(F)$ will be embedded.
Moreover, for each such shrub, we have also reserved a `private' set $R_F \subseteq V_i$ of the right size $|V_1(F)|$.

Now, for each shrub $F \in \mathcal{F}^i_A \cup \mathcal{F}^i_B$, we would like to show the existence of another index $\ell \in V(\Gamma)$, such that $F$ can be embedded in $G[V_i, V_\ell]$.
In fact, we will show that we can find so many such $\ell$ that some of the graphs $G[V_i,V_\ell]$ can be used to embed $F$ even if we have previously embedded other shrubs and we need to avoid using the space used by these shrubs.

We introduce some notation that we use to work with partial functions from one set to another.
Let $Z, X' \subseteq X$, $Y$ be an arbitrary set and $\varphi: X' \rightarrow Y$ be a function. We say that $\varphi$ is a \emph{partial function} from $X$ to $Y$,
and by $\varphi(Z)$, we mean $\varphi(Z \cap X')$, i.e. the set of images of the elements in $Z$ which do have their image defined by $\varphi$.

Now, recall that $\varphi_0$ is the embedding we have already defined on $W_A \cup W_B$.
We say that a partial function $\varphi : V(T) \rightarrow V(G)$ is a \emph{partial embedding of $T$} if
\stepcounter{propcounter}
\begin{enumerate}[(\Alph{propcounter}\arabic*),topsep=0.7em, itemsep=0.5em]
	\item $\varphi$ extends $\varphi_0$,
	\item $\varphi$ is defined precisely on $W_A \cup W_B$ and $\bigcup_{F \in \mathcal{F}} V(F)$, where $\mathcal{F} \subseteq \mathcal{F}_A \cup \mathcal{F}_B$; and
	\item for each $i \in [N]$, we have $\varphi(V_1(\mathcal{F}^i_A) \cup V_1(\mathcal{F}^i_B)) \subseteq V_i$.
	\item for each $i \in [N]$, we have $\varphi(V_2(\mathcal{F}^i_A) \cup V_2(\mathcal{F}^i_B)) \cap V_i = \emptyset$
\end{enumerate}
Thus, a partial embedding is an embedding defined on the seeds $W_A \cup W_B$ and a subset of the shrubs in $\mathcal{F}_A \cup \mathcal{F}_B$.
The last property means that if a shrub satisfies $F \in \mathcal{F}^{i}_A \cup \mathcal{F}^{i}_B$ ---that is, $F$ is a shrub allocated to the $i$th cluster--- and its image is defined by $\varphi$, then we naturally must have $\varphi(V_1(F)) \subseteq V_i$.
In the next step, we will try to extend partial embeddings by including one more shrub at a time.
 
Given a partial embedding $\varphi$, ${\AB} \in \{A, B\}$ and $1 \leq i \leq N$, 
we will say that $\ell \in [N]$ is a \emph{target index for $(\varphi, \AB, i)$} if 
\stepcounter{propcounter}
\begin{enumerate}[(\Alph{propcounter}\arabic*),topsep=0.7em, itemsep=0.5em]
	\item \label{it:embed-Y-Xell} $|\varphi(\mathcal{F}^i_{\AB}) \cap V_\ell| < (1 + \eta'/2)(1 - \eta') \frac{\gamma_{\AB}}{1 + \gamma_{\AB}} \sigma_{\AB}(\oriented{i \ell}) |V_\ell \setminus (Q \cup U)|$, and
	\item \label{it:embed-Y-V(t)} $\frac{(\eta')^2}{2} \frac{\gamma_{\AB}}{1 + \gamma_{\AB}} \sigma_{\AB}(\oriented{i \ell}) |V_\ell \setminus (Q \cup U)| \geq \rho k$.
\end{enumerate}

Intuitively, if $\ell$ is a target index for $(\varphi, \AB, i)$; then we can try to extend $\varphi$ by embedding a shrub $F \in \mathcal{F}^i_{\AB}$ in $G[V_i, V_\ell]$.
Now, we prove that target indexes always exist.

\begin{claim}[Finding a target index] \label{claim:embedding-suitablecluster}
	Let $\varphi$ be a partial embedding,
	and let $\AB \in \{A,B\}$ and $1 \leq i \leq N$ be such that $\mathcal{F}^i_{\AB} \neq \emptyset$.
	Then a target index exists for $(\varphi, \AB, i)$.
\end{claim}

\begin{proofclaim}
	First, we shall observe that actually most of the clusters $V_\ell$ satisfy~\ref{it:embed-Y-V(t)}.
	Intuitively, if~\ref{it:embed-Y-V(t)} fails for a given $1 \leq \ell \leq N$, the skew-matching $\sigma_{\AB}$ will have very little weight on the pair $\oriented{i \ell}$.
	More precisely, we will show the following.
	Denote by $\mathcal Y \subseteq V(\Gamma)$ the set of all $\ell$ which do not satisfy~\ref{it:embed-Y-V(t)}.
	We claim that
	\begin{align} \label{eq:Y}
		\sum_{\ell \in \mathcal Y}\frac{\gamma_{\AB}}{1+\gamma_{\AB}}\sigma_{\AB}(\oriented{i \ell})|V_\ell \setminus (Q\cup U)| \leq \frac{\eta'}{200} \sum_{\ell \in [N]}\frac{\gamma_{\AB}}{1+\gamma_{\AB}}\sigma_{\AB}(\oriented{i \ell})|V_\ell \setminus (Q\cup U)|.
	\end{align}
	Let $X_\mathcal{Y}$ be the left-hand side of \eqref{eq:Y}.
	By definition, we have
	\begin{align*}
		\frac{\gamma_L}{1 + \gamma_L} \sigma_L(\oriented{i \ell}) |V_\ell \setminus (Q \cup U)| < \frac{2 \rho k}{(\eta')^2}
	\end{align*}
	for any $\ell \in \mathcal{Y}$,
	and therefore, we have
	\begin{align*}
		X_\mathcal{Y}
	 < \frac{2\rho kN}{(\eta')^2}
		< \frac{\eta (\eta')^2}{400}\frac{n}{N},
	\end{align*}
	where we have used that $k \leq n$ and $\rho < (\eta')^4 \eta / (800 N^2)$ in the second inequality.
	From \eqref{equation:clustersize}, we can deduce $n/(2N) \leq |V_{i} \setminus (Q \cup U)|$,
	and therefore
	\begin{align*}
		X_\mathcal{Y}
		 \leq \frac{\eta (\eta')^2}{200}|V_{i} \setminus (Q \cup U)|
		\leq \frac{\eta \eta'}{200} \sigma^1_L(i) |V_{i} \setminus (Q \cup U)|,
	\end{align*}
		where in the last inequality, we used that $\mathcal{F}_L^i \neq \emptyset$ (by assumption) together with \ref{item:grouping-etaprime}.
		Next, using \eqref{equation:gammalarge} first and then \eqref{equation:sigma1Ai}--\eqref{equation:sigma1Bi}, we arrive at
	\begin{align*}
		X_\mathcal{Y}
		\leq \frac{\eta'}{200} \gamma_L \sigma^1_L(i) |V_{i} \setminus (Q \cup U)| 
		 \leq \frac{\eta'}{200} \sum_{\ell \in [N]}\frac{\gamma_{\AB}}{1+\gamma_{\AB}}\sigma_{\AB}(\oriented{i \ell})|V_\ell \setminus (Q\cup U)|,
	\end{align*}
	so indeed \eqref{eq:Y} holds.
	
	Now, we argue that there must exist $\ell \notin \mathcal{Y}$ which satisfies \ref{it:embed-Y-Xell}, i.e. such that the inequality $|\varphi(\mathcal{F}^i_{\AB}) \cap V_\ell| < (1 + \eta'/2)(1 - \eta') \frac{\gamma_{\AB}}{1 + \gamma_{\AB}} \sigma_{\AB}(\oriented{i \ell}) |V_\ell \setminus (Q \cup U)|$ holds.
	For this, we will find a lower bound for the sum
	\[ S := \sum_{\ell \not\in \mathcal Y}\left(1+\frac{\eta'}{2} \right)(1-\eta')\frac{\gamma_{\AB}}{1+\gamma_{\AB}}\sigma_{\AB}(\oriented{i \ell})|V_\ell \setminus (Q\cup U)|. \]
	Indeed, we have
	\begin{align*}
		S
		& = \left(1+\frac{\eta'}{2} \right)(1-\eta') \frac{\gamma_{\AB}}{1+\gamma_{\AB}} \sum_{\ell \not\in \mathcal Y}\sigma_{\AB}(\oriented{i \ell})|V_\ell \setminus (Q\cup U)| \\
		& = \left(1+\frac{\eta'}{2} \right)(1-\eta') \frac{\gamma_{\AB}}{1+\gamma_{\AB}} \left( \sum_{\ell \in [N]} \sigma_{\AB}(\oriented{i \ell})|V_\ell \setminus (Q\cup U)| - \sum_{\ell \in \mathcal{Y}} \sigma_{\AB}(\oriented{i \ell})|V_\ell \setminus (Q\cup U)| \right) \\
		& \overset{\eqref{eq:Y}}{\geq} \left(1+\frac{\eta'}{2} \right)(1-\eta') \left( 1 - \frac{\eta'}{200} \right) \sum_{\ell \in [N]} \frac{\gamma_{\AB}}{1+\gamma_{\AB}} \sigma_{\AB}(\oriented{i \ell})|V_\ell \setminus (Q\cup U)| \\
		& \overset{\text{\ref{it:grouping-v2}}}{\geq}
		\left(1+\frac
		{\eta'}{2}\right) \left(1-\frac{\eta'}{200}\right) |V_2( \mathcal{F}_{\AB}^i)|\ge  \left(1+\frac{\eta'}{4}\right)|V_2( \mathcal{F}_{\AB}^i)|>|V_2(\mathcal F_{\AB}^i)|.
	\end{align*}
	Thus, $S > |V_2(\mathcal F_{\AB}^i)|$.
	To finish, we use that the vertices in $V_1(\mathcal F_{\AB}^i)$ are mapped only to $V_i$ via $\varphi$ (because $\varphi$ is a partial embedding), and thus \[\sum_{\ell\notin \mathcal{Y}} |\varphi(\mathcal F_{\AB}^i)\cap V_\ell| \leq \sum_{\ell\in N_\Gamma(i)}|\varphi(\mathcal F_{\AB}^i)\cap V_\ell|=   |\varphi(\mathcal F_{\AB}^i)\setminus V_i|\le |V_2(\mathcal F_{\AB}^i)| < S.\]
	Hence, there is an $\ell$ satisfying~\ref{it:embed-Y-Xell} and \ref{it:embed-Y-V(t)}, as required.
\end{proofclaim}

Given a partial embedding $\varphi$, $\AB \in \{A,B\}$ and $1 \leq i \leq N$, let $\mathcal{L}_{\varphi, \AB, i}$ be the set of all target indexes for $(\varphi, \AB, i)$.
In this notation, \Cref{claim:embedding-suitablecluster} ensures that $\mathcal{L}_{\varphi, \AB, i}$ is non-empty for $\varphi$, $\AB$ and $i$ whenever $\mathcal{F}^i_{\AB}$ is nonempty.
For technical reasons, we will need to compare the sets $\mathcal{L}_{\varphi, \AB, i}$ and $\mathcal{L}_{\varphi', \AB, i}$ for two different partial embeddings $\varphi, \varphi'$.
What we need is the following easy fact, which says that the set of target indexes is ``decreasing''.

\begin{claim} \label{claim:embedding-nestedtargets}
	Let $\AB \in \{A,B\}$ and $1 \leq i \leq N$,
	and let $\varphi, \varphi'$ be two partial embeddings.
	If $\varphi'$ extends $\varphi$, then $\mathcal{L}_{\varphi', \AB, i} \subseteq \mathcal{L}_{\varphi, \AB, i}$.
\end{claim}

\begin{proofclaim}
	Given $\ell \in \mathcal{L}_{\varphi', \AB, i}$, we need to show that $\ell \in \mathcal{L}_{\varphi, \AB, i}$.
	We need to check that \ref{it:embed-Y-Xell}--\ref{it:embed-Y-V(t)} hold for $\ell$ and $\varphi$.
	Since \ref{it:embed-Y-V(t)} depends only on $i, L, \ell$ and not on $\varphi, \varphi'$; we only need to check that \ref{it:embed-Y-Xell} holds.
	Since $\varphi'$ extends $\varphi$, we have that $|\varphi(\mathcal{F}^i_{\AB}) \cap V_\ell| \leq |\varphi'(\mathcal{F}^i_{\AB}) \cap V_\ell|$.
	Combined with the fact that \ref{it:embed-Y-Xell} holds for $\ell, \varphi', L, i$,
	this shows that \ref{it:embed-Y-Xell} holds for $\ell, \varphi, L, i$, as required.
\end{proofclaim}

\medskip

\subsection{Proof of \Cref{lemma:treeembedding}: First Embedding Phase}
After all this preparation, it is finally the turn to embed the shrubs in $\mathcal{F}_A \cup \mathcal{F}_B$.
We will fix an enumeration of all the shrubs
and define our embedding iteratively, incorporating one shrub at a time, respecting this order.
In a first attempt, we shall try to use \Cref{lem:embed-embedding the shrubs} for each shrub while avoiding the set $Q$ for all but a constant number of vertices (which will correspond to the adventitious roots).
It turns out that this procedure can fail for some shrubs; those will not be embedded in this step but instead deferred to the Second Embedding Phase
(as we will show later, the number of vertices of the shrub which will fail to be embedded in the First Embedding Phase will be tiny).

From now on, we will work with a fixed enumeration of the shrubs in $\mathcal{F}_A \cup \mathcal{F}_B$.
Let $J := | \mathcal{F}_A| + | \mathcal{F}_B|$ be the total number of shrubs. 
Let $F_1, F_2, \dotsc, F_J$ be an enumeration of the shrubs in $\mathcal{F}_A \cup \mathcal{F}_B$,
chosen so that each tree in $\mathcal{F}^1_A \cup \mathcal{F}^1_B$ appears before each tree in $\mathcal{F}^2_A \cup \mathcal{F}^2_B$ in the ordering, and so on.
Formally, we will have that for every $1 \leq i_1 < i_2 \leq N$;
	if $F_{j_1} \in \mathcal{F}^{i_1}_A \cup \mathcal{F}^{i_1}_B$
	and $F_{j_2} \in \mathcal{F}^{i_2}_A \cup \mathcal{F}^{i_2}_B$,
	then $j_1 < j_2$.
	
Recall that for each shrub $F \in \mathcal{F}_A \cup \mathcal{F}_B$, we have defined sets $X_F$ (in Step 3) and $R_F$ (in Step 4), and we have identified a root $r_F \in V(F)$ and a seed $s_F \in W_A \cup W_B$ (in Step 4).
If $1 \leq j \leq J$ is such that $F = F_j$, from now on, we will call these objects $X_j, R_j, r_j, s_j$, respectively. 

Before continuing with the embedding, we will formally obtain ``rooted shrubs'' from the shrubs $F_1, \dotsc, F_J$ to be able to apply \Cref{lem:embed-embedding the shrubs}. 
To do so, we just need to specify the roots and adventitious roots in each $F_j$, which we do as follows.
We have already (at the beginning of Step 4) identified the seed $s_j \in W_A \cup W_B$ and the root $r_j \in V(F_j)$ of $F_j$.
$T$ is a tree and therefore every vertex from $W_A\cup W_B$ can have at most one neighbour in $F_j$. 
By \ref{item:fp-twoseeds}, there are at most two neighbours of $W_A \cup W_B$ in $V(F_j)$.
If there is no neighbour of $W_A \cup W_B$ in $V(F_j)$ apart from $r_j$, we set $x_j = \emptyset$, i.e. $F_j$ has no adventitious root.
Otherwise, there exists $x_j \neq r_j$ in $V(F_j)$ which is a neighbour of $W_A \cup W_B$.
By \ref{item:fp-distance}, $x_j$ and $r_j$ have distance at least $4$ in $F_j$.
Thus, we can set $x_j$ as the adventitious root of $F_j$.
In all cases, $(F_j, r_j, x_j)$ defines a valid rooted shrub, and thus $\{(F_j, r_j, x_j)\}_{j=1}^J$ is a  family of vertex-disjoint rooted shrubs.

Now, we describe our embedding process in detail.
We begin by describing two certain invariants that will be useful to track during the construction of the embeddings.
We will say that a partial embedding $\varphi$ is \emph{reasonable} if
\stepcounter{propcounter}
\begin{enumerate}[(\Alph{propcounter}\arabic*)]
	\item \label{item:reasonableL} for all distinct $i, \ell \in [N]$, and all $L \in \{A, B\}$, \[|\varphi(\mathcal{F}^i_L) \cap V_\ell| < \frac{\gamma_{\AB}}{1 + \gamma_{\AB}} \sigma_{\AB}(\oriented{i \ell}) |V_\ell \setminus (Q \cup U)|.\]   
\end{enumerate}
In what follows, we will work with reasonable embeddings only.

The second invariant we need concerns the sets of vertices we would like to avoid using when embedding $F_j$.
For $1 \leq j \leq J$, we will say that $\varphi$ is a \emph{$j$-partial embedding} if it is defined only for (a subset of) the trees in $F_1, \dotsc, F_{j}$ (and none of the trees $F_{j+1}, \dotsc, F_{J}$).
In our setting, a $(j-1)$-partial embedding $\varphi$ is given, and we want to find space to embed $F_j$.
We do not want to use vertices in $Q$ (those are reserved for the adventitious roots and the Second Embedding Phase), we do not want to use vertices in $R_k$ for $k > j$, and obviously, we do not want to use vertices already used by $\varphi$.
This defines the set of \emph{forbidden vertices for $(\varphi, j)$ as}
\[ \forbidden_{\varphi, j} := Q \cup \im(\varphi) \cup \bigcup_{j+1 \leq k \leq J} R_k. \]
We also set $\forbidden_{\varphi, J+1} := Q \cup \im(\varphi)$.

The importance of reasonable $(j-1)$-partial embeddings $\varphi$ is that they leave sufficient space in each cluster while avoiding the forbidden vertices $\forbidden_{\varphi, j}$.
We express this as a claim.

\begin{claim} \label{claim:reasonablespace}
	Let $1 \leq j \leq J+1$ and let $\varphi$ be a reasonable $(j-1)$-partial embedding.
	Then, for any $1 \leq i \leq N$, we have
	$|V_i \setminus \forbidden_{\varphi, j}| \geq |V_i \cap U|$.
\end{claim}

\begin{proofclaim}
	We first estimate the vertices of $V_i$ already used by the $V_1$-parts of embedded shrubs.
	Since for each $i \in [N]$ we have
	$\varphi(V_1(\mathcal{F}^i_A)\cup V_1(\mathcal{F}^i_B)) \subseteq V_i$,
	it follows that
	\begin{align}
		|V_i \cap \varphi(V_1(\mathcal{F}_A \cup \mathcal{F}_B))|
		= \sum_{L \in \{A,B\}} \sum_{k \in [N]} |V_i \cap \varphi(V_1(\mathcal{F}^k_L))|
		= \sum_{L \in \{A,B\}} |V_i \cap \varphi(V_1(\mathcal{F}^i_L))|. \label{equation:ViV1bound-1}
	\end{align}
	Moreover, since $\varphi$ is a $(j-1)$-partial embedding, only shrubs among $F_1,\dots,F_{j-1}$ contribute to the left-hand side above; while for $k \ge j+1$, the set $R_k$ contributes $|V_1(F_k)|$ whenever $F_k \in \mathcal{F}^i_A \cup \mathcal{F}^i_B$.
	Therefore,
	\begin{equation}
		|V_i \cap \varphi(V_1(\mathcal{F}_A \cup \mathcal{F}_B))|
		+ \sum_{j+1 \le k \le J} |V_i \cap R_k|
		\le |V_1(\mathcal{F}^i_A)| + |V_1(\mathcal{F}^i_B)|.
		\label{equation:ViV1bound}
	\end{equation}
	
	On the other hand, we have
	\begin{align}
		|V_i \cap \varphi(V_2(\mathcal{F}_A \cup \mathcal{F}_B))|
		&= \sum_{L \in \{A, B\}} \sum_{k \in [N]\setminus\{i\}} |V_i \cap \varphi(V_2(\mathcal{F}^k_L))| \nonumber\\
		&\le \sum_{L \in \{A, B\}} \sum_{k \in [N]\setminus\{i\}} |V_i \cap \varphi(\mathcal{F}^k_L)| \nonumber\\
		&\overset{\text{\ref{item:reasonableL}}}{\le}
		\sum_{L \in \{A, B\}} \sum_{k \in [N]\setminus\{i\}}
		\frac{\gamma_L}{1 + \gamma_L}\sigma_L(\oriented{k i}) |V_i \setminus (Q \cup U)| \nonumber\\
		&\le |V_i \setminus (Q \cup U)|(\sigma^2_A(i) + \sigma^2_B(i)),
		\label{equation:ViV2bound}
	\end{align}
	where we used the definition of $\sigma^2_A$ and $\sigma^2_B$ in the last step.
	
	Putting it all together, we have
	\begin{align*}
		& |V_i \setminus \forbidden_{\varphi, j}|
		=
		\left|V_i \setminus \left(Q\cup \im(\varphi)\cup \bigcup_{j+1 \leq k \leq J} R_k \right)\right| \\
		&\ge |V_i \setminus Q|
		- |V_i \cap \varphi(V_1(\mathcal{F}_A \cup \mathcal{F}_B))|
		- |V_i \cap \varphi(V_2(\mathcal{F}_A \cup \mathcal{F}_B))|
		- \sum_{j+1 \leq k \leq J}|V_i \cap R_k| \\
		&\overset{\eqref{equation:ViV1bound}}{\ge}
		|V_i \setminus Q|
		- |V_i \cap \varphi(V_2(\mathcal{F}_A \cup \mathcal{F}_B))|
		- (|V_1(\mathcal{F}^i_A)| + |V_1(\mathcal{F}^i_B)|) \\
		&\overset{\text{\ref{it:grouping-A}}}{\ge}
		|V_i \setminus Q|
		- |V_i \cap \varphi(V_2(\mathcal{F}_A \cup \mathcal{F}_B))|
		- |V_i \setminus (Q \cup U)|(\sigma^1_A(i) + \sigma^1_B(i)) \\
		&\overset{\eqref{equation:ViV2bound}}{\ge}
		|V_i \setminus Q|
		- |V_i \setminus (Q \cup U)|(\sigma^1_A(i) + \sigma^1_B(i) + \sigma^2_A(i) + \sigma^2_B(i)) \\
		&\ge |V_i \setminus Q| - |V_i \setminus (Q \cup U)| \ge |V_i \cap U|,
	\end{align*}
	where in the second to last we used that $\sigma_A,\sigma_B$ are disjoint skew-matchings, and in the last inequality we used that $U$ and $Q$ are disjoint.
	This proves the claim.
\end{proofclaim}

The importance of target indices is that they will allow us to extend reasonable embeddings to reasonable embeddings, as shown in the following claim.

\begin{claim} \label{claim:reasonabletoreasonable}
	Let $\varphi$ be a reasonable partial embedding that does not embed the shrub $F_j$.
	Let $i, \AB$ and $\ell$ be such that $F_j \in \mathcal{F}^i_{\AB}$, and that $\ell$ is a target index for $(\varphi, \AB, i)$.
	Suppose $\varphi'$ extends $\varphi$ by embedding $V_1(F_j)$ in $V_i$ and $V_2(F_j)$ in $V_\ell$.
	Then $\varphi'$ is a reasonable partial embedding.
\end{claim}

\begin{proofclaim}
	To see that \ref{item:reasonableL} holds for $\varphi'$, it is enough to check the property for $i$, $\ell$, and $L$, since $\varphi$ is reasonable and the embedding only changed because of $F_j$.
	Thus, we have
	\begin{align*}
		|\varphi'(\mathcal{F}^i_L) \cap V_\ell|
		& \leq |\varphi(\mathcal{F}^i_L) \cap V_\ell| + |V_2(F_j)|
		 \leq |\varphi(\mathcal{F}^i_L) \cap V_\ell| + |V(F_j)|
		 \overset{\text{\ref{item:fp-smallshrubs}}}{\leq} |\varphi(\mathcal{F}^i_L) \cap V_\ell| + \rho k \\
		& \leq \left( (1 + \eta'/2)(1 - \eta') + \frac{\eta'^2}{2} \right) \frac{\gamma_{\AB}}{1 + \gamma_{\AB}} \sigma_{\AB}(\oriented{i \ell}) |V_\ell \setminus (Q \cup U)|, \\
		& < \frac{\gamma_{\AB}}{1 + \gamma_{\AB}} \sigma_{\AB}(\oriented{i \ell}) |V_\ell \setminus (Q \cup U)|,
	\end{align*}
	where we used both \ref{it:embed-Y-Xell}--\ref{it:embed-Y-V(t)} of the definition of the target index in the second to last inequality.
	This proves that $\varphi'$ is reasonable.
\end{proofclaim}
As we said before, an embedding of a considered shrub can fail.
Now, we can precisely describe the condition that ensures that the embedding of a shrub is successful or not.
Suppose $1 \leq j \leq J$ is given, and let $\varphi$ be a $(j-1)$-partial embedding.
Let $1 \leq i \leq N$ and $L \in \{A, B\}$ be such that $F_j \in \mathcal{F}^i_L$.
We say that a target index $\ell \in \mathcal{L}_{\varphi, \AB, i}$ is \emph{successful for $F_j$} if there exists $v \in R_j$ such that
\begin{equation}
	\deg_G \left(v, V_{\ell} \setminus \forbidden_{\varphi, j} \right) \geq 2 \eps |V_\ell|, \label{equation:successful}
\end{equation}

We can now state the following crucial claim. This claim ensures that if a successful target index exists, then a reasonable $(j-1)$-partial embedding
can be extended to a reasonable $j$-partial embedding by including $F_j$.

\begin{claim} \label{claim:embedding-successfulcluster}
	Suppose $1 \leq j \leq J$ is given, that $1 \leq i \leq N$ and $L \in \{A, B\}$ are such that $F_j \in \mathcal{F}^i_L$.
	Suppose $\varphi$ is a reasonable $(j-1)$-partial embedding such that $R_j \cap \im(\varphi) = \emptyset$,
	only adventitious roots are embedded in $Q$,
	 and that $\ell$ is a successful target index for $(\varphi, \AB, i)$.
	Then, there exists a reasonable $j$-partial embedding $\varphi'$, that extends $\varphi$ by embedding $F_j$ in $G[V_i, V_\ell]$.
	Moreover, $\varphi'(r_j) \in R_j$, $\varphi'(x_j) \in Q \cap V_i$, and $V(F_j)\setminus \{r_j, x_j\}$ is mapped to $(V_i \cup V_\ell)\setminus \forbidden_{\varphi, j}$.
\end{claim}

To be clear, if $x_j = \emptyset$, the `moreover' part just claims that $\varphi'(r_j) \in R_j$ and $V(F_j) \setminus \{r_j\}$ is mapped to $(V_i \cup V_\ell) \setminus \forbidden_{\varphi, j}$.

\begin{proofclaim}
	If $x_j = \emptyset$, define $P_j := (Q \cap V_i) \setminus \im(\varphi)$;
	otherwise $x_j$ is an adventitious root of $(F_j, r_j, x_j)$.
	Therefore, by \ref{item:fp-twoseeds}, $x_j$ has a unique neighbour, say $y_j$, in $W_A \cup W_B$.
	We define $P_j := (Q \cap V_i \cap N(\varphi_0(y_j))) \setminus \im(\varphi)$.
	Define $\tilde U := (V_i \cup V_\ell)\setminus \forbidden_{\varphi, j}$.
	Since $\ell$ is successful, we can choose a vertex $v_j \in R_j$ that satisfies \eqref{equation:successful}, by the assumption that $R_j \cap \im(\varphi) = \emptyset$, we know that $v_j$ is available to use in the embedding.
	We wish to apply \Cref{lem:embed-embedding the shrubs} with
	\begin{center}
		\begin{tabular}{c|c|c|c|c|c|c|c}
			object & $V_i$ & $V_\ell$ & $\tilde U$ & $P_j$ & $v_j$ & $(F_j, r_j, x_j)$ & $\eta/4$ \\
			\hline
			in place of & $X$ & $Y$ & $U$ & $P$ & $v$ & $(T,r,x)$ & $\tilde{\eta}$
		\end{tabular}
	\end{center}
	We will check that the required \ref{item:shrubembed-XY}--\ref{item:shrubembed-Tsmol} hold.
	
	We begin by checking \ref{item:shrubembed-XY}, which means that we need to prove that $|V_i \cap \tilde{U}| \geq \eta m / 4$ and $|V_\ell \cap \tilde{U}| \geq \eta m / 4$ hold.
	We have that $|V_i \cap \tilde{U}| = |V_i \setminus \forbidden_{\varphi, j}|$, so from \Cref{claim:reasonablespace}, we deduce that $|V_i \cap \tilde{U}| \geq |V_i \cap U|$.
	From \eqref{equation:viUviR}, we have that $|V_i \cap U| = \eta m/4$, so we conclude that $|V_i \cap \tilde{U}| \geq \eta m/4$.
	The proof of $|V_\ell \cap \tilde{U}| \geq \eta m/4$ is identical.
	
	To check \ref{item:shrubembed-P}, we must show that $|P_j| \geq 2 \eps m$ holds.
	We might suppose we are in the case where the adventitious root $x_j$ exists and we have set $P_j = (Q \cap V_i \cap N(\varphi_0(y_j))) \setminus \im(\varphi)$, the other case is more straightforward.
	By assumption, only adventitious roots are contained in $\im(\varphi) \cap Q$.
	Since only predecessors of $W_A \cup W_B$ can be adventitious roots, we deduce that $|\im(\varphi) \cap Q| \leq |W_A \cup W_B| \leq 5 \eps |V_i|$, we used \eqref{eq:W_AcupW_B} and $\rho \leq \eps$ in the last step.
	On the other hand, we have that $y_j \in N_{T}(V(F_j)) \cap (W_A \cup W_B)$
	thus $i \in I_{y_j}$ by \ref{item:grouping-neighbours},
	and therefore $|Q \cap V_i \cap N(\varphi_0(y_j))| \geq 7 \eps |V_i|$ follows from \ref{it:embed-typical-A} or \ref{it:embed-typical-B}.
	Hence, we have $|P_j| \geq |P \cap V_i \cap N(\varphi_0(y_j))| - |\!\im(\varphi) \cap Q| \geq 7  \eps |V_i| - 5 \eps |V_i| = 2 \eps |V_i|$, so \ref{item:shrubembed-P} holds.
	
	Finally, \ref{item:shrubembed-degv} holds because $v_j$ was chosen to satisfy~\eqref{equation:successful};
	and \ref{item:shrubembed-Tsmol} follows from the fact that we are working with a $\rho |V(T)|$-fine partition 
	Thus, we can use \Cref{lem:embed-embedding the shrubs} as intended.
	Define $\varphi'$ as $\varphi$ extended by the inclusion of $F_j$.
	By construction, it is simple to check that indeed $\varphi'(r_j) \in R_j$, $\varphi'(x_j) \in Q \cap V_i$, and $V(F_j)\setminus \{r_j, x_j\}$ is mapped to $(V_i \cup V_\ell)\setminus \forbidden_{\varphi, j}$, as required.
	Finally, $\varphi'$ is a reasonable $j$-partial embedding by \Cref{claim:reasonabletoreasonable}.
\end{proofclaim}

Now, we can iteratively define $j$-partial embeddings for all $0 \leq j \leq J$.
We will also define an increasing family $\mathcal{T}_j$ of trees for all $0 \leq j \leq J$, $\mathcal{T}_j$ will correspond to the set of \emph{postponed shrubs} among $F_1, \dotsc, F_j$ whose embedding failed and have been postponed to the Second Embedding Phase.

Recall that $\varphi_0$ corresponds to the already-defined embedding of $T[W_A \cup W_B]$ done in Step 2.
Clearly, $\varphi_0$ is reasonable and $(Q \cup \bigcup_{k=1}^J R_k) \cap \im(\varphi_0) = \emptyset$.
Also, let $\mathcal{T}_0 = \emptyset$.

Next, suppose we are given $1 \leq j \leq J$, the set $\mathcal{T}_{j-1}$ and a reasonable $(j-1)$-partial embedding $\varphi_{j-1}$ such that $\bigcup_{k=j}^J R_k \cap \im(\varphi_{j-1}) = \emptyset$ and $Q \cap \im (\varphi_{j-1})$ contains only adventitious roots.
In the $j$-th round of the embedding process
we will define a reasonable $j$-partial embedding $\varphi_j$ which extends $\varphi_{j-1}$, and $\mathcal{T}_j \supseteq \mathcal{T}_{j-1}$, as follows.
Let $\AB \in \{A,B\}$ and $1 \leq i \leq N$ such that $F_j \in \mathcal{F}^i_{\AB}$.
We consider the set $\mathcal{L}_{\varphi_{j-1}, \AB, i}$ of target indexes for $(\varphi_{j-1}, \AB, i)$.
We know by \Cref{claim:embedding-suitablecluster} that $\mathcal{L}_{\varphi_{j-1}, \AB, i}$ is non-empty.
There are two possibilities: either there exists a successful index $\ell \in \mathcal{L}_{\varphi_{j-1}, \AB, i}$ or not.
If there exists a successful $\ell$, then by \cref{claim:embedding-successfulcluster} there exists a reasonable $j$-partial embedding $\varphi_j$ which extends $\varphi_{j-1}$, and embeds $F_j$.
Moreover, $\varphi_{j}(F_j) \cap \forbidden_{\varphi_{j-1}, j} \subseteq \{\varphi_j(x_j),\varphi_j(r_j)\}$, $\varphi_j(x_j)\subseteq Q$, $ \varphi_j(r_j) \subseteq  R_j$, $R_j \subseteq \forbidden_{\varphi_{j-1}, j} \setminus Q$ and $R_k \cap (Q\cup R_j) = \emptyset$ for all $k \geq j+1$. 

Therefore, we have that $\bigcup_{k={j+1}}^J R_k \cap \im(\varphi_{j}) = \emptyset$, and $Q \cap \im (\varphi_{j})$ contains only adventitious roots.
In this case, $F_j$ was not postponed and we set $\mathcal{T}_j := \mathcal{T}_{j-1}$.
Otherwise, if there is no successful $\ell \in \mathcal{L}_{\varphi_{j-1}, \AB, i}$,
then, we set $\varphi_j := \varphi_{j-1}$ and $\mathcal{T}_j := \mathcal{T}_{j-1} \cup \{F_j\}$.
This process finishes with a final $J$-partial embedding $\varphi^\star := \varphi_J$ and a final set of postponed shrubs $\mathcal{T}^\star := \mathcal{T}_J$, by construction each shrub not in $\mathcal{T}^\star$ is embedded by $\varphi^\star$.

Here, we record a crucial property of how we defined the embeddings.
Every shrub $F_j \in \mathcal{T}^\star$ was postponed only if there was no successful index $\ell$ among all its target indexes.
Thus, we have
\begin{enumerate}[(\Alph{propcounter}\arabic*), resume]
	\item \label{item:nottypical} for every $F_j \in \mathcal{T}^\star$ such that $F_j \in \mathcal{F}^i_{\AB}$, every target index $\ell \in \mathcal{L}_{\varphi_{j-1}, \AB, i}$, and for every $v \in R_j$, we have
	\[ \deg_G \left(v, V_{\ell} \setminus \forbidden_{\varphi_{j-1}, j} \right) < 2 \eps |V_\ell|. \]
\end{enumerate}

\subsection{Proof of \Cref{lemma:treeembedding}: Second Embedding Phase}
To conclude the embedding of the whole tree $T$,
it only remains to extend the embedding $\varphi^\star$ by defining the embedding of every postponed shrub in $\mathcal{T}^\star$.
For each $1 \leq i \leq N$, let $\mathcal{T}^\star_{i, A} = \mathcal{T}^\star \cap \mathcal{F}^i_A$,  $\mathcal{T}^\star_{i, B} = \mathcal{T}^\star \cap \mathcal{F}^i_B$,
and $\mathcal{T}^\star_{i} = \mathcal{T}^\star_{i, A} \cup \mathcal{T}^\star_{i, B}$.
We will proceed in rounds, extending $\varphi^\star$ one shrub at a time.
We do this by incorporating each $F_j \in \mathcal{T}^\star$ in increasing order of~$j$.

The following claim implies that each set $\mathcal{T}^\star_{i}$ of postponed trees is tiny.

\begin{claim} \label{claim:postponedshrubsatypical}
	For each $1 \leq i \leq N$ and $\AB \in \{A, B\}$, we have $|V_1(\mathcal{T}^\star_{i, \AB})| \leq \eps |V_i|$.
\end{claim}

\begin{proofclaim}
	To see this, we need to carefully track what happens when a shrub $F$ is postponed and ends up belonging in $\mathcal{T}^\star_{i, \AB}$.
	
	Let $j_{\max}$ be the maximum value of $j$ such that $F_j \in \mathcal{T}^\star_{i, \AB}$, and let $\ell_{\max}$ be any target index in $\mathcal{L}_{\varphi_{j_{\max} - 1}, L, i}$.
	We clearly have $\ell_{\max} \neq i$.
	Now, let $j$ be arbitrary such that $F_j \in \mathcal{T}^\star_{i, \AB}$.
	By the choice of $j_{\max}$, we have that $j \leq j_{\max}$ and that $\varphi_{j_{\max}}$ extends $\varphi_{j}$,
	therefore by \Cref{claim:embedding-nestedtargets}, we have that $\mathcal{L}_{\varphi_{j_{\max} - 1}, L, i} \subseteq \mathcal{L}_{\varphi_{j - 1}, L, i}$.
	This implies that $\ell_{\max}$ is a target index for $(\varphi_{j - 1},L,i)$.
	Therefore, by \ref{item:nottypical}, we deduce that for any $v \in R_j$,
	\begin{equation}
		\deg_G \left(v, V_{\ell_{\max}} \setminus \forbidden_{\varphi_{j-1}, j} \right) < 2 \eps |V_{\ell_{\max}}|.
		\label{equation:atypicalstepj}
	\end{equation}
	Now, we make the following crucial observation:
	Because of the way we ordered the shrubs, during the definition of the partial embeddings $\varphi_j, \varphi_{j+1}, \dotsc, \varphi_{j_{\max}}$, we have only considered trees which are in $\mathcal{F}^i_A \cup \mathcal{F}^i_B$.
	This means, by \ref{item:embedanchors-rightcluster}, that $R_j \cup R_{j+1} \cup \dotsb \cup R_{j_{\max}}$ are all contained in $V_i$, and thus are disjoint from $V_{\ell_{\max}}$.
	Together with $\im(\varphi_{j-1}) \subseteq \im(\varphi_{j_{\max}})$, this implies that
	\[
	 V_{\ell_{\max}}\setminus \forbidden_{\varphi_{j_{\max}-1}, j_{\max}} \subseteq V_{\ell_{\max}}\setminus\forbidden_{\varphi_{j-1}, j}.
	\]
	Combined with \eqref{equation:atypicalstepj}, we obtain that for each $v \in R_j$,
	\begin{equation}
		\deg_G \left(v, V_{\ell_{\max}} \setminus \forbidden_{\varphi_{j_{\max}-1},j_{\max}} \right) < 2 \eps |V_{\ell_{\max}}|.
		\label{equation:atypicalstepmax}
	\end{equation}
	Let $R_{i, L} := \bigcup_{j} R_j$, where the union is taken over all $j$ such that the shrub $F_j$ belongs in $\mathcal{T}^\star_{i, L}$.
	Since $j$ was arbitrary in the argument above, \eqref{equation:atypicalstepmax} holds for each $v \in R_{i, L}$.
	
	\Cref{claim:reasonablespace} implies that $|V_{\ell_{\max}} \setminus \forbidden_{\varphi_{j_{\max}-1},j_{\max}}| \geq |V_i \cap U| = \eta |V_i|/4$.
	Since $\varepsilon \ll \dred \ll \eta$, we deduce that
	\[ 2 \eps |V_{\ell_{\max}}| < (\dred - \eps) \left| V_{\ell_{\max}} \setminus \forbidden_{\varphi_{j_{\max}-1},j_{\max}} \right|. \]
	Using the terminology of \Cref{lemma:typicalvertex}, this and \eqref{equation:atypicalstepmax} implies that every vertex in $R_{i, L}$ is `atypical' with respect to $V_{\ell_{\max}} \setminus \forbidden_{\varphi_{j_{\max}-1},j_{\max}}$.
	Thus, \Cref{lemma:typicalvertex} implies that $|R_{i, L}| \leq \eps |V_i|$,
	and we can conclude that 
	\[ 
	|V_1(\mathcal{T}^\star_{i, L})| = \sum_{j : F_j \in \mathcal{T}^\star_{i, L}} |V_1(F_j)| = \sum_{j : F_j \in \mathcal{T}^\star_{i, L}} |R_j| = |R_{i, L}|, 
	\]
	where we used \ref{item:sizeofrJ} in the second to last equality.
\end{proofclaim}

Now, using this claim, we can define the embedding of the shrubs, as we explained before.
The following claim takes care of one round of this process by embedding a single extra shrub.

\begin{claim} \label{claim:finalembedpostponedshrub}
	Let $i, \AB$ be such that $F_j \in \mathcal{T}^\star_{i, \AB}$,
	and suppose $\varphi$ is a reasonable partial embedding of $T$  which extends $\varphi^\star$ and does not embed $F_j$.
	Moreover, suppose that $\im(\varphi) \cap Q \cap V_i$ contains only adventitious roots and roots of trees in $\mathcal{T}^\star_{i}$;
	and that $\ell$ is a target index for $(\varphi, i, \AB)$.
	Then there is a reasonable partial embedding $\varphi'$ of $T$  which extends $\varphi$ by embedding $F_j$.
	Moreover, $r_j, x_j$ are mapped to $Q \cap V_i$, and $V(F_j) \setminus \{ r_j, x_j \}$ is mapped to $(V_i \cup V_\ell) \setminus (Q \cup \im (\varphi))$.
\end{claim}

\begin{proofclaim}
	Define $\tilde{U} := (V_i \cup V_\ell) \setminus (Q \cup \im(\varphi))$.
	Note that $\varphi$ is a reasonable $J$-partial embedding and therefore $\forbidden_{\varphi, J+1} = Q \cup \im(\varphi)$ and $\tilde{U} = (V_i \cup V_\ell) \setminus \forbidden_{\varphi, J+1}$.
	Then \Cref{claim:reasonablespace} implies that
	$|V_i \cap \tilde{U}| \geq |V_i \cap U|$ and $|V_\ell \cap \tilde{U}| \geq |V_\ell \cap U|$.
	
	Recall that $s_j \in W_A \cup W_B$ is the seed of $F_j$ and is already embedded by $\varphi$ (since it extends $\varphi^\star$).
	We define $P'_j := (N(\varphi(s_j)) \cap Q \cap V_i) \setminus \im( \varphi )$ as the set of vertices where we will look for a vertex to embed $r_j$.
	By assumption, $\im(\varphi) \cap Q \cap V_i$ contains only adventitious roots and roots of trees in $\mathcal{T}^\star_i$ and adventitious roots embedded in the First Embedding Phase, which amounts to at most $|W_A\cup W_B|$.
	Then \Cref{claim:postponedshrubsatypical} implies that 
	\[ 
	|\im(\varphi) \cap Q \cap V_i| \leq |V_1(\mathcal{T}^\star_{i, A})| + |V_1(\mathcal{T}^\star_{i, B})|+|W_A\cup W_B| \leq 3 \eps |V_i|.
	 \]

	On the other hand, we have that $s_j \in N_T(V(F_j)) \cap (W_A \cup W_B)$, thus $i \in I_{s_j}$ by \ref{item:grouping-neighbours},
	and therefore $|Q \cap V_i \cap N(\varphi(s_j))| \geq 7 \eps |V_i|$ follows from \ref{it:embed-typical-A} or \ref{it:embed-typical-B}.
	Putting all together, we have that $|P'_j| \geq 7 \varepsilon |V_i| - 3 \eps |V_i| \geq 4 \eps |V_i|$.
	Since $|V_\ell \cap \tilde{U}| \geq |V_\ell \cap U| = \frac{\eta}{4}|V_i| \geq \eps |V_i|$, \Cref{lemma:typicalvertex} implies that all but at most $\eps |V_i|$ vertices in $V_i$ have degree at least $(\dred - \eps)|V_\ell \cap \tilde{U}| \geq 2 \eps |V_i|$ into $V_\ell \cap \tilde{U}$, where the inequality follows from \eqref{eq:param-embedding}.
	In particular, we can select a vertex $v_j \in P'_j$ such that $\deg_G(v_j, V_\ell \cap \tilde{U}) \geq 2 \eps |V_i|$.
	
	Now, we define the set $P_j$ where the adventitious root, if it exists, will be embedded.
	The only extra care needed here is that we cannot use $v_j$.
	If $x_j = \emptyset$, define $P_j := (Q \cap V_i) \setminus (\im(\varphi) \cup \{ v_j \})$;
	otherwise $x_j$ is an adventitious root of $(F_j, r_j, x_j)$.
	Therefore, by \ref{item:fp-twoseeds}, $x_j$ has a unique neighbour, say $y_j$, in $W_A \cup W_B$.
	We define $P_j := (Q \cap V_i \cap N_G(\varphi(y_j))) \setminus (\im(\varphi) \cup \{ v_j \})$.

	We will apply \Cref{lem:embed-embedding the shrubs} with
	\begin{center}
		\begin{tabular}{c|c|c|c|c|c|c|c}
			object & $V_i$ & $V_\ell$ & $\tilde U$ & $P_j$ & $v_j$ & $(F_j, r_j, x_j)$ & $\eta/4$ \\
			\hline
			in place of & $X$ & $Y$ & $U$ & $P$ & $v$ & $(T,r,x)$ & $\tilde{\eta}$
		\end{tabular}
	\end{center}
	We will quickly check that the required \ref{item:shrubembed-XY}--\ref{item:shrubembed-Tsmol} hold.
	We have already checked \ref{item:shrubembed-XY}.
	The proof of \ref{item:shrubembed-P} follows along the same lines we have used to check the lower bound on $|P'_j|$, so we omit it.
	The choice of $v_j$ was done precisely to ensure that \ref{item:shrubembed-degv} holds, and \ref{item:shrubembed-Tsmol} follows by \ref{item:fp-smallshrubs} since we are working with a $\rho |V(T)|$-fine partition.
	Thus, we can use \Cref{lem:embed-embedding the shrubs} as intended.
	
	Define $\varphi'$ as $\varphi$ extended by the inclusion of $F_j$.
	By construction, we have that $\varphi'(r_j) = v_j \in Q \cap V_i$, if $x_j \neq \emptyset$ then $\varphi'(x_j) \in Q \cap V_i$, and $V(F_j)\setminus \{r_j, x_j\}$ is mapped to $(V_i \cup V_\ell)\setminus (Q \cup \im(\varphi))$, as required.
	Finally, $\varphi'$ is reasonable by \Cref{claim:reasonabletoreasonable}.
\end{proofclaim}
Given the last claim, we can conclude as follows.
Let $J^\star := |\mathcal{T}^\star|$, i.e. the number of postponed trees.
We set $\varphi^\star_0 := \varphi^\star$, which, of course, is a partial embedding of $T$ that extends $\varphi^\star$.
By construction, $\varphi^\star$ only embeds adventitious roots in $Q$.

Next, assume that $1 \leq k \leq J^\star$ is given and that we have constructed a reasonable partial embedding $\varphi^\star_{k-1}$ which extends $\varphi^\star$ and embeds precisely the first $k-1$ shrubs in $\mathcal{T}^\star$.
Assume also that, for each $1 \leq i \leq N$, $\im(\varphi^\star_{k-1}) \cap Q \cap V_i$ contains only adventitious roots and roots of trees in $\mathcal{T}^\star_i$.
Let $F_j$ be the $k$th remaining shrub in $\mathcal{T}^\star$.
Let $i, \AB$ be such that $F_j \in \mathcal{T}^\star_{i, \AB}$.
Let $\ell$ be a target index for $(\varphi^\star_{k-1}, L, i)$, this exists by \Cref{claim:embedding-suitablecluster}.
We apply \Cref{claim:finalembedpostponedshrub} to find an 
reasonable partial embedding $\varphi^\star_{k}$ which extends $\varphi^\star_{k-1}$ by embedding $F_j$, and moreover $\varphi^\star_{k}(V(F_j)) \cap (Q \cap V_i) \subseteq \{ \varphi^\star_k(r_j), \varphi^\star_k(x_j) \}$. 
This implies that $\varphi^\star_k$ embeds precisely the first $k$ shrubs of $\mathcal{T}^\star$ and also, for each $1 \leq i \leq N$, $\im(\varphi^\star_{k}) \cap Q \cap V_i$ contains only adventitious roots and roots of trees in $\mathcal{T}^\star_i$; so we obtained a reasonable partial embedding which allows us to continue this process.

At the end of this process, we obtain a final embedding $\varphi^\star_{J^\star}$.
By construction, $\varphi^\star_{J^\star}$ extends $\varphi^\star$ and also embeds all shrubs in $\mathcal{T}^\star$. Thus, it is an embedding of the whole tree $T$.
This (finally!) finishes the proof of \Cref{lemma:treeembedding}. \hfill \qed

\section{Concluding remarks}\label{sec:conclusion}

\subsection{Comparison with the proof proposed by Ajtai, Koml\'os, Simonovits, and Szemer\'edi}\label{ssec:comparision}

This comparison is based on personal communication between the second author and Mikl\'os Simonovits and we focus here only on their proof when specialised to the setting of \Cref{cor:dense-approx E-S}, i.e., to the so-called  \emph{approximate dense case}. Their full result is much stronger than \Cref{cor:dense-approx E-S}.

The major difference comes from the basic properties we may assume about the cluster graph. While in the present paper we may only assume that the cluster graph has minimum degree\footnote{For simplicity of explanation, in this subsection all degrees in the cluster graph are scaled by a factor corresponding to the size of a cluster.} slightly more than $k/2$ and maximum degree slightly larger than $k$, in the proof by Ajtai, Koml\'os, Simonovits, and Szemer\'edi, they may assume in addition an average degree slightly above $k$. 

In some parts, they actually exploit only the existence of one cluster of degree slightly more than $k$ together with the minimal degree condition, i.e., they use the same properties of the cluster graph as we do. 
This leads to some similarities between the two proofs (in the part treated in \Cref{ssection:strucprop-step1} and \Cref{ssection:strucprop-step2}). 
However, then they do use the additional property of the average degree of the cluster graph.
Missing this additional property in our setting, we need to fight harder to obtain a suitable structure, treated by more case distinctions (represented by the cases analysed in \Cref{ssec:skew-matching-case}--\Cref{ssection:strucprop-last}).

\subsection{Further variations on Erd\H{o}s--Sós}
\label{ssec:stateoftheart}

We believe that our approach, combined with Simonovits' Stability Method and ad hoc analysis of the close-to-extremal graphs, should provide a proof of the Erd\H{o}s--S\'os conjecture (\Cref{conj:E-S}) in the context of dense graphs.
We recall that the best-known exact results for \Cref{conj:E-S} are in the dense setting by Besomi, Pavez-Signé and Stein~\cite{BPS2021}, then extended to sparse setting by Pokrovsky~\cite{Pokrovskiy2024b}; both of those results assume the tree is sufficiently large and that $\Delta(T) = O(1)$.
A recent result by Reed and Stein~\cite{Reed-Stein-2025} proves the conjecture for large trees whose number of vertices is very close to the number of vertices of the host graph, without maximum degree restrictions on the trees; an extension was recently announced in~\cite{BruceTalk}.

Regarding steps towards proving \Cref{conj:KPR}, replacing the minimum-degree requirement by $k/2$ and the large-degree requirement by $k$ could probably also be dealt with using Simonovits' Stability Method, as well.
However, requiring only a sublinear number of vertices of high degree (instead of $\Omega(n)$ such vertices) would require a novel approach.

It is possible that our techniques would also be of help in studying related tree-embedding conjectures which combine minimum and maximum degree conditions; this is in contrast to \Cref{conj:KPR} which requires many vertices of large degree.
We summarise some of those conjectures, starting with the following conjecture by Besomi, Pavez-Signé and Stein~\cite[Conjecture 1.1]{BPS2019}
\begin{conjecture}[Besomi, Pavez-Sign\'e, Stein] \label{conjecture:bps-alpha}
Let $k\in \mathbb{N}$, let $\alpha \in \left(0,\frac13\right)$, and let $G$ be a graph with $\delta(G)\ge (1+\alpha)\frac k2$ and $\Delta(G)\ge 2(1-\alpha)k$. Then $G$ contains each tree with $k$ edges.
\end{conjecture}
Hyde and Reed~\cite{HydeReed2023} proved a relaxation of the $\alpha = 0$ case, where the maximum degree condition is replaced with $\Delta(G) \geq f(k)$, for some larger function of $k$.
Such a relaxation is in some sense necessary, since if $\delta(G) = \lfloor (k-1)/2\rfloor$ then one in fact needs $\Delta(G)$ to be at least quadratic in $k$, as shown by some examples~\cite[\S 3]{HydeReed2023}.

The maximum degree condition in the previous conjecture is not suspected to be best-possible for $\alpha = 1/3$.
In that case, there is a previous conjecture by Havet, Reed, Stein and Wood~\cite[Conjecture 1.1]{HRSW2020}.

\begin{conjecture}[Havet, Reed, Stein, Wood] \label{conjecture:hrsw}
	Let $k \in \mathbb{N}$, let $G$ be a graph of minimum degree at least $\lfloor 2k/3 \rfloor$ and maximum degree at least $k$. Then every tree with $k$ edges is a subgraph of $G$.
\end{conjecture}

Partial results by replacing the minimum or maximum degree conditions can be found in~\cite{HRSW2020}; and a proof in the particular case of spanning trees (i.e. $k = n-1$) was obtained by Reed and Stein~\cite{Reed-Stein-2023, Reed-Stein-2023b}.
Finally, Pokrovskiy, Versteegen and Williams~\cite{PokrovskiyVersteegenWilliams2025} recently made progress both on \Cref{conjecture:bps-alpha} and \Cref{conjecture:hrsw} for large, bounded-degree trees.

\addtocontents{toc}{\protect\setcounter{tocdepth}{1}}
\subsection*{Acknowledgment}\label{sec:aknowledgement}
We would like to thank Jan Hladk\'y for drawing our attention on the fact that \Cref{thm:approDense-MinMax} implies the approximate Erd\H os-S\'os Conjecture for dense graphs (\Cref{cor:dense-approx E-S}).
The fourth author thanks Giovanne Santos for discussions.

An extended abstract of this work appeared in the proceedings of EUROCOMB `23~\cite{DPRS23}.
Part of the work leading to this result was done while all the authors were affiliated with The Czech Academy of Sciences, Institute of Computer Science and were supported by the long-term strategic development financing of the Institute of Computer Science (RVO: 67985807) and by the Czech Science Foundation, grant number 19-08740S. 
N. Sanhueza-Matamala was supported by ANID-FONDECYT Iniciaci\'on Nº11220269 grant and ANID-FONDECYT Regular Nº1251121 grant.

\printbibliography

\end{document}